**С.В. Курапов**
**М.В. Давидовский**

# АЛГОРИТМИЧЕСКИЕ МЕТОДЫ КОНЕЧНЫХ ДИСКРЕТНЫХ СТРУКТУР

# АВТОМОРФИЗМ НЕСЕПАРАБЕЛЬНЫХ ГРАФОВ

на правах рукописи

**Курапов С.В.**

**Давидовский М.В.**

Рассматривается задача построение группы автоморфизмов графа. Автоморфизм графа есть отображение множества вершин на себя, сохраняющее смежность. Множество таких автоморфизмов образует вершинную группу графа или просто группу графа. Основой для построения группы автоморфизмов графа является понятие орбиты. Построение орбиты тесно связано с количественной оценкой вершины или ребра графа, называемой весом. Для определения веса элемента используются инварианты графа, построенные на спектре реберных разрезов графа и спектре реберных циклов. Вес элементов графа позволяет выявлять образующие циклы и формировать орбиты. Приведены примеры построения группы автоморфизмов для некоторых видов графов.

Для научных работников, преподавателей, студентов и аспирантов высших учебных заведений, специализирующихся в области прикладной математики и информатики.

# Содержание

# Введение

Из задуманной серии работ «Алгоритмические методы конечных дискретных структур», описывающих алгоритмические методы для построения алгебраических структур графа, в 2021 году типография Запорожского университета опубликовала работу «Алгоритмические методы конечных дискретных структур. Изоморфизм несепарабельных графов». Данная работа является 1-ой частью описания векторных инвариантов графа для определения изоморфизма. Задача определения автоморфизма графа тесно связанна с задачей определения изоморфизма графов применением, для своего решения, одинаковых графовых структур. 2-ая часть работы носит название «Алгоритмические методы конечных дискретных структур. Автоморфизм несепарабельных графов», где основой описания автоморфизма графов служит явление симметрической устойчивости весов вершин изоморфных графов и выбор опорного цикла состоящего из изометрических циклов.

Автоморфизм графа есть отображение множества вершин на себя, сохраняющее смежность. Множество таких автоморфизмов образует вершинную группу графа или просто группу графа. Группа подстановок на множестве ребер называется реберной группой графа, которая тесно связана с вершинной:

Два графа $G(V_g, E_g)$ и $H(V_h, E_h)$ изоморфны, если между их множествами вершин $V_g$ и $V_h$ существует взаимно однозначное соответствие, сохраняющее смежность [8,9,10,17,31,34]. Однако возможно установить изоморфизм между вершинами самого графа.

Изоморфизм графов устанавливается методом сравнения его инвариантов. Такими инвариантами графа являются : $n$ - количество вершин графа G, $m$ - количество ребер графа G, вектор локальных степеней графа G, цифровой инвариант реберного графа L(G), инвариант спектра реберных разрезов графа G, инвариант реберных циклов графа G, интегральный инвариант графа G и другие характеристики графа. Таким образом, возникает задача построения математических моделей определения группы автоморфизмов графа, использующая инварианты для определения изоморфизма.

Существенную роль при построении группы автоморфизма графа является теорема Фрухта утверждающая: что каждая конечная группа изоморфна группе автоморфизмов конечного неориентированного графа.

Множество всех перестановок вершин графа G образует группу [2,3,4,7,12,32,33], называемую группой автоморфизма графа G. Обозначим ее как *Aut*(G)[30]. Если мы сопоставляем некоторую вершину и ее образ, в который она переходит при автоморфизме графа, то можно рассматривать такое отображение как перестановку на множестве вершин графа V, а *Aut*(G) как группу перестановок на V.

Вершинная группа графа G индуцирует группу перестановок *Aut*(G), называемую *вершинной группой графа* G, она действует на множестве вершин V(G). Реберная группа графа G индуцирует группу перестановок Γ₁(G), называемую *реберной группой графа* G, она действует на множестве ребер E(G).

Для иллюстрации различия групп *Aut*(G) и Γ₁(G) рассмотрим граф К₄ – e, показанный на рис. 1. Вершинная группа *Aut*(К₄ – e) состоит из четырех перестановок:

$p_0 = (v_1)(v_2)(v_3)(v_4);$
$p_1 = (v_1)(v_2,v_4)(v_3);$
$p_2 = (v_1,v_3)(v_2)(v_4);$
$p_3 = (v_1,v_3)(v_2,v_4).$

Реберная группа перестановок также состоит из четырех перестановок. В реберной группе индуцируется перестановка в которой ребро e₄ остается на месте, e₁ меняется с e₂, а e₃ на e₅. Таким образом, реберная группа Γ₁(К₄ – e) состоит из следующих перестановок, индуцируемых указанными выше элементами реберной группы:

$q_0 = (e_1)(e_2)(e_3)(e_4)(e_5);$
$q_1 = (e_1,e_3)(e_2,e_5)(e_4);$
$q_2 = (e_1,e_2)(e_3,e_5)(e_4);$
$q_3 = (e_1,e_5)(e_2,e_3)(e_4).$

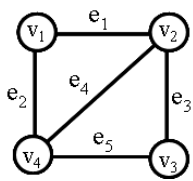

Рис. 1. Граф К₄ – e.

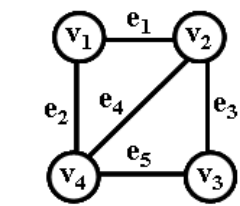 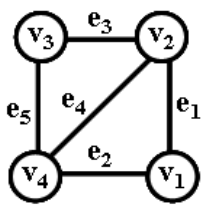 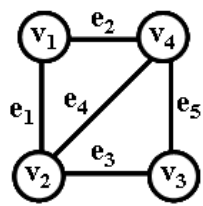 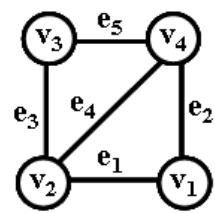

Рис. 2. Перестановка $p_0=(v_1)(v_2)(v_3)(v_4)$.    Рис. 3. Перестановка $p_1=(v_1)(v_2,v_4)(v_3)$.    Рис. 4. Перестановка $p_2=(v_1,v_3)(v_2)(v_4)$.    Рис. 5. Перестановка $p_3 = (v_1,v_3)(v_2,v_4)$.

Понятно, что реберная и вершинная группы графа К₄ – e в какой то степени изоморфны.

Но они, конечно, не могут быть идентичными, так как степень группы $\Gamma_1(K_4 - e)$ равна 5, а степень группы $Aut(K_4–e)$ равна 4 [34]. Очевидно, что для построения группы автоморфизмов графа следует рассматривать как вершинные группы, так и реберные группы сохраняющие расположение элементов матрицы инцидентности графа.

Как можно представить автоморфизм графа? Это можно представить в виде одинаковых рисунков графа [19] или в виде одинакового расположения элементов матрицы инциденций или матрицы смежностей (см. рис. 2 – рис.5).

Под рисунком графа будем понимать геометрическое представление топологического рисунка графа с учетом вращения его вершин [18-26].

Под сохранением расположения элементов матрицы инцидентности будем понимать адекватный вид матрицы с точностью до перестановки строк и столбцов (см. рис. 2 – рис.5).

Например, на рис. 2 – рис. 5 представлены перестановки графа $K_4$–e, рисунки графа, а так же вид расположения элементов матриц инцидентности и смежностей в зависимости от перестановки[11].

Для дальнейшего рассмотрения процесса определения группы автоморфизма нам понадобится векторные инварианты графа. Определим состав интегрального инварианта графа G необходимого для решения задачи определения группы автоморфизмов графа[16,17]:

Состав интегрального инварианта графа G:

- матрица инциденций графа B(G);
- $n$ – количество вершин графа G;
- $m$ – количество ребер графа G;
- множество изометрических циклов $C_\tau$ графа G;
- $k_s$ – количество уровней в спектре реберных разрезов графа G (max = 2);
- $k_c$ – количество уровней в спектре реберных циклов графа G (max = 1);
- $\xi_w(G)$ – кортеж весов ребер для спектра реберных разрезов графа;
- $\zeta_w(G)$ – кортеж весов вершин для спектра реберных разрезов графа.
- $\xi_\tau(G)$ – суммарный кортеж весов ребер для спектра реберных циклов графа;
- $\zeta_\tau(G)$ – кортеж весов вершин для спектра реберных циклов графа.
- $F_w(G)$ & $F_\tau(G)$ – инвариант спектра реберных разрезов графа G;
- $F_w(G)$ & $F_\tau(G)$ – инвариант спектра реберных циклов графа G;
- $F_w(G)$ & $F_w(G)$ & $F_\tau(G)$ & $F_\tau(G)$ – интегральный инвариант графа G.

Для графа $K_4$–e кортеж весов ребер имеет вид $\xi_w(K_4–e) = <7,7,7,4,7>$, а кортеж весов вершин $\zeta_w(K_4–e) = <14,18,14,18>$. Для удовлетворения условия существования изоморфизма внутри самого графа, очевидно, нужно переставлять местами только элементы равного веса. Выделим множества элементов равного веса. Для ребер: $\{e_1,e_2,e_3,e_5\}$ с весом равным 7 и $\{e_4\}$ с весом равным 4. Для вершин $\{v_1,v_3\}$ с весом равным 14, и $\{v_2,v_4\}$ с весом равным 18.

Запишем перестановки вершин с учетом весов:

$p_0 = (v_1)(v_2)(v_3)(v_4) \to (14)(18)(14)(18);$
$p_1 = (v_1)(v_2,v_4)(v_3) \to (14)(18,18)(14);$
$p_2 = (v_1,v_3)(v_2)(v_4) \to (14,14)(18)(18);$
$p_3 = (v_1,v_3)(v_2,v_4) \to (14,14)(18,18).$

Запишем перестановки ребер с учетом весов:

$q_0 = (e_1)(e_2)(e_3)(e_4)(e_5) \to (7)(7)(7)(4)(7);$
$q_1 = (e_1,e_3)(e_2,e_5)(e_4) \to (7,7)(7,7)(4);$
$q_2 = (e_1,e_2)(e_3,e_5)(e_4) \to (7,7)(7,7)(4);$
$q_3 = (e_1,e_5)(e_2,e_3)(e_4) \to (7,7)(7,7)(4).$

Таким образом, вводятся новые параметры графа – вес ребра и вес вершины, позволяющие проводить сравнительный анализ структур графа.

Количество перестановок в группе вершин и количество перестановок в группе ребер равны между собой. Поэтому для проведения практических вычислений удобно и достаточно пользоваться перестановками вершин[34].

Пусть $P$ группа перестановок на множестве элементов $M=\{1,2,…,n\}$. Подмножество $O \subset M$ называется орбитой группы $P$, если а) $(a)p \in O$ для любого $p \in P$ и любого $a \in O$, т.е. действие перестановок из $P$ на элементы $O$ не выводит за пределы $O$ [27];

б) любые два элемента из $O$ можно перевести друг в друга некоторой перестановкой из $P$.

Для определения орбит необходимо ответить на следующие вопросы:

1) какие вершины входят в орбиту;

2) сколько орбит имеет группа P;

3) какова длина каждой из орбит.

Ответы на поставленные вопросы можно получить, используя понятие веса вершины. Например, для графа $K_4$–е существуют только две орбиты, состоящие из 2-х элементов.

В сущности, применение к графу автоморфизма означает перенумерацию его вершин, причем отношение смежности должно сохраняться. Так что орбиты группы $Aut(G)$ - это просто классы «одинаковых» вершин. Здесь под «одинаковостью» следует понимать равенство весов элементов.

Более того, известно, что проблема распознавания принадлежности двух вершин произвольного графа одной орбите его группы автоморфизмов и проблема изоморфизма графов эквивалентны в том смысле, что любой алгоритм, эффективно решающий одну из этих проблем, может быть преобразован в эффективный алгоритм для другой. Построению группы автоморфизмов графа посвящено довольно большое количество публикаций [28-30,32,33]. К сожалению, эффективных алгоритмов для решения этих двух проблем, до сих пор не опубликовано.

Построение спектра реберных разрезов и спектра реберных циклов графа позволяют определить интегральный инвариант графа. Кортеж вершинных весов ставит в соответствие каждой вершине числовую характеристику называемую весом вершины, что позволяет определить множество вершин принадлежащих орбите.

Следует обратить внимание на тот факт, что направление вращения вершин графа для топологического рисунка графа в процессе перестановок может меняться.

С целью изучения процесса построение орбит с учетом существования интегрального инварианта графа будем рассматривать различные виды несепарабельных графов.

# Глава 1. АВТОМОРФИЗМ ГРАФА
## 1.1. Орбиты графа

Рассмотрим важное понятие орбиты группы подстановок. Пусть $\Gamma$ — произвольная группа подстановок на множестве V. Определим на V бинарное отношение ~, положив $u \sim v$ для $u,v \in V$ тогда и только тогда, когда в $\Gamma$ существует такая подстановка $s$, что $s(u)=v$. Очевидно, что отношение ~ является отношением эквивалентности и, следовательно, множество V разбивается на классы эквивалентных элементов: все элементы, входящие в один класс, переводятся подстановками из группы $\Gamma$ друг в друга, а элементы из разных классов друг в друга не переводятся. Эти классы называются *орбитами* группы $\Gamma$ [11].

При определении группы автоморфизмов графа, центральным вопросом является вопрос - выделения орбит группы подстановок. Это процесс соответствует построению орбит из множества вершин графа. Таким образом, *орбита графа* — это подмножество вершин равного веса, в которые может перейти вершина в результате действия автоморфизма.

Построение спектра реберных разрезов графа позволяет определить векторный интегральный инвариант графа. Кортеж весов вершин графа $\zeta_w(G)$ ставит в соответствие каждой вершине графа G числовую характеристику, называемую весом вершины. Кортеж реберных весов $\xi_w(G)$ также ставит в соответствие каждому ребру числовую характеристику, называемую весом ребра. Если в графе имеется вершины с одинаковым весом, то они принадлежат одной орбите.

Рассмотрим вопрос выделения орбит графа и нахождения группы автоморфизмов графа на следующем примере.

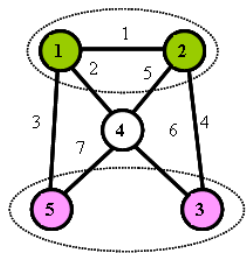

Рис. 1.1. Граф $G_1$.

***Пример 1.1.*** Найдем группу автоморфизмов графа $G_1$, изображенного на рисунке 1.1:

Вначале определим кортеж весов ребер $\xi_w(G_1)$ и кортеж вершинных весов $\zeta_w(G_1)$:

{ $\xi_w(G_1)$ = <8,9,8,8,9,7,7>, $\zeta_w(G_1)$ = <25,25,15,32,15>.

Выделяем в графе три орбиты: $o_1 = \{v_1, v_2\}$; $o_2 = \{v_3, v_5\}$; $o_3 = \{v_4\}$.

Соответственно ребра можно записать с учетом веса ребер и вершин в виде: $e_1$ = [8,25,25]; $e_2$ = [9,25,32]; $e_3$ = [8,15,25]; $e_4$ = [8,15,25]; $e_5$ = [9,25,32]; $e_6$ = [7,15,32]; $e_7$ = [7,15,32]. Первый элемент в тройке характеризует вес ребра, второй и третьей веса инцидентных вершин.

Распределим ребра по подмножествам с равным весом:

E = {[e_1],[e_2,e_5],[e_3,e_4],[e_6,e_7]}

Будем рассматривать все возможные перестановки вершин в орбитах, меняя местами только вершины в орбитах $o_1, o_2$, оставляя на месте вершину орбиты $o_3$.

Это перестановки $o_1$, $o_2$, $o_1 o_2$.

Совместная перестановка вершин в орбитах $o_1 o_2$ приводит к совпадению элементов матриц смежностей с точностью до перестановки строк и столбцов. А вот перестановка вершин в отдельно взятых орбитах по отдельности, не приводит к цели.

Автоморфизмы группы графа $G_1$, представлены двумя подстановками на множестве его вершин: одной тривиальной перестановкой $p_0$ и другой нетривиальной $p_1$.

$$Aut(G_1) = \begin{cases} p_0 = (1)(2)(3)(4)(5) - \textit{тождественный автоморфизм}, \\ p_1 = (12)(4)(35). \end{cases}$$

Стоит заметить, что при автоморфизме степени вершин в орбитах совпадают. У первой и второй вершин степень равняется трем, у третьей и пятой – двум, у четвертой – четырем.

### 1.2. Орбиты графа и перестановки

В состав интегрального инварианта графа входят два параметра: кортеж весов ребер $\xi_w(G_1)$ и кортеж весов вершин $\zeta_w(G_1)$. Дополним определение понятия орбиты. ***Будем считать орбитами графа, подмножество вершин имеющих равный вес в кортеже весов вершин графа.***

Таким образом, перестановка вершин будет осуществляться не хаотически, а только в рамках вершин с равным весом. Естественно, что процесс построения группы автоморфизмов графа облегчается, когда существует приемлемый рисунок графа.

***Пример 1.2.*** На следующем рисунке показаны орбиты вершин графа $G_2$.

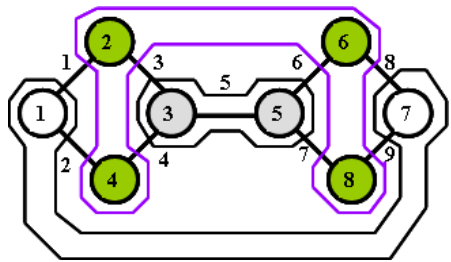

Рис. 1.2. Подмножества вершин с весами графа $G_2$.

Выделим из векторного интегрального инварианта графа $G_2$ представленного на рис. 1.2, кортеж весов вершин $\zeta_w(G_2)$ = <30,35,52,35,52,35,30,35>. Образуем подмножества вершин с равным весом V = {[$v_1,v_7$][$v_2,v_4,v_6,v_8$][$v_3,v_5$]}. Таким образом, множество вершин графа V разбито на орбиты.

Распределим множество ребер по эквивалентным классам E={[$e_1,e_2,e_8,e_9$],[$e_3,e_4,e_6,e_7$],[$e_5$]}.

Займемся перестановкой местами вершин графа. Будем переставлять местами вершины только в рамках орбит. Перестановку местами вершин представленных на рис. 1.3 будем считать начальной (нулевой, тривиальной) перестановкой. Выделим маршруты длиной 6 между вершинами $v_1$ и $v_7$:

$<v_1,v_2,v_3,v_5,v_6,v_7>$; $<v_1,v_2,v_3,v_5,v_8,v_7>$; $<v_1,v_4,v_3,v_5,v_6,v_7>$; $<v_1,v_4,v_3,v_5,v_8,v_7>$;

$<v_1,v_2,v_3,v_5,v_8,v_7>$; $<v_1,v_4,v_3,v_5,v_6,v_7>$;

Как видно, перестановка вершин $v_1$ и $v_7$ автоматически производит перестановку вершин $v_3$ и $v_5$. Выделим маршруты между парами вершин $(v_2,v_4),(v_2,v_6),(v_2,v_8),(v_4,v_6),(v_4,v_8)$ и $(v_6,v_8)$. Маршруты длиной 4 $<v_2,v_3,v_5,v_6>,<v_2,v_3,v_5,v_8>,<v_4,v_3,v_5,v_6>,<v_4,v_3,v_5,v_8>$ входят в маршруты между вершинами $v_1$ и $v_7$, поэтому могут отдельно не рассматриваться. А вот маршруты $<v_2,v_3,v_4>$ и $<v_6,v_5,v_8>$ длиной 3, не входят в маршруты между вершинами $v_1$ и $v_7$, поэтому четверку вершин $[v_2,v_4,v_6,v_8]$ следует разбить на две пары вершин $[v_2,v_4]$ и $[v_6,v_8]$.

Будем осуществлять перестановки вершин.

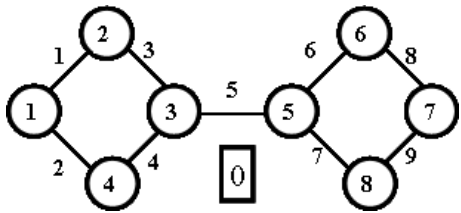 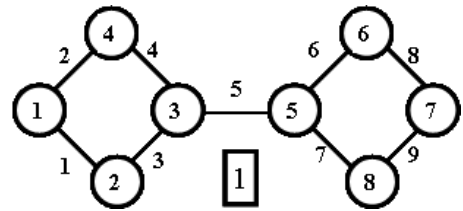

Рис. 1.3. 0 – элемент группы.     Рис. 1.4. 1-ый элемент группы.

$$p_0 = \begin{pmatrix} 1 & 2 & 3 & 4 & 5 & 6 & 7 & 8 \\ 1 & 2 & 3 & 4 & 5 & 6 & 7 & 8 \end{pmatrix} = (1)(2)(3)(4)(5)(6)(7)(8)$$

Для рисунка 1.4 определим следующую перестановку, меняя местами вершины $v_2$ и $v_4$:

$$p_1 = \begin{pmatrix} 1 & 2 & 3 & 4 & 5 & 6 & 7 & 8 \\ 1 & 4 & 3 & 2 & 5 & 6 & 7 & 8 \end{pmatrix} = (1)(2\ 4)(3)(5)(6)(7)(8).$$

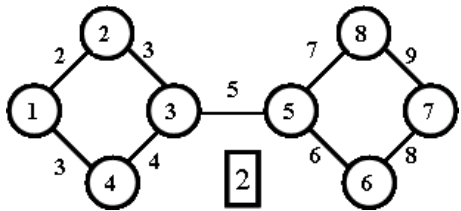 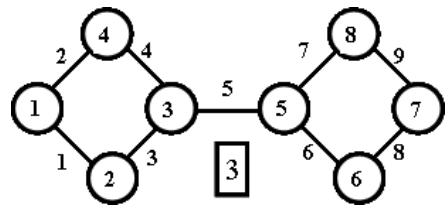

Рис. 1.5. 2-ой – элемент группы.     Рис. 1.6. 3-ий элемент группы.

Для рисунка 1.5 определим следующую перестановку, меняя местами вершины $v_6$ и $v_8$:

$$p_2 = \begin{pmatrix} 1 & 2 & 3 & 4 & 5 & 6 & 7 & 8 \\ 1 & 2 & 3 & 4 & 5 & 8 & 7 & 6 \end{pmatrix} = (1)(2)(3)(4)(5)(6\ 8)(7).$$

Для рисунка 1.6 определим следующую перестановку, одновременно меняя местами вершины $v_6$ и $v_8$, а также $v_2$ и $v_4$:

$$p_3 = \begin{pmatrix} 1 & 2 & 3 & 4 & 5 & 6 & 7 & 8 \\ 1 & 4 & 3 & 2 & 5 & 8 & 7 & 6 \end{pmatrix} = (1)(2\,4)(3)(5)(6\,8)(7).$$

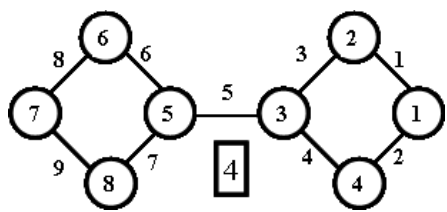 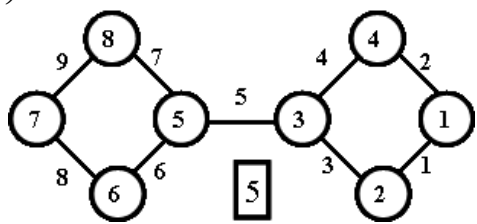

Рис. 1.7. 4-ый элемент группы.    Рис. 1.8. 5-ый – элемент группы.

Для рисунка 1.7 определим следующую перестановку, одновременно меняя местами пару вершин $(v_1,v_7)$, а также $(v_3,v_5)$, $(v_2,v_6)$, $(v_4,v_8)$.

$$p_4 = \begin{pmatrix} 1 & 2 & 3 & 4 & 5 & 6 & 7 & 8 \\ 7 & 6 & 5 & 8 & 3 & 2 & 1 & 4 \end{pmatrix} = (1\,7)(2\,6)(3\,5)(4\,8).$$

Для рисунка 1.8 определим следующую перестановку, одновременно меняя местами пару вершин $(v_1,v_7)$, а также пары $(v_3,v_5)$, $(v_2,v_8)$, $(v_4,v_6)$.

$$p_5 = \begin{pmatrix} 1 & 2 & 3 & 4 & 5 & 6 & 7 & 8 \\ 7 & 8 & 5 & 6 & 3 & 4 & 1 & 2 \end{pmatrix} = (1\,7)(2\,8)(3\,5)(4\,6).$$

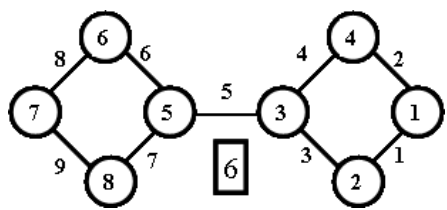 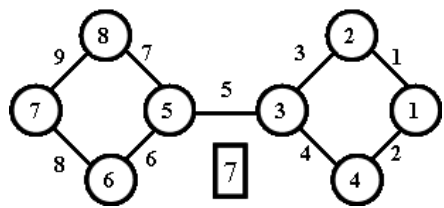

Рис. 1.9. 6-ой элемент группы.    Рис. 1.10. 7-ой элемент группы.

Для рисунка 1.9 определим следующую перестановку, одновременно меняя местами пару вершин $(v_1,v_7)$, а также пару $(v_3,v_5)$, и вершины $(v_2,v_6,v_4,v_8)$.

$$p_6 = \begin{pmatrix} 1 & 2 & 3 & 4 & 5 & 6 & 7 & 8 \\ 7 & 6 & 5 & 8 & 3 & 4 & 1 & 2 \end{pmatrix} = (1\,7)(2\,6\,4\,8)(3\,5).$$

Для рисунка 1.10 определим следующую перестановку, одновременно меняя местами пару вершин $(v_1,v_7)$, а также вершины $(v_3,v_5)$ и $(v_2,v_8,v_4,v_6)$.

$$p_7 = \begin{pmatrix} 1 & 2 & 3 & 4 & 5 & 6 & 7 & 8 \\ 7 & 8 & 5 & 6 & 3 & 2 & 1 & 4 \end{pmatrix} = (1\,7)(2\,8\,4\,6)(3\,5).$$

Граф $G_2$ рассматриваемый в данном примере, имеет только одну связь между вершинами одной орбиты - ребро $(v_3,v_5)$. Перестановка местами вершин $(v_1,v_7)$ определяет и перестановку местами вершин $(v_2,v_4,v_6,v_8)$, так как эти вершины связаны одним ребром с вершинами $(v_1,v_7)$. В свою очередь, вершины $(v_2,v_4,v_6,v_8)$ связаны одним ребром и с вершинами $(v_3,v_5)$. Поэтому, одна только перестановка вершин $(v_1,v_7)$ приводит к перестановки вершин $(v_2,v_4,v_6,v_8)$ и к перестановки вершин $(v_3,v_5)$. Это перестановка ($p_4,p_5,p_6,p_7$). Если не переставлять местами вершины ($v_1,v_7$), то можно кроме начальной подстановки определить еще три. Это перестановка ($p_0,p_1,p_2,p_3$) для вершин ($v_2,v_4$) и ($v_6,v_8$) по отдельности, а также вариант их совместной перестановки. Следует обратить внимание, на тот факт, что перестановка вершин в орбите возможна только в случае существования определенного маршрута в графе.

Так как возможно только переставлять вершины в рамке одной орбиты, то в графе $G_2$ можно выделить четыре орбиты: $o_1 = (v_3,v_5)$ с весом 52; $o_2 = (v_1,v_7)$ с весом 30; $o_3 = (v_2,v_4)$ с весом 35; $o_4 = (v_6,v_8)$ с весом 35.

Составим таблицу Кэли умножения группы автоморфизмов графа $G_2$ [11]:

$Aut(G)=$

| ○ | $p_0$ | $p_1$ | $p_2$ | $p_3$ | $p_4$ | $p_5$ | $p_6$ | $p_7$ |
|---|---|---|---|---|---|---|---|---|
| $p_0$ | $p_0$ | $p_1$ | $p_2$ | $p_3$ | $p_4$ | $p_5$ | $p_6$ | $p_7$ |
| $p_1$ | $p_1$ | $p_0$ | $p_3$ | $p_2$ | $p_7$ | $p_6$ | $p_5$ | $p_4$ |
| $p_2$ | $p_2$ | $p_3$ | $p_0$ | $p_1$ | $p_6$ | $p_7$ | $p_4$ | $p_5$ |
| $p_3$ | $p_3$ | $p_2$ | $p_1$ | $p_0$ | $p_5$ | $p_4$ | $p_7$ | $p_6$ |
| $p_4$ | $p_4$ | $p_6$ | $p_7$ | $p_5$ | $p_0$ | $p_3$ | $p_1$ | $p_2$ |
| $p_5$ | $p_5$ | $p_7$ | $p_6$ | $p_4$ | $p_3$ | $p_0$ | $p_2$ | $p_1$ |
| $p_6$ | $p_6$ | $p_4$ | $p_5$ | $p_7$ | $p_2$ | $p_0$ | $p_3$ | $p_1$ |
| $p_7$ | $p_7$ | $p_5$ | $p_4$ | $p_6$ | $p_1$ | $p_2$ | $p_0$ | $p_3$ |

***Пример 1.3.*** Рассмотрим несепарабельный граф, представленный на рис. 1.11.

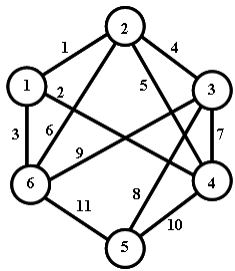

Рис. 1.11. Граф $G_3$.

Определим веса вершин, используя инвариант, построенный на спектре реберных разрезов.

Определим группу автоморфизмов графа $G_3$: Кортеж весов вершин: $\zeta_w(G_3)$ = <10,16,16,14,10,14>. Выделим подмножества вершин с равным весом: $o_1 = \{v_1,v_5\}$ с весом 10, $o_2 = \{v_2,v_3\}$ с весом 16, $o_3 = \{v_4,v_6\}$ с весом 14 [16].

$$Aut(G_3) = \begin{cases} p_0 = (1)(2)(3)(4)(5)(6); \\ p_1 = (1\ 5)(2\ 3)(4)(6); \\ p_2 = (1)(2)(3)(4\ 6)(5); \\ p_3 = (1\ 5)(2\ 3)(4\ 6). \end{cases}$$

Составим таблицу умножения Кэли для этой группы автоморфизмов:

$Aut(G_3)=$

| ∘ | $p_0$ | $p_1$ | $p_2$ | $p_3$ |
|---|---|---|---|---|
| $p_0$ | $p_0$ | $p_1$ | $p_2$ | $p_3$ |
| $p_1$ | $p_1$ | $p_0$ | $p_3$ | $p_2$ |
| $p_2$ | $p_2$ | $p_3$ | $p_0$ | $p_1$ |
| $p_3$ | $p_3$ | $p_2$ | $p_1$ | $p_0$ |

Если сформировать орбиты $o_1 = \{v_1, v_5\}$, $o_2 = \{v_2, v_3\}$, $o_3 = \{v_4, v_6\}$. То для определения всех перестановок нужно перебрать 7 вариантов перестановок $(o_1)$, $(o_2)$, $(o_3)$, $(o_1, o_2)$, $(o_1, o_3)$, $(o_2, o_3)$, $(o_1, o_2, o_3)$ определяемых суммированием сочетаний $C_3^1 + C_3^2 + C_3^3 = 3+3+1$.

Таким образом, для построения автоморфизмов графа следует переставлять местами строки и столбцы матрицы смежностей. Однако, строить матрицы и располагать элементы по строкам и столбцам дело довольно трудоемкое и хлопотное.

Рассмотрим следующий алгоритм, основанный на представления элементов матрицы однострочными структурными числами [1]. Перестановка называется транспозицией, если переставляются местами только два элемента множества, тогда как остальные элементы остаются на своих местах. Любую перестановку элементов множества М можно осуществить посредством нескольких транспозиций. Поэтому будем рассматривать орбиты, состоящие только из двух вершин.

Запишем матрицу смежностей графа в виде однострочных структурных чисел характеризующих номера столбцов и строк. Однострочное структурное число описывается множеством вершин смежных к заданной вершине. Для графа $G_3$ структурное число, состоящее из однострочных структурных чисел, имеет вид:

```
1 :  2  4  6
2 :  1  3  4  6
3 :  2  4  5  6
4 :  1  2  3  5
5 :  3  4  6
6 :  1  2  3  5
```

Осуществим перестановку вершин в орбите $o_1$. Будем менять местами строку 1 и строку 5. Затем заменим цифры 1 и 5 в структурных числах.

```
5 :  3  4  6           5₁:  3  4  6            5 :  3  4  6
2 :  1  3  4  6        2 :  1₅ 3  4  6         2 :  3  4  1  6
3 :  2  4  5  6        3 :  2  4  5₁ 6         3 :  5  2  4  6
4 :  1  2  3  5        4 :  1₅ 2  3  5₁        4 :  1  2  3  5
1 :  2  4  6           1₅:  2  4  6            1 :  2  4  6
6 :  1  2  3  5        6 :  1₅ 2  3  5₁        6 :  1  2  3  5
```

Произведем построчное сравнение соответствующих записей исходного и вычисленного структурных чисел. Если перестановка вершин в графе существует, то записи в соответствующих строках структурных чисел должны совпадать, как элементы равных подмножеств.

В этом примере, записи строк совпадают, следовательно, перестановка корректна.

Осуществим перестановку вершин в орбите $o_2$. Будем менять местами строку 2 и строку 3. Затем заменим цифры 2 и 3 в структурных числах.

| 1 : 2 4 6      | 1 : $2_3$ 4 6            | 1 : 3 4 6     |
| 3 : 2 4 5 6    | $3_2$: $2_3$ 4 5 6       | 2 : 3 4 5 6   |
| 2 : 1 3 4 6    | $2_3$: 1 $3_2$ 4 6       | 3 : 1 2 4 6   |
| 4 : 1 2 3 5    | 4 : 1 $2_3$ $3_2$ 5      | 4 : 1 2 3 5   |
| 5 : 3 4 6      | 5 : $3_2$ 4 6            | 5 : 2 4 6     |
| 6 : 1 2 3 5    | 6 : 1 $2_3$ $3_2$ 5      | 6 : 1 2 3 5   |

Записи в соответствующих строках содержат разные элементы, перестановка не корректна.

Осуществим перестановку вершин в орбите $o_3$. Будем менять местами строку 4 и строку 6. Затем заменим цифры 4 и 6 в структурных числах.

| 1 : 2 4 6      | 1 : 2 $4_6$ $6_4$        | 1 : 2 4 6     |
| 2 : 1 3 4 6    | 2 : 1 3 $4_6$ $6_4$      | 2 : 1 3 4 6   |
| 3 : 2 4 5 6    | 3 : 2 $4_6$ 5 $6_4$      | 3 : 2 4 5 6   |
| 6 : 1 2 3 5    | $6_4$: 1 2 3 5           | 4 : 1 2 3 5   |
| 5 : 3 4 6      | 5 : 3 $4_6$ $6_4$        | 5 : 3 4 6     |
| 4 : 1 2 3 5    | $4_6$: 1 2 3 5           | 6 : 1 2 3 5   |

Записи в соответствующих строках равны, следовательно, перестановка корректна.

Осуществим перестановку вершин в орбитах $o_1 o_2$. Будем менять местами строку 1 и строку 5, строку 2 и строку 3. Затем заменим цифры 1 и 5, 2 и 3 в структурных числах.

| 5 : 3 4 6      | $5_1$: $3_2$ 4 6               | 1 : 2 4 6     |
| 3 : 2 4 5 6    | $3_2$: $2_3$ 4 $5_1$ 6         | 2 : 1 3 4 6   |
| 2 : 1 3 4 6    | $2_3$: $1_5$ $3_2$ 4 6         | 3 : 2 4 5 6   |
| 4 : 1 2 3 5    | 4 : $1_5$ $2_3$ $3_2$ $5_1$    | 4 : 1 2 3 5   |
| 1 : 2 4 6      | $1_5$: $2_3$ 4 6               | 5 : 3 4 6     |
| 6 : 1 2 3 5    | 6 : $1_5$ $2_3$ $3_2$ $5_1$    | 6 : 1 2 3 5   |

Записи совпадают, следовательно, перестановке корректна.

Осуществим перестановку вершин в орбитах $o_1 o_3$. Будем менять местами строку 1 и строку 5, строку 4 и строку 6. Затем заменим цифры 1 и 5, 4 и 6 в структурных числах.

| 5 : 3 4 6      | $5_1$: 3 $4_6$ $6_4$           | 1 : 3 4 6     |
| 2 : 1 3 4 6    | 2 : $1_5$ 3 $4_6$ $6_4$        | 2 : 3 4 5 6   |
| 3 : 2 4 5 6    | 3 : 2 $4_6$ $5_1$ $6_4$        | 3 : 1 2 4 6   |
| 6 : 1 2 3 5    | 6 : $1_5$ 2 3 $5_1$            | 4 : 1 2 3 5   |
| 1 : 2 4 6      | $1_5$: 2 $4_6$ $6_4$           | 5 : 2 4 6     |
| 4 : 1 2 3 5    | 4 : $1_5$ 2 3 $5_1$            | 6 : 1 2 3 5   |

Перестановка корректна.

Осуществим перестановку вершин в орбитах $o_2 o_3$. Будем менять местами строку 2 и строку 3, строку 4 и строку 6. Затем заменим цифры 2 и 3, 4 и 6 в структурных числах.

| 1 : 2 4 6      | 1 : $2_3$ $4_6$ $6_4$          | 1 : 3 4 6     |
| 3 : 2 4 5 6 | **3₂: 2₃ 4₆ 5 6₄** | 2 : 3 4 5 6 |
| 2 : 1 3 4 6 | **2₃: 1 3₂ 4₆ 6₄** | 3 : 1 2 4 6 |
| 6 : 1 2 3 5 | **6₄: 1 2₃ 3₂ 5** | 4 : 1 2 3 5 |
| 5 : 3 4 6 | 5 : **3₂ 4₆ 6₄** | 5 : 2 4 6 |
| 4 : 1 2 3 5 | **4₆: 1 2₃ 3₂ 5** | 6 : 1 2 3 5 |

Перестановка не корректна.

Осуществим перестановку вершин в орбитах o₁o₂o₃. Будем менять местами строку 2 и строку 3, строку 4 и строку 6, строку 5 и строку 1. Затем заменим цифры 2 и 3, 4 и 6, 1 и 5 в структурных числах.

| | | |
|---|---|---|
| 5 : 3 4 6 | **5₁: 3₂ 4₆ 6₄** | 1 : 2 4 6 |
| 3 : 2 4 5 6 | **3₂: 2₃ 4₆ 5₁ 6₄** | 2 : 1 3 4 6 |
| 2 : 1 3 4 6 | **2₃: 1₅ 3₂ 4₆ 6₄** | 3 : 2 4 5 6 |
| 6 : 1 2 3 5 | **6₄: 1₅ 2₃ 3₂ 5₁** | 4 : 1 2 3 5 |
| 1 : 2 4 6 | **1₅: 2₃ 4₆ 6₄** | 5 : 3 4 6 |
| 4 : 1 2 3 5 | **4₆: 1₅ 2₃ 3₂ 5₁** | 6 : 1 2 3 5 |

Перестановки вершин в орбитах o₃,o₁o₂, o₁o₂o₃ порождает изоморфизм матриц смежностей графа $G_3$. Таким образом, если структурное число, полученное в результате преобразований, совпадает с первоначальным структурным числом, то матрицы смежности изоморфны

***Пример 1.4.*** Рассмотрим граф $G_4$ представленный на рис. 1.12,а. Спектр реберных разрезов графа:

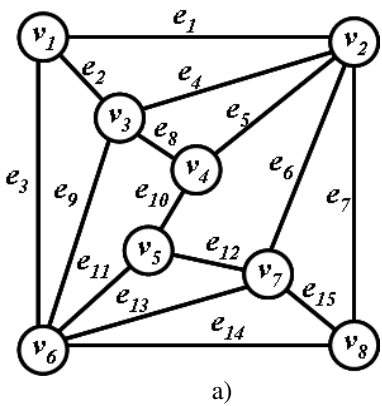
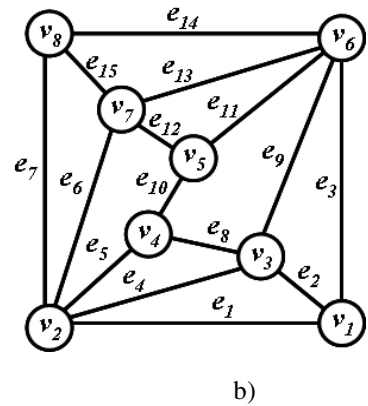

а) b)

Рис. 1.12. Автоморфизмы графа $G_4$.

Кортеж весов ребер: <5,4,5,8,5,9,5,5,9,6,5,5,8,5,4>;

Кортеж весов вершин: <14,32,26,16,16,32,26,14>.

Обозначим орбиты: $o_1 = \{v_1,v_8\}$, $o_2 = \{v_2,v_6\}$, $o_3 = \{v_3,v_7\}$, $o_4 = \{v_4,v_5\}$. Таким образом, образуется четыре орбиты. Очевидно, что нужно переставлять только вершины, принадлежащие одной орбите. Для определения автоморфизма графа нужно перебрать все возможные сочетания орбит, количество которых можно выразить формулой

$$k = C_n^1 + C_n^2 + C_n^3 + ... + C_n^n \tag{1.1}$$

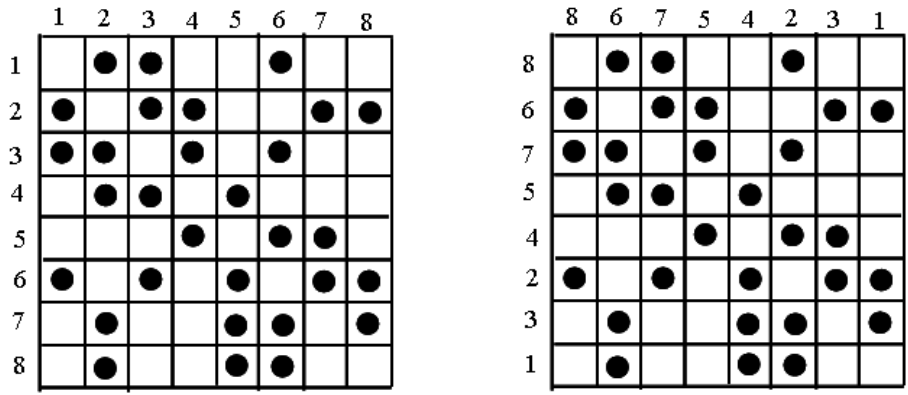

Рис. 1.13. Матрица смежностей графа $G_4$.

Перебрав все сочетания орбит определим, что возможна только одна нетривиальная перестановка вершин $p_1$.

Так как все орбиты имеют разные веса и в каждой орбите по два элемента, то нетривиальную перестановку можно определить, выделив DFS-дерево алгоритмом поиска в глубину. Затем произвести замену вершин в орбитах.

Например, построим дерево:

$(v_1,v_3)+(v_3,v_4)+(v_4,v_5)+(v_5,v_7)+(v_7,v_8)+(v_8,v_6)+(v_8,v_2);$

Поменяем местами вершины в орбитах:

$(v_8,v_7)+(v_7,v_5)+(v_5,v_4)+(v_4,v_3)+(v_3,v_1)+(v_1,v_2)+(v_1,v_6).$

В результате получаем нетривиальную подстановку $p_1$:

$$p_1 = \begin{pmatrix} 1 & 2 & 3 & 4 & 5 & 6 & 7 & 8 \\ 8 & 6 & 7 & 5 & 4 & 2 & 3 & 1 \end{pmatrix}.$$

Результат перестановки $p_1$ представлен на рис. 1.12,b. Автоморфизмы графа $G_4$ состоят из двух перестановок $p_0$ и $p_1$.

$$Aut(G_4) = \begin{cases} p_0 = (1)(2)(3)(4)(5)(6)(7)(8)(9)(10); \\ p_1 = (1\ 8)(2\ 6)(3\ 7)(4\ 5). \end{cases}$$

***Пример 1.5.*** Рассмотрим граф $G_5$ представленный на рис. 1.14,a. Спектр реберных разрезов:

Кортеж весов ребер : <5,7,8,6,5,8,4,10,5,5,4,7,6,5,5>;
Кортеж весов вершин : <20,24,32,14,32,14,20,24>.

Выделим орбиты в графе $o_1 = \{v_1,v_7\}$, $o_2 = \{v_2,v_8\}$, $o_3 = \{v_3,v_5\}$, $o_4 = \{v_4,v_6\}$.

Для определения автоморфизмов графа, нужно перебрать все возможные сочетания орбит, количество которых можно определить по формуле (1.1). Для этого графа возможна только одна нетривиальная перестановка вершин. Так как все орбиты имеют разные веса и в каждой орбите по два элемента, то нетривиальную перестановку можно определить, выделив DFS-дерево алгоритмом поиска в глубину.

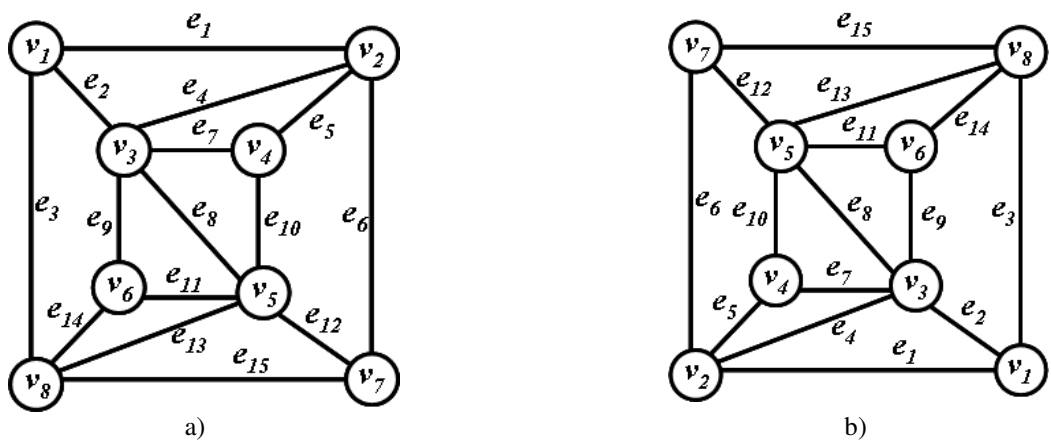

Рис. 1.14. Автоморфизм графа $G_5$.

Построим DFS-дерево :

$(v_1,v_2)+(v_2,v_3)+(v_3,v_4)+(v_4,v_5)+(v_5,v_6)+(v_6,v_8)+(v_8,v_7)$;

Поменяем местами вершины в орбитах:

$(v_7,v_8)+(v_8,v_5)+(v_5,v_6)+(v_6,v_3)+(v_3,v_4)+(v_4,v_2)+(v_2,v_1)$.

В результате получаем нетривиальную подстановку $p_1$:

$$p_1 = \begin{pmatrix} 1 & 2 & 3 & 4 & 5 & 6 & 7 & 8 \\ 7 & 8 & 5 & 6 & 3 & 4 & 1 & 2 \end{pmatrix}.$$

Результат перестановки представлен на рис. 1.14,b. Автоморфизм графа $G_5$ состоит из двух перестановок $p_0$ и $p_1$.

$$Aut(G_5) = \begin{cases} p_0 = (1)(2)(3)(4)(5)(6)(7)(8)(9)(10); \\ p_1 = (1\ 7)(2\ 8(3\ 5(4\ 6). \end{cases}$$

***Пример 1.6***. Рассмотрим полный граф $K_9$ представленный на рис. 1.16.

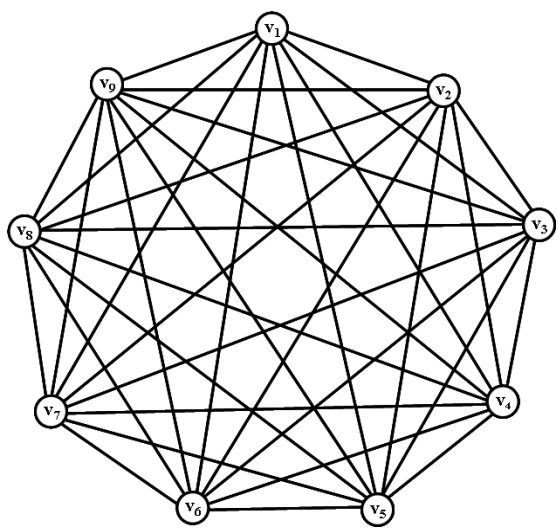

Рис. 1.15. Полный граф $K_9$.

Группой автоморфизма полного графа $K_n$ является симметрическая группа $S_n$. Число элементов симметрической группы для конечного множества равно числу перестановок элементов, то есть [факториалу](link) мощности $|S_n| = n!$.

*Пример 1.7.* Будем рассматривать рисунок полного двудольного графа $K_{4,4}$.

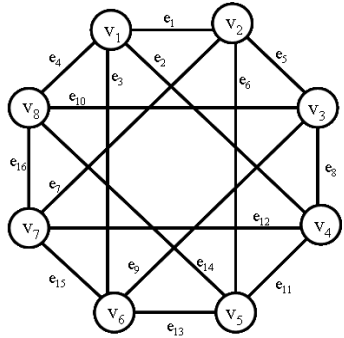

Рис. 1.16. Граф $G_{4,4}$.

Кортеж весов рёбер графа $\xi_w(K_{4,4}) = (16 \times 6)$. Кортеж весов вершин графа $\zeta_w(K_{4,4}) = \{24,24,24,24,24,24,24,24\}$.

Выделим множество изометрических циклов:

$c_1 = \{e_1,e_2,e_5,e_8\} = \{v_1,v_2,v_3,v_4\}$;          $c_2 = \{e_1,e_2,e_6,e_{11}\} = \{v_1,v_2,v_4,v_5\}$;
$c_3 = \{e_1,e_2,e_7,e_{12}\} = \{v_1,v_2,v_4,v_7\}$;       $c_4 = \{e_1,e_3,e_5,e_9\} = \{v_1,v_2,v_3,v_6\}$;
$c_5 = \{e_1,e_3,e_6,e_{13}\} = \{v_1,v_2,v_5,v_6\}$;       $c_6 = \{e_1,e_3,e_7,e_{15}\} = \{v_1,v_2,v_6,v_7\}$;
$c_7 = \{e_1,e_4,e_5,e_{10}\} = \{v_1,v_2,v_3,v_8\}$;       $c_8 = \{e_1,e_4,e_6,e_{14}\} = \{v_1,v_2,v_5,v_8\}$;
$c_9 = \{e_1,e_4,e_7,e_{16}\} = \{v_1,v_2,v_7,v_8\}$;       $c_{10} = \{e_2,e_3,e_8,e_9\} = \{v_1,v_3,v_4,v_6\}$;
$c_{11} = \{e_2,e_3,e_{11},e_{13}\} = \{v_1,v_4,v_5,v_6\}$; $c_{12} = \{e_2,e_3,e_{12},e_{15}\} = \{v_1,v_4,v_6,v_7\}$;
$c_{13} = \{e_2,e_4,e_8,e_{10}\} = \{v_1,v_3,v_4,v_8\}$;    $c_{14} = \{e_2,e_4,e_{11},e_{14}\} = \{v_1,v_4,v_5,v_8\}$;
$c_{15} = \{e_2,e_4,e_{12},e_{16}\} = \{v_1,v_4,v_7,v_8\}$; $c_{16} = \{e_3,e_4,e_9,e_{10}\} = \{v_1,v_3,v_6,v_8\}$;
$c_{17} = \{e_3,e_4,e_{13},e_{14}\} = \{v_1,v_5,v_6,v_8\}$; $c_{18} = \{e_3,e_4,e_{15},e_{16}\} = \{v_1,v_6,v_7,v_8\}$;
$c_{19} = \{e_5,e_6,e_8,e_{11}\} = \{v_2,v_3,v_4,v_5\}$;    $c_{20} = \{e_5,e_6,e_9,e_{13}\} = \{v_2,v_3,v_5,v_6\}$;
$c_{21} = \{e_5,e_6,e_{10},e_{14}\} = \{v_2,v_3,v_5,v_8\}$; $c_{22} = \{e_5,e_7,e_8,e_{12}\} = \{v_2,v_3,v_4,v_7\}$;
$c_{23} = \{e_5,e_7,e_9,e_{15}\} = \{v_2,v_3,v_6,v_7\}$;    $c_{24} = \{e_5,e_7,e_{10},e_{16}\} = \{v_2,v_3,v_7,v_8\}$;
$c_{25} = \{e_6,e_7,e_{11},e_{12}\} = \{v_2,v_4,v_5,v_7\}$; $c_{26} = \{e_6,e_7,e_{13},e_{15}\} = \{v_2,v_5,v_6,v_7\}$;
$c_{27} = \{e_6,e_7,e_{14},e_{16}\} = \{v_2,v_5,v_7,v_8\}$; $c_{28} = \{e_8,e_9,e_{11},e_{13}\} = \{v_3,v_4,v_5,v_6\}$;
$c_{29} = \{e_8,e_9,e_{12},e_{15}\} = \{v_3,v_4,v_6,v_7\}$; $c_{30} = \{e_8,e_{10},e_{11},e_{14}\} = \{v_3,v_4,v_5,v_8\}$;
$c_{31} = \{e_8,e_{10},e_{12},e_{16}\} = \{v_3,v_4,v_7,v_8\}$; $c_{32} = \{e_9,e_{10},e_{13},e_{14}\} = \{v_3,v_5,v_6,v_8\}$;
$c_{33} = \{e_9,e_{10},e_{15},e_{16}\} = \{v_3,v_6,v_7,v_8\}$; $c_{34} = \{e_{11},e_{12},e_{13},e_{15}\} = \{v_4,v_5,v_6,v_7\}$;
$c_{35} = \{e_{11},e_{12},e_{14},e_{16}\} = \{v_4,v_5,v_7,v_8\}$; $c_{36} = \{e_{13},e_{14},e_{15},e_{16}\} = \{v_5,v_6,v_7,v_8\}$.

Каждое объединение следующих циклов порождает (накрывает) всё множество вершин графа:

$c_1 \cup c_{36} = \{v_1,v_2,v_3,v_4\} \cup \{v_5,v_6,v_7,v_8\} = \{v_1,v_2,v_3,v_4,v_5,v_6,v_7,v_8\}$;
$c_2 \cup c_{33} = \{v_1,v_2,v_4,v_5\} \cup \{v_3,v_6,v_7,v_8\} = \{v_1,v_2,v_3,v_4,v_5,v_6,v_7,v_8\}$;
$c_3 \cup c_{32} = \{v_1,v_2,v_4,v_7\} \cup \{v_3,v_5,v_6,v_8\} = \{v_1,v_2,v_3,v_4,v_5,v_6,v_7,v_8\}$;
$c_4 \cup c_{35} = \{v_1,v_2,v_3,v_6\} \cup \{v_4,v_5,v_7,v_8\} = \{v_1,v_2,v_3,v_4,v_5,v_6,v_7,v_8\}$;
$c_5 \cup c_{31} = \{v_1,v_2,v_5,v_6\} \cup \{v_3,v_4,v_7,v_8\} = \{v_1,v_2,v_3,v_4,v_5,v_6,v_7,v_8\}$;
$c_6 \cup c_{30} = \{v_1,v_2,v_6,v_7\} \cup \{v_3,v_4,v_5,v_8\} = \{v_1,v_2,v_3,v_4,v_5,v_6,v_7,v_8\}$;
$c_7 \cup c_{34} = \{v_1,v_2,v_3,v_8\} \cup \{v_4,v_5,v_6,v_7\} = \{v_1,v_2,v_3,v_4,v_5,v_6,v_7,v_8\}$;
$c_8 \cup c_{29} = \{v_1,v_2,v_5,v_8\} \cup \{v_3,v_4,v_6,v_7\} = \{v_1,v_2,v_3,v_4,v_5,v_6,v_7,v_8\}$;
$c_9 \cup c_{28} = \{v_1,v_2,v_7,v_8\} \cup \{v_3,v_4,v_5,v_6\} = \{v_1,v_2,v_3,v_4,v_5,v_6,v_7,v_8\}$;
$c_{10} \cup c_{27} = \{v_1,v_3,v_4,v_6\} \cup \{v_2,v_5,v_7,v_8\} = \{v_1,v_2,v_3,v_4,v_5,v_6,v_7,v_8\}$;

$c_{11} \cup c_{24} = \{v_1,v_4,v_5,v_6\} \cup \{v_2,v_3,v_7,v_8\} = \{v_1,v_2,v_3,v_4,v_5,v_6,v_7,v_8\}$;
$c_{12} \cup c_{21} = \{v_1,v_4,v_6,v_7\} \cup \{v_2,v_3,v_5,v_8\} = \{v_1,v_2,v_3,v_4,v_5,v_6,v_7,v_8\}$;
$c_{13} \cup c_{26} = \{v_1,v_3,v_4,v_8\} \cup \{v_2,v_5,v_6,v_7\} = \{v_1,v_2,v_3,v_4,v_5,v_6,v_7,v_8\}$;
$c_{14} \cup c_{23} = \{v_1,v_4,v_5,v_8\} \cup \{v_2,v_3,v_6,v_7\} = \{v_1,v_2,v_3,v_4,v_5,v_6,v_7,v_8\}$;
$c_{15} \cup c_{20} = \{v_1,v_4,v_7,v_8\} \cup \{v_2,v_3,v_5,v_6\} = \{v_1,v_2,v_3,v_4,v_5,v_6,v_7,v_8\}$;
$c_{16} \cup c_{25} = \{v_1,v_3,v_6,v_8\} \cup \{v_2,v_4,v_5,v_7\} = \{v_1,v_2,v_3,v_4,v_5,v_6,v_7,v_8\}$;
$c_{17} \cup c_{22} = \{v_1,v_5,v_6,v_8\} \cup \{v_2,v_3,v_4,v_7\} = \{v_1,v_2,v_3,v_4,v_5,v_6,v_7,v_8\}$;
$c_{18} \cup c_{19} = \{v_1,v_6,v_7,v_8\} \cup \{v_2,v_3,v_4,v_5\} = \{v_1,v_2,v_3,v_4,v_5,v_6,v_7,v_8\}$;

Характерной особенностью двудольных графов является существование только изометрических циклов четной длины. Известно, что вершины двудольных графов раскрашиваются в два цвета. Поэтому схему перестановок в двудольном графе можно представить в виде:

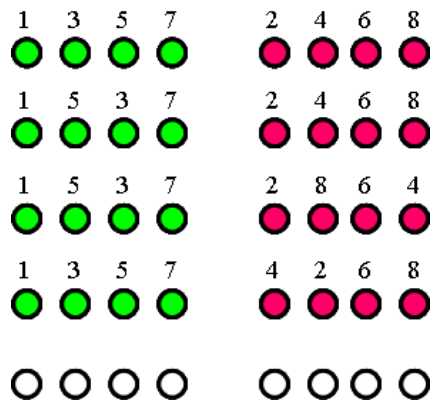

Рис. 1.17. Перестановки вершин в двудольном графе.

Вершины двудольного графа можно разбить на два несвязанных подмножества вершин называемых частями (долями). Вершины одной части соединяются с вершинами другой части, но не соединяются между собой. Поэтому перестановки вершин возможно только в одной части графа.

Рисунок полного двудольного графа $K_{4,4}$ представлен на рис. 1.17. Множество всех перестановок в графе определяются как перестановки всех несвязанных вершин в одной части на перестановку всех несвязанных вершин внутри второй части $|K_{n,n}| = \frac{n}{2}! \times \frac{n}{2}!$. В нашем случае $4! \times 4! = 24 \times 24 = 576$.

**Комментарии**

Применение векторных инвариантов графа в задаче определения группы автоморфизмов графа, позволяет выявить орбиты графа и определить их количество. Подмножества вершин, имеющие равный вес определяют состав и структуру каждой орбиты. В свою очередь состав орбиты может быть разбит на несколько транспозиций вершин с одинаковыми весами в зависимости от структуры графа.

В целях экономии вычислительного процесса, для определения симметричной устойчивости весов вершин, следует рассматривать кортеж весов вершин для спектра реберных разрезов графа. Так как вычислительная сложность построения кортежа весов вершин для спектра реберных разрезов меньше, чем вычислительная сложность других кортежей весов вершин.

Представлен алгоритм для определения изоморфизма матриц смежностей графа в зависимости от перестановки транспозиций, основанный на методах алгебры структурных чисел.

В качестве иллюстраций мы применяли различные выды графов. Однако, не следует забывать, что понятие векторного инварианта относится только к несепарабельным графам.

Приведены примеры вычисления групп автоморфизмов для графов, имеющих различные веса вершин.

## Глава 2. Геометрические и топологические методы построения перестановок
### 2.1. Геометрическое представление группы перестановок

Будем рассматривать преобразования обычного трехмерного (евклидова) пространства R.

**Определение 2.1.** Движением называется такое преобразование вещественного евклидова пространства R, при котором расстояния между точками не меняются: если точка P переходит в P′ а точка Q — в Q', то расстояние P'Q' равно PQ.

Такие преобразования образуют, очевидно, группу, называемую группой движений евклидова пространства, или евклидовой группой [5].

В группе движений пространства R выделим множество тех движений, при которых некоторая фиксированная точка O (начало координат) остается неподвижной, т. е. переходит сама в себя. Такие движения тоже, очевидно, образуют группу—подгруппу группы движений, называемую центроевклидовой, или полной ортогональной группой.

Каждому движению с неподвижной точкой O отвечает определенное преобразование соответствующего векторного пространства: если точка P переходит в P', то вектор $\overline{OP}$ переходит в вектор $\overline{OP'}$. При этом преобразовании длины векторов не меняются; легко видеть, что не меняются также и углы -между векторами, т. е. рассматриваемое преобразование сохраняет скалярное произведение векторов.

Рассмотрим группу вращений правильного тетраэдра ABCD. Он переходит в себя при следующих нетождественных поворотах:

а) При поворотах вокруг каждой из осей типа AP (рис. 2.1,а), соединяющих вершину тетраэдра с центром противолежащей грани, на углы $\frac{2\pi}{3}$ и $\frac{4\pi}{3}$. Всего таких вращений имеется $4\times 2=8$.

б) При поворотах на угол $\pi$ вокруг каждой из трех прямых типа MN (см. рис. 2.1,b), соединяющих середины противоположных ребер. Так как MN перпендикулярна BC, MN перпендикулярна AD, BM = MC и AN = ND, то при повороте вокруг прямой MN на угол $\pi$ точка B перейдет в C, C — в B, A — в D и D — в A. Всего, вместе с тождественным поворотом, мы имеем $1 + 8 + 3 = 12$ поворотов, при которых тетраэдр переходит в себя. Для удобства записи, мы иногда будем применять цифровое обозначение вершин, буквенное обозначение (ABCD) заменим цифровым обозначением (1 2 3 4). Этим 12-ти поворотам отвечают такие перестановки вершин:

$$\delta_1 = \begin{pmatrix} 1 & 2 & 3 & 4 \\ 1 & 2 & 3 & 4 \end{pmatrix}, \delta_2 = \begin{pmatrix} 1 & 2 & 3 & 4 \\ 1 & 3 & 4 & 2 \end{pmatrix}, \delta_3 = \begin{pmatrix} 1 & 2 & 3 & 4 \\ 1 & 4 & 2 & 3 \end{pmatrix}, \delta_4 = \begin{pmatrix} 1 & 2 & 3 & 4 \\ 3 & 2 & 4 & 1 \end{pmatrix},$$

$$\delta_5 = \begin{pmatrix} 1 & 2 & 3 & 4 \\ 4 & 2 & 1 & 3 \end{pmatrix}, \delta_6 = \begin{pmatrix} 1 & 2 & 3 & 4 \\ 2 & 4 & 3 & 1 \end{pmatrix}, \delta_7 = \begin{pmatrix} 1 & 2 & 3 & 4 \\ 4 & 1 & 3 & 2 \end{pmatrix}, \delta_8 = \begin{pmatrix} 1 & 2 & 3 & 4 \\ 2 & 3 & 1 & 4 \end{pmatrix},$$

$$\delta_9 = \begin{pmatrix} 1 & 2 & 3 & 4 \\ 3 & 1 & 2 & 4 \end{pmatrix}, \delta_{10} = \begin{pmatrix} 1 & 2 & 3 & 4 \\ 2 & 1 & 4 & 3 \end{pmatrix}, \delta_{11} = \begin{pmatrix} 1 & 2 & 3 & 4 \\ 3 & 4 & 1 & 2 \end{pmatrix}, \delta_{12} = \begin{pmatrix} 1 & 2 & 3 & 4 \\ 4 & 3 & 2 & 1 \end{pmatrix}.$$

Все эти подстановки — четные, и значит, соответствующие им повороты действительно образуют группу T, изоморфную знакопеременной подгруппе $A_4$ симметрической группы $S_4$.

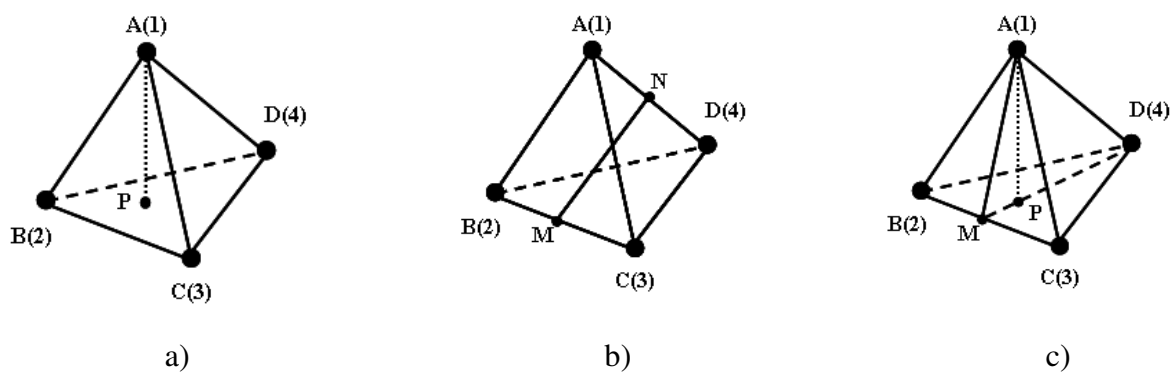

Рис. 2.1. Тетраэдр.

Рассмотрим группу вращений куба O. Легко видеть, что куб переходит в себя при следующих нетождественных вращениях:

а) При трех поворотах на углы $\pi/2$, $\pi$ и $3\pi/2$ вокруг каждой из трех прямых типа MN (рис. 2.2,a), соединяющих центры противоположных граней. Всего таких поворотов $3 \times 3 = 9$.

б) При двух поворотах вокруг каждой из четырех диагоналей на углы $2\pi/3$ и $4\pi/3$ (правильный треугольник ACD' при этом переходит в себя). Всего таких поворотов $2 \times 4 = 8$ (см. рис. 2.2,b) в) При шести поворотах на угол $\pi$ — вокруг каждой из прямых типа PQ (рис. 2.2,c), соединяющих середины противоположных ребер.

Всего, вместе с тождественным преобразованием, мы нашли $1+9 + 8 + 6 = 24$ поворота при которых куб переходит в себя.

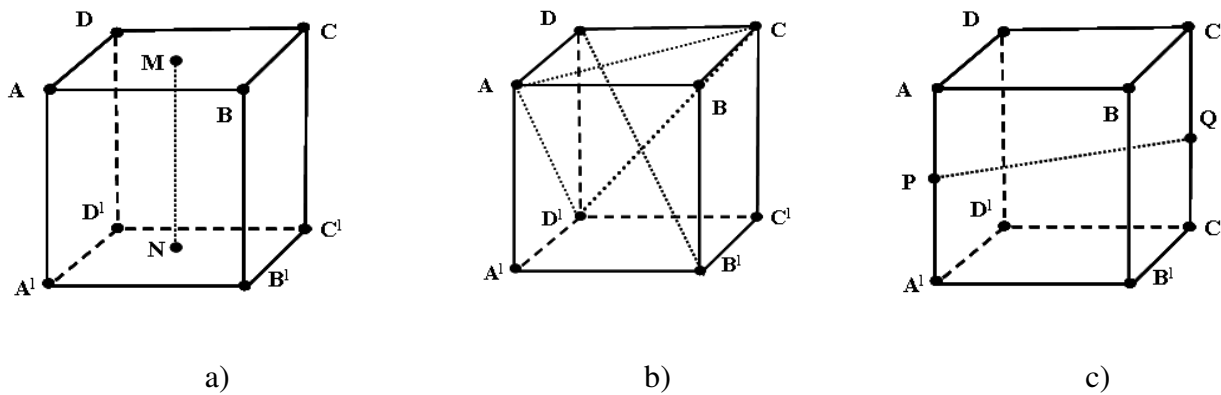

Рис. 2.3.. Куб (гексаэдр)

Для удобства записи перестановок в кубе, мы иногда будем применять цифровое обозначение вершин вместо буквенного (ABCD A′B′C′D′) → (1 2 3 4 5 6 7 8). Этим 24-м поворотам отвечают такие перестановки вершин:

для перестановок вершин согласно рис. 2,3,a

$$\xi_1 = \begin{pmatrix} 1 & 2 & 3 & 4 & 5 & 6 & 7 & 8 \\ 1 & 2 & 3 & 4 & 5 & 6 & 7 & 8 \end{pmatrix}, \xi_2 = \begin{pmatrix} 1 & 2 & 3 & 4 & 5 & 6 & 7 & 8 \\ 4 & 1 & 2 & 3 & 8 & 5 & 6 & 7 \end{pmatrix},$$

$$\xi_3 = \begin{pmatrix} 1 & 2 & 3 & 4 & 5 & 6 & 7 & 8 \\ 3 & 4 & 1 & 2 & 7 & 8 & 5 & 6 \end{pmatrix}, \xi_4 = \begin{pmatrix} 1 & 2 & 3 & 4 & 5 & 6 & 7 & 8 \\ 2 & 3 & 4 & 1 & 6 & 7 & 8 & 5 \end{pmatrix},$$

$$\xi_5 = \begin{pmatrix} 1 & 2 & 3 & 4 & 5 & 6 & 7 & 8 \\ 2 & 6 & 7 & 3 & 1 & 5 & 8 & 4 \end{pmatrix}, \xi_6 = \begin{pmatrix} 1 & 2 & 3 & 4 & 5 & 6 & 7 & 8 \\ 6 & 5 & 8 & 7 & 2 & 1 & 4 & 3 \end{pmatrix},$$

$$\xi_7 = \begin{pmatrix} 1 & 2 & 3 & 4 & 5 & 6 & 7 & 8 \\ 5 & 1 & 4 & 8 & 6 & 2 & 3 & 7 \end{pmatrix}, \xi_8 = \begin{pmatrix} 1 & 2 & 3 & 4 & 5 & 6 & 7 & 8 \\ 4 & 3 & 7 & 8 & 1 & 2 & 6 & 5 \end{pmatrix},$$

$$\xi_9 = \begin{pmatrix} 1 & 2 & 3 & 4 & 5 & 6 & 7 & 8 \\ 8 & 7 & 6 & 5 & 4 & 3 & 2 & 1 \end{pmatrix}, \xi_{10} = \begin{pmatrix} 1 & 2 & 3 & 4 & 5 & 6 & 7 & 8 \\ 5 & 6 & 2 & 1 & 8 & 7 & 3 & 4 \end{pmatrix},$$

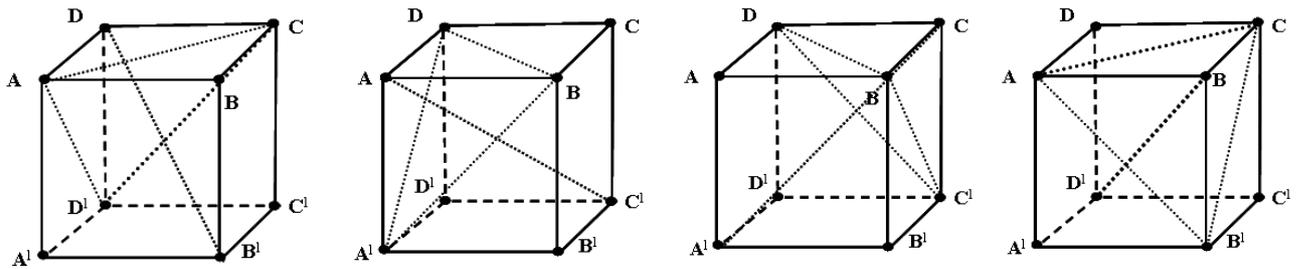

Рис. 2.4. Вращение вокруг диагоналей.

для перестановок вершин согласно рис. 2.3,b

$$\xi_{11} = \begin{pmatrix} 1 & 2 & 3 & 4 & 5 & 6 & 7 & 8 \\ 3 & 7 & 8 & 4 & 2 & 6 & 5 & 1 \end{pmatrix}, \xi_{12} = \begin{pmatrix} 1 & 2 & 3 & 4 & 5 & 6 & 7 & 8 \\ 8 & 5 & 1 & 4 & 7 & 6 & 2 & 3 \end{pmatrix},$$

$$\xi_{13} = \begin{pmatrix} 1 & 2 & 3 & 4 & 5 & 6 & 7 & 8 \\ 1 & 5 & 6 & 2 & 4 & 8 & 7 & 3 \end{pmatrix}, \xi_{14} = \begin{pmatrix} 1 & 2 & 3 & 4 & 5 & 6 & 7 & 8 \\ 1 & 4 & 8 & 5 & 2 & 3 & 7 & 6 \end{pmatrix},$$

$$\xi_{15} = \begin{pmatrix} 1 & 2 & 3 & 4 & 5 & 6 & 7 & 8 \\ 8 & 4 & 3 & 7 & 5 & 1 & 2 & 6 \end{pmatrix}, \xi_{16} = \begin{pmatrix} 1 & 2 & 3 & 4 & 5 & 6 & 7 & 8 \\ 6 & 7 & 3 & 2 & 5 & 8 & 4 & 1 \end{pmatrix},$$

$$\xi_{17} = \begin{pmatrix} 1 & 2 & 3 & 4 & 5 & 6 & 7 & 8 \\ 6 & 2 & 1 & 5 & 7 & 3 & 4 & 8 \end{pmatrix}, \xi_{18} = \begin{pmatrix} 1 & 2 & 3 & 4 & 5 & 6 & 7 & 8 \\ 3 & 2 & 6 & 7 & 4 & 1 & 5 & 8 \end{pmatrix},$$

для перестановок вершин согласно рис. 2.3,c

$$\xi_{19} = \begin{pmatrix} 1 & 2 & 3 & 4 & 5 & 6 & 7 & 8 \\ 7 & 6 & 5 & 8 & 3 & 2 & 1 & 4 \end{pmatrix}, \xi_{20} = \begin{pmatrix} 1 & 2 & 3 & 4 & 5 & 6 & 7 & 8 \\ 5 & 8 & 7 & 6 & 1 & 4 & 3 & 2 \end{pmatrix},$$

$$\xi_{21} = \begin{pmatrix} 1 & 2 & 3 & 4 & 5 & 6 & 7 & 8 \\ 7 & 8 & 4 & 3 & 6 & 5 & 1 & 2 \end{pmatrix}, \xi_{22} = \begin{pmatrix} 1 & 2 & 3 & 4 & 5 & 6 & 7 & 8 \\ 2 & 1 & 5 & 6 & 3 & 4 & 8 & 7 \end{pmatrix},$$

$$\xi_{23} = \begin{pmatrix} 1 & 2 & 3 & 4 & 5 & 6 & 7 & 8 \\ 4 & 8 & 5 & 1 & 3 & 7 & 6 & 2 \end{pmatrix}, \xi_{24} = \begin{pmatrix} 1 & 2 & 3 & 4 & 5 & 6 & 7 & 8 \\ 7 & 3 & 2 & 6 & 8 & 4 & 1 & 6 \end{pmatrix}.$$

Рассмотрим группу симметрий правильного тетраэдра $T_a$. Правильный тетраэдр имеет шесть плоскостей симметрии. К 12 вращениям, при которых тетраэдр переходит в себя (и которые отвечают, четным подстановкам его вершин), добавим одну из симметрии, например, симметрию относительно плоскости ADM '(рис. 2.1,c) — ей соответствует нечетная перестановка $\begin{pmatrix} 1 & 2 & 3 & 4 \\ 1 & 3 & 2 & 4 \end{pmatrix}$ вершин тетраэдра. Если умножить эту симметрию на каждый из 12 поворотов, при которых тетраэдр переходит в себя, мы получим еще 12 преобразований, отвечающих нечетным подстановкам вершин.

$$\delta_{13} = \begin{pmatrix} 1 & 2 & 3 & 4 \\ 1 & 3 & 2 & 4 \end{pmatrix}, \delta_{14} = \begin{pmatrix} 1 & 2 & 3 & 4 \\ 1 & 4 & 3 & 2 \end{pmatrix}, \delta_{15} = \begin{pmatrix} 1 & 2 & 3 & 4 \\ 1 & 2 & 4 & 3 \end{pmatrix}, \delta_{16} = \begin{pmatrix} 1 & 2 & 3 & 4 \\ 3 & 4 & 2 & 1 \end{pmatrix},$$

$$\delta_{17} = \begin{pmatrix} 1 & 2 & 3 & 4 \\ 4 & 1 & 2 & 3 \end{pmatrix}, \delta_{18} = \begin{pmatrix} 1 & 2 & 3 & 4 \\ 2 & 3 & 4 & 1 \end{pmatrix}, \delta_{19} = \begin{pmatrix} 1 & 2 & 3 & 4 \\ 4 & 3 & 1 & 2 \end{pmatrix}, \delta_{20} = \begin{pmatrix} 1 & 2 & 3 & 4 \\ 2 & 1 & 3 & 4 \end{pmatrix},$$

$$\delta_{21} = \begin{pmatrix} 1 & 2 & 3 & 4 \\ 3 & 2 & 1 & 4 \end{pmatrix}, \delta_{22} = \begin{pmatrix} 1 & 2 & 3 & 4 \\ 2 & 4 & 1 & 3 \end{pmatrix}, \delta_{23} = \begin{pmatrix} 1 & 2 & 3 & 4 \\ 3 & 1 & 4 & 2 \end{pmatrix}, \delta_{24} = \begin{pmatrix} 1 & 2 & 3 & 4 \\ 4 & 2 & 3 & 1 \end{pmatrix}.$$

Будем рассматривать группу симметрии куба $O_h$. Куб имеет 9 плоскостей симметрии (и центр симметрии): три плоскости симметрии — такие, как PQMN на рис. 2.5,а, и шесть диагональных плоскостей — таких, как DBB'D' на рис. 2.5,б.

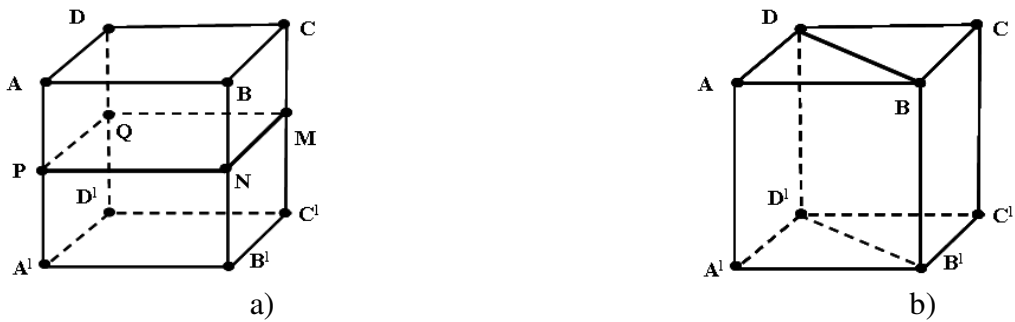

а)          b)

Рис. 2.5. Симметрия куба.

Осуществим симметричное преобразование относительно плоскости PQMN

$$\xi = \begin{pmatrix} 1 & 2 & 3 & 4 & 5 & 6 & 7 & 8 \\ 5 & 6 & 7 & 8 & 1 & 2 & 3 & 4 \end{pmatrix}.$$

Умножим полученное преобразование на 24 вращений куба:

$$\xi_{25} = \begin{pmatrix} 1 & 2 & 3 & 4 & 5 & 6 & 7 & 8 \\ 5 & 6 & 7 & 8 & 1 & 2 & 3 & 4 \end{pmatrix}, \xi_{26} = \begin{pmatrix} 1 & 2 & 3 & 4 & 5 & 6 & 7 & 8 \\ 8 & 5 & 6 & 7 & 4 & 1 & 2 & 3 \end{pmatrix},$$

$$\xi_{27} = \begin{pmatrix} 1 & 2 & 3 & 4 & 5 & 6 & 7 & 8 \\ 7 & 8 & 5 & 6 & 3 & 4 & 1 & 2 \end{pmatrix}, \xi_{28} = \begin{pmatrix} 1 & 2 & 3 & 4 & 5 & 6 & 7 & 8 \\ 6 & 7 & 8 & 5 & 2 & 3 & 4 & 1 \end{pmatrix},$$

$$\xi_{29} = \begin{pmatrix} 1 & 2 & 3 & 4 & 5 & 6 & 7 & 8 \\ 1 & 5 & 8 & 4 & 2 & 6 & 7 & 3 \end{pmatrix}, \xi_{30} = \begin{pmatrix} 1 & 2 & 3 & 4 & 5 & 6 & 7 & 8 \\ 2 & 1 & 4 & 3 & 6 & 5 & 8 & 7 \end{pmatrix},$$

$$\xi_{31} = \begin{pmatrix} 1 & 2 & 3 & 4 & 5 & 6 & 7 & 8 \\ 6 & 2 & 3 & 7 & 5 & 1 & 4 & 8 \end{pmatrix}, \xi_{32} = \begin{pmatrix} 1 & 2 & 3 & 4 & 5 & 6 & 7 & 8 \\ 1 & 2 & 6 & 5 & 4 & 3 & 7 & 8 \end{pmatrix},$$

$$\xi_{33} = \begin{pmatrix} 1 & 2 & 3 & 4 & 5 & 6 & 7 & 8 \\ 4 & 3 & 2 & 1 & 8 & 7 & 6 & 5 \end{pmatrix}, \xi_{34} = \begin{pmatrix} 1 & 2 & 3 & 4 & 5 & 6 & 7 & 8 \\ 8 & 7 & 3 & 4 & 5 & 6 & 2 & 1 \end{pmatrix},$$

$$\xi_{35} = \begin{pmatrix} 1 & 2 & 3 & 4 & 5 & 6 & 7 & 8 \\ 2 & 6 & 5 & 1 & 3 & 7 & 8 & 4 \end{pmatrix}, \xi_{36} = \begin{pmatrix} 1 & 2 & 3 & 4 & 5 & 6 & 7 & 8 \\ 7 & 6 & 2 & 3 & 8 & 5 & 1 & 4 \end{pmatrix},$$

$$\xi_{37} = \begin{pmatrix} 1 & 2 & 3 & 4 & 5 & 6 & 7 & 8 \\ 4 & 8 & 7 & 3 & 1 & 5 & 6 & 2 \end{pmatrix}, \xi_{38} = \begin{pmatrix} 1 & 2 & 3 & 4 & 5 & 6 & 7 & 8 \\ 2 & 3 & 7 & 6 & 1 & 4 & 8 & 5 \end{pmatrix},$$

$$\xi_{39} = \begin{pmatrix} 1 & 2 & 3 & 4 & 5 & 6 & 7 & 8 \\ 5 & 1 & 2 & 6 & 8 & 4 & 3 & 7 \end{pmatrix}, \xi_{40} = \begin{pmatrix} 1 & 2 & 3 & 4 & 5 & 6 & 7 & 8 \\ 5 & 8 & 4 & 1 & 6 & 7 & 3 & 2 \end{pmatrix},$$

$$\xi_{41} = \begin{pmatrix} 1 & 2 & 3 & 4 & 5 & 6 & 7 & 8 \\ 7 & 3 & 4 & 8 & 6 & 2 & 1 & 5 \end{pmatrix}, \xi_{42} = \begin{pmatrix} 1 & 2 & 3 & 4 & 5 & 6 & 7 & 8 \\ 4 & 1 & 5 & 8 & 3 & 2 & 6 & 7 \end{pmatrix},$$

$$\xi_{43} = \begin{pmatrix} 1 & 2 & 3 & 4 & 5 & 6 & 7 & 8 \\ 3 & 2 & 1 & 4 & 7 & 6 & 5 & 8 \end{pmatrix}, \xi_{44} = \begin{pmatrix} 1 & 2 & 3 & 4 & 5 & 6 & 7 & 8 \\ 1 & 4 & 3 & 2 & 5 & 8 & 7 & 6 \end{pmatrix},$$

$$\xi_{45} = \begin{pmatrix} 1 & 2 & 3 & 4 & 5 & 6 & 7 & 8 \\ 6 & 5 & 1 & 2 & 7 & 8 & 4 & 3 \end{pmatrix}, \xi_{46} = \begin{pmatrix} 1 & 2 & 3 & 4 & 5 & 6 & 7 & 8 \\ 3 & 4 & 8 & 7 & 2 & 1 & 5 & 6 \end{pmatrix},$$

$$\xi_{47} = \begin{pmatrix} 1 & 2 & 3 & 4 & 5 & 6 & 7 & 8 \\ 3 & 7 & 6 & 2 & 4 & 8 & 5 & 1 \end{pmatrix}, \xi_{48} = \begin{pmatrix} 1 & 2 & 3 & 4 & 5 & 6 & 7 & 8 \\ 8 & 4 & 1 & 6 & 7 & 3 & 2 & 6 \end{pmatrix}.$$

Правильным многогранником называется такой (выпуклый) многогранник, все грани которого — равные между собой правильные многоугольники и все многогранные углы которого равны между собой. В древности было известно, что существует всего пять правильных многогранников: правильный тетраэдр ограниченный четырьмя правильными треугольными гранями, и имеющий 6 ребер и 4 вершины; куб (или правильный гексаэдр), ограниченный шестью квадратными гранями и имеющий 12 ребер и 8 вершин; правильный октаэдр, ограниченный восемью треугольными гранями и имеющий 12 ребер и 6 вершин (рис.2.6,a); правильный додекаэдр, ограниченный двенадцатью пятиугольными гранями и имеющий 30 ребер и 20 вершин (рис. 2.6,b); и, наконец, правильный икосаэдр, ограниченный двадцатью треугольными гранями и имеющий 30 ребер и 12 вершин (рис. 2.6,c)[5].

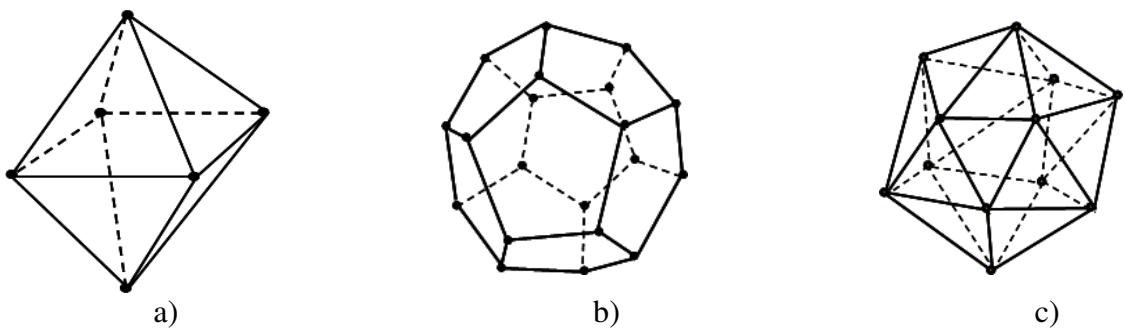

Рис. 2.6. Многогранники.

Куб и правильный октаэдр в определенном смысле двойственны друг другу: если соединить центры граней куба, как указано на рис. 2.7,,а, то получится правильный октаэдр, и, наоборот, если соединить центры граней октаэдра, то получится куб (рис. 2.7,b). Поэтому группа вращений октаэдра изоморфна группе вращений куба (и обозначается эта группа буквой О), а группа симметрии октаэдра — группе симметрии куба $O_h$.

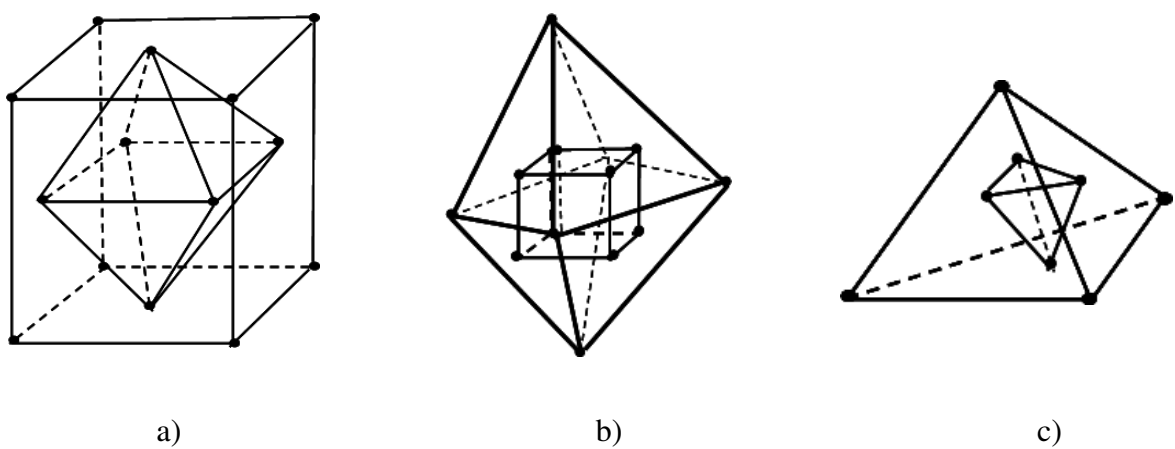

Рис. 2.7. Свойство многогранников.

В описанном смысле правильный тетраэдр двойствен сам себе (рис. 2.7,в)[5].

### 2.2. Диэдральная группа

В предыдущем разделе мы рассмотрели геометрические методы построения перестановок. Рассмотрим процесс построения перестановок топологическими методами. С этой целью, сравним результаты построения перестановок геометрическими и топологическими методами на примере пяти правильных платоновых тел.

В математике, диэдральной группой ($D_n$, группа диэдра) называется группа симметрии правильного многоугольника, включающая как вращения, так и осевые симметрии. Диэдральные группы являются простейшими примерами конечных групп и играют важную роль в теории групп, геометрии и химии [5].

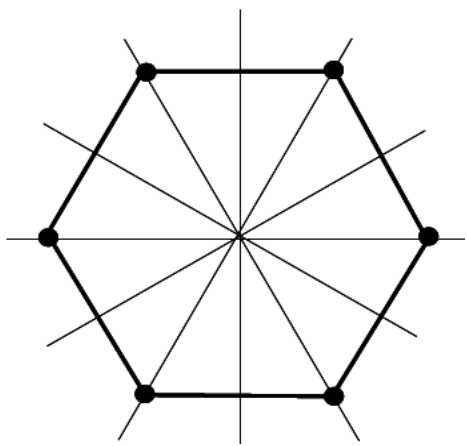

Рис. 2.8. Правильный многоугольник.

Правильный *n*-угольник (см. рис. 2.8) имеет различных симметрий: поворотов и осевых отражений, образующих диэдральную группу . Если нечётно, каждая ось симметрии проходит через середину одной из сторон и противоположную вершину. Если чётно, имеется осей симметрии, соединяющих середины противоположных сторон и осей, соединяющих противоположные вершины. В любом случае, имеется осей симметрии и элементов в группе симметрий. Отражение относительно одной оси, а затем относительно другой, приводит к вращению на удвоенный угол между осями.

Следует заметить, что многоугольнику можно поставить в соответствие простой цикл графа.

В качестве примера рассмотрим граф цикла длиной три.

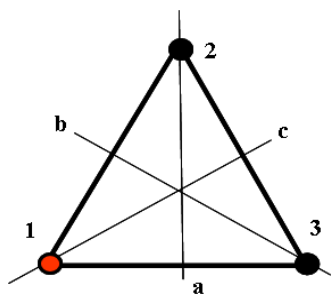 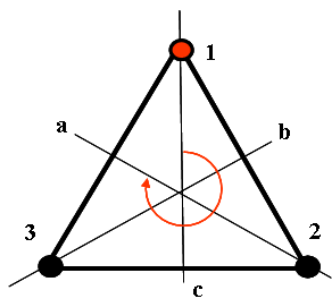 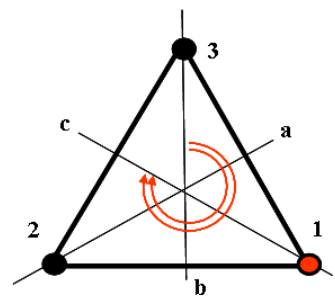

Рис. 2.9. Перестановка (1 2 3).    Рис. 2.10. Перестановка (3 1 2).    Рис. 2.11. Перестановка (2 3 1).

Классически перестановки обозначаются следующим образом: $p_1 = \begin{pmatrix} 1\ 2\ 3 \\ 3\ 1\ 2 \end{pmatrix}$.

Будем пользоваться сокращенной записью, записывая только нижнюю строчку (3 1 2).

На рис. 2.9 представлена начальная (нулевая) перестановка вершин треугольника $p_0$ = (1 2 3). На рис. 2.10 представлена перестановка вершин треугольника, определяемая вращением вокруг вертикальной оси на угол 120 градусов по часовой стрелке $p_1$ = (3 1 2). На рис. 2.11 представлена перестановка вершин треугольника, определяемая вращением вокруг вертикальной оси на угол 240 градусов по часовой стрелке $p_2$ = (2 3 1).

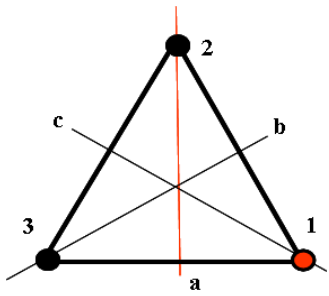 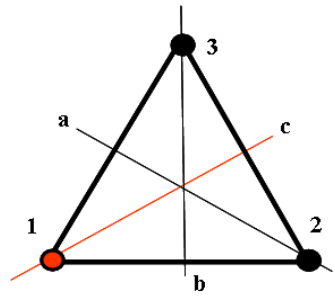 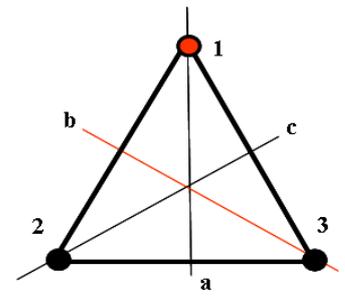

Рис. 2.12. Перестановка (3 2 1).  Рис. 2.13. Перестановка (1 3 2).  Рис. 2.14. Перестановка (2 1 3).

На рис. 2.12 представлена перестановка вершин треугольника, определяемая вращением вокруг оси *a*, $p_3$ = (3 2 1). На рис. 2.13 представлена перестановка вершин треугольника, определяемая вращением вокруг оси *c*, $p_4$ = (1 3 2). На рис. 2.14 представлена перестановка вершин треугольника, определяемая вращением вокруг оси *b*, $p_5$ = (2 1 3).

### 2.3. Граф ТЕТРАЭДРа

Рассмотрим автоморфизм графа тетраэдра К$_4$. Перед тем как рассмотреть свойства графа тетраэдра, рассмотрим многоугольник с центральной точкой (см. рис. 2.15).

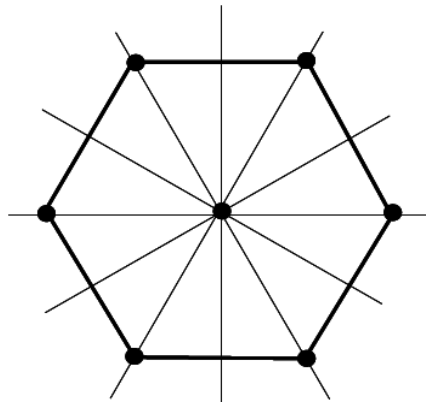

Рис. 2.15. Правильный многоугольник с центральной точкой.

**Определение 2.2**. Диэдральной группой с центральной точкой ($D_{(n-1)+1}$, группа диэдра с центральной точкой) будем называть группу симметрии правильного многоугольника с центральной точкой, включающая как вращения с центральной точкой, так и осевые симметрии с центральной точкой.

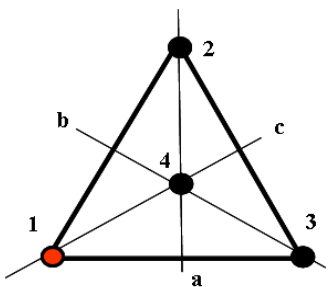 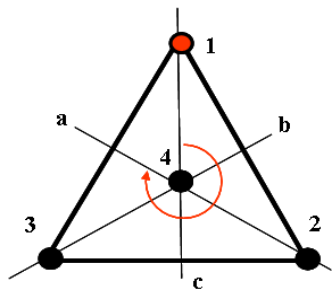 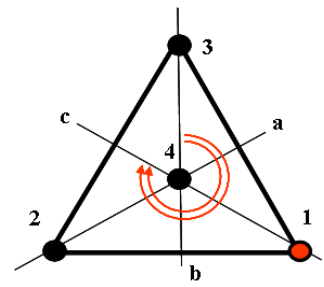

Рис. 2.16. Перестановка (2 3 1 4).  Рис. 2.17. Перестановка (2 3 1 4).  Рис. 2.18. Перестановка (2 3 1 4).

На рис. 2.16 – 2.18 представлены перестановки образованные вращением вокруг центральной точки.

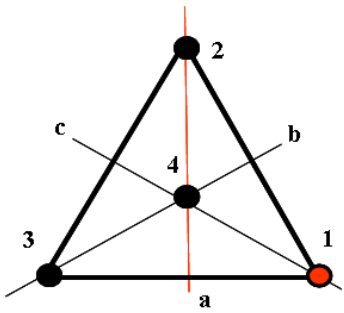 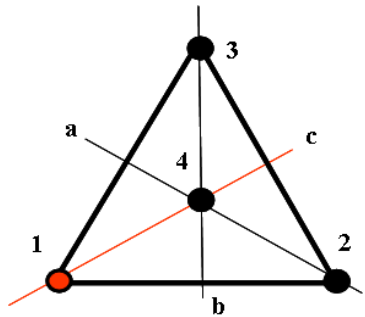 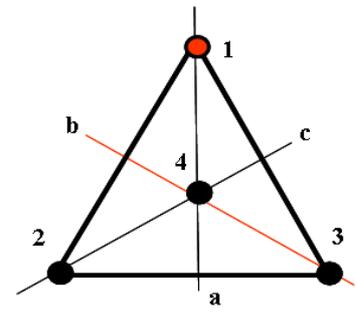

Рис. 2.19. Перестановка (2 3 1 4).   Рис. 2.20. Перестановка (2 3 1 4).   Рис. 2.21. Перестановка (2 3 1 4).

На рис. 2.19 – 2.21 представлены перестановки образованные отображением вокруг оси.

Очевидно, что количество перестановок в группе $D_{3+1}$, равно количеству перестановок в диэдральной группе $D_3$. А вот количество вершин в графе на единицу больше, чем у диэдра.

В качестве примера, рассмотрим полный граф $K_4$.

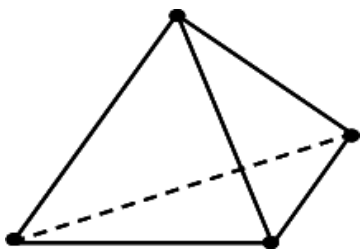 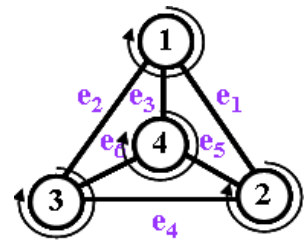

Рис. 2.22. Тетраэдр - полный граф $K_4$ и его топологический рисунок.

Выделим множество изометрических циклов в графе $K_4$.

$c_1 = \{e_1,e_2,e_5\} = \{v_1,v_2,v_4\};$          $c_2 = \{e_4,e_5,e_6\} = \{v_2,v_3,v_4\};$
$c_3 = \{e_2,e_3,e_6\} = \{v_1,v_3,v_4\};$          $c_4 = \{e_1,e_3,e_4\} = \{v_1,v_2,v_3\}.$

Очевидно, что основой для группы автоморфизма графа $K_4$ будет служить диэдральная группа с центральной неподвижной точкой.

**Определение 2.3**. Будем называть *образующим циклом* – внешний цикл, проходящий по вершинам обода многоугольника (геометрический рисунок), соответствующий кольцевому суммированию изометрических циклов графа.

В тетраэдре каждый изометрический цикл может быть выбран в качестве образующего цикла, так как он индуцирует правильный геометрический треугольник, соответствующий грани. Количество граней в тетраэдре соответствует количеству изометрических циклов графа тетраэдра. И тогда можно определить группу автоморфизмов графа как симметричные перестановки вершин относительно заданных точек двумерного пространства [5].

В качестве образующего цикла, для порождения перестановок, выберем изометрический цикл $c_4$. Перестановки, индуцируемые циклом $c_4$, представлены на рис. 2.16 – рис. 2.21.

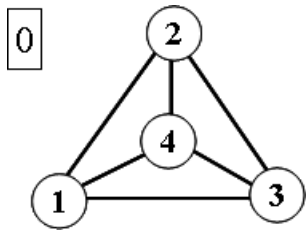 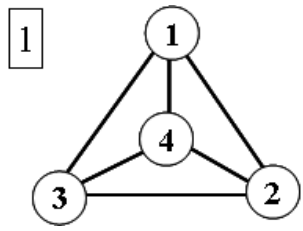 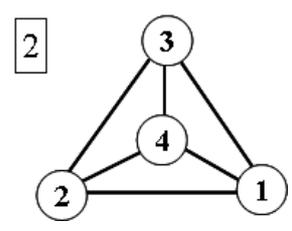

Рис. 2.23. Перестановка $p_0$.

$$\begin{pmatrix} 1 & 2 & 3 & 4 \\ 1 & 2 & 3 & 4 \end{pmatrix}$$

Рис. 2.24. Перестановка $p_1$.

$$\begin{pmatrix} 1 & 2 & 3 & 4 \\ 3 & 1 & 2 & 4 \end{pmatrix}$$

Рис. 2.25. Перестанова $p_2$.

$$\begin{pmatrix} 1 & 2 & 3 & 4 \\ 2 & 3 & 1 & 4 \end{pmatrix}$$

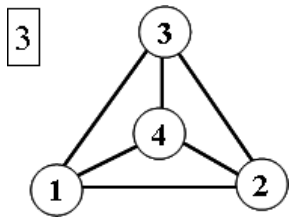 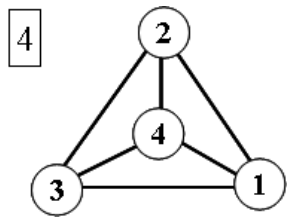 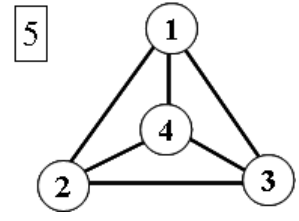

Рис. 2.26. Перестановка $p_3$.

$$\begin{pmatrix} 1 & 2 & 3 & 4 \\ 1 & 3 & 2 & 4 \end{pmatrix}$$

Рис. 2.27. Перестановка $p_4$

$$\begin{pmatrix} 1 & 2 & 3 & 4 \\ 3 & 2 & 1 & 4 \end{pmatrix}$$

Рис. 2.28. Перестановка $p_5$.

$$\begin{pmatrix} 1 & 2 & 3 & 4 \\ 2 & 1 & 3 & 4 \end{pmatrix}$$

В качестве опорного цикла выберем изометрический цикл $c_1$.

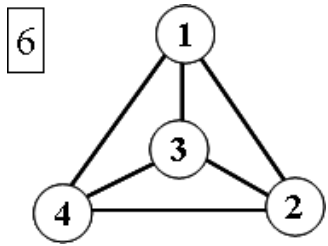 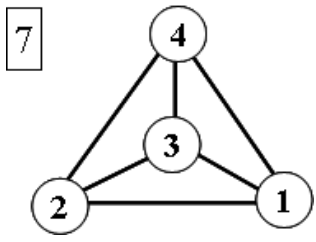 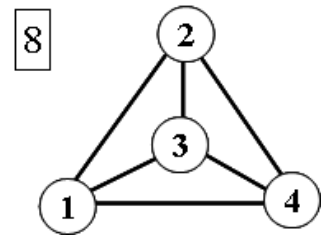

Рис. 2.29. Перестановка $p_6$.

$$\begin{pmatrix} 1 & 2 & 3 & 4 \\ 1 & 2 & 4 & 3 \end{pmatrix}$$

Рис. 2.30. Перестановка $p_7$.

$$\begin{pmatrix} 1 & 2 & 3 & 4 \\ 4 & 1 & 2 & 3 \end{pmatrix}$$

Рис. 2.31. Перестановка $p_8$.

$$\begin{pmatrix} 1 & 2 & 3 & 4 \\ 2 & 4 & 1 & 3 \end{pmatrix}$$

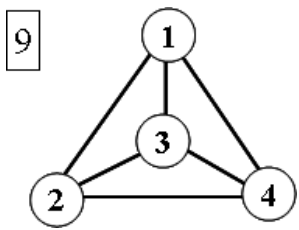 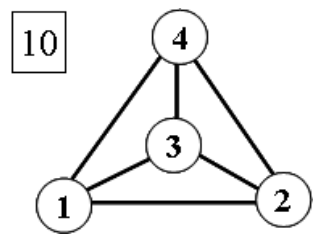 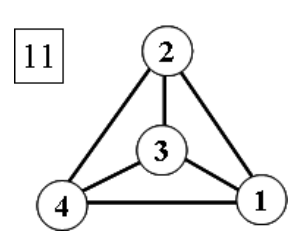

Рис. 2.32. Перестановка $p_9$.

$$\begin{pmatrix} 1 & 2 & 3 & 4 \\ 1 & 4 & 2 & 3 \end{pmatrix}$$

Рис. 2.33. Перестановка $p_{10}$.

$$\begin{pmatrix} 1 & 2 & 3 & 4 \\ 4 & 2 & 1 & 3 \end{pmatrix}$$

Рис. 2.34. Перестановка $p_{11}$.

$$\begin{pmatrix} 1 & 2 & 3 & 4 \\ 2 & 1 & 4 & 3 \end{pmatrix}$$

Перестановки, индуцируемые циклом $c_1$ представлены на рис. 2.22 – рис. 2.27.

В качестве образующего цикла выберем изометрический цикл $c_3$.

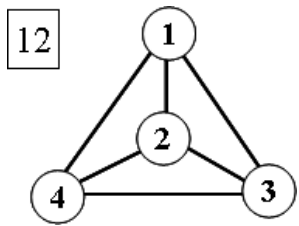 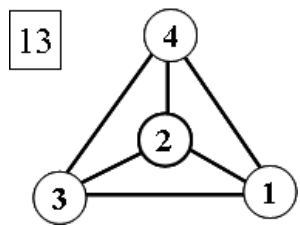 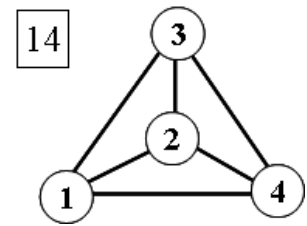

Рис. 2.35. Перестановка $p_{12}$.  Рис. 2.36. Перестановка $p_{13}$.  Рис. 2.37. Перестановка $p_{14}$.

$$\begin{pmatrix} 1 & 2 & 3 & 4 \\ 1 & 3 & 4 & 2 \end{pmatrix} \qquad \begin{pmatrix} 1 & 2 & 3 & 4 \\ 4 & 1 & 3 & 2 \end{pmatrix} \qquad \begin{pmatrix} 1 & 2 & 3 & 4 \\ 3 & 4 & 1 & 2 \end{pmatrix}$$

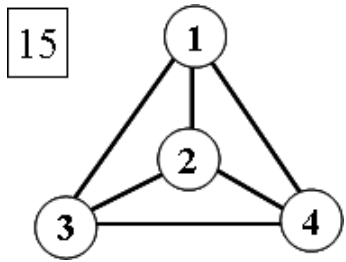 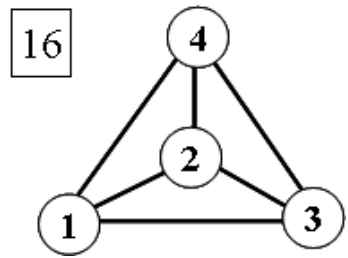 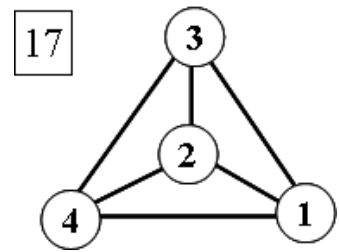

Рис. 2.38. Перестановка $p_{15}$.  Рис. 2.39. Перестановка $p_{16}$.  Рис. 2.40. Перестановка $p_{17}$.

$$\begin{pmatrix} 1 & 2 & 3 & 4 \\ 1 & 4 & 3 & 2 \end{pmatrix} \qquad \begin{pmatrix} 1 & 2 & 3 & 4 \\ 4 & 3 & 1 & 2 \end{pmatrix} \qquad \begin{pmatrix} 1 & 2 & 3 & 4 \\ 3 & 1 & 4 & 2 \end{pmatrix}$$

Образующий цикл соответствует грани тетраэдра. Количество перестановок для образующего цикла определяется как удвоенное произведение количества сторон многоугольника. Поэтому количество всех перестановок вершин для образующего цикла равно 2×3 = 6.

Перестановки, индуцируемые циклом $c_3$, представлены на рис. 2.35 – рис. 2.40.

В качестве образующего цикла выберем изометрический цикл $c_2$.

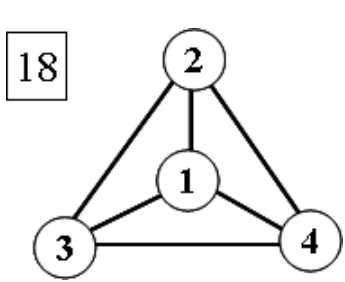 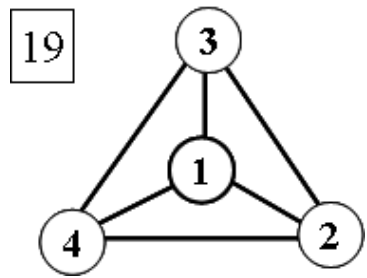 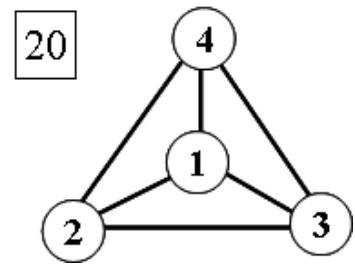

Рис. 2.41. Перестановка $p_{18}$.  Рис. 2.42. Перестановка $p_{19}$.  Рис. 2.43. Перестановка $p_{20}$.

$$\begin{pmatrix} 1 & 2 & 3 & 4 \\ 2 & 4 & 3 & 1 \end{pmatrix} \qquad \begin{pmatrix} 1 & 2 & 3 & 4 \\ 3 & 2 & 4 & 1 \end{pmatrix} \qquad \begin{pmatrix} 1 & 2 & 3 & 4 \\ 4 & 3 & 2 & 1 \end{pmatrix}$$

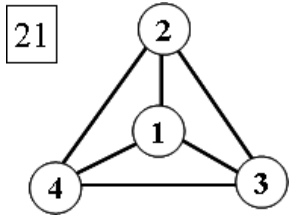 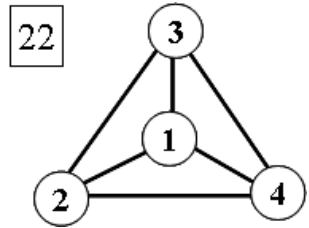 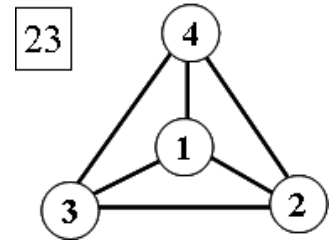

Рис. 2.44. Перестановка $p_{21}$.  Рис. 2.45. Перестановка $p_{22}$.  Рис. 2.46. Перестановка $p_{23}$.

$$\begin{pmatrix} 1 & 2 & 3 & 4 \\ 2 & 3 & 4 & 1 \end{pmatrix} \qquad \begin{pmatrix} 1 & 2 & 3 & 4 \\ 3 & 4 & 2 & 1 \end{pmatrix} \qquad \begin{pmatrix} 1 & 2 & 3 & 4 \\ 4 & 2 & 3 & 1 \end{pmatrix}$$

Перестановки, индуцируемые циклом $c_2$ представлены на рис. 2.41 – рис. 2.46.

Множество перестановок будем записывать в виде кортежа состоящего из нижней последовательности элементов в подстановке или в виде цикла перестановки:

$p_0 = <1,2,3,4> = (1)(2)(3)(4);$
$p_1 = <3,1,2,4> = (1\ 3\ 2)(4);$
$p_2 = <2,3,1,4> = (1\ 2\ 3)(4);$
$p_3 = <1,3,2,4> = (1)(2\ 3)(4);$
$p_4 = <3,2,1,4> = (1\ 3)(2)(4);$
$p_5 = <2,1,3,4> = (1\ 2)(3)(4);$
$p_6 = <1,2,4,3> = (1)(2)(3\ 4);$
$p_7 = <4,1,2,3> = (1\ 4\ 3\ 2);$
$p_8 = <2,4,1,3> = (1\ 2\ 4\ 3);$
$p_9 = <1,4,2,3> = (1)(2\ 4\ 3);$
$p_{10} = <4,2,1,3> = (1\ 4\ 3)(2);$
$p_{11} = <2,1,4,3> = (1\ 2)(3\ 4);$
$p_{12} = <1,3,4,2> = (1)(2\ 3\ 4);$
$p_{13} = <4,1,3,2> = (1\ 4\ 2)(3);$
$p_{14} = <3,4,1,2> = (1\ 3)(2\ 4);$
$p_{15} = <1,4,3,2> = (1)(2\ 4)(3);$
$p_{16} = <4,3,1,2> = (1\ 4\ 2\ 3);$
$p_{17} = <3,1,4,2> = (1\ 3\ 4\ 2);$
$p_{18} = <2,4,3,1> = (1\ 2\ 4)(3);$
$p_{19} = <3,2,4,1> = (1\ 3\ 4)(2);$
$p_{20} = <4,3,2,1> = (1\ 4)(2\ 3);$
$p_{21} = <2,3,4,1> = (1\ 2\ 3\ 4);$
$p_{22} = <3,4,2,1> = (1\ 3\ 2\ 4);$
$p_{23} = <4,2,3,1> = (1\ 4)(2)(3).$

Количество перестановок (порядок) для данного графа тетраэдра, определяется произведением количества образующих циклов, на удвоенное количество сторон многоугольника порожденного образующим циклом. В нашем случае $4\times(2\times3) = 4\times6$ и равно 24. Представим произведения элементов группы в виде таблицы Кэли состоящей из подгрупп перестановок [4], где размер квадратного блока для подгруппы равен удвоенной длине образующего цикла:

$Aut(K_4) =$

|        | $P_{11}$ | $P_{12}$ | $P_{13}$ | $P_{14}$ |
|--------|----------|----------|----------|----------|
|        | $P_{21}$ | $P_{22}$ | $P_{23}$ | $P_{24}$ |
|        | $P_{31}$ | $P_{32}$ | $P_{33}$ | $P_{34}$ |
|        | $P_{41}$ | $P_{42}$ | $P_{43}$ | $P_{44}$ |

$P_{11} =$

|       | $p_0$ | $p_1$ | $p_2$ | $p_3$ | $p_4$ | $p_5$ |
|-------|-------|-------|-------|-------|-------|-------|
| $p_0$ | $p_0$ | $p_1$ | $p_2$ | $p_3$ | $p_4$ | $p_5$ |
| $p_1$ | $p_1$ | $p_2$ | $p_0$ | $p_5$ | $p_3$ | $p_4$ |
| $p_2$ | $p_2$ | $p_0$ | $p_1$ | $p_4$ | $p_5$ | $p_3$ |
| $p_3$ | $p_3$ | $p_4$ | $p_5$ | $p_0$ | $p_1$ | $p_2$ |
| $p_4$ | $p_4$ | $p_5$ | $p_3$ | $p_2$ | $p_0$ | $p_1$ |
| $p_5$ | $p_5$ | $p_3$ | $p_4$ | $p_1$ | $p_2$ | $p_0$ |

$P_{12} =$

|       | $p_6$ | $p_7$ | $p_8$ | $p_9$ | $p_{10}$ | $p_{11}$ |
|-------|-------|-------|-------|-------|----------|----------|
| $p_0$ | $p_6$ | $p_7$ | $p_8$ | $p_9$ | $p_{10}$ | $p_{11}$ |
| $p_1$ | $p_7$ | $p_8$ | $p_6$ | $p_{11}$ | $p_9$ | $p_{10}$ |
| $p_2$ | $p_8$ | $p_6$ | $p_7$ | $p_{10}$ | $p_{11}$ | $p_9$ |
| $p_3$ | $p_9$ | $p_{10}$ | $p_{11}$ | $p_6$ | $p_7$ | $p_8$ |
| $p_4$ | $p_{10}$ | $p_{11}$ | $p_9$ | $p_8$ | $p_6$ | $p_7$ |
| $p_5$ | $p_{11}$ | $p_9$ | $p_{10}$ | $p_7$ | $p_8$ | $p_6$ |

$P_{13} =$

|       | $p_{12}$ | $p_{13}$ | $p_{14}$ | $p_{15}$ | $p_{16}$ | $p_{17}$ |
|-------|----------|----------|----------|----------|----------|----------|
| $p_0$ | $p_{12}$ | $p_{13}$ | $p_{14}$ | $p_{15}$ | $p_{16}$ | $p_{17}$ |
| $p_1$ | $p_{13}$ | $p_{14}$ | $p_{12}$ | $p_{17}$ | $p_{15}$ | $p_{16}$ |
| $p_2$ | $p_{14}$ | $p_{12}$ | $p_{13}$ | $p_{16}$ | $p_{17}$ | $p_{15}$ |
| $p_3$ | $p_{15}$ | $p_{16}$ | $p_{17}$ | $p_{12}$ | $p_{13}$ | $p_{14}$ |
| $p_4$ | $p_{16}$ | $p_{17}$ | $p_{15}$ | $p_{14}$ | $p_{12}$ | $p_{13}$ |
| $p_5$ | $p_{17}$ | $p_{15}$ | $p_{16}$ | $p_{13}$ | $p_{14}$ | $p_{12}$ |

$P_{14} =$

|       | $p_{18}$ | $p_{19}$ | $p_{20}$ | $p_{21}$ | $p_{22}$ | $p_{23}$ |
|-------|----------|----------|----------|----------|----------|----------|
| $p_0$ | $p_{18}$ | $p_{19}$ | $p_{20}$ | $p_{21}$ | $p_{22}$ | $p_{23}$ |
| $p_1$ | $p_{19}$ | $p_{20}$ | $p_{18}$ | $p_{23}$ | $p_{21}$ | $p_{22}$ |
| $p_2$ | $p_{20}$ | $p_{18}$ | $p_{19}$ | $p_{22}$ | $p_{23}$ | $p_{21}$ |
| $p_3$ | $p_{21}$ | $p_{22}$ | $p_{23}$ | $p_{18}$ | $p_{19}$ | $p_{20}$ |
| $p_4$ | $p_{22}$ | $p_{23}$ | $p_{21}$ | $p_{20}$ | $p_{18}$ | $p_{19}$ |
| $p_5$ | $p_{23}$ | $p_{21}$ | $p_{22}$ | $p_{19}$ | $p_{20}$ | $p_{18}$ |

$P_{21} =$

|          | $p_0$ | $p_1$ | $p_2$ | $p_3$ | $p_4$ | $p_5$ |
|----------|-------|-------|-------|-------|-------|-------|
| $p_6$    | $p_6$ | $p_{17}$ | $p_{21}$ | $p_{12}$ | $p_{19}$ | $p_{11}$ |
| $p_7$    | $p_7$ | $p_{16}$ | $p_{23}$ | $p_{13}$ | $p_{20}$ | $p_{10}$ |
| $p_8$    | $p_8$ | $p_{15}$ | $p_{22}$ | $p_{14}$ | $p_{18}$ | $p_9$ |
| $p_9$    | $p_9$ | $p_{14}$ | $p_{18}$ | $p_{15}$ | $p_{22}$ | $p_8$ |
| $p_{10}$ | $p_{10}$ | $p_{13}$ | $p_{20}$ | $p_{16}$ | $p_{23}$ | $p_7$ |
| $p_{11}$ | $p_{11}$ | $p_{12}$ | $p_{19}$ | $p_{17}$ | $p_{21}$ | $p_6$ |

$P_{22} =$

|          | $p_6$ | $p_7$ | $p_8$ | $p_9$ | $p_{10}$ | $p_{11}$ |
|----------|-------|-------|-------|-------|----------|----------|
| $p_6$    | $p_0$ | $p_{13}$ | $p_{18}$ | $p_{15}$ | $p_{23}$ | $p_5$ |
| $p_7$    | $p_1$ | $p_{14}$ | $p_{19}$ | $p_{17}$ | $p_{22}$ | $p_4$ |
| $p_8$    | $p_2$ | $p_{12}$ | $p_{20}$ | $p_{16}$ | $p_{21}$ | $p_3$ |
| $p_9$    | $p_3$ | $p_{16}$ | $p_{21}$ | $p_{12}$ | $p_{20}$ | $p_2$ |
| $p_{10}$ | $p_4$ | $p_{17}$ | $p_{22}$ | $p_{14}$ | $p_{19}$ | $p_1$ |
| $p_{11}$ | $p_5$ | $p_{15}$ | $p_{23}$ | $p_{13}$ | $p_{18}$ | $p_0$ |

$P_{23} =$

|          | $p_{12}$ | $p_{13}$ | $p_{14}$ | $p_{15}$ | $p_{16}$ | $p_{17}$ |
|----------|----------|----------|----------|----------|----------|----------|
| $p_6$    | $p_3$ | $p_7$ | $p_{22}$ | $p_9$ | $p_{20}$ | $p_1$ |
| $p_7$    | $p_5$ | $p_8$ | $p_{21}$ | $p_{11}$ | $p_{18}$ | $p_2$ |
| $p_8$    | $p_4$ | $p_6$ | $p_{23}$ | $p_{10}$ | $p_{19}$ | $p_0$ |
| $p_9$    | $p_0$ | $p_{10}$ | $p_{19}$ | $p_6$ | $p_{23}$ | $p_4$ |
| $p_{10}$ | $p_2$ | $p_{11}$ | $p_{18}$ | $p_8$ | $p_{21}$ | $p_5$ |
| $p_{11}$ | $p_1$ | $p_9$ | $p_{20}$ | $p_7$ | $p_{22}$ | $p_3$ |

$P_{24} =$

|          | $p_{18}$ | $p_{19}$ | $p_{20}$ | $p_{21}$ | $p_{22}$ | $p_{23}$ |
|----------|----------|----------|----------|----------|----------|----------|
| $p_6$    | $p_8$ | $p_4$ | $p_{16}$ | $p_2$ | $p_{14}$ | $p_{10}$ |
| $p_7$    | $p_6$ | $p_3$ | $p_{15}$ | $p_0$ | $p_{12}$ | $p_9$ |
| $p_8$    | $p_7$ | $p_5$ | $p_{17}$ | $p_1$ | $p_{13}$ | $p_{11}$ |
| $p_9$    | $p_{11}$ | $p_1$ | $p_{13}$ | $p_5$ | $p_{17}$ | $p_7$ |
| $p_{10}$ | $p_9$ | $p_0$ | $p_{12}$ | $p_3$ | $p_{15}$ | $p_6$ |
| $p_{11}$ | $p_{10}$ | $p_2$ | $p_{14}$ | $p_4$ | $p_{16}$ | $p_8$ |

$P_{31}=$

| | $p_0$ | $p_1$ | $p_2$ | $p_3$ | $p_4$ | $p_5$ |
|---|---|---|---|---|---|---|
| $p_{12}$ | $p_{12}$ | $p_{19}$ | $p_{11}$ | $p_6$ | $p_{17}$ | $p_{21}$ |
| $p_{13}$ | $p_{13}$ | $p_{20}$ | $p_{10}$ | $p_7$ | $p_{16}$ | $p_{23}$ |
| $p_{14}$ | $p_{14}$ | $p_{18}$ | $p_9$ | $p_8$ | $p_{15}$ | $p_{22}$ |
| $p_{15}$ | $p_{15}$ | $p_{22}$ | $p_8$ | $p_9$ | $p_{14}$ | $p_{18}$ |
| $p_{16}$ | $p_{16}$ | $p_{23}$ | $p_7$ | $p_{10}$ | $p_{13}$ | $p_{20}$ |
| $p_{17}$ | $p_{17}$ | $p_{21}$ | $p_6$ | $p_{11}$ | $p_{12}$ | $p_1$ |

$P_{32}=$

| | $p_6$ | $p_7$ | $p_8$ | $p_9$ | $p_{10}$ | $p_{11}$ |
|---|---|---|---|---|---|---|
| $p_{12}$ | $p_{15}$ | $p_{23}$ | $p_5$ | $p_0$ | $p_{13}$ | $p_{18}$ |
| $p_{13}$ | $p_{17}$ | $p_{22}$ | $p_4$ | $p_1$ | $p_{14}$ | $p_{19}$ |
| $p_{14}$ | $p_{16}$ | $p_{21}$ | $p_3$ | $p_2$ | $p_{12}$ | $p_{20}$ |
| $p_{15}$ | $p_{12}$ | $p_{20}$ | $p_2$ | $p_3$ | $p_{16}$ | $p_{21}$ |
| $p_{16}$ | $p_{14}$ | $p_{19}$ | $p_1$ | $p_4$ | $p_{17}$ | $p_{22}$ |
| $p_{17}$ | $p_{13}$ | $p_{18}$ | $p_0$ | $p_5$ | $p_{15}$ | $p_{23}$ |

$P_{33}=$

| | $p_{12}$ | $p_{13}$ | $p_{14}$ | $p_{15}$ | $p_{16}$ | $p_{17}$ |
|---|---|---|---|---|---|---|
| $p_{12}$ | $p_9$ | $p_{20}$ | $p_1$ | $p_3$ | $p_7$ | $p_{22}$ |
| $p_{13}$ | $p_{11}$ | $p_{18}$ | $p_2$ | $p_5$ | $p_8$ | $p_{21}$ |
| $p_{14}$ | $p_{10}$ | $p_{19}$ | $p_0$ | $p_4$ | $p_6$ | $p_{23}$ |
| $p_{15}$ | $p_6$ | $p_{23}$ | $p_4$ | $p_0$ | $p_{10}$ | $p_{19}$ |
| $p_{16}$ | $p_8$ | $p_{21}$ | $p_5$ | $p_2$ | $p_{11}$ | $p_{18}$ |
| $p_{17}$ | $p_7$ | $p_{22}$ | $p_3$ | $p_1$ | $p_9$ | $p_{20}$ |

$P_{34}=$

| | $p_{18}$ | $p_{19}$ | $p_{20}$ | $p_{21}$ | $p_{22}$ | $p_{23}$ |
|---|---|---|---|---|---|---|
| $p_{12}$ | $p_2$ | $p_{14}$ | $p_{10}$ | $p_8$ | $p_4$ | $p_{16}$ |
| $p_{13}$ | $p_0$ | $p_{12}$ | $p_9$ | $p_6$ | $p_3$ | $p_{15}$ |
| $p_{14}$ | $p_1$ | $p_{13}$ | $p_{11}$ | $p_7$ | $p_5$ | $p_{17}$ |
| $p_{15}$ | $p_5$ | $p_{17}$ | $p_7$ | $p_{11}$ | $p_1$ | $p_{13}$ |
| $p_{16}$ | $p_3$ | $p_{15}$ | $p_6$ | $p_9$ | $p_0$ | $p_{12}$ |
| $p_{17}$ | $p_4$ | $p_{16}$ | $p_8$ | $p_{10}$ | $p_2$ | $p_{14}$ |

$P_{41}=$

| | $p_0$ | $p_1$ | $p_2$ | $p_3$ | $p_4$ | $p_5$ |
|---|---|---|---|---|---|---|
| $p_{18}$ | $p_{18}$ | $p_9$ | $p_{14}$ | $p_{22}$ | $p_8$ | $p_{15}$ |
| $p_{19}$ | $p_{19}$ | $p_{11}$ | $p_{12}$ | $p_{21}$ | $p_6$ | $p_{17}$ |
| $p_{20}$ | $p_{20}$ | $p_{10}$ | $p_{13}$ | $p_{23}$ | $p_7$ | $p_{16}$ |
| $p_{21}$ | $p_{21}$ | $p_6$ | $p_{17}$ | $p_{19}$ | $p_{11}$ | $p_{12}$ |
| $p_{22}$ | $p_{22}$ | $p_8$ | $p_{15}$ | $p_{18}$ | $p_9$ | $p_{14}$ |
| $p_{23}$ | $p_{23}$ | $p_7$ | $p_{16}$ | $p_{20}$ | $p_{10}$ | $p_{13}$ |

$P_{42}=$

| | $p_6$ | $p_7$ | $p_8$ | $p_9$ | $p_{10}$ | $p_{11}$ |
|---|---|---|---|---|---|---|
| $p_{18}$ | $p_{21}$ | $p_3$ | $p_{16}$ | $p_{20}$ | $p_2$ | $p_{12}$ |
| $p_{19}$ | $p_{23}$ | $p_5$ | $p_{15}$ | $p_{18}$ | $p_0$ | $p_{13}$ |
| $p_{20}$ | $p_{22}$ | $p_4$ | $p_{17}$ | $p_{19}$ | $p_1$ | $p_{14}$ |
| $p_{21}$ | $p_{18}$ | $p_0$ | $p_{13}$ | $p_{23}$ | $p_5$ | $p_{15}$ |
| $p_{22}$ | $p_{20}$ | $p_2$ | $p_{12}$ | $p_{21}$ | $p_3$ | $p_{16}$ |
| $p_{23}$ | $p_{19}$ | $p_1$ | $p_{14}$ | $p_{22}$ | $p_4$ | $p_{17}$ |

$P_{43}=$

| | $p_{12}$ | $p_{13}$ | $p_{14}$ | $p_{15}$ | $p_{16}$ | $p_{17}$ |
|---|---|---|---|---|---|---|
| $p_{18}$ | $p_{19}$ | $p_0$ | $p_{10}$ | $p_{23}$ | $p_4$ | $p_6$ |
| $p_{19}$ | $p_{20}$ | $p_1$ | $p_9$ | $p_{22}$ | $p_3$ | $p_7$ |
| $p_{20}$ | $p_{18}$ | $p_2$ | $p_{11}$ | $p_{21}$ | $p_5$ | $p_8$ |
| $p_{21}$ | $p_{22}$ | $p_3$ | $p_7$ | $p_{20}$ | $p_1$ | $p_9$ |
| $p_{22}$ | $p_{23}$ | $p_4$ | $p_6$ | $p_{19}$ | $p_0$ | $p_{10}$ |
| $p_{23}$ | $p_{21}$ | $p_5$ | $p_8$ | $p_{18}$ | $p_2$ | $p_{11}$ |

$P_{44}=$

| | $p_{18}$ | $p_{19}$ | $p_{20}$ | $p_{21}$ | $p_{22}$ | $p_{23}$ |
|---|---|---|---|---|---|---|
| $p_{18}$ | $p_{13}$ | $p_{11}$ | $p_1$ | $p_{17}$ | $p_7$ | $p_5$ |
| $p_{19}$ | $p_{14}$ | $p_{10}$ | $p_2$ | $p_{16}$ | $p_8$ | $p_4$ |
| $p_{20}$ | $p_{12}$ | $p_9$ | $p_0$ | $p_{15}$ | $p_6$ | $p_3$ |
| $p_{21}$ | $p_{16}$ | $p_8$ | $p_4$ | $p_{14}$ | $p_{10}$ | $p_2$ |
| $p_{22}$ | $p_{17}$ | $p_7$ | $p_5$ | $p_{13}$ | $p_{11}$ | $p_1$ |
| $p_{23}$ | $p_{15}$ | $p_6$ | $p_3$ | $p_{12}$ | $p_9$ | $p_0$ |

Группа всех симметрий тетраэдра изоморфна группе $S_4$, симметрической группе перестановок четырёх элементов. Группа автоморфизма графа $K_4$ содержит подгруппу, состоящую из четырех элементов:

$p_0$ = ‹1,2,3,4› = (1)(2)(3)(4);
$p_{11}$ = ‹2,1,4,3› = (1 2)(3 4);

$p_{14}$ = <3,4,1,2> = (1 3)(2 4);
$p_{20}$ = <4,3,2,1> = (1 4)(2 3).

Данная подгруппа называемой группой Клейна четвертого порядка.

Так как граф тетраэдра планарен, то кольцевая сумма изометрических циклов есть пустое множество.

Веса вершин векторного инварианта для спектра реберных разрезов графа тетраэдра равны между собой $\zeta_w(K_4)$ =<12,12,12,12>. Веса ребер также равны между собой $F_{es}$ = <4,4,4,4,4,4>. Орбита графа характеризуется множеством вершин образующего изометрического цикла. Например, изометрический цикл $c_1=\{e_1,e_2,e_5\}$ множество вершин $\{v_1,v_2,v_3\}$ которого, образует орбиту графа К$_4$ [15].

## 2.4. Граф ГЕКТАЭДРа (КУБа)

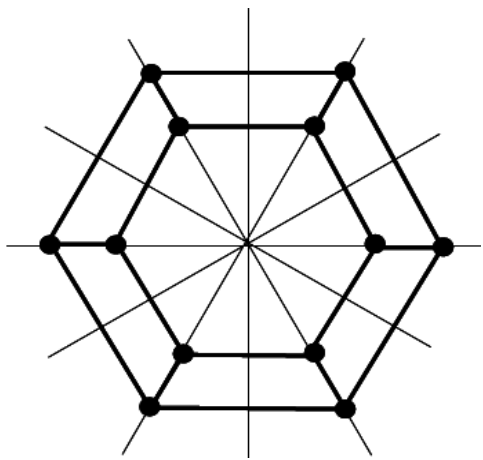

Рис. 2.47. Многоугольник с двумя связанными контурами.

Будем рассматривать следующее симметричное построение многоугольника с 2n вершинами (см. рис. 2.47).

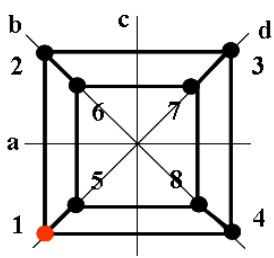 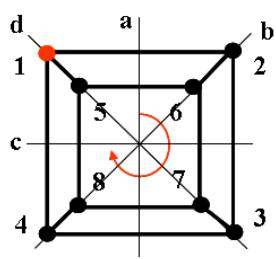 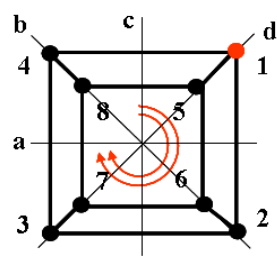 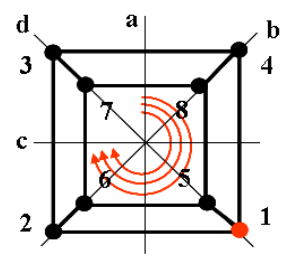

Рис. 2.48. Перестановка (1 2 3 4 5 6 7 8).    Рис. 2.49. Перестановка (4 1 2 3 8 5 6 7).    Рис. 2.50. Перестановка (3 4 1 2 7 8 5 6).    Рис. 2.51. Перестановка (2 3 4 1 6 7 8 5).

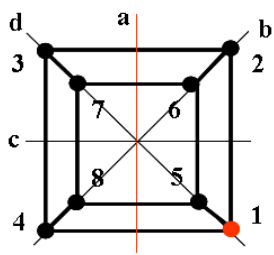 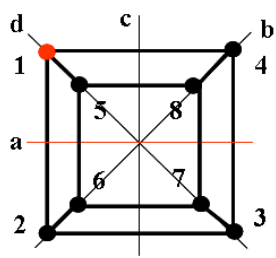 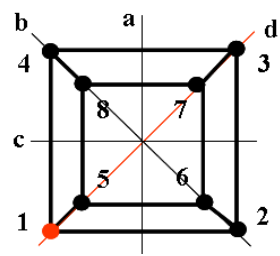 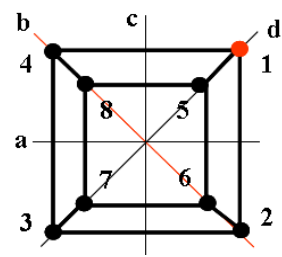

Рис. 2.52. Перестановка (4 3 2 1 8 7 6 5).    Рис. 2.53. Перестановка (2 1 4 3 6 5 8 7).    Рис. 2.54. Перестановка (1 4 3 2 5 8 7 6).    Рис. 2.55. Перестановка (3 4 1 2 7 8 5 6).

Представленный правильный *n*-угольник с двумя связанными контурами имеет различных симметрий: поворотов и осевых отражений, образующих диэдральную группу . Если нечётно, каждая ось симметрии проходит через середину одной из сторон и противоположную вершину. Если чётно, имеется осей симметрии, соединяющих середины противоположных сторон и осей, соединяющих противоположные вершины. В любом случае, имеется осей симметрии и элементов в группе симметрий. Отражение относительно одной оси, а затем относительно другой, приводит к вращению на удвоенный угол между осями.

Рассмотрим граф куба (см. рис. 2.56).

Рис. 2.56. Граф куба $G_7$ и его топологический рисунок.

Определим орбиты из векторного интегрального инварианта графа для спектра реберных разрезов. Вектор весов вершин $F(\zeta_w(G_7)) = 8\times36$. Все вершины и ребра графа имеют равные веса.

Множество изометрических циклов запишем в следующем виде:

$c_1 = \{e_1,e_2,e_4,e_6\} \rightarrow \{v_1,v_2,v_3,v_4\}$;  $c_2 = \{e_4,e_5,e_7,e_{11}\} \rightarrow \{v_2,v_3,v_6,v_7\}$;
$c_3 = \{e_2,e_3,e_8,e_{10}\} \rightarrow \{v_1,v_4,v_5,v_8\}$;  $c_4 = \{e_9,e_{10},e_{11},e_{12}\} \rightarrow \{v_5,v_6,v_7,v_8\}$;
$c_5 = \{e_6,e_7,e_8,e_{12}\} \rightarrow \{v_3,v_4,v_7,v_8\}$;  $c_6 = \{e_1,e_3,e_5,e_9\} \rightarrow \{v_1,v_2,v_5,v_6\}$.

Рис. 2.57. Перестановка $p_0$.
$\begin{pmatrix} 1 & 2 & 3 & 4 & 5 & 6 & 7 & 8 \\ 1 & 2 & 3 & 4 & 5 & 6 & 7 & 8 \end{pmatrix}$

Рис. 2.58. Перестановка $p_1$.
$\begin{pmatrix} 1 & 2 & 3 & 4 & 5 & 6 & 7 & 8 \\ 4 & 1 & 2 & 3 & 8 & 5 & 6 & 7 \end{pmatrix}$

Рис. 2.59. Перестановка $p_2$.
$\begin{pmatrix} 1 & 2 & 3 & 4 & 5 & 6 & 7 & 8 \\ 3 & 4 & 1 & 2 & 7 & 8 & 5 & 6 \end{pmatrix}$

Рис. 2.60. Перестановка $p_3$.
$\begin{pmatrix} 1 & 2 & 3 & 4 & 5 & 6 & 7 & 8 \\ 2 & 3 & 4 & 1 & 6 & 7 & 8 & 5 \end{pmatrix}$

Так как все изометрические циклы имеют в своем составе равные веса вершин и ребер, то их можно представлять как образующие изометрические циклы и орбиты. Будем считать, что каждый образующий цикл графа $G_1$ индуцирует геометрический рисунок правильного многоугольника. В данном случае это квадрат. Количество образующих изометрических циклов равно шести. Количество перестановок определяется симметрической группой квадрата [6]. Будем рассматривать перестановки для каждого многоугольника индуцированного образующим циклом графа. Всего таких перестановок будет $6\times8 = 48$, которые и определяют все перестановки вершин.

В качестве образующего цикла выбираем цикл $c_1$.

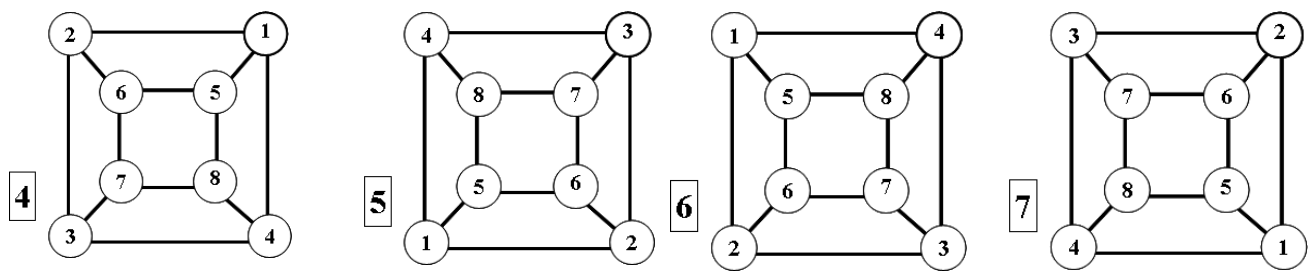

Рис. 2.61.  Рис. 2.62.  Рис. 2.63.  Рис. 2.64.
Перестановка $p_4$.  Перестановка $p_5$.  Перестановка $p_6$.  Перестановка $p_7$.

$\begin{pmatrix} 1 & 2 & 3 & 4 & 5 & 6 & 7 & 8 \\ 2 & 1 & 4 & 3 & 6 & 5 & 8 & 7 \end{pmatrix}$ $\begin{pmatrix} 1 & 2 & 3 & 4 & 5 & 6 & 7 & 8 \\ 4 & 3 & 2 & 1 & 8 & 7 & 6 & 5 \end{pmatrix}$ $\begin{pmatrix} 1 & 2 & 3 & 4 & 5 & 6 & 7 & 8 \\ 1 & 4 & 3 & 2 & 5 & 8 & 7 & 6 \end{pmatrix}$ $\begin{pmatrix} 1 & 2 & 3 & 4 & 5 & 6 & 7 & 8 \\ 3 & 2 & 1 & 4 & 7 & 6 & 5 & 8 \end{pmatrix}$

Перестановки определяемые образующими циклами $c_2, c_3, c_4, c_5$ предлагается проделать в качестве самостоятельных упражнений.

Выберем в качестве образующего цикла – цикл $c_6$. Построим перестановки определяемые циклом $c_6$.

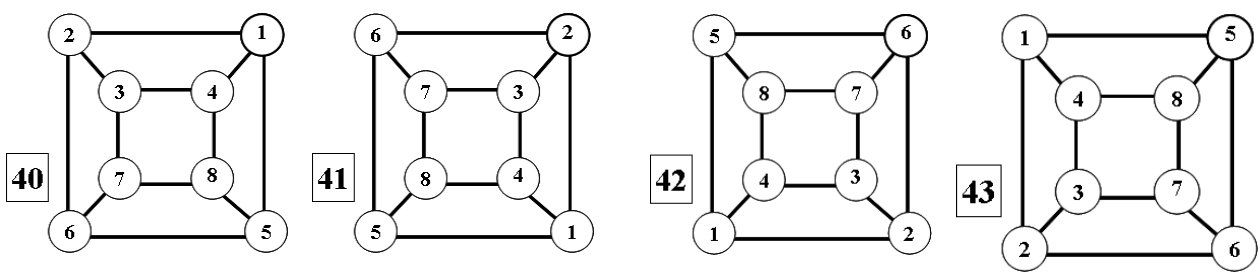

Рис. 2658.  Рис. 2.66.  Рис. 2.67.  Рис. 2.68.
Перестановка $p_{40}$.  Перестановка $p_{41}$.  Перестановка $p_{42}$.  Перестановка $p_{43}$.

$\begin{pmatrix} 1 & 2 & 3 & 4 & 5 & 6 & 7 & 8 \\ 2 & 1 & 5 & 6 & 3 & 4 & 8 & 7 \end{pmatrix}$ $\begin{pmatrix} 1 & 2 & 3 & 4 & 5 & 6 & 7 & 8 \\ 6 & 2 & 1 & 5 & 7 & 3 & 4 & 8 \end{pmatrix}$ $\begin{pmatrix} 1 & 2 & 3 & 4 & 5 & 6 & 7 & 8 \\ 5 & 6 & 2 & 1 & 8 & 7 & 3 & 4 \end{pmatrix}$ $\begin{pmatrix} 1 & 2 & 3 & 4 & 5 & 6 & 7 & 8 \\ 1 & 5 & 6 & 2 & 4 & 8 & 7 & 3 \end{pmatrix}$

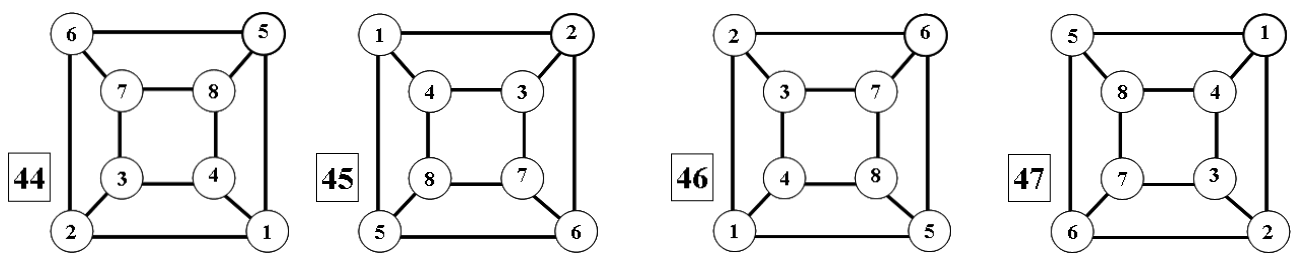

Рис. 2.69.  Рис. 2.70.  Рис. 2.71..  Рис. 2.72.
Перестановка $p_{44}$.  Перестановка $p_{45}$.  Перестановка $p_{46}$.  Перестановка $p_{47}$.

$\begin{pmatrix} 1 & 2 & 3 & 4 & 5 & 6 & 7 & 8 \\ 6 & 5 & 1 & 2 & 7 & 8 & 4 & 3 \end{pmatrix}$ $\begin{pmatrix} 1 & 2 & 3 & 4 & 5 & 6 & 7 & 8 \\ 1 & 2 & 6 & 5 & 4 & 3 & 7 & 8 \end{pmatrix}$ $\begin{pmatrix} 1 & 2 & 3 & 4 & 5 & 6 & 7 & 8 \\ 2 & 6 & 5 & 1 & 3 & 7 & 8 & 4 \end{pmatrix}$ $\begin{pmatrix} 1 & 2 & 3 & 4 & 5 & 6 & 7 & 8 \\ 5 & 1 & 2 & 6 & 8 & 4 & 3 & 7 \end{pmatrix}$

Запишем множество перестановок графа $G_1$ индуцированные различными рисунками:

$p_0 = <1,2,3,4,5,6,7,8> = (1)(2)(3)(4)(5)(6)(7)(8)$;

$p_1 = \langle 4,1,2,3,8,5,6,7 \rangle = (1\ 4\ 3\ 2)(5\ 8\ 7\ 6)$;
$p_2 = \langle 3,4,1,2,7,8,5,6 \rangle = (1\ 3)(2\ 4)(5\ 7)(6\ 8)$;
$p_3 = \langle 2,3,4,1,6,7,8,5 \rangle = (1\ 2\ 3\ 4)(5\ 6\ 7\ 8)$;
$p_4 = \langle 2,1,4,3,6,5,8,7 \rangle = (1\ 2)(3\ 4)(5\ 6)(7\ 8)$;
$p_5 = \langle 4,3,2,1,8,7,6,5 \rangle = (1\ 4)(2\ 3)(5\ 8)(6\ 7)$;
$p_6 = \langle 1,4,3,2,5,8,7,6 \rangle = (1)(2\ 4)(3)(5)(6\ 8)(7)$;
$p_7 = \langle 3,2,1,4,7,6,5,8 \rangle = (1\ 3)(2)(4)(5\ 7)(6)(8)$;
$p_8 = \langle 2,6,7,3,1,5,8,4 \rangle = (1\ 2\ 6\ 5)(3\ 7\ 8\ 4)$;
$p_9 = \langle 3,2,6,7,4,1,5,8 \rangle = (1\ 3\ 6)(2)(4\ 7\ 5)(8)$;
$p_{10} = \langle 7,3,2,6,8,4,1,5 \rangle = (1\ 7)(2\ 3)(4\ 6)(5\ 8)$;
$p_{11} = \langle 6,7,3,2,5,8,4,1 \rangle = (1\ 6\ 8)(2\ 7\ 4)(3)(5)$;
$p_{12} = \langle 6,2,3,7,5,1,4,8 \rangle = (1\ 6)(2)(3)(4\ 7)(5)(8)$;
$p_{13} = \langle 3,7,6,2,4,8,5,1 \rangle = (1\ 3\ 6\ 8)(2\ 7\ 5\ 4)$;
$p_{14} = \langle 7,6,2,3,8,5,1,4 \rangle = (1\ 7)(2\ 6\ 5\ 8\ 4\ 3)$;
$p_{15} = \langle 2,3,7,6,1,4,8,5 \rangle = (1\ 2\ 3\ 7\ 8\ 5)(4\ 6)$;
$p_{16} = \langle 1,4,8,5,2,3,7,6 \rangle = (1)(2\ 4\ 5)(3\ 8\ 6)(7)$;
$p_{17} = \langle 5,1,4,8,6,2,3,7 \rangle = (1\ 5\ 6\ 2)(3\ 4\ 8\ 7)$;
$p_{18} = \langle 8,5,1,4,7,6,2,3 \rangle = (1\ 8\ 3)(2\ 5\ 7)(4)(6)$;
$p_{19} = \langle 4,8,5,1,3,7,6,2 \rangle = (1\ 4)(2\ 8)(3\ 5)(6\ 7)$;
$p_{20} = \langle 5,8,4,1,6,7,3,2 \rangle = (1\ 5\ 6\ 7\ 3\ 4)(2\ 8)$;
$p_{21} = \langle 4,1,5,8,3,2,6,7 \rangle = (1\ 4\ 8\ 7\ 6\ 2)(3\ 5)$;
$p_{22} = \langle 8,4,1,5,7,3,2,6 \rangle = (1\ 8\ 6\ 3)(2\ 4\ 5\ 7)$;
$p_{23} = \langle 1,5,8,4,2,6,7,3 \rangle = (1)(2\ 5)(3\ 8)(4)(6\ 7)$;
$p_{24} = \langle 8,7,6,5,4,3,2,1 \rangle = (1\ 8)(2\ 7)(3\ 6)(4\ 5)$;
$p_{25} = \langle 5,8,7,6,1,4,3,2 \rangle = (1\ 5)(2\ 8)(3\ 7)(4\ 6)$;
$p_{26} = \langle 6,5,8,7,2,1,4,3 \rangle = (1\ 6)(2\ 5)(3\ 8)(4\ 7)$;
$p_{27} = \langle 7,6,5,8,3,2,1,4 \rangle = (1\ 7)(2\ 6)(3\ 5)(4\ 8)$;
$p_{28} = \langle 5,6,7,8,1,2,3,4 \rangle = (1\ 5)(2\ 6)(3\ 7)(4\ 8)$;
$p_{29} = \langle 7,8,5,6,3,4,1,2 \rangle = (1\ 7)(2\ 8)(3\ 5)(4\ 6)$;
$p_{30} = \langle 8,5,6,7,4,1,2,3 \rangle = (1\ 8\ 3\ 6)(2\ 5\ 4\ 7)$;
$p_{31} = \langle 6,7,8,5,2,3,4,1 \rangle = (1\ 6\ 3\ 8)(2\ 7\ 4\ 5)$;
$p_{32} = \langle 4,3,7,8,1,2,6,5 \rangle = (1\ 4\ 8\ 5)(2\ 3\ 7\ 6)$;
$p_{33} = \langle 8,4,3,7,5,1,2,6 \rangle = (1\ 8\ 6)(2\ 4\ 7)(3)(5)$;
$p_{34} = \langle 7,8,4,3,6,5,1,2 \rangle = (1\ 7)(2\ 8)(3\ 4)(5\ 6)$;
$p_{35} = \langle 3,7,8,4,2,6,5,1 \rangle = (1\ 3\ 8)(2\ 7\ 5)(4)(6)$;
$p_{36} = \langle 8,7,3,4,5,6,2,1 \rangle = (1\ 8)(2\ 7)(3\ 4)(5)(6)$;
$p_{37} = \langle 3,4,8,7,2,1,5,6 \rangle = (1\ 3\ 8\ 6)(2\ 4\ 7\ 5)$;
$p_{38} = \langle 4,8,7,3,1,5,6,2 \rangle = (1\ 4\ 3\ 7\ 6\ 5)(2\ 8)$;
$p_{39} = \langle 7,3,4,8,6,2,1,5 \rangle = (1\ 7)(2\ 3\ 4\ 8\ 5\ 6)$;
$p_{40} = \langle 2,1,5,6,3,4,8,7 \rangle = (1\ 2)(3\ 5)(4\ 6)(7\ 8)$;
$p_{41} = \langle 6,2,1,5,7,3,4,8 \rangle = (1\ 6\ 3)(2)(4\ 5\ 7)(8)$;
$p_{42} = \langle 5,6,2,1,8,7,3,4 \rangle = (1\ 5\ 8\ 4)(2\ 6\ 7\ 3)$;
$p_{43} = \langle 1,5,6,2,4,8,7,3 \rangle = (1)(2\ 5\ 4)(3\ 6\ 8)(7)$;
$p_{44} = \langle 6,5,1,2,7,8,4,3 \rangle = (1\ 6\ 8\ 3)(2\ 5\ 7\ 4)$;
$p_{45} = \langle 1,2,6,5,4,3,7,8 \rangle = (1)(2)(3\ 6)(4\ 5)(7)(8)$;
$p_{46} = \langle 2,6,5,1,3,7,8,4 \rangle = (1\ 2\ 6\ 7\ 8\ 4)(3\ 5)$;
$p_{47} = \langle 5,1,2,6,8,4,3,7 \rangle = (1\ 5\ 8\ 7\ 3\ 2)(4\ 6)$.

Граф $G_7$ является планарным графом, кольцевая сумма изометрических циклов есть пустое множество $\sum_{i=1}^{6} c_i = \varnothing$.

Разобьем таблицу умножения перестановок на блоки. Каждый блок соответствует подгруппе и характеризует образующий цикл графа. С другой стороны – это орбита графа:

$$Aut(G_1) = \begin{bmatrix} P_{11} & P_{12} & P_{13} & P_{14} & P_{15} & P_{16} \\ P_{21} & P_{22} & P_{23} & P_{24} & P_{25} & P_{26} \\ P_{31} & P_{32} & P_{33} & P_{34} & P_{35} & P_{36} \\ P_{41} & P_{42} & P_{43} & P_{44} & P_{45} & P_{46} \\ P_{51} & P_{52} & P_{53} & P_{54} & P_{55} & P_{56} \\ P_{61} & P_{62} & P_{63} & P_{64} & P_{65} & P_{66} \end{bmatrix}$$

$P_{11}=$

|       | $p_0$ | $p_1$ | $p_2$ | $p_3$ | $p_4$ | $p_5$ | $p_6$ | $p_7$ |
|-------|-------|-------|-------|-------|-------|-------|-------|-------|
| $p_0$ | $p_0$ | $p_1$ | $p_2$ | $p_3$ | $p_4$ | $p_5$ | $p_6$ | $p_7$ |
| $p_1$ | $p_1$ | $p_2$ | $p_3$ | $p_0$ | $p_7$ | $p_6$ | $p_4$ | $p_5$ |
| $p_2$ | $p_2$ | $p_3$ | $p_0$ | $p_1$ | $p_5$ | $p_4$ | $p_7$ | $p_6$ |
| $p_3$ | $p_3$ | $p_0$ | $p_1$ | $p_2$ | $p_6$ | $p_7$ | $p_5$ | $p_4$ |
| $p_4$ | $p_4$ | $p_6$ | $p_5$ | $p_7$ | $p_0$ | $p_2$ | $p_1$ | $p_3$ |
| $p_5$ | $p_5$ | $p_7$ | $p_4$ | $p_6$ | $p_2$ | $p_0$ | $p_3$ | $p_1$ |
| $p_6$ | $p_6$ | $p_5$ | $p_7$ | $p_4$ | $p_3$ | $p_1$ | $p_0$ | $p_2$ |
| $p_7$ | $p_7$ | $p_4$ | $p_6$ | $p_5$ | $p_1$ | $p_3$ | $p_2$ | $p_0$ |

$P_{12}=$

|       | $p_8$ | $p_9$ | $p_{10}$ | $p_{11}$ | $p_{12}$ | $p_{13}$ | $p_{14}$ | $p_{15}$ |
|-------|-------|-------|----------|----------|----------|----------|----------|----------|
| $p_0$ | $p_8$ | $p_9$ | $p_{10}$ | $p_{11}$ | $p_{12}$ | $p_{13}$ | $p_{14}$ | $p_{15}$ |
| $p_1$ | $p_9$ | $p_{10}$ | $p_{11}$ | $p_8$ | $p_{14}$ | $p_{15}$ | $p_{13}$ | $p_{12}$ |
| $p_2$ | $p_{10}$ | $p_{11}$ | $p_8$ | $p_9$ | $p_{13}$ | $p_{12}$ | $p_{15}$ | $p_{14}$ |
| $p_3$ | $p_{11}$ | $p_8$ | $p_9$ | $p_{10}$ | $p_{15}$ | $p_{14}$ | $p_{12}$ | $p_{13}$ |
| $p_4$ | $p_{12}$ | $p_{15}$ | $p_{13}$ | $p_{14}$ | $p_8$ | $p_{10}$ | $p_{11}$ | $p_9$ |
| $p_5$ | $p_{13}$ | $p_{14}$ | $p_{12}$ | $p_{15}$ | $p_{10}$ | $p_8$ | $p_9$ | $p_{11}$ |
| $p_6$ | $p_{15}$ | $p_{13}$ | $p_{14}$ | $p_{12}$ | $p_{11}$ | $p_9$ | $p_{10}$ | $p_8$ |
| $p_7$ | $p_{14}$ | $p_{12}$ | $p_{15}$ | $p_{13}$ | $p_9$ | $p_{11}$ | $p_8$ | $p_{10}$ |

$P_{13}=$

|       | $p_{16}$ | $p_{17}$ | $p_{18}$ | $p_{19}$ | $p_{20}$ | $p_{21}$ | $p_{22}$ | $p_{23}$ |
|-------|----------|----------|----------|----------|----------|----------|----------|----------|
| $p_0$ | $p_{16}$ | $p_{17}$ | $p_{18}$ | $p_{19}$ | $p_{20}$ | $p_{21}$ | $p_{22}$ | $p_{23}$ |
| $p_1$ | $p_{17}$ | $p_{18}$ | $p_{19}$ | $p_{16}$ | $p_{23}$ | $p_{22}$ | $p_{20}$ | $p_{21}$ |
| $p_2$ | $p_{18}$ | $p_{19}$ | $p_{16}$ | $p_{17}$ | $p_{21}$ | $p_{20}$ | $p_{23}$ | $p_{22}$ |
| $p_3$ | $p_{19}$ | $p_{16}$ | $p_{17}$ | $p_{18}$ | $p_{22}$ | $p_{23}$ | $p_{21}$ | $p_{20}$ |
| $p_4$ | $p_{21}$ | $p_{23}$ | $p_{20}$ | $p_{22}$ | $p_{18}$ | $p_{16}$ | $p_{19}$ | $p_{17}$ |
| $p_5$ | $p_{20}$ | $p_{22}$ | $p_{21}$ | $p_{23}$ | $p_{16}$ | $p_{18}$ | $p_{17}$ | $p_{19}$ |
| $p_6$ | $p_{23}$ | $p_{20}$ | $p_{22}$ | $p_{21}$ | $p_{17}$ | $p_{19}$ | $p_{18}$ | $p_{16}$ |
| $p_7$ | $p_{22}$ | $p_{21}$ | $p_{23}$ | $p_{20}$ | $p_{19}$ | $p_{17}$ | $p_{16}$ | $p_{18}$ |

$P_{14}=$

|       | $p_{24}$ | $p_{25}$ | $p_{26}$ | $p_{27}$ | $p_{28}$ | $p_{29}$ | $p_{30}$ | $p_{31}$ |
|-------|----------|----------|----------|----------|----------|----------|----------|----------|
| $p_0$ | $p_{24}$ | $p_{25}$ | $p_{26}$ | $p_{27}$ | $p_{28}$ | $p_{29}$ | $p_{30}$ | $p_{31}$ |
| $p_1$ | $p_{25}$ | $p_{26}$ | $p_{27}$ | $p_{24}$ | $p_{30}$ | $p_{31}$ | $p_{29}$ | $p_{28}$ |
| $p_2$ | $p_{26}$ | $p_{27}$ | $p_{24}$ | $p_{25}$ | $p_{29}$ | $p_{28}$ | $p_{31}$ | $p_{30}$ |
| $p_3$ | $p_{27}$ | $p_{24}$ | $p_{25}$ | $p_{26}$ | $p_{31}$ | $p_{30}$ | $p_{28}$ | $p_{29}$ |
| $p_4$ | $p_{29}$ | $p_{30}$ | $p_{28}$ | $p_{31}$ | $p_{26}$ | $p_{24}$ | $p_{25}$ | $p_{27}$ |
| $p_5$ | $p_{28}$ | $p_{31}$ | $p_{29}$ | $p_{30}$ | $p_{24}$ | $p_{26}$ | $p_{27}$ | $p_{25}$ |
| $p_6$ | $p_{30}$ | $p_{28}$ | $p_{31}$ | $p_{29}$ | $p_{25}$ | $p_{27}$ | $p_{24}$ | $p_{26}$ |
| $p_7$ | $p_{31}$ | $p_{29}$ | $p_{30}$ | $p_{28}$ | $p_{27}$ | $p_{25}$ | $p_{26}$ | $p_{24}$ |

$P_{15}=$

|  | $p_{32}$ | $p_{33}$ | $p_{34}$ | $p_{35}$ | $p_{36}$ | $p_{37}$ | $p_{38}$ | $p_{39}$ |
|---|---|---|---|---|---|---|---|---|
| $p_0$ | $p_{32}$ | $p_{33}$ | $p_{34}$ | $p_{35}$ | $p_{36}$ | $p_{37}$ | $p_{38}$ | $p_{39}$ |
| $p_1$ | $p_{33}$ | $p_{34}$ | $p_{35}$ | $p_{32}$ | $p_{38}$ | $p_{39}$ | $p_{37}$ | $p_{36}$ |
| $p_2$ | $p_{34}$ | $p_{35}$ | $p_{32}$ | $p_{33}$ | $p_{37}$ | $p_{36}$ | $p_{39}$ | $p_{38}$ |
| $p_3$ | $p_{35}$ | $p_{32}$ | $p_{33}$ | $p_{34}$ | $p_{39}$ | $p_{38}$ | $p_{36}$ | $p_{37}$ |
| $p_4$ | $p_{37}$ | $p_{38}$ | $p_{36}$ | $p_{39}$ | $p_{34}$ | $p_{32}$ | $p_{33}$ | $p_{35}$ |
| $p_5$ | $p_{36}$ | $p_{39}$ | $p_{37}$ | $p_{38}$ | $p_{32}$ | $p_{34}$ | $p_{35}$ | $p_{33}$ |
| $p_6$ | $p_{38}$ | $p_{36}$ | $p_{39}$ | $p_{37}$ | $p_{33}$ | $p_{35}$ | $p_{32}$ | $p_{34}$ |
| $p_7$ | $p_{39}$ | $p_{37}$ | $p_{38}$ | $p_{36}$ | $p_{35}$ | $p_{33}$ | $p_{34}$ | $p_{32}$ |

$P_{16}=$

|  | $p_{40}$ | $p_{41}$ | $p_{42}$ | $p_{43}$ | $p_{44}$ | $p_{45}$ | $p_{46}$ | $p_{47}$ |
|---|---|---|---|---|---|---|---|---|
| $p_0$ | $p_{40}$ | $p_{41}$ | $p_{42}$ | $p_{43}$ | $p_{44}$ | $p_{45}$ | $p_{46}$ | $p_{47}$ |
| $p_1$ | $p_{41}$ | $p_{42}$ | $p_{43}$ | $p_{40}$ | $p_{46}$ | $p_{47}$ | $p_{45}$ | $p_{44}$ |
| $p_2$ | $p_{42}$ | $p_{43}$ | $p_{40}$ | $p_{41}$ | $p_{45}$ | $p_{44}$ | $p_{47}$ | $p_{46}$ |
| $p_3$ | $p_{43}$ | $p_{40}$ | $p_{41}$ | $p_{42}$ | $p_{47}$ | $p_{46}$ | $p_{44}$ | $p_{45}$ |
| $p_4$ | $p_{45}$ | $p_{46}$ | $p_{44}$ | $p_{47}$ | $p_{42}$ | $p_{40}$ | $p_{41}$ | $p_{43}$ |
| $p_5$ | $p_{44}$ | $p_{47}$ | $p_{45}$ | $p_{46}$ | $p_{40}$ | $p_{42}$ | $p_{43}$ | $p_{41}$ |
| $p_6$ | $p_{46}$ | $p_{44}$ | $p_{47}$ | $p_{45}$ | $p_{41}$ | $p_{43}$ | $p_{40}$ | $p_{42}$ |
| $p_7$ | $p_{47}$ | $p_{45}$ | $p_{46}$ | $p_{44}$ | $p_{43}$ | $p_{41}$ | $p_{42}$ | $p_{40}$ |

$P_{21}=$

|  | $p_0$ | $p_1$ | $p_2$ | $p_3$ | $p_4$ | $p_5$ | $p_6$ | $p_7$ |
|---|---|---|---|---|---|---|---|---|
| $p_8$ | $p_8$ | $p_{43}$ | $p_{19}$ | $p_{35}$ | $p_{23}$ | $p_{13}$ | $p_{38}$ | $p_{46}$ |
| $p_9$ | $p_9$ | $p_{40}$ | $p_{16}$ | $p_{32}$ | $p_{21}$ | $p_{15}$ | $p_{37}$ | $p_{45}$ |
| $p_{10}$ | $p_{10}$ | $p_{41}$ | $p_{17}$ | $p_{33}$ | $p_{22}$ | $p_{12}$ | $p_{39}$ | $p_{47}$ |
| $p_{11}$ | $p_{11}$ | $p_{42}$ | $p_{18}$ | $p_{34}$ | $p_{20}$ | $p_{14}$ | $p_{36}$ | $p_{44}$ |
| $p_{12}$ | $p_{12}$ | $p_{47}$ | $p_{22}$ | $p_{39}$ | $p_{17}$ | $p_{10}$ | $p_{33}$ | $p_{41}$ |
| $p_{13}$ | $p_{13}$ | $p_{46}$ | $p_{23}$ | $p_{38}$ | $p_{19}$ | $p_8$ | $p_{35}$ | $p_{43}$ |
| $p_{14}$ | $p_{14}$ | $p_{44}$ | $p_{20}$ | $p_{36}$ | $p_{18}$ | $p_{11}$ | $p_{34}$ | $p_{42}$ |
| $p_{15}$ | $p_{15}$ | $p_{45}$ | $p_{21}$ | $p_{37}$ | $p_{16}$ | $p_9$ | $p_{32}$ | $p_{40}$ |

$P_{22}=$

|  | $p_8$ | $p_9$ | $p_{10}$ | $p_{11}$ | $p_{12}$ | $p_{13}$ | $p_{14}$ | $p_{15}$ |
|---|---|---|---|---|---|---|---|---|
| $p_8$ | $p_{26}$ | $p_{40}$ | $p_2$ | $p_{34}$ | $p_4$ | $p_{29}$ | $p_{44}$ | $p_{37}$ |
| $p_9$ | $p_{27}$ | $p_{41}$ | $p_3$ | $p_{35}$ | $p_7$ | $p_{31}$ | $p_{46}$ | $p_{39}$ |
| $p_{10}$ | $p_{24}$ | $p_{42}$ | $p_0$ | $p_{32}$ | $p_5$ | $p_{28}$ | $p_{45}$ | $p_{36}$ |
| $p_{11}$ | $p_{25}$ | $p_{43}$ | $p_1$ | $p_{33}$ | $p_6$ | $p_{30}$ | $p_{47}$ | $p_{38}$ |
| $p_{12}$ | $p_{28}$ | $p_{45}$ | $p_5$ | $p_{36}$ | $p_0$ | $p_{24}$ | $p_{42}$ | $p_{32}$ |
| $p_{13}$ | $p_{29}$ | $p_{44}$ | $p_4$ | $p_{37}$ | $p_2$ | $p_{26}$ | $p_{40}$ | $p_{34}$ |
| $p_{14}$ | $p_{30}$ | $p_{47}$ | $p_6$ | $p_{38}$ | $p_1$ | $p_{25}$ | $p_{43}$ | $p_{33}$ |
| $p_{15}$ | $p_{31}$ | $p_{46}$ | $p_7$ | $p_{39}$ | $p_3$ | $p_{27}$ | $p_{41}$ | $p_{35}$ |

$P_{23}=$

|  | $p_{16}$ | $p_{17}$ | $p_{18}$ | $p_{19}$ | $p_{20}$ | $p_{21}$ | $p_{22}$ | $p_{23}$ |
|---|---|---|---|---|---|---|---|---|
| $p_8$ | $p_{32}$ | $p_0$ | $p_{42}$ | $p_{24}$ | $p_{36}$ | $p_{45}$ | $p_5$ | $p_{28}$ |
| $p_9$ | $p_{33}$ | $p_1$ | $p_{43}$ | $p_{25}$ | $p_{38}$ | $p_{47}$ | $p_6$ | $p_{30}$ |
| $p_{10}$ | $p_{34}$ | $p_2$ | $p_{40}$ | $p_{26}$ | $p_{37}$ | $p_{44}$ | $p_4$ | $p_{29}$ |
| $p_{11}$ | $p_{35}$ | $p_3$ | $p_{41}$ | $p_{27}$ | $p_{39}$ | $p_{46}$ | $p_7$ | $p_{31}$ |
| $p_{12}$ | $p_{37}$ | $p_4$ | $p_{44}$ | $p_{29}$ | $p_{34}$ | $p_{40}$ | $p_2$ | $p_{26}$ |
| $p_{13}$ | $p_{36}$ | $p_5$ | $p_{45}$ | $p_{28}$ | $p_{32}$ | $p_{42}$ | $p_0$ | $p_{24}$ |
| $p_{14}$ | $p_{39}$ | $p_7$ | $p_{46}$ | $p_{31}$ | $p_{35}$ | $p_{41}$ | $p_3$ | $p_{27}$ |
| $p_{15}$ | $p_{38}$ | $p_6$ | $p_{47}$ | $p_{30}$ | $p_{33}$ | $p_{43}$ | $p_1$ | $p_{25}$ |

$P_{24}=$

|  | $p_{24}$ | $p_{25}$ | $p_{26}$ | $p_{27}$ | $p_{28}$ | $p_{29}$ | $p_{30}$ | $p_{31}$ |
|---|---|---|---|---|---|---|---|---|
| $p_8$ | $p_{10}$ | $p_{33}$ | $p_{17}$ | $p_{41}$ | $p_{12}$ | $p_{22}$ | $p_{47}$ | $p_{39}$ |
| $p_9$ | $p_{11}$ | $p_{34}$ | $p_{18}$ | $p_{42}$ | $p_{14}$ | $p_{20}$ | $p_{44}$ | $p_{36}$ |
| $p_{10}$ | $p_8$ | $p_{35}$ | $p_{19}$ | $p_{43}$ | $p_{13}$ | $p_{23}$ | $p_{46}$ | $p_{38}$ |
| $p_{11}$ | $p_9$ | $p_{32}$ | $p_{16}$ | $p_{40}$ | $p_{15}$ | $p_{21}$ | $p_{45}$ | $p_{37}$ |
| $p_{12}$ | $p_{13}$ | $p_{38}$ | $p_{23}$ | $p_{46}$ | $p_8$ | $p_{19}$ | $p_{43}$ | $p_{35}$ |
| $p_{13}$ | $p_{12}$ | $p_{39}$ | $p_{22}$ | $p_{47}$ | $p_{10}$ | $p_{17}$ | $p_{41}$ | $p_{33}$ |
| $p_{14}$ | $p_{15}$ | $p_{37}$ | $p_{21}$ | $p_{45}$ | $p_9$ | $p_{16}$ | $p_{40}$ | $p_{32}$ |
| $p_{15}$ | $p_{14}$ | $p_{36}$ | $p_{20}$ | $p_{44}$ | $p_{11}$ | $p_{18}$ | $p_{42}$ | $p_{34}$ |

$P_{25}=$

|  | $p_{32}$ | $p_{33}$ | $p_{34}$ | $p_{35}$ | $p_{36}$ | $p_{37}$ | $p_{38}$ | $p_{39}$ |
|---|---|---|---|---|---|---|---|---|
| $p_8$ | $p_9$ | $p_1$ | $p_{18}$ | $p_{27}$ | $p_{14}$ | $p_{21}$ | $p_{30}$ | $p_7$ |
| $p_9$ | $p_{10}$ | $p_2$ | $p_{19}$ | $p_{24}$ | $p_{13}$ | $p_{22}$ | $p_{29}$ | $p_5$ |
| $p_{10}$ | $p_{11}$ | $p_3$ | $p_{16}$ | $p_{25}$ | $p_{15}$ | $p_{20}$ | $p_{31}$ | $p_6$ |
| $p_{11}$ | $p_8$ | $p_0$ | $p_{17}$ | $p_{26}$ | $p_{12}$ | $p_{23}$ | $p_{28}$ | $p_4$ |
| $p_{12}$ | $p_{15}$ | $p_6$ | $p_{20}$ | $p_{31}$ | $p_{11}$ | $p_{16}$ | $p_{25}$ | $p_3$ |
| $p_{13}$ | $p_{14}$ | $p_7$ | $p_{21}$ | $p_{30}$ | $p_9$ | $p_{18}$ | $p_{27}$ | $p_1$ |
| $p_{14}$ | $p_{12}$ | $p_4$ | $p_{23}$ | $p_{28}$ | $p_8$ | $p_{17}$ | $p_{26}$ | $p_0$ |
| $p_{15}$ | $p_{13}$ | $p_5$ | $p_{22}$ | $p_{29}$ | $p_{10}$ | $p_{19}$ | $p_{24}$ | $p_2$ |

$P_{26}=$

|  | $p_{40}$ | $p_{41}$ | $p_{42}$ | $p_{43}$ | $p_{44}$ | $p_{45}$ | $p_{46}$ | $p_{47}$ |
|---|---|---|---|---|---|---|---|---|
| $p_8$ | $p_{16}$ | $p_3$ | $p_{11}$ | $p_{25}$ | $p_{20}$ | $p_{15}$ | $p_{31}$ | $p_6$ |
| $p_9$ | $p_{17}$ | $p_0$ | $p_8$ | $p_{26}$ | $p_{23}$ | $p_{12}$ | $p_{28}$ | $p_4$ |
| $p_{10}$ | $p_{18}$ | $p_1$ | $p_9$ | $p_{27}$ | $p_{21}$ | $p_{14}$ | $p_{30}$ | $p_7$ |
| $p_{11}$ | $p_{19}$ | $p_2$ | $p_{10}$ | $p_{24}$ | $p_{22}$ | $p_{13}$ | $p_{29}$ | $p_5$ |
| $p_{12}$ | $p_{21}$ | $p_7$ | $p_{14}$ | $p_{30}$ | $p_{18}$ | $p_9$ | $p_{27}$ | $p_1$ |
| $p_{13}$ | $p_{20}$ | $p_6$ | $p_{15}$ | $p_{31}$ | $p_{16}$ | $p_{11}$ | $p_{25}$ | $p_3$ |
| $p_{14}$ | $p_{22}$ | $p_5$ | $p_{13}$ | $p_{29}$ | $p_{19}$ | $p_{10}$ | $p_{24}$ | $p_2$ |
| $p_{15}$ | $p_{23}$ | $p_4$ | $p_{12}$ | $p_{28}$ | $p_{17}$ | $p_8$ | $p_{26}$ | $p_0$ |

$P_{31}=$

|  | $p_0$ | $p_1$ | $p_2$ | $p_3$ | $p_4$ | $p_5$ | $p_6$ | $p_7$ |
|---|---|---|---|---|---|---|---|---|
| $p_{16}$ | $p_{16}$ | $p_{32}$ | $p_9$ | $p_{40}$ | $p_{15}$ | $p_{21}$ | $p_{45}$ | $p_{37}$ |
| $p_{17}$ | $p_{17}$ | $p_{33}$ | $p_{10}$ | $p_{41}$ | $p_{12}$ | $p_{22}$ | $p_{47}$ | $p_{39}$ |
| $p_{18}$ | $p_{18}$ | $p_{34}$ | $p_{11}$ | $p_{42}$ | $p_{14}$ | $p_{20}$ | $p_{44}$ | $p_{36}$ |
| $p_{19}$ | $p_{19}$ | $p_{35}$ | $p_8$ | $p_{43}$ | $p_{13}$ | $p_{23}$ | $p_{46}$ | $p_{38}$ |
| $p_{20}$ | $p_{20}$ | $p_{36}$ | $p_{14}$ | $p_{44}$ | $p_{11}$ | $p_{18}$ | $p_{42}$ | $p_{34}$ |
| $p_{21}$ | $p_{21}$ | $p_{37}$ | $p_{15}$ | $p_{45}$ | $p_9$ | $p_{16}$ | $p_{40}$ | $p_{32}$ |
| $p_{22}$ | $p_{22}$ | $p_{39}$ | $p_{12}$ | $p_{47}$ | $p_{10}$ | $p_{17}$ | $p_{41}$ | $p_{33}$ |
| $p_{23}$ | $p_{23}$ | $p_{38}$ | $p_{13}$ | $p_{46}$ | $p_8$ | $p_{19}$ | $p_{43}$ | $p_{35}$ |

$P_{32}=$

|  | $p_8$ | $p_9$ | $p_{10}$ | $p_{11}$ | $p_{12}$ | $p_{13}$ | $p_{14}$ | $p_{15}$ |
|---|---|---|---|---|---|---|---|---|
| $p_{16}$ | $p_3$ | $p_{35}$ | $p_{27}$ | $p_{41}$ | $p_{31}$ | $p_7$ | $p_{39}$ | $p_{46}$ |
| $p_{17}$ | $p_0$ | $p_{32}$ | $p_{24}$ | $p_{42}$ | $p_{28}$ | $p_5$ | $p_{36}$ | $p_{45}$ |
| $p_{18}$ | $p_1$ | $p_{33}$ | $p_{25}$ | $p_{43}$ | $p_{30}$ | $p_6$ | $p_{38}$ | $p_{47}$ |
| $p_{19}$ | $p_2$ | $p_{34}$ | $p_{26}$ | $p_{40}$ | $p_{29}$ | $p_4$ | $p_{37}$ | $p_{44}$ |
| $p_{20}$ | $p_6$ | $p_{38}$ | $p_{30}$ | $p_{47}$ | $p_{25}$ | $p_1$ | $p_{33}$ | $p_{43}$ |
| $p_{21}$ | $p_7$ | $p_{39}$ | $p_{31}$ | $p_{46}$ | $p_{27}$ | $p_3$ | $p_{35}$ | $p_{41}$ |
| $p_{22}$ | $p_5$ | $p_{36}$ | $p_{28}$ | $p_{45}$ | $p_{24}$ | $p_0$ | $p_{32}$ | $p_{42}$ |
| $p_{23}$ | $p_4$ | $p_{37}$ | $p_{29}$ | $p_{44}$ | $p_{26}$ | $p_2$ | $p_{34}$ | $p_{40}$ |

$P_{33}=$

|  | $p_{16}$ | $p_{17}$ | $p_{18}$ | $p_{19}$ | $p_{20}$ | $p_{21}$ | $p_{22}$ | $p_{23}$ |
|---|---|---|---|---|---|---|---|---|
| $p_{16}$ | $p_{43}$ | $p_{25}$ | $p_{33}$ | $p_{1}$ | $p_{47}$ | $p_{38}$ | $p_{30}$ | $p_{6}$ |
| $p_{17}$ | $p_{40}$ | $p_{26}$ | $p_{34}$ | $p_{2}$ | $p_{44}$ | $p_{37}$ | $p_{29}$ | $p_{4}$ |
| $p_{18}$ | $p_{41}$ | $p_{27}$ | $p_{35}$ | $p_{3}$ | $p_{46}$ | $p_{39}$ | $p_{31}$ | $p_{7}$ |
| $p_{19}$ | $p_{42}$ | $p_{24}$ | $p_{32}$ | $p_{0}$ | $p_{45}$ | $p_{36}$ | $p_{28}$ | $p_{5}$ |
| $p_{20}$ | $p_{46}$ | $p_{31}$ | $p_{39}$ | $p_{7}$ | $p_{41}$ | $p_{35}$ | $p_{27}$ | $p_{3}$ |
| $p_{21}$ | $p_{47}$ | $p_{30}$ | $p_{38}$ | $p_{6}$ | $p_{43}$ | $p_{33}$ | $p_{25}$ | $p_{1}$ |
| $p_{22}$ | $p_{44}$ | $p_{29}$ | $p_{37}$ | $p_{4}$ | $p_{40}$ | $p_{34}$ | $p_{26}$ | $p_{2}$ |
| $p_{23}$ | $p_{45}$ | $p_{28}$ | $p_{36}$ | $p_{5}$ | $p_{42}$ | $p_{32}$ | $p_{24}$ | $p_{0}$ |

$P_{34}=$

|  | $p_{24}$ | $p_{25}$ | $p_{26}$ | $p_{27}$ | $p_{28}$ | $p_{29}$ | $p_{30}$ | $p_{31}$ |
|---|---|---|---|---|---|---|---|---|
| $p_{16}$ | $p_{18}$ | $p_{42}$ | $p_{11}$ | $p_{34}$ | $p_{20}$ | $p_{14}$ | $p_{36}$ | $p_{44}$ |
| $p_{17}$ | $p_{19}$ | $p_{43}$ | $p_{8}$ | $p_{35}$ | $p_{23}$ | $p_{13}$ | $p_{38}$ | $p_{46}$ |
| $p_{18}$ | $p_{16}$ | $p_{40}$ | $p_{9}$ | $p_{32}$ | $p_{21}$ | $p_{15}$ | $p_{37}$ | $p_{45}$ |
| $p_{19}$ | $p_{17}$ | $p_{41}$ | $p_{10}$ | $p_{33}$ | $p_{22}$ | $p_{12}$ | $p_{39}$ | $p_{47}$ |
| $p_{20}$ | $p_{21}$ | $p_{45}$ | $p_{15}$ | $p_{37}$ | $p_{16}$ | $p_{9}$ | $p_{32}$ | $p_{40}$ |
| $p_{21}$ | $p_{20}$ | $p_{44}$ | $p_{14}$ | $p_{36}$ | $p_{18}$ | $p_{11}$ | $p_{34}$ | $p_{42}$ |
| $p_{22}$ | $p_{23}$ | $p_{46}$ | $p_{13}$ | $p_{38}$ | $p_{19}$ | $p_{8}$ | $p_{35}$ | $p_{43}$ |
| $p_{23}$ | $p_{22}$ | $p_{47}$ | $p_{12}$ | $p_{39}$ | $p_{17}$ | $p_{10}$ | $p_{33}$ | $p_{41}$ |

$P_{35}=$

|  | $p_{32}$ | $p_{33}$ | $p_{34}$ | $p_{35}$ | $p_{36}$ | $p_{37}$ | $p_{38}$ | $p_{39}$ |
|---|---|---|---|---|---|---|---|---|
| $p_{16}$ | $p_{19}$ | $p_{24}$ | $p_{10}$ | $p_{2}$ | $p_{22}$ | $p_{13}$ | $p_{5}$ | $p_{29}$ |
| $p_{17}$ | $p_{16}$ | $p_{25}$ | $p_{11}$ | $p_{3}$ | $p_{20}$ | $p_{15}$ | $p_{6}$ | $p_{31}$ |
| $p_{18}$ | $p_{17}$ | $p_{26}$ | $p_{8}$ | $p_{0}$ | $p_{23}$ | $p_{12}$ | $p_{4}$ | $p_{28}$ |
| $p_{19}$ | $p_{18}$ | $p_{27}$ | $p_{9}$ | $p_{1}$ | $p_{21}$ | $p_{14}$ | $p_{7}$ | $p_{30}$ |
| $p_{20}$ | $p_{23}$ | $p_{28}$ | $p_{12}$ | $p_{4}$ | $p_{17}$ | $p_{8}$ | $p_{0}$ | $p_{26}$ |
| $p_{21}$ | $p_{22}$ | $p_{29}$ | $p_{13}$ | $p_{5}$ | $p_{19}$ | $p_{10}$ | $p_{2}$ | $p_{24}$ |
| $p_{22}$ | $p_{20}$ | $p_{31}$ | $p_{15}$ | $p_{6}$ | $p_{16}$ | $p_{11}$ | $p_{3}$ | $p_{25}$ |
| $p_{23}$ | $p_{21}$ | $p_{30}$ | $p_{14}$ | $p_{7}$ | $p_{18}$ | $p_{9}$ | $p_{1}$ | $p_{27}$ |

$P_{36}=$

|  | $p_{40}$ | $p_{41}$ | $p_{42}$ | $p_{43}$ | $p_{44}$ | $p_{45}$ | $p_{46}$ | $p_{47}$ |
|---|---|---|---|---|---|---|---|---|
| $p_{16}$ | $p_{8}$ | $p_{26}$ | $p_{17}$ | $p_{0}$ | $p_{12}$ | $p_{23}$ | $p_{4}$ | $p_{28}$ |
| $p_{17}$ | $p_{9}$ | $p_{27}$ | $p_{18}$ | $p_{1}$ | $p_{14}$ | $p_{21}$ | $p_{7}$ | $p_{30}$ |
| $p_{18}$ | $p_{10}$ | $p_{24}$ | $p_{19}$ | $p_{2}$ | $p_{13}$ | $p_{22}$ | $p_{5}$ | $p_{29}$ |
| $p_{19}$ | $p_{11}$ | $p_{25}$ | $p_{16}$ | $p_{3}$ | $p_{15}$ | $p_{20}$ | $p_{6}$ | $p_{31}$ |
| $p_{20}$ | $p_{13}$ | $p_{29}$ | $p_{22}$ | $p_{5}$ | $p_{10}$ | $p_{19}$ | $p_{2}$ | $p_{24}$ |
| $p_{21}$ | $p_{12}$ | $p_{28}$ | $p_{23}$ | $p_{4}$ | $p_{8}$ | $p_{17}$ | $p_{0}$ | $p_{26}$ |
| $p_{22}$ | $p_{14}$ | $p_{30}$ | $p_{21}$ | $p_{7}$ | $p_{9}$ | $p_{18}$ | $p_{1}$ | $p_{27}$ |
| $p_{23}$ | $p_{15}$ | $p_{31}$ | $p_{20}$ | $p_{6}$ | $p_{11}$ | $p_{16}$ | $p_{3}$ | $p_{25}$ |

$P_{41}=$

|  | $p_{0}$ | $p_{1}$ | $p_{2}$ | $p_{3}$ | $p_{4}$ | $p_{5}$ | $p_{6}$ | $p_{7}$ |
|---|---|---|---|---|---|---|---|---|
| $p_{24}$ | $p_{24}$ | $p_{27}$ | $p_{26}$ | $p_{25}$ | $p_{29}$ | $p_{28}$ | $p_{31}$ | $p_{30}$ |
| $p_{25}$ | $p_{25}$ | $p_{24}$ | $p_{27}$ | $p_{26}$ | $p_{31}$ | $p_{30}$ | $p_{28}$ | $p_{29}$ |
| $p_{26}$ | $p_{26}$ | $p_{25}$ | $p_{24}$ | $p_{27}$ | $p_{28}$ | $p_{29}$ | $p_{30}$ | $p_{31}$ |
| $p_{27}$ | $p_{27}$ | $p_{26}$ | $p_{25}$ | $p_{24}$ | $p_{30}$ | $p_{31}$ | $p_{29}$ | $p_{28}$ |
| $p_{28}$ | $p_{28}$ | $p_{30}$ | $p_{29}$ | $p_{31}$ | $p_{26}$ | $p_{24}$ | $p_{25}$ | $p_{27}$ |
| $p_{29}$ | $p_{29}$ | $p_{31}$ | $p_{28}$ | $p_{30}$ | $p_{24}$ | $p_{26}$ | $p_{27}$ | $p_{25}$ |
| $p_{30}$ | $p_{30}$ | $p_{29}$ | $p_{31}$ | $p_{28}$ | $p_{27}$ | $p_{25}$ | $p_{26}$ | $p_{24}$ |
| $p_{31}$ | $p_{31}$ | $p_{28}$ | $p_{30}$ | $p_{29}$ | $p_{25}$ | $p_{27}$ | $p_{24}$ | $p_{26}$ |

$$P_{42} = \begin{array}{c|cccccccc} & p_8 & p_9 & p_{10} & p_{11} & p_{12} & p_{13} & p_{14} & p_{15} \\ \hline p_{24} & p_{19} & p_{18} & p_{17} & p_{16} & p_{22} & p_{23} & p_{21} & p_{20} \\ p_{25} & p_{16} & p_{19} & p_{18} & p_{17} & p_{20} & p_{21} & p_{22} & p_{23} \\ p_{26} & p_{17} & p_{16} & p_{19} & p_{18} & p_{23} & p_{22} & p_{20} & p_{21} \\ p_{27} & p_{18} & p_{17} & p_{16} & p_{19} & p_{21} & p_{20} & p_{23} & p_{22} \\ p_{28} & p_{23} & p_{21} & p_{22} & p_{20} & p_{17} & p_{19} & p_{18} & p_{16} \\ p_{29} & p_{22} & p_{20} & p_{23} & p_{21} & p_{19} & p_{17} & p_{16} & p_{18} \\ p_{30} & p_{21} & p_{22} & p_{20} & p_{23} & p_{18} & p_{16} & p_{19} & p_{17} \\ p_{31} & p_{20} & p_{23} & p_{21} & p_{22} & p_{16} & p_{18} & p_{17} & p_{19} \end{array}$$

$$P_{43} = \begin{array}{c|cccccccc} & p_{16} & p_{17} & p_{18} & p_{19} & p_{20} & p_{21} & p_{22} & p_{23} \\ \hline p_{24} & p_{11} & p_{10} & p_9 & p_8 & p_{15} & p_{14} & p_{12} & p_{13} \\ p_{25} & p_8 & p_{11} & p_{10} & p_9 & p_{12} & p_{13} & p_{14} & p_{15} \\ p_{26} & p_9 & p_8 & p_{11} & p_{10} & p_{14} & p_{15} & p_{13} & p_{12} \\ p_{27} & p_{10} & p_9 & p_8 & p_{11} & p_{13} & p_{12} & p_{15} & p_{14} \\ p_{28} & p_{15} & p_{12} & p_{14} & p_{13} & p_{11} & p_9 & p_{10} & p_8 \\ p_{29} & p_{14} & p_{13} & p_{15} & p_{12} & p_9 & p_{11} & p_8 & p_{10} \\ p_{30} & p_{12} & p_{14} & p_{13} & p_{15} & p_8 & p_{10} & p_{11} & p_9 \\ p_{31} & p_{13} & p_{15} & p_{12} & p_{14} & p_{10} & p_8 & p_9 & p_{11} \end{array}$$

$$P_{44} = \begin{array}{c|cccccccc} & p_{24} & p_{25} & p_{26} & p_{27} & p_{28} & p_{29} & p_{30} & p_{31} \\ \hline p_{24} & p_0 & p_3 & p_2 & p_1 & p_5 & p_4 & p_7 & p_6 \\ p_{25} & p_1 & p_0 & p_3 & p_2 & p_6 & p_7 & p_5 & p_4 \\ p_{26} & p_2 & p_1 & p_0 & p_3 & p_4 & p_5 & p_6 & p_7 \\ p_{27} & p_3 & p_2 & p_1 & p_0 & p_7 & p_6 & p_4 & p_5 \\ p_{28} & p_5 & p_6 & p_4 & p_7 & p_0 & p_2 & p_1 & p_3 \\ p_{29} & p_4 & p_7 & p_5 & p_6 & p_2 & p_0 & p_3 & p_1 \\ p_{30} & p_6 & p_4 & p_7 & p_5 & p_1 & p_3 & p_2 & p_0 \\ p_{31} & p_7 & p_5 & p_6 & p_4 & p_3 & p_1 & p_0 & p_2 \end{array}$$

$$P_{45} = \begin{array}{c|cccccccc} & p_{32} & p_{33} & p_{34} & p_{35} & p_{36} & p_{37} & p_{38} & p_{39} \\ \hline p_{24} & p_{42} & p_{41} & p_{40} & p_{43} & p_{45} & p_{44} & p_{46} & p_{47} \\ p_{25} & p_{43} & p_{42} & p_{41} & p_{40} & p_{47} & p_{46} & p_{45} & p_{44} \\ p_{26} & p_{40} & p_{43} & p_{42} & p_{41} & p_{44} & p_{45} & p_{47} & p_{46} \\ p_{27} & p_{41} & p_{40} & p_{43} & p_{42} & p_{42} & p_{40} & p_{43} & p_{41} \\ p_{28} & p_{45} & p_{47} & p_{44} & p_{46} & p_{17} & p_8 & p_0 & p_{26} \\ p_{29} & p_{44} & p_{46} & p_{45} & p_{47} & p_{40} & p_{42} & p_{41} & p_{43} \\ p_{30} & p_{47} & p_{44} & p_{46} & p_{45} & p_{43} & p_{41} & p_{40} & p_{42} \\ p_{31} & p_{46} & p_{45} & p_{47} & p_{44} & p_{41} & p_{43} & p_{42} & p_{40} \end{array}$$

$$P_{46} = \begin{array}{c|cccccccc} & p_{40} & p_{41} & p_{42} & p_{43} & p_{44} & p_{45} & p_{46} & p_{47} \\ \hline p_{24} & p_{34} & p_{33} & p_{32} & p_{35} & p_{37} & p_{36} & p_{38} & p_{39} \\ p_{25} & p_{35} & p_{34} & p_{33} & p_{32} & p_{39} & p_{38} & p_{37} & p_{36} \\ p_{26} & p_{32} & p_{35} & p_{34} & p_{33} & p_{36} & p_{37} & p_{39} & p_{38} \\ p_{27} & p_{33} & p_{32} & p_{35} & p_{34} & p_{38} & p_{39} & p_{36} & p_{37} \\ p_{28} & p_{37} & p_{39} & p_{36} & p_{38} & p_{34} & p_{32} & p_{35} & p_{33} \\ p_{29} & p_{36} & p_{38} & p_{37} & p_{39} & p_{32} & p_{34} & p_{33} & p_{35} \\ p_{30} & p_{39} & p_{36} & p_{38} & p_{37} & p_{35} & p_{33} & p_{32} & p_{34} \\ p_{31} & p_{38} & p_{37} & p_{39} & p_{36} & p_{33} & p_{35} & p_{34} & p_{32} \end{array}$$

$$P_{51} = \begin{array}{c|cccccccc} & p_0 & p_1 & p_2 & p_3 & p_4 & p_5 & p_6 & p_7 \\ \hline p_{32} & p_{32} & p_9 & p_{40} & p_{16} & p_{37} & p_{45} & p_{15} & p_{21} \\ p_{33} & p_{33} & p_{10} & p_{41} & p_{17} & p_{39} & p_{47} & p_{12} & p_{22} \\ p_{34} & p_{34} & p_{11} & p_{42} & p_{18} & p_{36} & p_{44} & p_{14} & p_{20} \\ p_{35} & p_{35} & p_8 & p_{43} & p_{19} & p_{38} & p_{46} & p_{13} & p_{23} \\ p_{36} & p_{36} & p_{14} & p_{44} & p_{20} & p_{34} & p_{42} & p_{11} & p_{18} \\ p_{37} & p_{37} & p_{15} & p_{45} & p_{21} & p_{32} & p_{40} & p_9 & p_{16} \\ p_{38} & p_{38} & p_{13} & p_{46} & p_{23} & p_{35} & p_{43} & p_8 & p_{19} \\ p_{39} & p_{39} & p_{12} & p_{47} & p_{22} & p_{33} & p_{41} & p_{10} & p_{17} \end{array}$$

$$P_{52} = \begin{array}{c|cccccccc} & p_8 & p_9 & p_{10} & p_{11} & p_{12} & p_{13} & p_{14} & p_{15} \\ \hline p_{32} & p_{35} & p_{27} & p_{41} & p_3 & p_{39} & p_{46} & p_7 & p_{31} \\ p_{33} & p_{32} & p_{24} & p_{42} & p_0 & p_{36} & p_{45} & p_5 & p_{28} \\ p_{34} & p_{33} & p_{25} & p_{43} & p_1 & p_{38} & p_{47} & p_6 & p_{30} \\ p_{35} & p_{34} & p_{26} & p_{40} & p_2 & p_{37} & p_{44} & p_4 & p_{29} \\ p_{36} & p_{38} & p_{30} & p_{47} & p_6 & p_{33} & p_{43} & p_1 & p_{25} \\ p_{37} & p_{39} & p_{31} & p_{46} & p_7 & p_{35} & p_{41} & p_3 & p_{27} \\ p_{38} & p_{37} & p_{29} & p_{44} & p_4 & p_{34} & p_{40} & p_2 & p_{26} \\ p_{39} & p_{36} & p_{28} & p_{45} & p_5 & p_{32} & p_{42} & p_0 & p_{24} \end{array}$$

$$P_{53} = \begin{array}{c|cccccccc} & p_{16} & p_{17} & p_{18} & p_{19} & p_{20} & p_{21} & p_{22} & p_{23} \\ \hline p_{32} & p_{25} & p_{33} & p_1 & p_{43} & p_6 & p_{30} & p_{47} & p_{38} \\ p_{33} & p_{26} & p_{34} & p_2 & p_{40} & p_4 & p_{29} & p_{44} & p_{37} \\ p_{34} & p_{27} & p_{35} & p_3 & p_{41} & p_7 & p_{31} & p_{46} & p_{39} \\ p_{35} & p_{24} & p_{32} & p_0 & p_{42} & p_5 & p_{28} & p_{45} & p_{36} \\ p_{36} & p_{31} & p_{39} & p_7 & p_{46} & p_3 & p_{27} & p_{41} & p_{35} \\ p_{37} & p_{30} & p_{38} & p_6 & p_{47} & p_1 & p_{25} & p_{43} & p_{33} \\ p_{38} & p_{28} & p_{36} & p_5 & p_{45} & p_0 & p_{24} & p_{42} & p_{32} \\ p_{39} & p_{29} & p_{37} & p_4 & p_{44} & p_2 & p_{26} & p_{40} & p_{34} \end{array}$$

$$P_{54} = \begin{array}{c|cccccccc} & p_{24} & p_{25} & p_{26} & p_{27} & p_{28} & p_{29} & p_{30} & p_{31} \\ \hline p_{32} & p_{42} & p_{11} & p_{34} & p_{18} & p_{36} & p_{44} & p_{14} & p_{20} \\ p_{33} & p_{43} & p_8 & p_{35} & p_{19} & p_{38} & p_{46} & p_{13} & p_{23} \\ p_{34} & p_{40} & p_9 & p_{32} & p_{16} & p_{37} & p_{45} & p_{15} & p_{21} \\ p_{35} & p_{41} & p_{10} & p_{33} & p_{17} & p_{39} & p_{47} & p_{12} & p_{22} \\ p_{36} & p_{45} & p_{15} & p_{37} & p_{21} & p_{32} & p_{40} & p_9 & p_{16} \\ p_{37} & p_{44} & p_{14} & p_{36} & p_{20} & p_{34} & p_{42} & p_{11} & p_{18} \\ p_{38} & p_{47} & p_{12} & p_{39} & p_{22} & p_{33} & p_{41} & p_{10} & p_{17} \\ p_{39} & p_{46} & p_{13} & p_{38} & p_{23} & p_{35} & p_{43} & p_8 & p_{19} \end{array}$$

$$P_{55} = \begin{array}{c|cccccccc} & p_{32} & p_{33} & p_{34} & p_{35} & p_{36} & p_{37} & p_{38} & p_{39} \\ \hline p_{32} & p_{24} & p_{10} & p_2 & p_{19} & p_5 & p_{29} & p_{13} & p_{22} \\ p_{33} & p_{25} & p_{11} & p_3 & p_{16} & p_6 & p_{31} & p_{15} & p_{20} \\ p_{34} & p_{26} & p_8 & p_0 & p_{17} & p_4 & p_{28} & p_{12} & p_{23} \\ p_{35} & p_{27} & p_9 & p_1 & p_{18} & p_7 & p_{30} & p_{14} & p_{21} \\ p_{36} & p_{28} & p_{12} & p_4 & p_{23} & p_0 & p_{26} & p_8 & p_{17} \\ p_{37} & p_{29} & p_{13} & p_5 & p_{22} & p_2 & p_{24} & p_{10} & p_{19} \\ p_{38} & p_{30} & p_{14} & p_7 & p_{21} & p_1 & p_{27} & p_9 & p_{18} \\ p_{39} & p_{31} & p_{15} & p_6 & p_{20} & p_3 & p_{25} & p_{11} & p_{16} \end{array}$$

$P_{56} =$

|  | $p_{40}$ | $p_{41}$ | $p_{42}$ | $p_{43}$ | $p_{44}$ | $p_{45}$ | $p_{46}$ | $p_{47}$ |
|---|---|---|---|---|---|---|---|---|
| $p_{32}$ | $p_{26}$ | $p_{17}$ | $p_0$ | $p_8$ | $p_4$ | $p_{28}$ | $p_{23}$ | $p_{12}$ |
| $p_{33}$ | $p_{27}$ | $p_{18}$ | $p_1$ | $p_9$ | $p_7$ | $p_{30}$ | $p_{21}$ | $p_{14}$ |
| $p_{34}$ | $p_{24}$ | $p_{19}$ | $p_2$ | $p_{10}$ | $p_5$ | $p_{29}$ | $p_{22}$ | $p_{13}$ |
| $p_{35}$ | $p_{25}$ | $p_{16}$ | $p_3$ | $p_{11}$ | $p_6$ | $p_{31}$ | $p_{20}$ | $p_{15}$ |
| $p_{36}$ | $p_{29}$ | $p_{22}$ | $p_5$ | $p_{13}$ | $p_2$ | $p_{24}$ | $p_{19}$ | $p_{10}$ |
| $p_{37}$ | $p_{28}$ | $p_{23}$ | $p_4$ | $p_{12}$ | $p_0$ | $p_{26}$ | $p_{17}$ | $p_8$ |
| $p_{38}$ | $p_{31}$ | $p_{20}$ | $p_6$ | $p_{15}$ | $p_3$ | $p_{25}$ | $p_{16}$ | $p_{11}$ |
| $p_{39}$ | $p_{30}$ | $p_{21}$ | $p_7$ | $p_{14}$ | $p_1$ | $p_{27}$ | $p_{18}$ | $p_9$ |

$P_{61} =$

|  | $p_0$ | $p_1$ | $p_2$ | $p_3$ | $p_4$ | $p_5$ | $p_6$ | $p_7$ |
|---|---|---|---|---|---|---|---|---|
| $p_{40}$ | $p_{40}$ | $p_{16}$ | $p_{32}$ | $p_9$ | $p_{45}$ | $p_{37}$ | $p_{21}$ | $p_{15}$ |
| $p_{41}$ | $p_{41}$ | $p_{17}$ | $p_{33}$ | $p_{10}$ | $p_{47}$ | $p_{39}$ | $p_{22}$ | $p_{12}$ |
| $p_{42}$ | $p_{42}$ | $p_{18}$ | $p_{34}$ | $p_{11}$ | $p_{44}$ | $p_{36}$ | $p_{20}$ | $p_{14}$ |
| $p_{43}$ | $p_{43}$ | $p_{19}$ | $p_{35}$ | $p_8$ | $p_{46}$ | $p_{38}$ | $p_{23}$ | $p_{13}$ |
| $p_{44}$ | $p_{44}$ | $p_{20}$ | $p_{36}$ | $p_{14}$ | $p_{42}$ | $p_{34}$ | $p_{18}$ | $p_{11}$ |
| $p_{45}$ | $p_{45}$ | $p_{21}$ | $p_{37}$ | $p_{15}$ | $p_{40}$ | $p_{32}$ | $p_{16}$ | $p_9$ |
| $p_{46}$ | $p_{46}$ | $p_{23}$ | $p_{38}$ | $p_{13}$ | $p_{43}$ | $p_{35}$ | $p_{19}$ | $p_8$ |
| $p_{47}$ | $p_{47}$ | $p_{22}$ | $p_{39}$ | $p_{12}$ | $p_{41}$ | $p_{33}$ | $p_{17}$ | $p_{10}$ |

$P_{62} =$

|  | $p_8$ | $p_9$ | $p_{10}$ | $p_{11}$ | $p_{12}$ | $p_{13}$ | $p_{14}$ | $p_{15}$ |
|---|---|---|---|---|---|---|---|---|
| $p_{40}$ | $p_{41}$ | $p_3$ | $p_{35}$ | $p_{27}$ | $p_{46}$ | $p_{39}$ | $p_{31}$ | $p_7$ |
| $p_{41}$ | $p_{42}$ | $p_0$ | $p_{32}$ | $p_{24}$ | $p_{45}$ | $p_{36}$ | $p_{28}$ | $p_5$ |
| $p_{42}$ | $p_{43}$ | $p_1$ | $p_{33}$ | $p_{25}$ | $p_{47}$ | $p_{38}$ | $p_{30}$ | $p_6$ |
| $p_{43}$ | $p_{40}$ | $p_2$ | $p_{34}$ | $p_{26}$ | $p_{44}$ | $p_{37}$ | $p_{29}$ | $p_4$ |
| $p_{44}$ | $p_{47}$ | $p_6$ | $p_{38}$ | $p_{30}$ | $p_{43}$ | $p_{33}$ | $p_{25}$ | $p_1$ |
| $p_{45}$ | $p_{46}$ | $p_7$ | $p_{39}$ | $p_{31}$ | $p_{41}$ | $p_{35}$ | $p_{27}$ | $p_3$ |
| $p_{46}$ | $p_{44}$ | $p_4$ | $p_{37}$ | $p_{29}$ | $p_{40}$ | $p_{34}$ | $p_{26}$ | $p_2$ |
| $p_{47}$ | $p_{45}$ | $p_5$ | $p_{36}$ | $p_{28}$ | $p_{42}$ | $p_{32}$ | $p_{24}$ | $p_0$ |

$P_{63} =$

|  | $p_{16}$ | $p_{17}$ | $p_{18}$ | $p_{19}$ | $p_{20}$ | $p_{21}$ | $p_{22}$ | $p_{23}$ |
|---|---|---|---|---|---|---|---|---|
| $p_{40}$ | $p_1$ | $p_{43}$ | $p_{25}$ | $p_{33}$ | $p_{30}$ | $p_6$ | $p_{38}$ | $p_{47}$ |
| $p_{41}$ | $p_2$ | $p_{40}$ | $p_{26}$ | $p_{34}$ | $p_{29}$ | $p_4$ | $p_{37}$ | $p_{44}$ |
| $p_{42}$ | $p_3$ | $p_{41}$ | $p_{27}$ | $p_{35}$ | $p_{31}$ | $p_7$ | $p_{39}$ | $p_{46}$ |
| $p_{43}$ | $p_0$ | $p_{42}$ | $p_{24}$ | $p_{32}$ | $p_{28}$ | $p_5$ | $p_{36}$ | $p_{45}$ |
| $p_{44}$ | $p_7$ | $p_{46}$ | $p_{31}$ | $p_{39}$ | $p_{27}$ | $p_3$ | $p_{35}$ | $p_{41}$ |
| $p_{45}$ | $p_6$ | $p_{47}$ | $p_{30}$ | $p_{38}$ | $p_{25}$ | $p_1$ | $p_{33}$ | $p_{43}$ |
| $p_{46}$ | $p_5$ | $p_{45}$ | $p_{28}$ | $p_{36}$ | $p_{24}$ | $p_0$ | $p_{32}$ | $p_{42}$ |
| $p_{47}$ | $p_4$ | $p_{44}$ | $p_{29}$ | $p_{37}$ | $p_{26}$ | $p_2$ | $p_{34}$ | $p_{40}$ |

$P_{64} =$

|  | $p_{24}$ | $p_{25}$ | $p_{26}$ | $p_{27}$ | $p_{28}$ | $p_{29}$ | $p_{30}$ | $p_{31}$ |
|---|---|---|---|---|---|---|---|---|
| $p_{40}$ | $p_{34}$ | $p_{18}$ | $p_{42}$ | $p_{11}$ | $p_{44}$ | $p_{36}$ | $p_{20}$ | $p_{14}$ |
| $p_{41}$ | $p_{35}$ | $p_{19}$ | $p_{43}$ | $p_8$ | $p_{46}$ | $p_{38}$ | $p_{23}$ | $p_{13}$ |
| $p_{42}$ | $p_{32}$ | $p_{16}$ | $p_{40}$ | $p_9$ | $p_{45}$ | $p_{37}$ | $p_{21}$ | $p_{15}$ |
| $p_{43}$ | $p_{33}$ | $p_{17}$ | $p_{41}$ | $p_{10}$ | $p_{47}$ | $p_{39}$ | $p_{22}$ | $p_{12}$ |
| $p_{44}$ | $p_{37}$ | $p_{21}$ | $p_{45}$ | $p_{15}$ | $p_{40}$ | $p_{32}$ | $p_{16}$ | $p_9$ |
| $p_{45}$ | $p_{36}$ | $p_{20}$ | $p_{44}$ | $p_{14}$ | $p_{42}$ | $p_{34}$ | $p_{18}$ | $p_{11}$ |
| $p_{46}$ | $p_{39}$ | $p_{22}$ | $p_{47}$ | $p_{12}$ | $p_{41}$ | $p_{33}$ | $p_{17}$ | $p_{10}$ |
| $p_{47}$ | $p_{38}$ | $p_{23}$ | $p_{46}$ | $p_{13}$ | $p_{43}$ | $p_{35}$ | $p_{19}$ | $p_8$ |

$P_{65}=$

|     | $p_{32}$ | $p_{33}$ | $p_{34}$ | $p_{35}$ | $p_{36}$ | $p_{37}$ | $p_{38}$ | $p_{39}$ |
|---|---|---|---|---|---|---|---|---|
| $p_{40}$ | $p_2$ | $p_{19}$ | $p_{24}$ | $p_{10}$ | $p_{29}$ | $p_5$ | $p_{22}$ | $p_{13}$ |
| $p_{41}$ | $p_3$ | $p_{16}$ | $p_{25}$ | $p_{11}$ | $p_{31}$ | $p_6$ | $p_{20}$ | $p_{15}$ |
| $p_{42}$ | $p_0$ | $p_{17}$ | $p_{26}$ | $p_8$ | $p_{28}$ | $p_4$ | $p_{23}$ | $p_{12}$ |
| $p_{43}$ | $p_1$ | $p_{18}$ | $p_{27}$ | $p_9$ | $p_{30}$ | $p_7$ | $p_{21}$ | $p_{14}$ |
| $p_{44}$ | $p_4$ | $p_{23}$ | $p_{28}$ | $p_{12}$ | $p_{26}$ | $p_0$ | $p_{17}$ | $p_8$ |
| $p_{45}$ | $p_5$ | $p_{22}$ | $p_{29}$ | $p_{13}$ | $p_{24}$ | $p_2$ | $p_{19}$ | $p_{10}$ |
| $p_{46}$ | $p_7$ | $p_{21}$ | $p_{30}$ | $p_{14}$ | $p_{27}$ | $p_1$ | $p_{18}$ | $p_9$ |
| $p_{47}$ | $p_6$ | $p_{20}$ | $p_{31}$ | $p_{15}$ | $p_{25}$ | $p_3$ | $p_{16}$ | $p_{11}$ |

$P_{66}=$

|     | $p_{40}$ | $p_{41}$ | $p_{42}$ | $p_{43}$ | $p_{44}$ | $p_{45}$ | $p_{46}$ | $p_{47}$ |
|---|---|---|---|---|---|---|---|---|
| $p_{40}$ | $p_0$ | $p_8$ | $p_{26}$ | $p_{17}$ | $p_{28}$ | $p_4$ | $p_{12}$ | $p_{23}$ |
| $p_{41}$ | $p_1$ | $p_9$ | $p_{27}$ | $p_{18}$ | $p_{30}$ | $p_7$ | $p_{14}$ | $p_{21}$ |
| $p_{42}$ | $p_2$ | $p_{10}$ | $p_{24}$ | $p_{19}$ | $p_{29}$ | $p_5$ | $p_{13}$ | $p_{22}$ |
| $p_{43}$ | $p_3$ | $p_{11}$ | $p_{25}$ | $p_{16}$ | $p_{31}$ | $p_6$ | $p_{15}$ | $p_{20}$ |
| $p_{44}$ | $p_5$ | $p_{13}$ | $p_{29}$ | $p_{22}$ | $p_{24}$ | $p_2$ | $p_{10}$ | $p_{19}$ |
| $p_{45}$ | $p_4$ | $p_{12}$ | $p_{28}$ | $p_{23}$ | $p_{26}$ | $p_0$ | $p_8$ | $p_{17}$ |
| $p_{46}$ | $p_6$ | $p_{15}$ | $p_{31}$ | $p_{20}$ | $p_{25}$ | $p_3$ | $p_{11}$ | $p_{16}$ |
| $p_{47}$ | $p_7$ | $p_{14}$ | $p_{30}$ | $p_{21}$ | $p_{27}$ | $p_1$ | $p_9$ | $p_{18}$ |

## 2.5. Граф ОКТАЭДРа

Рассмотрим автоморфизм графа октаэдра $G_8$.

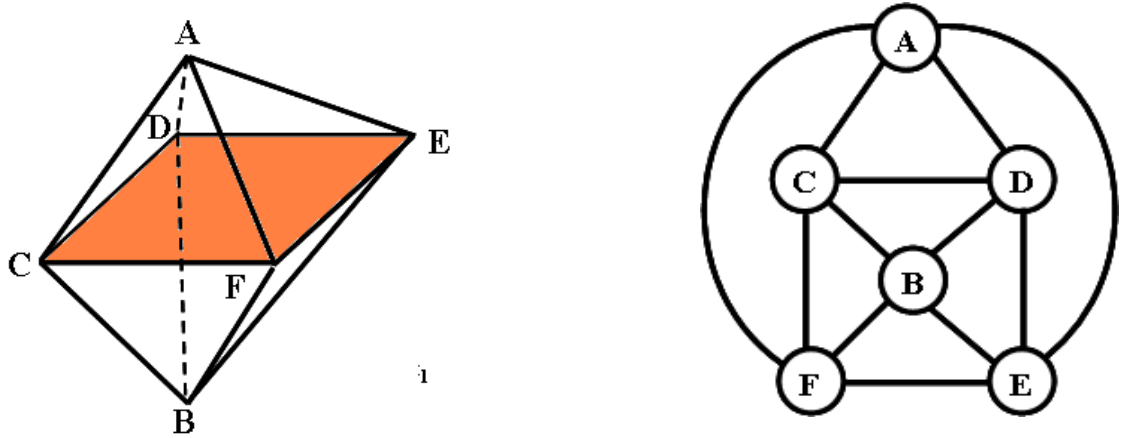

Рис. 2.73. Граф октаэдра $G_8$.

Количество вершин графа = 6. Количество ребер графа = 12

Количество изометрических циклов = 11

Смежность графа:

вершина 1:  2  4  5  6
вершина 2:  1  3  5  6
вершина 3:  2  4  5  6
вершина 4:  1  3  5  6
вершина 5:  1  2  3  4
вершина 6:  1  2  3  4

Инцидентность графа:

ребро $e_1$: $(v_1,v_2)$ или $(v_2,v_1)$;
ребро $e_2$: $(v_1,v_4)$ или $(v_4,v_1)$;
ребро $e_3$: $(v_1,v_5)$ или $(v_5,v_1)$;
ребро $e_4$: $(v_1,v_6)$ или $(v_6,v_1)$;
ребро $e_5$: $(v_2,v_3)$ или $(v_3,v_2)$;
ребро $e_6$: $(v_2,v_5)$ или $(v_5,v_2)$;
ребро $e_7$: $(v_2,v_6)$ или $(v_6,v_2)$;
ребро $e_8$: $(v_3,v_4)$ или $(v_4,v_3)$;
ребро $e_9$: $(v_3,v_5)$ или $(v_5,v_3)$;
ребро $e_{10}$: $(v_3,v_6)$ или $(v_6,v_3)$;
ребро $e_{11}$: $(v_4,v_5)$ или $(v_5,v_4)$;
ребро $e_{12}$: $(v_4,v_6)$ или $(v_6,v_4)$.

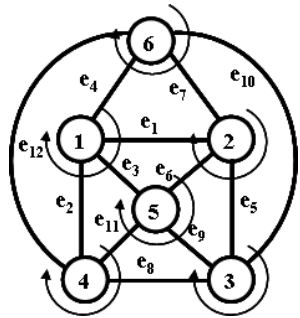

Кортеж весов ребер для векторного инварианта спектра реберных разрезов $\xi_w(G_8)$ = <12×14>, кортеж весов вершин $\zeta_w(G_8)$ = <6×56>.

Множество изометрических циклов:

$c_1 = \{e_1,e_3,e_6\}$;         $c_2 = \{e_2,e_3,e_{11}\}$;       $c_3 = \{e_8,e_9,e_{11}\}$;
$c_4 = \{e_5,e_6,e_9\}$;         $c_5 = \{e_2,e_4,e_{12}\}$;       $c_6 = \{e_5,e_7,e_{10}\}$;
$c_7 = \{e_1,e_4,e_7\}$;         $c_8 = \{e_8,e_{10},e_{12}\}$;    $c_9 = \{e_1,e_2,e_5,e_8\}$;

$c_{10} = \{e_3,e_4,e_9,e_{10}\}; \qquad c_{11} = \{e_6,e_7,e_{11},e_{12}\}.$

В данном случае, выбираем в качестве образующих циклов – циклы $c_9,c_{10},c_{11}$ (см. рис. 2.67). Следует заметить, что данные изометрические циклы могут быть образованы кольцевым суммированием других изометрических циклов.

$c_9 = c_1 \oplus c_2 \oplus c_3 \oplus c_4 = c_5 \oplus c_6 \oplus c_7 \oplus c_8 = \{e_1,e_2,e_5,e_8\};$
$c_{10} = c_1 \oplus c_2 \oplus c_5 \oplus c_7 = c_3 \oplus c_4 \oplus c_6 \oplus c_8 = \{e_6,e_7,e_{11},e_{12}\};$
$c_{11} = c_1 \oplus c_4 \oplus c_6 \oplus c_7 = c_2 \oplus c_3 \oplus c_5 \oplus c_8 = \{e_3,e_4,e_9,e_{10}\}.$

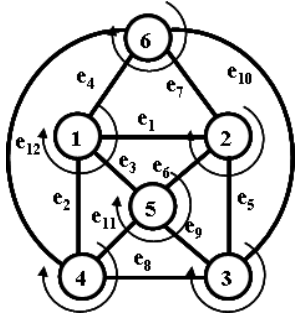

Рис. 2.74. Топологический рисунок графа октаэдра.

Здесь, каждый простой четырехугольный цикл имеет в центре различные вершины. Например, цикл $\{e_1,e_2,e_5,e_8\}$ имеет в середине либо вершину $v_5$, либо вершину $v_6$. Таким образом, получаем три образующихся циклов, но шесть топологических рисунков [19]. Так как вершины A и B могут меняться местами (перемещаться в пространстве $R^3$). Определим следующие перестановки. Будем переставлять вершины в многоугольнике, затем переставим вершины A,B.

В качестве образующего цикла выберем цикл $c_9$. Осуществим все 8 перестановок, расположив вершину 5 в центре.

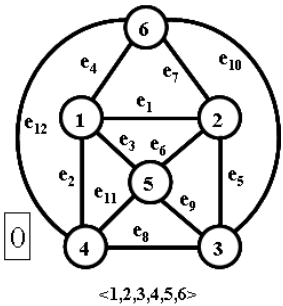 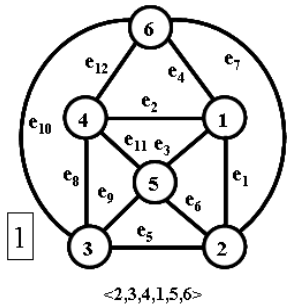 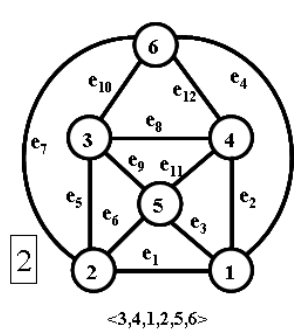 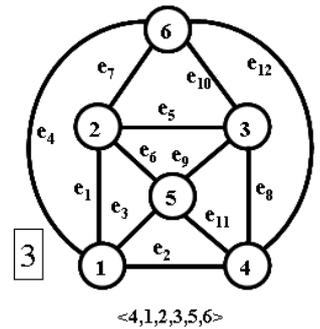

| Рис. 2.75. Перестановка <1,2,3,4,5,6>. | Рис. 2.76. Перестановка <2,3,4,1,5,6>. | Рис. 2.77. Перестановка <3,4,1,2,5,6>. | Рис. 2.78. Перестановка <4,1,2,3,5,6>. |

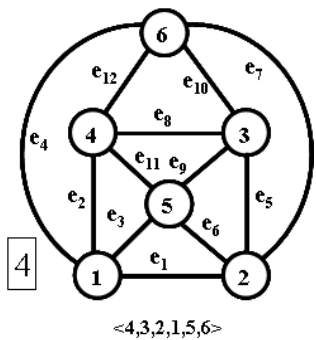 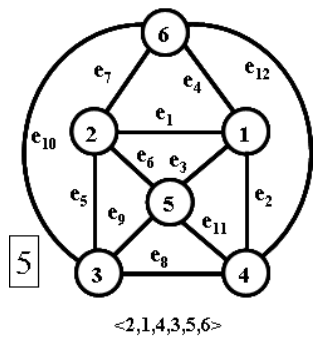 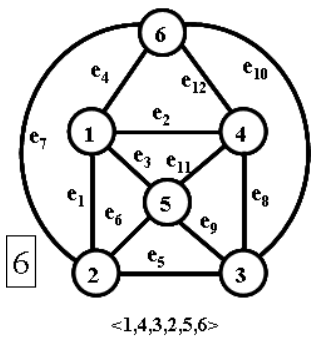 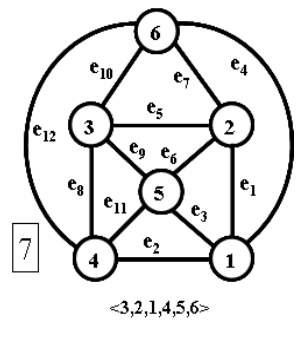

Рис. 2.79. Перестановка <4,3,2,1,5,6>.  
Рис. 2.80. Перестановка <2,1,4,3,5,6>.  
Рис. 2.81. Перестановка <1,4,3,2,5,6>.  
Рис. 2.82. Перестановка <3,2,1,4,5,6>.

В качестве образующего цикла выберем цикл $c_{10}$. Осуществим все 8 перестановок, расположив вершину 1 в центре.

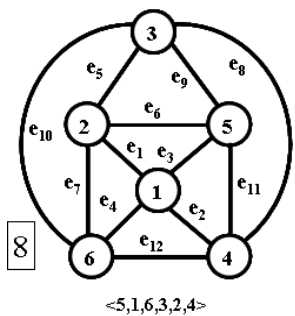 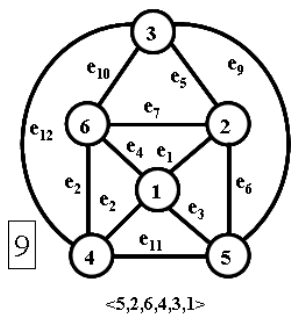 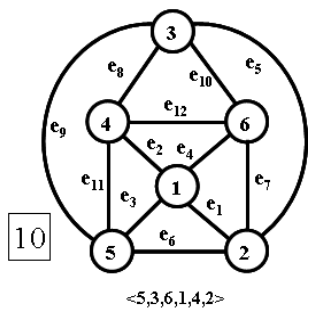 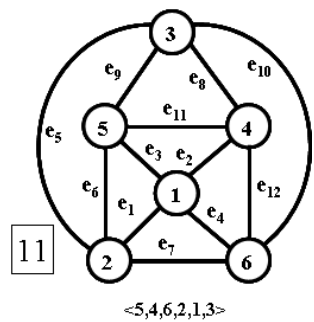

Рис. 2.83. Перестановка <5,1,6,3,2,4>.  
Рис. 2.84. Перестановка <5,2,6,4,3,1>.  
Рис. 2.85. Перестановка <5,3,6,1,4,2>.  
Рис. 2.86. Перестановка <5,4,6,2,1,3>.

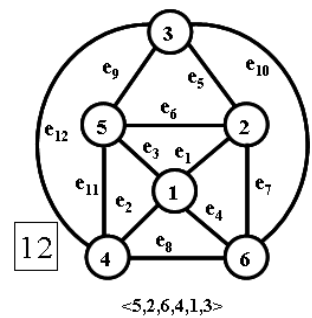 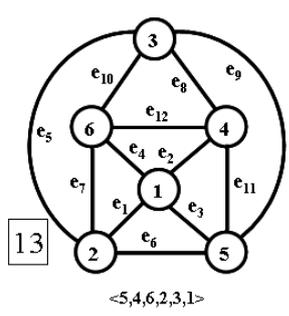 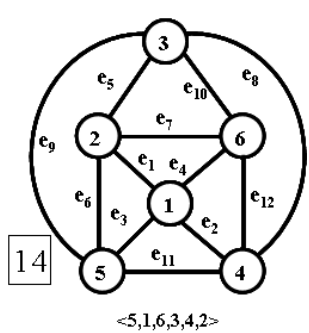 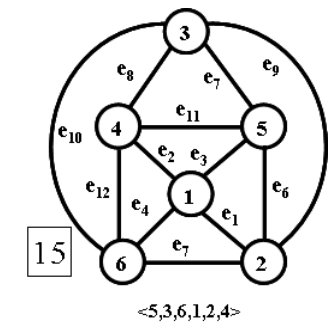

Рис. 2.87. Перестановка <5,2,6,4,1,3>.  
Рис. 2.88. Перестановка <5,4,6,2,3,1>.  
Рис. 2.89. Перестановка <5,1,6,3,4,2>.  
Рис. 2.90. Перестановка <5,3,6,1,2,4>.

В качестве образующего цикла выберем цикл $c_{11}$. Осуществим все 8 перестановок, расположив вершину 2 в центре.

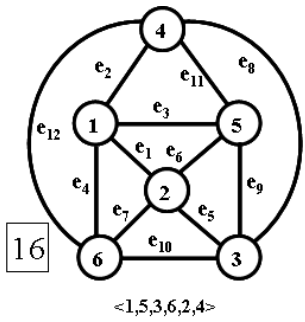 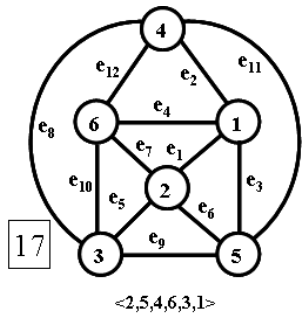 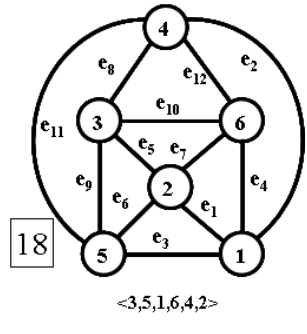 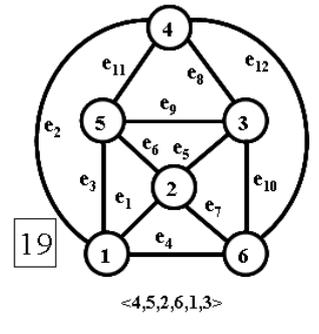

Рис. 2.91. Перестановка <1,5,3,6,2,4>. | Рис. 2.92. Перестановка <2,5,4,6,3,1>. | Рис. 2.93. Перестановка <3,5,1,6,4,2>. | Рис. 2.94. Перестановка <4,5,2,6,1,3>.

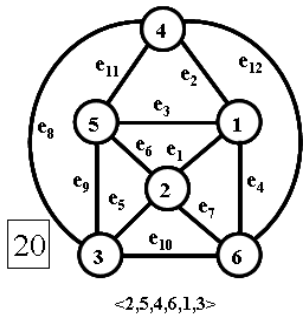 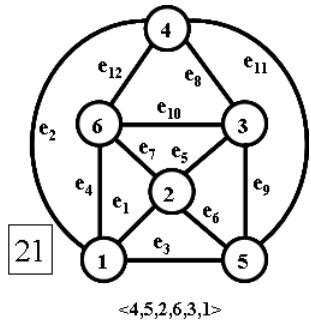 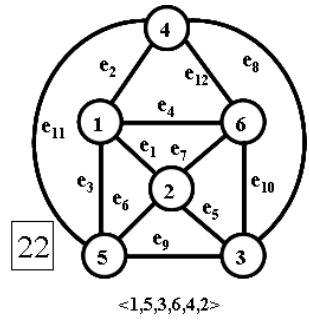 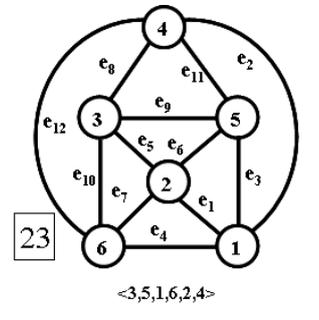

Рис. 2.95. Перестановка <2,5,4,6,1,3>. | Рис. 2.96. Перестановка <4,5,2,6,3,1>. | Рис. 2.97. Перестановка <1,5,3,6,4,2>. | Рис. 2.98. Перестановка <3,5,1,6,2,4>.

В качестве образующего цикла выберем цикл $c_9$. Осуществим все 8 перестановок, расположив вершину 6 в центре.

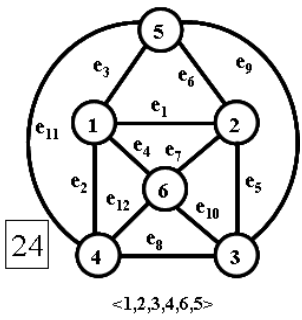 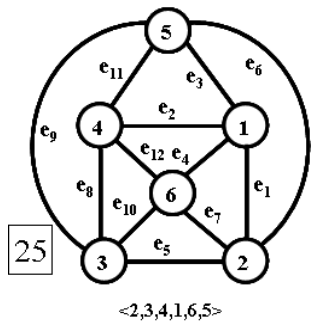 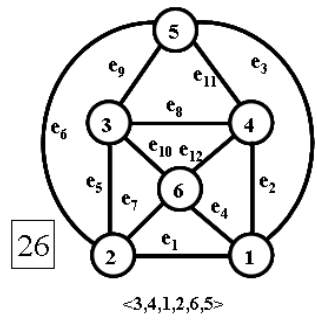 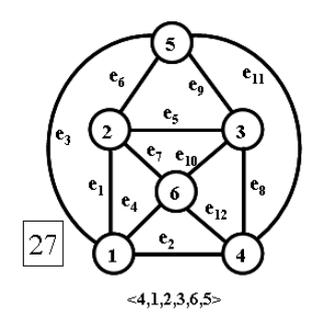

Рис. 2.99. Перестановка <1,2,3,4,6,5>. | Рис. 2.100. Перестановка <2,3,4,1,6,5>. | Рис. 2.101. Перестановка <3,4,1,2,6,5>. | Рис. 2.102. Перестановка <4,1,2,3,6,5>.

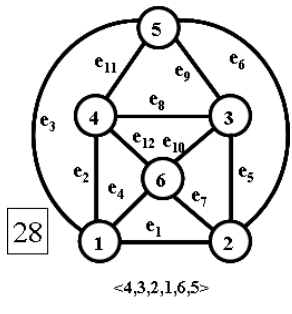 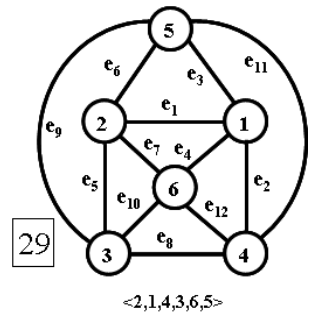 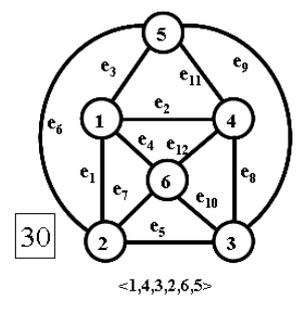 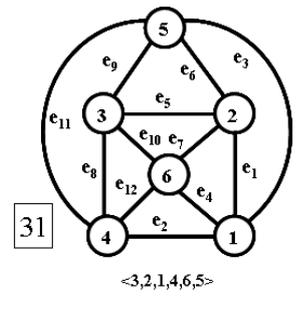

|28 <4,3,2,1,6,5> | 29 <2,1,4,3,6,5> | 30 <1,4,3,2,6,5> | 31 <3,2,1,4,6,5> |

Рис. 2.103. Перестановка <4,3,2,1,6,5>.  
Рис. 2.104. Перестановка <2,1,4,3,6,5>.  
Рис. 2.105. Перестановка <1,4,3,2,6,5>.  
Рис. 2.106. Перестановка <3,2,1,4,6,5>.

В качестве образующего цикла выберем цикл $c_{10}$. Осуществим все 8 перестановок, расположив вершину 3 в центре.

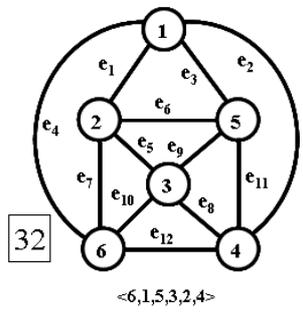 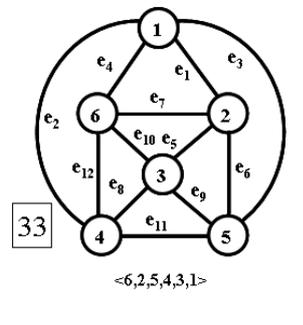 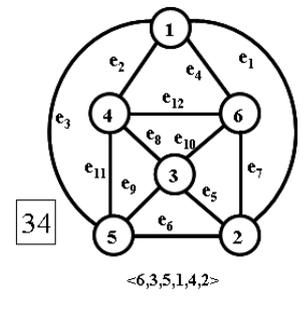 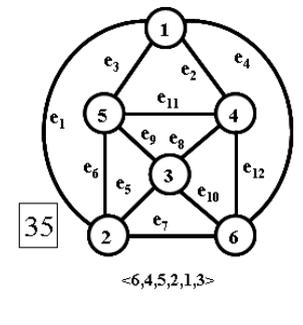

|32 <6,1,5,3,2,4> | 33 <6,2,5,4,3,1> | 34 <6,3,5,1,4,2> | 35 <6,4,5,2,1,3> |

Рис. 2.107. Перестановка <6,1,5,3,2,4>.  
Рис. 2.108. Перестановка <6,2,5,4,3,1>.  
Рис. 2.109. Перестановка <6,3,5,1,4,2>.  
Рис. 2.110. Перестановка <6,4,5,2,1,3>.

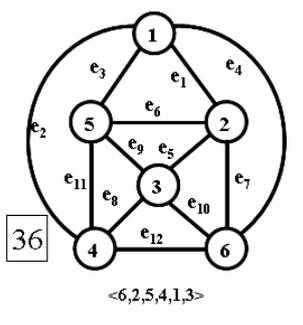 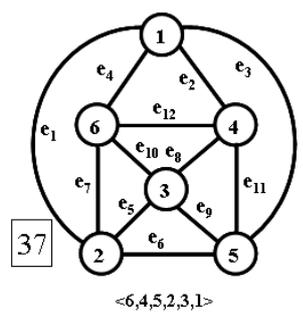 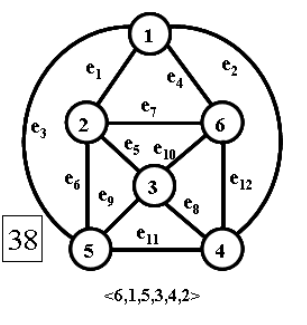 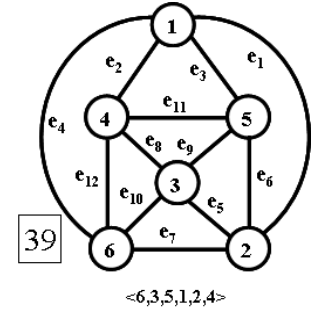

|36 <6,2,5,4,1,3> | 37 <6,4,5,2,3,1> | 38 <6,1,5,3,4,2> | 39 <6,3,5,1,2,4> |

Рис. 2.111. Перестановка <6,2,5,4,1,3>.  
Рис. 2.112. Перестановка <6,4,5,2,3,1>.  
Рис. 2.113. Перестановка <6,1,5,3,4,2>.  
Рис. 2.114. Перестановка <6,3,5,1,2,4>.

В качестве образующего цикла выберем цикл $c_{11}$. Осуществим все 8 перестановок, расположив вершину 4 в центре.

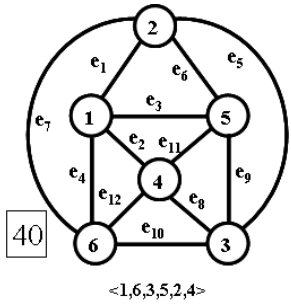 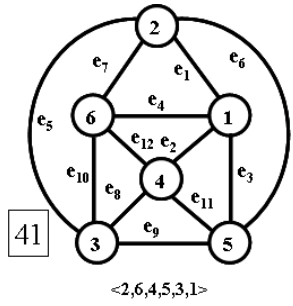 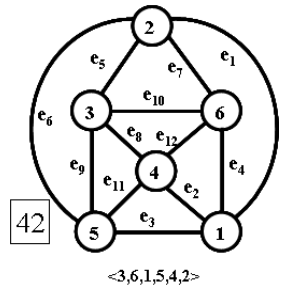 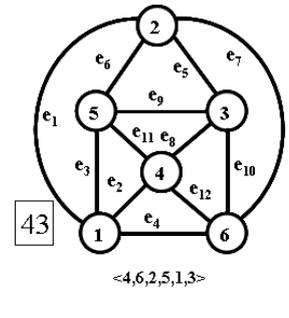

Рис. 2.115. Перестановка <1,6,3,5,2,4>.

Рис. 2.116. Перестановка <2,6,4,5,3,1>.

Рис. 2.117. Перестановка <3,6,1,5,4,2>.

Рис. 2.118. Перестановка <4,6,2,5,1,3>.

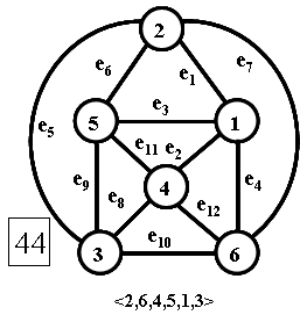 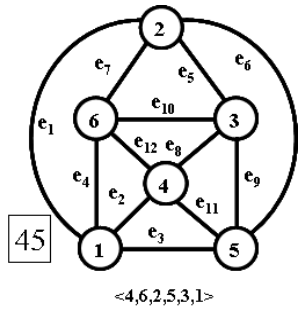 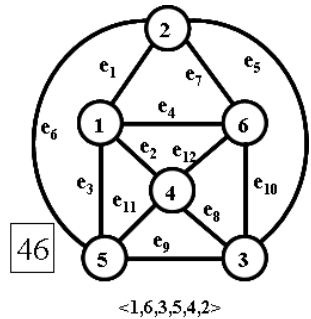 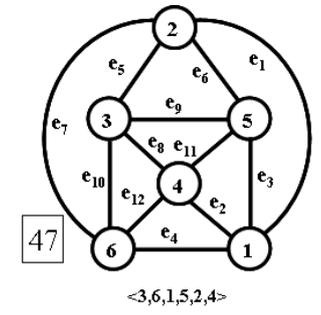

Рис. 2.119. Перестановка <2,6,4,6,1,3>.

Рис. 2.120. Перестановка <4,6,2,5,3,1>.

Рис. 2.121. Перестановка <1,6,3,5,4,2>.

Рис. 2.122. Перестановка <3,6,1,5,2,4>.

Граф $G_2$ является планарным графом, кольцевая сумма изометрических циклов есть пустое множество $\sum_{i=1}^{11} c_i = \varnothing$.

Таким образом, выделено 6×8 = 48 перестановок:

$p_0$ = <1,2,3,4,5,6> = (1)(2)(3)(4)(5)(6);
$p_1$ = <2,3,4,1,5,6> = (1 2 3 4)(5)(6);
$p_2$ = <3,4,1,2,5,6> = (1 3)(2 4)(5)(6);
$p_3$ = <4,1,2,3,5,6> = (1 4 3 2)(5)(6);
$p_4$ = <4,3,2,1,5,6> = (1 4)(2 3)(5)(6);
$p_5$ = <2,1,4,3,5,6> = (1 2)(3 4)(5)(6);
$p_6$ = <1,4,3,2,5,6> = (1)(2 4)(3)(5)(6);
$p_7$ = <3,2,1,4,5,6> = (1 3)(2)(4)(5)(6);
$p_8$ = <5,1,6,3,2,4> = (1 5 2)(3 6 4);
$p_9$ = <5,2,6,4,3,1> = (1 5 3 6)(2)(4);
$p_{10}$ = <5,3,6,1,4,2> = (1 5 4)(2 4 6);
$p_{11}$ = <5,4,6,2,1,3> = (1 5)(2 4)(3 6);
$p_{12}$ = <5,2,6,4,1,3> = (1 5)(2)(3 6)(4);
$p_{13}$ = <5,4,6,2,3,1> = (1 5 3 6)(2 4);
$p_{14}$ = <5,1,6,3,4,2> = (1 5 4 3 6 2);
$p_{15}$ = <5,3,6,1,2,4> = (1 5 2 3 6 4);
$p_{16}$ = <1,5,3,6,2,4> = (1)(2 5)(3)(4 6);

$p_{17}= \langle 2,5,4,6,3,1 \rangle = (1\ 2\ 5\ 3\ 4\ 6)$;
$p_{18}= \langle 3,5,1,6,4,2 \rangle = (1\ 3)(2\ 5\ 4\ 6)$;
$p_{19}= \langle 4,5,2,6,1,3 \rangle = (1\ 4\ 6\ 3\ 2\ 5)$;
$p_{20}= \langle 2,5,4,6,1,3 \rangle = (1\ 2\ 5)(3\ 4\ 6)$;
$p_{21}= \langle 4,5,2,6,3,1 \rangle = (1\ 4\ 6)(2\ 5\ 3)$;
$p_{22}= \langle 1,5,3,6,4,2 \rangle = (1)(2\ 5\ 4\ 6)(3)$;
$p_{23}= \langle 3,5,1,6,2,4 \rangle = (1\ 3)(2\ 5)(4\ 6)$;
$p_{24}= \langle 1,2,3,4,6,5 \rangle = (1)(2)(3)(4)(5\ 6)$;
$p_{25}= \langle 2,3,4,1,6,5 \rangle = (1\ 2\ 3\ 4)(5\ 6)$;
$p_{26}= \langle 3,4,1,2,6,5 \rangle = (1\ 3)(2\ 4)(5\ 6)$;
$p_{27}= \langle 4,1,2,3,6,5 \rangle = (1\ 4\ 3\ 2)(5\ 6)$;
$p_{28}= \langle 4,3,2,1,6,5 \rangle = (1\ 4)(2\ 3)(5\ 6)$;
$p_{29}= \langle 2,1,4,3,6,5 \rangle = (1\ 2)(3\ 4)(5\ 6)$;
$p_{30}= \langle 1,4,3,2,6,5 \rangle = (1)(2\ 4)(3)(5\ 6)$;
$p_{32}= \langle 3,2,1,4,6,5 \rangle = (1\ 3)(2)(4)(5\ 6)$;
$p_{32}= \langle 6,1,5,3,2,4 \rangle = (1\ 6\ 4\ 3\ 5\ 2)$;
$p_{33}= \langle 6,2,5,4,3,1 \rangle = (1\ 6)(2)(3\ 5)(4)$;
$p_{34}= \langle 6,3,5,1,4,2 \rangle = (1\ 6\ 2\ 3\ 5\ 4)$;
$p_{35}= \langle 6,4,5,2,1,3 \rangle = (1\ 6\ 3\ 5)(2\ 4)$;
$p_{36}= \langle 6,2,5,4,1,3 \rangle = (1\ 6\ 3\ 5)(2\ 4)$;
$p_{37}= \langle 6,4,5,2,3,1 \rangle = (1\ 6)(2\ 4)(3\ 5)$;
$p_{38}= \langle 6,1,5,3,4,2 \rangle = (1\ 6\ 2)(3\ 5\ 4)$;
$p_{39}= \langle 6,3,5,1,2,4 \rangle = (1\ 6\ 4)(2\ 3\ 5)$;
$p_{40}= \langle 1,6,3,5,2,4 \rangle = (1)(2\ 6\ 4\ 5)(3)$;
$p_{41}= \langle 2,6,4,5,3,1 \rangle = (1\ 2\ 6)(3\ 4\ 5)$;
$p_{42}= \langle 3,6,1,5,4,2 \rangle = (1\ 3)(2\ 6)(4\ 5)$;
$p_{43}= \langle 4,6,2,5,1,3 \rangle = (1\ 4\ 5)(2\ 6\ 3)$;
$p_{44}= \langle 2,6,4,5,1,3 \rangle = (1\ 2\ 6\ 3\ 4\ 5)$;
$p_{45}= \langle 4,6,2,5,3,1 \rangle = (1\ 4\ 5\ 3\ 2\ 6)$;
$p_{46}= \langle 1,6,3,5,4,2 \rangle = (1)(2\ 6)(3)(4\ 5)$;
$p_{47}= \langle 3,6,1,5,2,4 \rangle = (1\ 3)(2\ 6\ 4\ 5)$.

Таблица Кэли для графа октаэдра:

$$Aut(G_2) = \begin{array}{|c|c|c|c|c|c|} \hline P_{11} & P_{12} & P_{13} & P_{14} & P_{15} & P_{16} \\ \hline P_{21} & P_{22} & P_{23} & P_{24} & P_{25} & P_{26} \\ \hline P_{31} & P_{32} & P_{33} & P_{34} & P_{35} & P_{36} \\ \hline P_{41} & P_{42} & P_{43} & P_{44} & P_{45} & P_{46} \\ \hline P_{51} & P_{52} & P_{53} & P_{54} & P_{55} & P_{56} \\ \hline P_{61} & P_{62} & P_{63} & P_{64} & P_{65} & P_{66} \\ \hline \end{array}$$

Таблица Кэли для первой строки, состоящей из блоков $P_{11} – P_{16}$:

$P_{11}=$

|       | $p_0$ | $p_1$ | $p_2$ | $p_3$ | $p_4$ | $p_5$ | $p_6$ | $p_7$ |
|-------|-------|-------|-------|-------|-------|-------|-------|-------|
| $p_0$ | $p_0$ | $p_1$ | $p_2$ | $p_3$ | $p_4$ | $p_5$ | $p_6$ | $p_7$ |
| $p_1$ | $p_1$ | $p_2$ | $p_3$ | $p_0$ | $p_7$ | $p_6$ | $p_4$ | $p_5$ |
| $p_2$ | $p_2$ | $p_3$ | $p_0$ | $p_1$ | $p_5$ | $p_4$ | $p_7$ | $p_6$ |
| $p_3$ | $p_3$ | $p_0$ | $p_1$ | $p_2$ | $p_6$ | $p_7$ | $p_5$ | $p_4$ |
| $p_4$ | $p_4$ | $p_6$ | $p_5$ | $p_7$ | $p_0$ | $p_2$ | $p_1$ | $p_3$ |
| $p_5$ | $p_5$ | $p_7$ | $p_4$ | $p_6$ | $p_2$ | $p_0$ | $p_3$ | $p_1$ |
| $p_6$ | $p_6$ | $p_5$ | $p_7$ | $p_4$ | $p_3$ | $p_1$ | $p_0$ | $p_2$ |
| $p_7$ | $p_7$ | $p_4$ | $p_6$ | $p_5$ | $p_1$ | $p_3$ | $p_2$ | $p_0$ |

$P_{12}=$

|     | $p_8$ | $p_9$ | $p_{10}$ | $p_{11}$ | $p_{12}$ | $p_{13}$ | $p_{14}$ | $p_{15}$ |
|-----|-------|-------|----------|----------|----------|----------|----------|----------|
| $p_0$ | $p_8$ | $p_9$ | $p_{10}$ | $p_{11}$ | $p_{12}$ | $p_{13}$ | $p_{14}$ | $p_{15}$ |
| $p_1$ | $p_{40}$ | $p_{41}$ | $p_{42}$ | $p_{43}$ | $p_{44}$ | $p_{45}$ | $p_{46}$ | $p_{47}$ |
| $p_2$ | $p_{39}$ | $p_{37}$ | $p_{38}$ | $p_{36}$ | $p_{35}$ | $p_{33}$ | $p_{345}$ | $p_{32}$ |
| $p_3$ | $p_{23}$ | $p_{21}$ | $p_{22}$ | $p_{20}$ | $p_{19}$ | $p_{17}$ | $p_{18}$ | $p_{16}$ |
| $p_4$ | $p_{47}$ | $p_{45}$ | $p_{46}$ | $p_{44}$ | $p_{43}$ | $p_{41}$ | $p_{42}$ | $p_{40}$ |
| $p_5$ | $p_{16}$ | $p_{17}$ | $p_{18}$ | $p_{19}$ | $p_{20}$ | $p_{21}$ | $p_{22}$ | $p_{23}$ |
| $p_6$ | $p_{15}$ | $p_{13}$ | $p_{14}$ | $p_{12}$ | $p_{11}$ | $p_9$ | $p_{10}$ | $p_8$ |
| $p_7$ | $p_{32}$ | $p_{33}$ | $p_{34}$ | $p_{35}$ | $p_{36}$ | $p_{37}$ | $p_{38}$ | $p_{39}$ |

$P_{13}=$

|     | $p_{16}$ | $p_{17}$ | $p_{18}$ | $p_{19}$ | $p_{20}$ | $p_{21}$ | $p_{22}$ | $p_{23}$ |
|-----|----------|----------|----------|----------|----------|----------|----------|----------|
| $p_0$ | $p_{16}$ | $p_{17}$ | $p_{18}$ | $p_{19}$ | $p_{20}$ | $p_{21}$ | $p_{22}$ | $p_{23}$ |
| $p_1$ | $p_{15}$ | $p_{13}$ | $p_{14}$ | $p_{12}$ | $p_{11}$ | $p_9$ | $p_{10}$ | $p_8$ |
| $p_2$ | $p_{47}$ | $p_{45}$ | $p_{46}$ | $p_{44}$ | $p_{43}$ | $p_{41}$ | $p_{42}$ | $p_{40}$ |
| $p_3$ | $p_{32}$ | $p_{33}$ | $p_{34}$ | $p_{35}$ | $p_{36}$ | $p_{37}$ | $p_{38}$ | $p_{39}$ |
| $p_4$ | $p_{39}$ | $p_{37}$ | $p_{38}$ | $p_{36}$ | $p_{35}$ | $p_{33}$ | $p_{34}$ | $p_{32}$ |
| $p_5$ | $p_8$ | $p_9$ | $p_{10}$ | $p_{11}$ | $p_{12}$ | $p_{13}$ | $p_{14}$ | $p_{15}$ |
| $p_6$ | $p_{40}$ | $p_{41}$ | $p_{42}$ | $p_{43}$ | $p_{44}$ | $p_{45}$ | $p_{46}$ | $p_{47}$ |
| $p_7$ | $p_{23}$ | $p_{21}$ | $p_{22}$ | $p_{20}$ | $p_{19}$ | $p_{17}$ | $p_{18}$ | $p_{16}$ |

$P_{14}=$

|     | $p_{24}$ | $p_{25}$ | $p_{26}$ | $p_{27}$ | $p_{28}$ | $p_{29}$ | $p_{30}$ | $p_{31}$ |
|-----|----------|----------|----------|----------|----------|----------|----------|----------|
| $p_0$ | $p_{24}$ | $p_{25}$ | $p_{26}$ | $p_{27}$ | $p_{28}$ | $p_{29}$ | $p_{30}$ | $p_{31}$ |
| $p_1$ | $p_{25}$ | $p_{26}$ | $p_{27}$ | $p_{24}$ | $p_{31}$ | $p_{30}$ | $p_{28}$ | $p_{29}$ |
| $p_2$ | $p_{26}$ | $p_{27}$ | $p_{24}$ | $p_{25}$ | $p_{29}$ | $p_{28}$ | $p_{31}$ | $p_{30}$ |
| $p_3$ | $p_{27}$ | $p_{24}$ | $p_{25}$ | $p_{26}$ | $p_{30}$ | $p_{31}$ | $p_{29}$ | $p_{28}$ |
| $p_4$ | $p_{28}$ | $p_{30}$ | $p_{29}$ | $p_{31}$ | $p_{24}$ | $p_{26}$ | $p_{25}$ | $p_{27}$ |
| $p_5$ | $p_{29}$ | $p_{31}$ | $p_{28}$ | $p_{30}$ | $p_{26}$ | $p_{24}$ | $p_{27}$ | $p_{25}$ |
| $p_6$ | $p_{30}$ | $p_{29}$ | $p_{31}$ | $p_{28}$ | $p_{27}$ | $p_{25}$ | $p_{24}$ | $p_{26}$ |
| $p_7$ | $p_{31}$ | $p_{28}$ | $p_{30}$ | $p_{29}$ | $p_{25}$ | $p_{27}$ | $p_{26}$ | $p_{24}$ |

$P_{15}=$

|     | $p_{32}$ | $p_{33}$ | $p_{34}$ | $p_{35}$ | $p_{36}$ | $p_{37}$ | $p_{38}$ | $p_{39}$ |
|-----|----------|----------|----------|----------|----------|----------|----------|----------|
| $p_0$ | $p_{32}$ | $p_{33}$ | $p_{34}$ | $p_{35}$ | $p_{36}$ | $p_{37}$ | $p_{38}$ | $p_{39}$ |
| $p_1$ | $p_{16}$ | $p_{17}$ | $p_{18}$ | $p_{19}$ | $p_{20}$ | $p_{21}$ | $p_{22}$ | $p_{23}$ |
| $p_2$ | $p_{15}$ | $p_{13}$ | $p_{14}$ | $p_{12}$ | $p_{11}$ | $p_9$ | $p_{10}$ | $p_8$ |
| $p_3$ | $p_{47}$ | $p_{45}$ | $p_{46}$ | $p_{44}$ | $p_{43}$ | $p_{41}$ | $p_{42}$ | $p_{40}$ |
| $p_4$ | $p_{23}$ | $p_{21}$ | $p_{22}$ | $p_{20}$ | $p_{19}$ | $p_{17}$ | $p_{18}$ | $p_{16}$ |
| $p_5$ | $p_{40}$ | $p_{41}$ | $p_{42}$ | $p_{43}$ | $p_{44}$ | $p_{45}$ | $p_{46}$ | $p_{47}$ |
| $p_6$ | $p_{39}$ | $p_{37}$ | $p_{38}$ | $p_{36}$ | $p_{35}$ | $p_{33}$ | $p_{34}$ | $p_{32}$ |
| $p_7$ | $p_8$ | $p_9$ | $p_{10}$ | $p_{11}$ | $p_{12}$ | $p_{13}$ | $p_{14}$ | $p_{15}$ |

$P_{16}=$

|     | $p_{40}$ | $p_{41}$ | $p_{42}$ | $p_{43}$ | $p_{44}$ | $p_{45}$ | $p_{46}$ | $p_{47}$ |
|-----|----------|----------|----------|----------|----------|----------|----------|----------|
| $p_0$ | $p_{40}$ | $p_{41}$ | $p_{42}$ | $p_{43}$ | $p_{44}$ | $p_{45}$ | $p_{46}$ | $p_{47}$ |
| $p_1$ | $p_{39}$ | $p_{37}$ | $p_{38}$ | $p_{36}$ | $p_{35}$ | $p_{33}$ | $p_{34}$ | $p_{32}$ |
| $p_2$ | $p_{23}$ | $p_{21}$ | $p_{22}$ | $p_{20}$ | $p_{19}$ | $p_{17}$ | $p_{18}$ | $p_{16}$ |
| $p_3$ | $p_8$ | $p_9$ | $p_{10}$ | $p_{11}$ | $p_{12}$ | $p_{13}$ | $p_{14}$ | $p_{15}$ |
| $p_4$ | $p_{15}$ | $p_{13}$ | $p_{14}$ | $p_{12}$ | $p_{11}$ | $p_9$ | $p_{10}$ | $p_8$ |
| $p_5$ | $p_{32}$ | $p_{33}$ | $p_{34}$ | $p_{35}$ | $p_{36}$ | $p_{37}$ | $p_{38}$ | $p_{39}$ |
| $p_6$ | $p_{16}$ | $p_{17}$ | $p_{18}$ | $p_{19}$ | $p_{20}$ | $p_{21}$ | $p_{22}$ | $p_{23}$ |
| $p_7$ | $p_{47}$ | $p_{45}$ | $p_{46}$ | $p_{44}$ | $p_{43}$ | $p_{41}$ | $p_{42}$ | $p_{40}$ |

## 2.6. Граф ДОДЕКАЭДРа

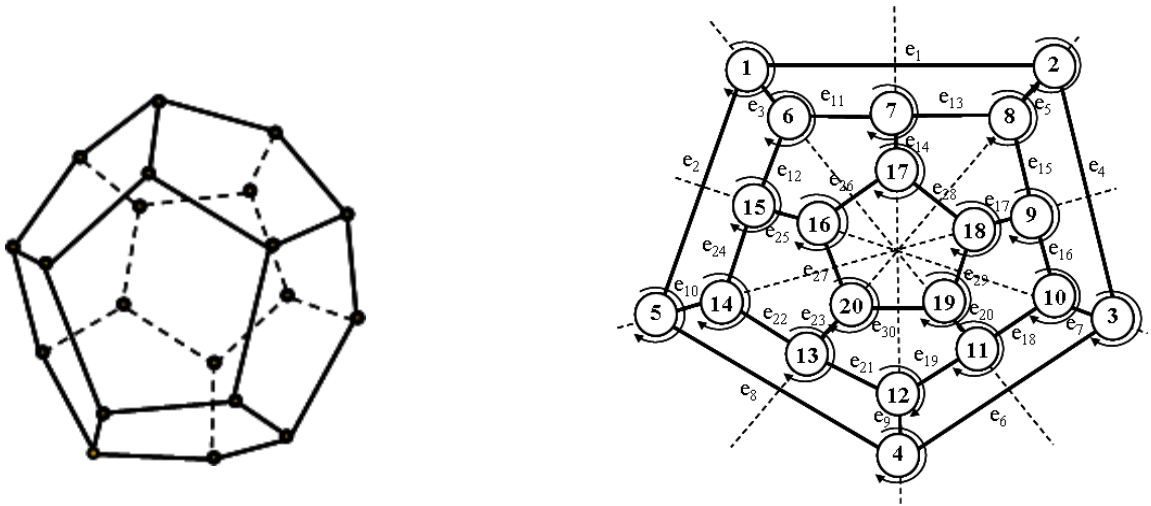

Рис. 2.123. Граф додекаэдра $G_9$ и его топологический рисунок.

Количество вершин графа = 20. Количество рёбер графа = 30.

Смежность графа:

вершина $v_1$: $v_2$ $v_5$ $v_6$;
вершина $v_2$: $v_1$ $v_3$ $v_8$;
вершина $v_3$: $v_2$ $v_4$ $v_{10}$;
вершина $v_4$: $v_3$ $v_5$ $v_{12}$;
вершина $v_5$: $v_1$ $v_4$ $v_{14}$;
вершина $v_6$: $v_1$ $v_7$ $v_{15}$;
вершина $v_7$: $v_6$ $v_8$ $v_{17}$;
вершина $v_8$: $v_2$ $v_7$ $v_9$;
вершина $v_9$: $v_8$ $v_{10}$ $v_{18}$;
вершина $v_{10}$: $v_3$ $v_9$ $v_{11}$;
вершина $v_{11}$: $v_{10}$ $v_{12}$ $v_{19}$;
вершина $v_{12}$: $v_4$ $v_{11}$ $v_{13}$;
вершина $v_{13}$: $v_{12}$ $v_{14}$ $v_{20}$;
вершина $v_{14}$: $v_5$ $v_{13}$ $v_{15}$;
вершина $v_{15}$: $v_6$ $v_{14}$ $v_{16}$;
вершина $v_{16}$: $v_{15}$ $v_{17}$ $v_{20}$;
вершина $v_{17}$: $v_7$ $v_{16}$ $v_{18}$;
вершина $v_{18}$: $v_9$ $v_{17}$ $v_{19}$;
вершина $v_{19}$: $v_{11}$ $v_{18}$ $v_{20}$;
вершина $v_{20}$: $v_{13}$ $v_{16}$ $v_{19}$.

Инцидентность графа:

ребро $e_1$: $(v_1,v_2)$ или $(v_2,v_1)$;     ребро $e_2$: $(v_1,v_5)$ или $(v_5,v_1)$;
ребро $e_3$: $(v_1,v_6)$ или $(v_6,v_1)$;     ребро $e_4$: $(v_2,v_3)$ или $(v_3,v_2)$;
ребро $e_5$: $(v_2,v_8)$ или $(v_8,v_2)$;     ребро $e_6$: $(v_3,v_4)$ или $(v_4,v_3)$;
ребро $e_7$: $(v_3,v_{10})$ или $(v_{10},v_3)$;     ребро $e_8$: $(v_4,v_5)$ или $(v_5,v_4)$;
ребро $e_9$: $(v_4,v_{12})$ или $(v_{12},v_4)$;     ребро $e_{10}$: $(v_5,v_{14})$ или $(v_{14},v_5)$;

ребро $e_{11}$: $(v_6,v_7)$ или $(v_7,v_6)$;  ребро $e_{12}$: $(v_6,v_{15})$ или $(v_{15},v_6)$;
ребро $e_{13}$: $(v_7,v_8)$ или $(v_8,v_7)$;  ребро $e_{14}$: $(v_7,v_{17})$ или $(v_{17},v_7)$;
ребро $e_{15}$: $(v_8,v_9)$ или $(v_9,v_8)$;  ребро $e_{16}$: $(v_9,v_{10})$ или $(v_{10},v_9)$;
ребро $e_{17}$: $(v_9,v_{18})$ или $(v_{18},v_9)$;  ребро $e_{18}$: $(v_{10},v_{11})$ или $(v_{11},v_{10})$;
ребро $e_{19}$: $(v_{11},v_{12})$ или $(v_{12},v_{11})$;  ребро $e_{20}$: $(v_{11},v_{19})$ или $(v_{19},v_{11})$;
ребро $e_{21}$: $(v_{12},v_{13})$ или $(v_{13},v_{12})$;  ребро $e_{22}$: $(v_{13},v_{14})$ или $(v_{14},v_{13})$;
ребро $e_{23}$: $(v_{13},v_{20})$ или $(v_{20},v_{13})$;  ребро $e_{24}$: $(v_{14},v_{15})$ или $(v_{15},v_{14})$;
ребро $e_{25}$: $(v_{15},v_{16})$ или $(v_{16},v_{15})$;  ребро $e_{26}$: $(v_{16},v_{17})$ или $(v_{17},v_{16})$;
ребро $e_{27}$: $(v_{16},v_{20})$ или $(v_{20},v_{16})$;  ребро $e_{28}$: $(v_{17},v_{18})$ или $(v_{18},v_{17})$;
ребро $e_{29}$: $(v_{18},v_{19})$ или $(v_{19},v_{18})$;  ребро $e_{30}$: $(v_{19},v_{20})$ или $(v_{20},v_{19})$.

Множество изометрических циклов графа:

$c_1 = \{e_1,e_2,e_4,e_6,e_8\} \to \{v_1,v_2,v_3,v_4,v_5\}$;
$c_2 = \{e_1,e_3,e_5,e_{11},e_{13}\} \to \{v_1,v_2,v_6,v_7,v_8\}$;
$c_3 = \{e_2,e_3,e_{10},e_{12},e_{24}\} \to \{v_1,v_5,v_6,v_{14},v_{15}\}$;
$c_4 = \{e_4,e_5,e_7,e_{15},e_{16}\} \to \{v_2,v_3,v_8,v_9,v_{10}\}$;
$c_5 = \{e_6,e_7,e_9,e_{18},e_{19}\} \to \{v_3,v_4,v_{10},v_{11},v_{12}\}$;
$c_6 = \{e_8,e_9,e_{10},e_{21},e_{22}\} \to \{v_4,v_5,v_{12},v_{13},v_{14}\}$;
$c_7 = \{e_{11},e_{12},e_{14},e_{25},e_{26}\} \to \{v_6,v_7,v_{15},v_{16},v_{17}\}$;
$c_8 = \{e_{13},e_{14},e_{15},e_{17},e_{28}\} \to \{v_7,v_8,v_9,v_{17},v_{18}\}$;
$c_9 = \{e_{16},e_{17},e_{18},e_{20},e_{29}\} \to \{v_9,v_{10},v_{11},v_{18},v_{19}\}$;
$c_{10} = \{e_{19},e_{20},e_{21},e_{23},e_{30}\} \to \{v_{11},v_{12},v_{13},v_{19},v_{20}\}$;
$c_{11} = \{e_{22},e_{23},e_{24},e_{25},e_{27}\} \to \{v_{13},v_{14},v_{15},v_{16},v_{20}\}$;
$c_{12} = \{e_{26},e_{27},e_{28},e_{29},e_{30}\} \to \{v_{16},v_{17},v_{18},v_{19},v_{20}\}$.

Вектор весов ребер: $F(\xi_w(G_9)) = (30 \times 16)$;

Вектор весов вершин: $F(\zeta_w(G_9)) = (20 \times 48)$

В качестве образующего цикла, выбираем цикл $c_1$.

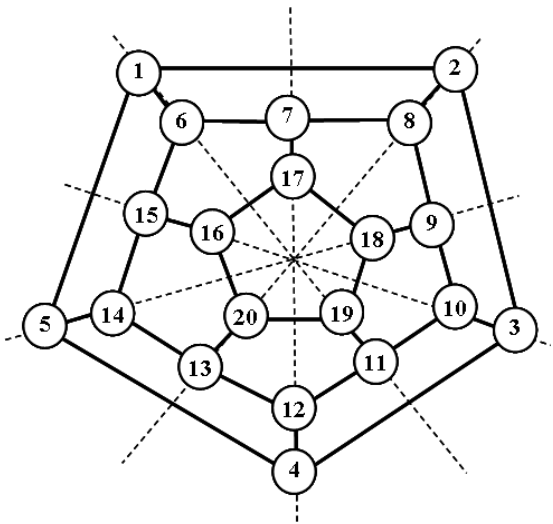

Рис. 2.124. Перестановка $p_0 =$
<1 2 3 4 5 6 7 8 9 10 11 12 13 14 15 16 17 18 19 20>.

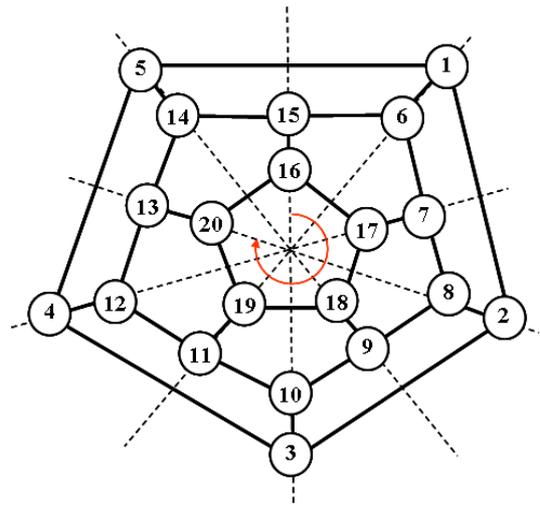

Рис. 2.125. Перестановка $p_1 =$
<5 1 2 3 4 14 15 6 7 8 9 10 11 12 13 20 16 17 18 19>.

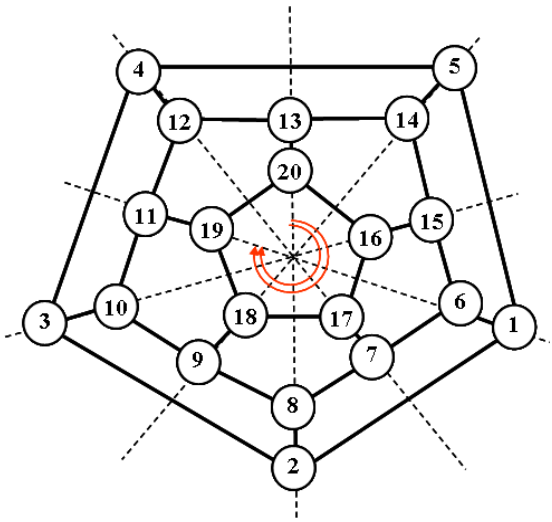
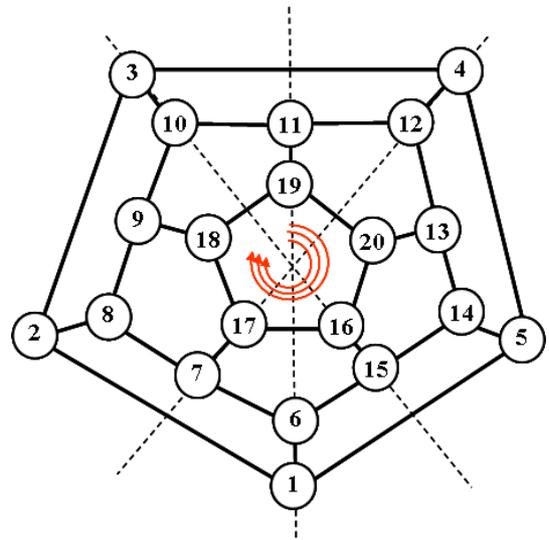

Рис. 2.126. Перестановка $p_2$ = 
<4 5 1 2 3 12 13 14 15 6 7 8 9 10 11 19 20 16 17 18>.

Рис. 2.127. Перестановка $p_3$ = 
<3 4 5 1 2 10 11 12 13 14 15 6 7 8 9 18 19 20 16 17>.

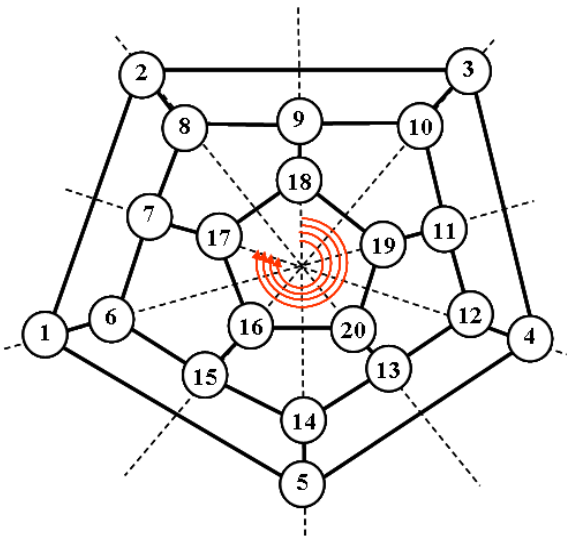
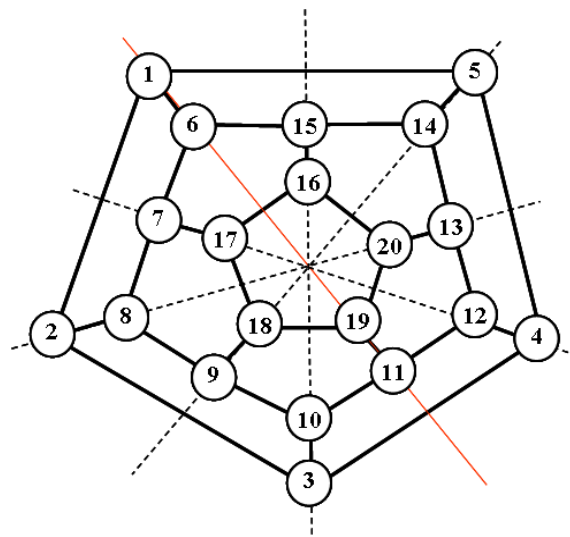

Рис. 2.128. Перестановка $p_4$ = 
<2 3 4 5 1 8 9 10 11 12 13 14 15 6 7 17 18 19 20 16>.

Рис. 2.129. Перестановка $p_5$ = 
<1 5 4 3 2 6 15 14 13 12 11 10 9 8 7 17 16 20 19 18>.

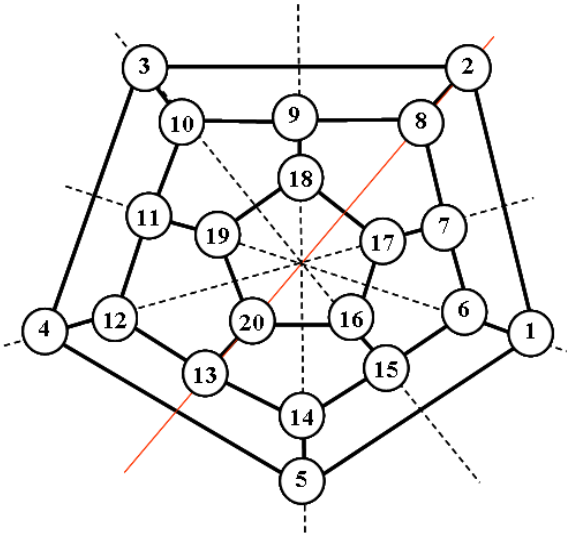
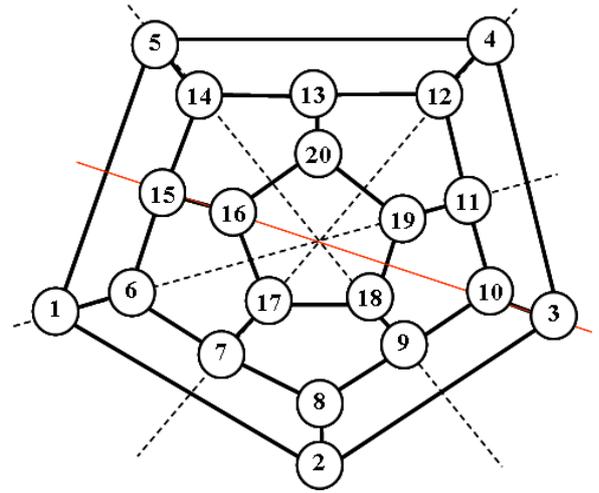

Рис. 2.130. Перестановка $p_6 =$
<3 2 1 5 4 10 9 8 7 6 15 14 13 12 11 19 18 17 16 20>.

Рис. 2.131. Перестановка $p_7 =$
<5 4 3 2 1 14 13 12 11 10 9 8 7 6 15 16 20 19 18 17>.

Количество образующих циклов в графе додекаэдра равно 12. Каждый образующий цикл порождает десять перестановок, то общее количество перестановок в графе додекаэдра определится как $12 \times 10 = 120$.

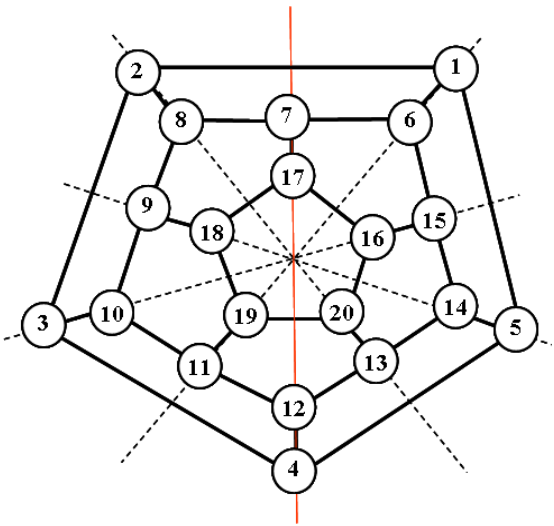
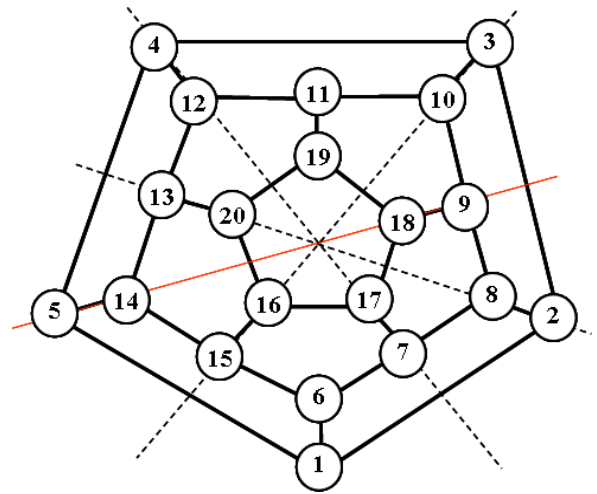

Рис. 2.132. Перестановка $p_8 =$
<2 1 5 4 3 8 7 6 15 14 13 12 11 10 9 18 17 16 20 19>.

Рис. 2.133. Перестановка $p_9 =$
<4 3 2 1 5 12 11 10 9 8 7 6 15 14 13 20 19 18 17 16>.

Первые десять перестановок графа додекаэдра индуцированные циклом $c_1$:

$p_0$ = <1 2 3 4 5 6 7 8 9 10 11 12 13 14 15 16 17 18 19 20> =
= (1)(2)(3)(4)(5)(6)(7)(8)(9)(10)(11)(12)(13)(14)(15)(16)(17)(18)(19)(20);
$p_1$ = <5 1 2 3 4 14 15 6 7 8 9 10 11 12 13 20 16 17 18 19> =
= (1 4 2 5 3)(6 14 12 10 8)(7 15 13 11 9)(16 20 19 18 17);
$p_2$ = <4 5 1 2 3 12 13 14 15 6 7 8 9 10 11 19 20 16 17 18> =
= (1 4 2 5 3)(6 12 8 14 10)(7 13 9 15 11)(16 19 17 20 18);

$p_3$ = <3 4 5 1 2 10 11 12 13 14 15 6 7 8 9 18 19 20 16 17> =
= (1 3 5 2 4)(6 10 14 8 12)(7 11 15 9 13)(16 18 20 17 19);
$p_4$ = <2 3 4 5 1 8 9 10 11 12 13 14 15 6 7 17 18 19 20 16> =
= (1 2 3 4 5)(6 8 10 12 14)(7 9 11 13 15)(16 17 18 19 20);
$p_5$ = <1 5 4 3 2 6 15 14 13 12 11 10 9 8 7 17 16 20 19 18> =
= (1)(2 5)(3 4)(6)(7 15)(8 14)(9 13)(10 12)(11)(16 17)(18 20)(19);
$p_6$ = <3 2 1 5 4 10 9 8 7 6 15 14 13 12 11 19 18 17 16 20> =
= (1 3)(2)(4 5)(6 10)(7 9)(8)(11 15)(12 14)(13)(16 19)(17 18)(20);
$p_7$ = <5 4 3 2 1 14 13 12 11 10 9 8 7 6 15 16 20 19 18 17> =
= (1 5)(2 4)(3)(6 14)(7 13)(8 12)(9 11)(10)(15)(16)(17 20)(18 19);
$p_8$ = <2 1 5 4 3 8 7 6 15 14 13 12 11 10 9 18 17 16 20 19> =
= (1 2)(3 5)(4)(6 8)(7)(9 15)(10 14)(11 13)(12)(16 18)(17)(19 20);
$p_9$ = <4 3 2 1 5 12 11 10 9 8 7 6 15 14 13 20 19 18 17 16> =
= (1 4)(2 3)(5)(6 12)(7 11)(8 10)(9)(13 15)(14)(16 20)(17 19)(18).

Остальные 110 перестановок можно получить в качестве упражнения используя вышеописанный метод.

Граф $G_9$ является планарным, кольцевая сумма изометрических циклов есть пустое множество $\sum_{i=1}^{12} c_i = \varnothing$.

Таблица Кэли для графа додекаэдра:

$Aut(G_9) = $

| $P_{1\,1}$ | $P_{1\,2}$ | $P_{1\,3}$ | $P_{1\,4}$ | … | $P_{1\,8}$ | $P_{1\,9}$ | $P_{1\,10}$ | $P_{1\,11}$ | $P_{1\,12}$ |
|---|---|---|---|---|---|---|---|---|---|
| $P_{2\,1}$ | $P_{2\,2}$ | $P_{2\,3}$ | $P_{2\,4}$ | … | $P_{2\,8}$ | $P_{2\,9}$ | $P_{2\,10}$ | $P_{2\,11}$ | $P_{2\,12}$ |
| $P_{3\,1}$ | $P_{3\,2}$ | $P_{3\,3}$ | $P_{3\,4}$ | … | $P_{3\,8}$ | $P_{3\,9}$ | $P_{3\,10}$ | $P_{3\,11}$ | $P_{3\,12}$ |
| $P_{4\,1}$ | $P_{4\,2}$ | $P_{4\,3}$ | $P_{4\,4}$ | … | $P_{4\,8}$ | $P_{4\,9}$ | $P_{4\,10}$ | $P_{4\,11}$ | $P_{4\,12}$ |
| $P_{5\,1}$ | $P_{5\,2}$ | $P_{5\,3}$ | $P_{5\,4}$ | … | $P_{5\,8}$ | $P_{5\,9}$ | $P_{5\,10}$ | $P_{5\,11}$ | $P_{5\,12}$ |
| $P_{6\,1}$ | $P_{6\,2}$ | $P_{6\,3}$ | $P_{6\,4}$ | … | $P_{6\,8}$ | $P_{6\,9}$ | $P_{6\,10}$ | $P_{6\,11}$ | $P_{6\,12}$ |
| $P_{7\,1}$ | $P_{7\,2}$ | $P_{7\,3}$ | $P_{7\,4}$ | … | $P_{7\,8}$ | $P_{7\,9}$ | $P_{7\,10}$ | $P_{7\,11}$ | $P_{7\,12}$ |
| $P_{8\,1}$ | $P_{8\,2}$ | $P_{8\,3}$ | $P_{8\,4}$ | … | $P_{8\,8}$ | $P_{8\,9}$ | $P_{8\,10}$ | $P_{8\,11}$ | $P_{8\,12}$ |
| $P_{9\,1}$ | $P_{9\,2}$ | $P_{9\,3}$ | $P_{9\,4}$ | … | $P_{9\,8}$ | $P_{9\,9}$ | $P_{9\,10}$ | $P_{9\,11}$ | $P_{9\,12}$ |
| $P_{10\,1}$ | $P_{10\,2}$ | $P_{10\,3}$ | $P_{10\,4}$ | … | $P_{10\,8}$ | $P_{10\,9}$ | $P_{10\,10}$ | $P_{10\,11}$ | $P_{10\,12}$ |
| $P_{11\,1}$ | $P_{11\,2}$ | $P_{11\,3}$ | $P_{11\,4}$ | … | $P_{11\,8}$ | $P_{11\,9}$ | $P_{11\,10}$ | $P_{11\,11}$ | $P_{11\,12}$ |
| $P_{12\,1}$ | $P_{12\,2}$ | $P_{12\,3}$ | $P_{12\,4}$ | … | $P_{12\,8}$ | $P_{12\,9}$ | $P_{12\,10}$ | $P_{12\,11}$ | $P_{12\,12}$ |

Таблица Кэли для первого блока $P_{1\,1}$:

$P_{1\,1} = $

|  | $p_0$ | $p_1$ | $p_2$ | $p_3$ | $p_4$ | $p_5$ | $p_6$ | $p_7$ | $p_8$ | $p_9$ |
|---|---|---|---|---|---|---|---|---|---|---|
| $p_0$ | $p_0$ | $p_1$ | $p_2$ | $p_3$ | $p_4$ | $p_5$ | $p_6$ | $p_7$ | $p_8$ | $p_9$ |
| $p_1$ | $p_1$ | $p_2$ | $p_3$ | $p_4$ | $p_0$ | $p_8$ | $p_9$ | $p_5$ | $p_6$ | $p_7$ |
| $p_2$ | $p_2$ | $p_3$ | $p_4$ | $p_0$ | $p_1$ | $p_6$ | $p_7$ | $p_8$ | $p_9$ | $p_5$ |
| $p_3$ | $p_3$ | $p_4$ | $p_0$ | $p_1$ | $p_2$ | $p_9$ | $p_5$ | $p_6$ | $p_7$ | $p_8$ |
| $p_4$ | $p_4$ | $p_0$ | $p_1$ | $p_2$ | $p_3$ | $p_7$ | $p_8$ | $p_9$ | $p_5$ | $p_6$ |
| $p_5$ | $p_5$ | $p_7$ | $p_9$ | $p_6$ | $p_8$ | $p_0$ | $p_3$ | $p_1$ | $p_4$ | $p_2$ |
| $p_6$ | $p_6$ | $p_8$ | $p_5$ | $p_7$ | $p_9$ | $p_2$ | $p_0$ | $p_3$ | $p_1$ | $p_4$ |
| $p_7$ | $p_7$ | $p_9$ | $p_6$ | $p_8$ | $p_5$ | $p_4$ | $p_2$ | $p_0$ | $p_3$ | $p_1$ |
| $p_8$ | $p_8$ | $p_5$ | $p_7$ | $p_9$ | $p_6$ | $p_1$ | $p_4$ | $p_2$ | $p_0$ | $p_3$ |
| $p_9$ | $p_9$ | $p_6$ | $p_8$ | $p_5$ | $p_7$ | $p_3$ | $p_1$ | $p_4$ | $p_2$ | $p_0$ |

## 2.7. Граф ИКОСАЭДРа

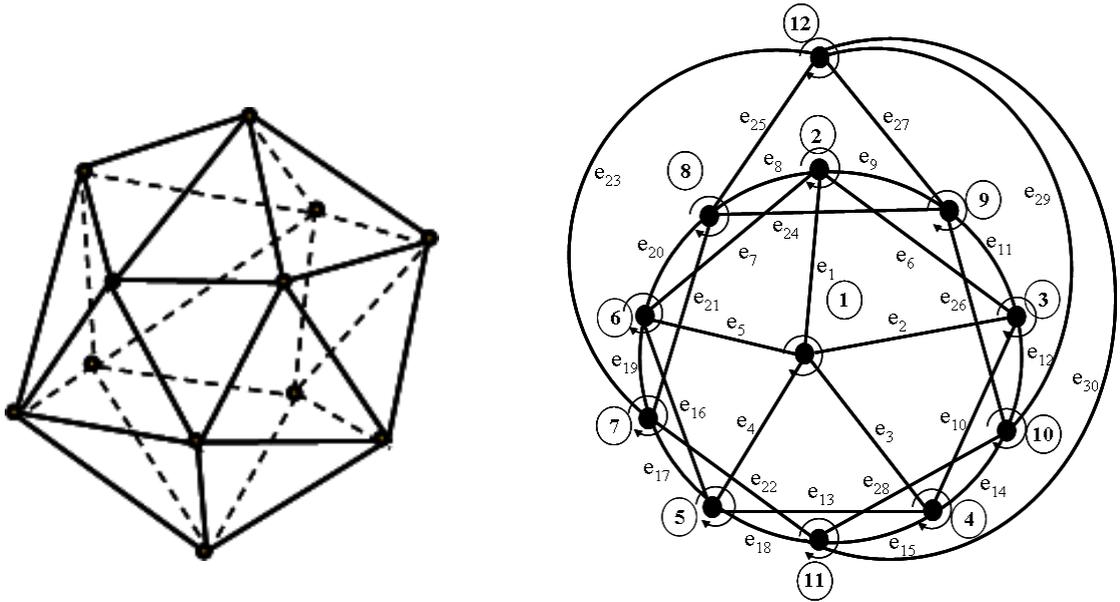

Рис. 2.134. Граф икосаэдра $G_{10}$ и его топологический рисунок.

Будем рассматривать граф икосаэдра.

Количество вершин графа = 12. Количество рёбер графа = 30

Количество изометрических циклов = 20

Смежность графа:

вершина $v_1$:  $v_2$  $v_3$  $v_4$  $v_5$  $v_6$;
вершина $v_2$:  $v_1$  $v_3$  $v_6$  $v_8$  $v_9$;
вершина $v_3$:  $v_1$  $v_2$  $v_4$  $v_9$  $v_{10}$;
вершина $v_4$:  $v_1$  $v_3$  $v_5$  $v_{10}$  $v_{11}$;
вершина $v_5$:  $v_1$  $v_4$  $v_6$  $v_7$  $v_{11}$;
вершина $v_6$:  $v_1$  $v_2$  $v_5$  $v_7$  $v_8$;
вершина $v_7$:  $v_5$  $v_6$  $v_8$  $v_{11}$  $v_{12}$;
вершина $v_8$:  $v_2$  $v_6$  $v_7$  $v_9$  $v_{12}$;
вершина $v_9$:  $v_2$  $v_3$  $v_8$  $v_{10}$  $v_{12}$;
вершина $v_{10}$:  $v_3$  $v_4$  $v_9$  $v_{11}$  $v_{12}$;
вершина $v_{11}$:  $v_4$  $v_5$  $v_7$  $v_{10}$  $v_{12}$;
вершина $v_{12}$:  $v_7$  $v_8$  $v_9$  $v_{10}$  $v_{11}$.

Инцидентность графа:

ребро $e_1$: $(v_1,v_2)$ или $(v_2,v_1)$;   ребро $e_2$: $(v_1,v_3)$ или $(v_3,v_1)$;
ребро $e_3$: $(v_1,v_4)$ или $(v_4,v_1)$;   ребро $e_4$: $(v_1,v_5)$ или $(v_5,v_1)$;
ребро $e_5$: $(v_1,v_6)$ или $(v_6,v_1)$;   ребро $e_6$: $(v_2,v_3)$ или $(v_3,v_2)$;
ребро $e_7$: $(v_2,v_6)$ или $(v_6,v_2)$;   ребро $e_8$: $(v_2,v_8)$ или $(v_8,v_2)$;
ребро $e_9$: $(v_2,v_9)$ или $(v_9,v_2)$;   ребро $e_{10}$: $(v_3,v_4)$ или $(v_4,v_3)$;
ребро $e_{11}$: $(v_3,v_9)$ или $(v_9,v_3)$;   ребро $e_{12}$: $(v_3,v_{10})$ или $(v_{10},v_3)$;
ребро $e_{13}$: $(v_4,v_5)$ или $(v_5,v_4)$;   ребро $e_{14}$: $(v_4,v_{10})$ или $(v_{10},v_4)$;
ребро $e_{15}$: $(v_4,v_{11})$ или $(v_{11},v_4)$;   ребро $e_{16}$: $(v_5,v_6)$ или $(v_6,v_5)$;

ребро $e_{17}$: $(v_5,v_7)$ или $(v_7,v_5)$; ребро $e_{18}$: $(v_5,v_{11})$ или $(v_{11},v_5)$;
ребро $e_{19}$: $(v_6,v_7)$ или $(v_7,v_6)$; ребро $e_{20}$: $(v_6,v_8)$ или $(v_8,v_6)$;
ребро $e_{21}$: $(v_7,v_8)$ или $(v_8,v_7)$; ребро $e_{22}$: $(v_7,v_{11})$ или $(v_{11},v_7)$;
ребро $e_{23}$: $(v_7,v_{12})$ или $(v_{12},v_7)$; ребро $e_{24}$: $(v_8,v_9)$ или $(v_9,v_8)$;
ребро $e_{25}$: $(v_8,v_{12})$ или $(v_{12},v_8)$; ребро $e_{26}$: $(v_9,v_{10})$ или $(v_{10},v_9)$;
ребро $e_{27}$: $(v_9,v_{12})$ или $(v_{12},v_9)$; ребро $e_{28}$: $(v_{10},v_{11})$ или $(v_{11},v_{10})$;
ребро $e_{29}$: $(v_{10},v_{12})$ или $(v_{12},v_{10})$; ребро $e_{30}$: $(v_{11},v_{12})$ или $(v_{12},v_{11})$.

Множество изометрических циклов графа:

$c_1 = \{e_1,e_2,e_6\} \to \{v_1,v_2,v_3\}$;  $c_2 = \{e_1,e_5,e_7\} \to \{v_1,v_2,v_6\}$;
$c_3 = \{e_2,e_3,e_{10}\} \to \{v_1,v_3,v_4\}$;  $c_4 = \{e_3,e_4,e_{13}\} \to \{v_1,v_4,v_5\}$;
$c_5 = \{e_4,e_5,e_{16}\} \to \{v_1,v_5,v_6\}$;  $c_6 = \{e_6,e_9,e_{11}\} \to \{v_2,v_3,v_9\}$;
$c_7 = \{e_7,e_8,e_{20}\} \to \{v_2,v_6,v_8\}$;  $c_8 = \{e_8,e_9,e_{24}\} \to \{v_2,v_8,v_9\}$;
$c_9 = \{e_{10},e_{12},e_{14}\} \to \{v_3,v_4,v_{10}\}$;  $c_{10} = \{e_{11},e_{12},e_{26}\} \to \{v_3,v_9,v_{10}\}$;
$c_{11} = \{e_{13},e_{15},e_{18}\} \to \{v_4,v_5,v_{11}\}$;  $c_{12} = \{e_{14},e_{15},e_{28}\} \to \{v_4,v_{10},v_{11}\}$;
$c_{13} = \{e_{16},e_{17},e_{19}\} \to \{v_5,v_6,v_7\}$;  $c_{14} = \{e_{17},e_{18},e_{22}\} \to \{v_5,v_7,v_{11}\}$;
$c_{15} = \{e_{19},e_{20},e_{21}\} \to \{v_6,v_7,v_8\}$;  $c_{16} = \{e_{21},e_{23},e_{25}\} \to \{v_7,v_8,v_{12}\}$;
$c_{17} = \{e_{22},e_{23},e_{30}\} \to \{v_7,v_{11},v_{12}\}$;  $c_{18} = \{e_{24},e_{25},e_{27}\} \to \{v_8,v_9,v_{12}\}$;
$c_{19} = \{e_{26},e_{27},e_{29}\} \to \{v_9,v_{10},v_{12}\}$;  $c_{20} = \{e_{28},e_{29},e_{30}\} \to \{v_{10},v_{11},v_{12}\}$.

Будем формировать образующие изометрические циклы. Начнем формирование с циклов окружающих соответствующие вершины:

$c_{v1} = c_1 \oplus c_2 \oplus c_3 \oplus c_4 \oplus c_5 = \{e_6,e_7,e_{10},e_{13},e_{16}\} \to \{v_2,v_3,v_4,v_5,v_6\}$;
$c_{v2} = c_1 \oplus c_2 \oplus c_6 \oplus c_7 \oplus c_8 = \{e_2,e_5,e_{11},e_{20},e_{24}\} \to \{v_1,v_3,v_6,v_8,v_9\}$;
$c_{v3} = c_1 \oplus c_3 \oplus c_6 \oplus c_9 \oplus c_{10} = \{e_1,e_3,e_9,e_{14},e_{26}\} \to \{v_1,v_2,v_4,v_9,v_{10}\}$;
$c_{v4} = c_3 \oplus c_4 \oplus c_9 \oplus c_{11} \oplus c_{12} = \{e_2,e_4,e_{12},e_{18},e_{28}\} \to \{v_1,v_3,v_5,v_{10},v_{11}\}$;
$c_{v5} = c_4 \oplus c_5 \oplus c_{11} \oplus c_{13} \oplus c_{14} = \{e_3,e_5,e_{15},e_{19},e_{22}\} \to \{v_1,v_4,v_6,v_7,v_{11}\}$;
$c_{v6} = c_2 \oplus c_5 \oplus c_7 \oplus c_{13} \oplus c_{15} = \{e_1,e_4,e_8,e_{17},e_{21}\} \to \{v_1,v_2,v_5,v_7,v_8\}$;
$c_{v7} = c_{13} \oplus c_{14} \oplus c_{15} \oplus c_{16} \oplus c_{17} = \{e_{16},e_{18},e_{20},e_{25},e_{30}\} \to \{v_5,v_6,v_8,v_{11},v_{12}\}$;
$c_{v8} = c_7 \oplus c_8 \oplus c_{15} \oplus c_{16} \oplus c_{18} = \{e_7,e_9,e_{19},e_{23},e_{27}\} \to \{v_2,v_6,v_7,v_9,v_{12}\}$;
$c_{v9} = c_6 \oplus c_8 \oplus c_{10} \oplus c_{18} \oplus c_{19} = \{e_6,e_8,e_{12},e_{25},e_{29}\} \to \{v_2,v_3,v_8,v_{10},v_{12}\}$;
$c_{v10} = c_9 \oplus c_{10} \oplus c_{12} \oplus c_{19} \oplus c_{20} = \{e_{10},e_{11},e_{15},e_{27},e_{30}\} \to \{v_3,v_4,v_9,v_{11},v_{12}\}$;
$c_{v11} = c_{11} \oplus c_{12} \oplus c_{14} \oplus c_{17} \oplus c_{20} = \{e_{13},e_{14},e_{17},e_{23},e_{29}\} \to \{v_4,v_5,v_7,v_{10},v_{12}\}$;
$c_{v12} = c_{16} \oplus c_{17} \oplus c_{18} \oplus c_{19} \oplus c_{20} = \{e_{21},e_{22},e_{24},e_{26},e_{28}\} \to \{v_7,v_8,v_9,v_{10},v_{11}\}$.

Инвариант интегрального инварианта спектра реберных разрезов состоит из:

кортежа весов ребер $\xi_w(G_{10}) = (30 \times 24)$, и кортежа весов вершин: $\zeta_w(G_{10}) = (12 \times 120)$.

С учетом симметрии, концевые вершины в топологическом рисунке будут располагаться парами. Составим список этих пар: $(v_1,v_{12})$, $(v_2,v_{11})$, $(v_3,v_7)$, $(v_4,v_8)$, $(v_5,v_9)$, $(v_6,v_{10})$. Нам необходимо создать топологический рисунок в котором присутствуют все вершины графа, и все изометрические циклы графа. И с учетом этого выбрать образующие циклы для построения группы автомормизма графа. Такие циклы должны быть образованы относительно пар вершин.

Рассмотрим формирование образующего цикла для пары вершин $(v_1,v_{12})$.

$c_{\text{обр}1} = c_{v1} \oplus c_6 \oplus c_7 \oplus c_9 \oplus c_{11} \oplus c_{13} = \{e_6,e_7,e_{10},e_{13},e_{16}\} \oplus \{e_6,e_9,e_{11}\} \oplus \{e_7,e_8,e_{20}\} \oplus$
$\oplus \{e_{13},e_{15},e_{18}\} \oplus \{e_{10},e_{12},e_{14}\} \oplus \{e_{16},e_{17},e_{19}\} = \{e_8,e_9,e_{11},e_{12},e_{14},e_{15},e_{17},e_{18},e_{19},e_{20}\}$;

$c_{обр1} = c_{v12} \oplus c_8 \oplus c_{10} \oplus c_{12} \oplus c_{14} \oplus c_{15} = \{e_{21}, e_{22}, e_{24}, e_{26}, e_{28}\} \oplus \{e_8, e_9, e_{24}\} \oplus \{e_{11}, e_{12}, e_{26}\} \oplus$
$\oplus \{e_{14}, e_{15}, e_{28}\} \oplus \{e_{17}, e_{18}, e_{22}\} \oplus \{e_{19}, e_{20}, e_{21}\} = \{e_8, e_9, e_{11}, e_{12}, e_{14}, e_{15}, e_{17}, e_{18}, e_{19}, e_{20}\}.$

Образующий цикл $c_{обр1}$ в векторной записи (см. рис. 2.135):
$c_{обр1} = (v_2, v_9) + (v_9, v_3) + (v_3, v_{10}) + (v_{10}, v_4) + (v_4, v_{11}) + (v_{11}, v_5) + (v_5, v_7) + (v_7, v_6) + (v_6, v_8) + (v_8, v_2).$

Построим образующий цикл для пары вершин $(v_2, v_{11})$.

$c_{обр2} = c_{v2} \oplus c_3 \oplus c_5 \oplus c_{10} \oplus c_{15} \oplus c_{18} = \{e_2, e_5, e_{11}, e_{20}, e_{24}\} \oplus \{e_2, e_3, e_{10}\} \oplus \{e_4, e_5, e_{16}\} \oplus$
$\oplus \{e_{11}, e_{12}, e_{26}\} \oplus \{e_{19}, e_{20}, e_{21}\} \oplus \{e_{24}, e_{25}, e_{27}\} = \{e_3, e_4, e_{10}, e_{12}, e_{16}, e_{19}, e_{21}, e_{25}, e_{26}, e_{27}\};$
$c_{обр2} = c_{v12} \oplus c_{11} \oplus c_{12} \oplus c_{14} \oplus c_{17} \oplus c_{20} = \{e_{13}, e_{14}, e_{17}, e_{23}, e_{29}\} \oplus \{e_3, e_4, e_{13}\} \oplus \{e_{10}, e_{12}, e_{14}\} \oplus$
$\oplus \{e_{16}, e_{17}, e_{19}\} \oplus \{e_{21}, e_{23}, e_{25}\} \oplus \{e_{26}, e_{27}, e_{29}\} = \{e_3, e_4, e_{10}, e_{12}, e_{16}, e_{19}, e_{21}, e_{25}, e_{26}, e_{27}\}.$

Образующий цикл $c_{обр2}$ в векторной записи (см. рис. 2.136):
$c_{обр2} = (v_1, v_4) + (v_4, v_3) + (v_3, v_{10}) + (v_{10}, v_9) + (v_9, v_{12}) + (v_{12}, v_8) + (v_8, v_7) + (v_7, v_6) + (v_6, v_5) + (v_5, v_1).$

Создадим образующий цикл для пары вершин $(v_3, v_7)$.

$c_{обр3} = c_{v3} \oplus c_2 \oplus c_4 \oplus c_8 \oplus c_{13} \oplus c_{19} = \{e_1, e_3, e_9, e_{14}, e_{26}\} \oplus \{e_1, e_5, e_7\} \oplus \{e_3, e_4, e_{13}\} \oplus$
$\oplus \{e_8, e_9, e_{24}\} \oplus \{e_{14}, e_{15}, e_{28}\} \oplus \{e_{26}, e_{27}, e_{29}\} = \{e_4, e_5, e_7, e_8, e_{13}, e_{15}, e_{24}, e_{27}, e_{28}, e_{29}\};$
$c_{обр3} = c_{v7} \oplus c_5 \oplus c_7 \oplus c_{11} \oplus c_{18} \oplus c_{20} = \{e_{16}, e_{18}, e_{20}, e_{25}, e_{30}\} \oplus \{e_4, e_5, e_{16}\} \oplus \{e_7, e_8, e_{20}\} \oplus$
$\oplus \{e_{13}, e_{15}, e_{18}\} \oplus \{e_{24}, e_{25}, e_{27}\} \oplus \{e_{28}, e_{29}, e_{30}\} = \{e_4, e_5, e_7, e_8, e_{13}, e_{15}, e_{24}, e_{27}, e_{28}, e_{29}\}.$

Образующий цикл $c_{обр3}$ в векторной записи (см. рис. 2.137):
$c_{обр3} = (v_1, v_6) + (v_6, v_2) + (v_2, v_8) + (v_8, v_9) + (v_9, v_{12}) + (v_{12}, v_{10}) + (v_{10}, v_{11}) + (v_{11}, v_4) + (v_4, v_5) + (v_5, v_1).$

Сформируем образующий цикл для пары вершин $(v_4, v_8)$.

$c_{обр4} = c_{v4} \oplus c_1 \oplus c_5 \oplus c_{10} \oplus c_{14} \oplus c_{20} = \{e_2, e_4, e_{12}, e_{18}, e_{28}\} \oplus \{e_1, e_2, e_6\} \oplus \{e_4, e_5, e_{16}\} \oplus$
$\oplus \{e_{11}, e_{12}, e_{26}\} \oplus \{e_{17}, e_{18}, e_{22}\} \oplus \{e_{28}, e_{29}, e_{30}\} = \{e_1, e_5, e_6, e_{11}, e_{16}, e_{17}, e_{22}, e_{26}, e_{29}, e_{30}\};$
$c_{обр4} = c_{v8} \oplus c_2 \oplus c_6 \oplus c_{13} \oplus c_{17} \oplus c_{19} = \{e_7, e_9, e_{19}, e_{23}, e_{27}\} \oplus \{e_1, e_5, e_7\} \oplus \{e_6, e_9, e_{11}\} \oplus$
$\oplus \{e_{16}, e_{17}, e_{19}\} \oplus \{e_{22}, e_{23}, e_{30}\} \oplus \{e_{26}, e_{27}, e_{29}\} = \{e_1, e_5, e_6, e_{11}, e_{16}, e_{17}, e_{22}, e_{26}, e_{29}, e_{30}\}.$

Образующий цикл $c_{обр4}$ в векторной записи (см. рис. 2.138):
$c_{обр4} = (v_1, v_2) + (v_2, v_3) + (v_3, v_9) + (v_9, v_{10}) + (v_{10}, v_{12}) + (v_{12}, v_{11}) + (v_{11}, v_7) + (v_7, v_5) + (v_5, v_6) + (v_6, v_1).$

Рассмотрим процесс формирования образующего цикла, для пары вершин $(v_5, v_9)$.

$c_{обр5} = c_{v5} \oplus c_2 \oplus c_3 \oplus c_{12} \oplus c_{15} \oplus c_{17} = \{e_2, e_4, e_{12}, e_{18}, e_{28}\} \oplus \{e_1, e_5, e_7\} \oplus \{e_2, e_3, e_{10}\} \oplus$
$\oplus \{e_{14}, e_{15}, e_{28}\} \oplus \{e_{19}, e_{20}, e_{21}\} \oplus \{e_{22}, e_{23}, e_{30}\} = \{e_1, e_2, e_7, e_{10}, e_{14}, e_{20}, e_{21}, e_{23}, e_{28}, e_{30}\};$
$c_{обр5} = c_{v9} \oplus c_1 \oplus c_7 \oplus c_9 \oplus c_{16} \oplus c_{20} = \{e_6, e_8, e_{12}, e_{25}, e_{29}\} \oplus \{e_1, e_2, e_6\} \oplus \{e_7, e_8, e_{20}\} \oplus$
$\oplus \{e_{10}, e_{12}, e_{14}\} \oplus \{e_{21}, e_{23}, e_{25}\} \oplus \{e_{28}, e_{29}, e_{30}\} = \{e_1, e_2, e_7, e_{10}, e_{14}, e_{20}, e_{21}, e_{23}, e_{28}, e_{30}\}.$

Образующий цикл $c_{обр5}$ в векторной записи (см. рис. 2.139):
$c_{обр5} = (v_1, v_3) + (v_3, v_4) + (v_4, v_{10}) + (v_{10}, v_{11}) + (v_{11}, v_{12}) + (v_{12}, v_7) + (v_7, v_8) + (v_8, v_6) + (v_6, v_2) + (v_2, v_1).$

Кольцевая сумма изометрических циклов формирует образующий цикл для пары вершин $(v_6, v_{10})$.

$c_{обр6} = c_{v6} \oplus c_1 \oplus c_4 \oplus c_8 \oplus c_{14} \oplus c_{16} = \{e_1, e_4, e_8, e_{17}, e_{21}\} \oplus \{e_1, e_2, e_6\} \oplus \{e_3, e_4, e_{13}\} \oplus$
$\oplus \{e_8, e_9, e_{24}\} \oplus \{e_{17}, e_{18}, e_{22}\} \oplus \{e_{21}, e_{23}, e_{25}\} = \{e_2, e_3, e_6, e_9, e_{13}, e_{18}, e_{22}, e_{23}, e_{24}, e_{25}\};$
$c_{обр6} = c_{v10} \oplus c_3 \oplus c_6 \oplus c_{11} \oplus c_{17} \oplus c_{18} = \{e_{10}, e_{11}, e_{15}, e_{27}, e_{30}\} \oplus \{e_2, e_3, e_{10}\} \oplus \{e_6, e_9, e_{11}\} \oplus$
$\oplus \{e_{13}, e_{15}, e_{18}\} \oplus \{e_{22}, e_{23}, e_{30}\} \oplus \{e_{24}, e_{25}, e_{27}\} = \{e_2, e_3, e_6, e_9, e_{13}, e_{18}, e_{22}, e_{23}, e_{24}, e_{25}\}.$

Образующий цикл $c_{обр2}$ в векторной записи (см. рис. 2.140):
$c_{обр6} = (v_1, v_3) + (v_3, v_2) + (v_2, v_9) + (v_9, v_8) + (v_8, v_{12}) + (v_{12}, v_7) + (v_7, v_{11}) + (v_{11}, v_5) + (v_5, v_4) + (v_4, v_1).$

Будем индуцировать перестановки вершин, выбрав в качестве симметричного изображения правильный многоугольник с десятью сторонами. В качестве образующего цикла выберем $c_{обр1}$ = {$e_8,e_9,e_{11},e_{12},e_{14},e_{15},e_{17},e_{18},e_{19},e_{20}$}. Для наглядности процесса порождения перестановок, поместим в середину вершины $v_1, v_{12}$.

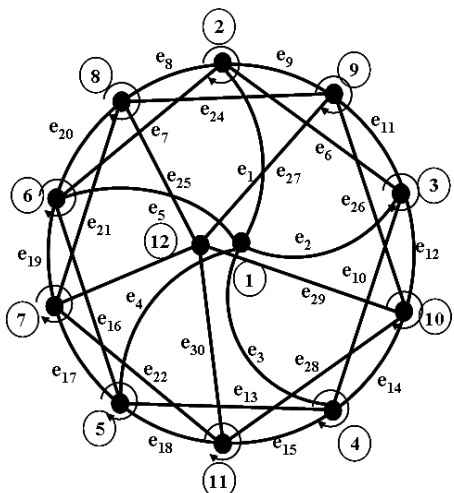
Рис. 2.135. 1-ый образующий цикл.

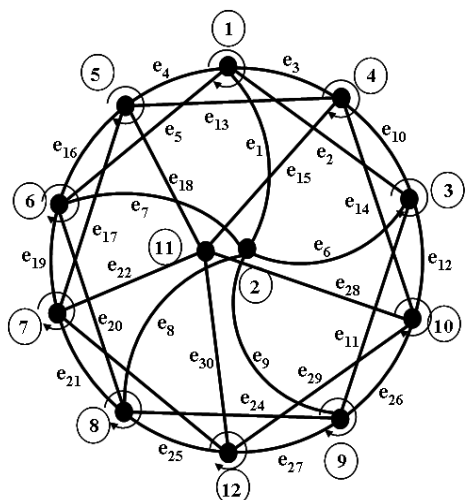
Рис. 2.136. 2-ой образующий цикл.

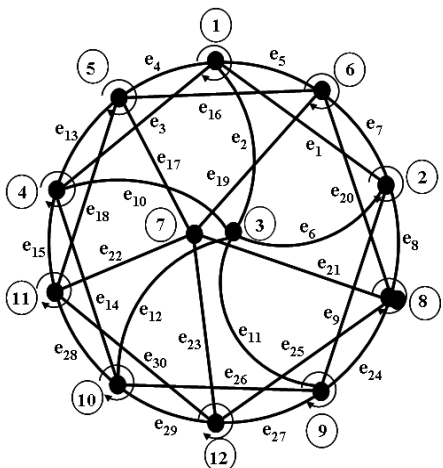
Рис. 2.137. 3-ий образующий цикл.

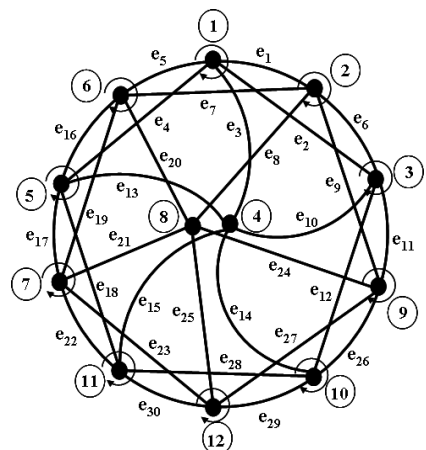
Рис. 2.138. 4-ый образующий цикл.

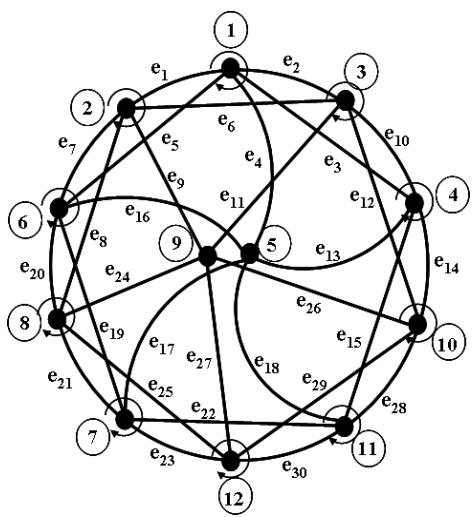
Рис. 2.139. 5-ый образующий цикл.

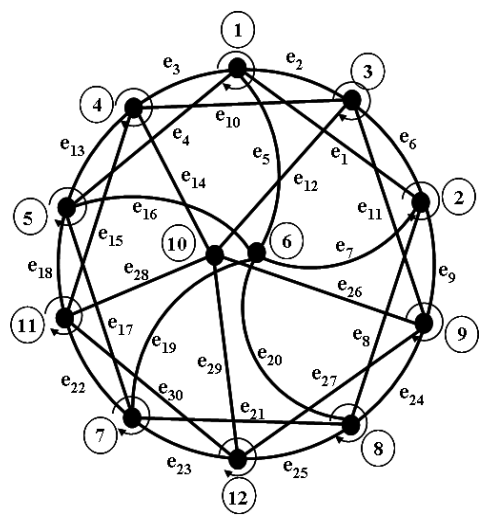
Рис. 2.140. 6-ой образующий цикл.

Будем рассматривать перестановки для заданного образующего цикла.

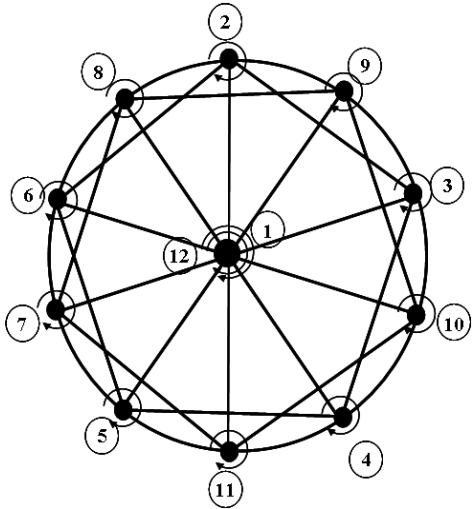

Рис. 2.141. Перестановка $p_0$.

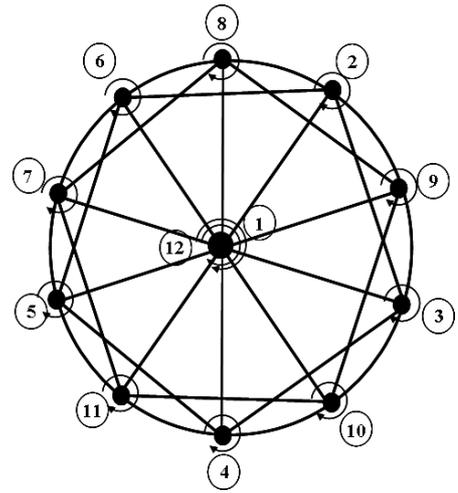

Рис. 2.142. Перестановка $p_1$.

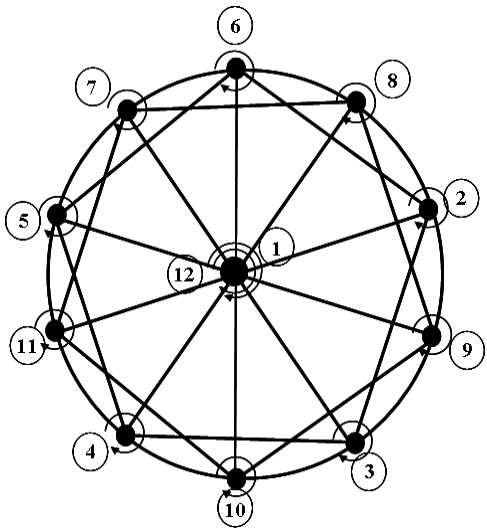

Рис. 2.143. Перестановка $p_2$.

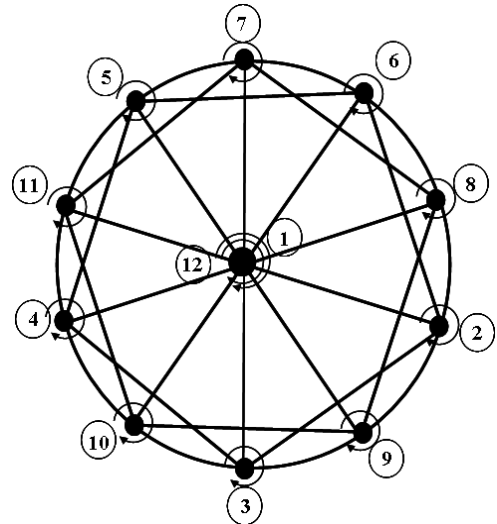

Рис. 2.144. Перестановка $p_3$.

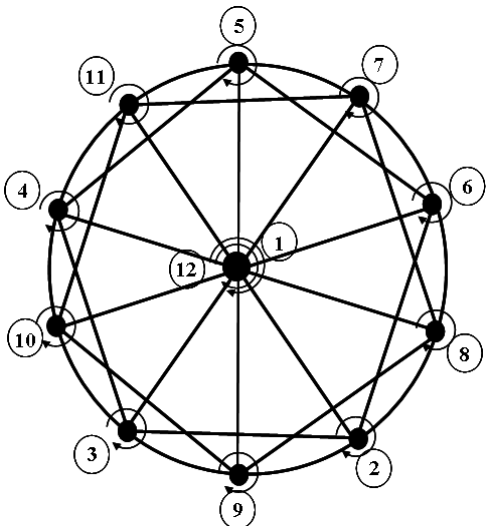

Рис. 2.145. Перестановка $p_4$.

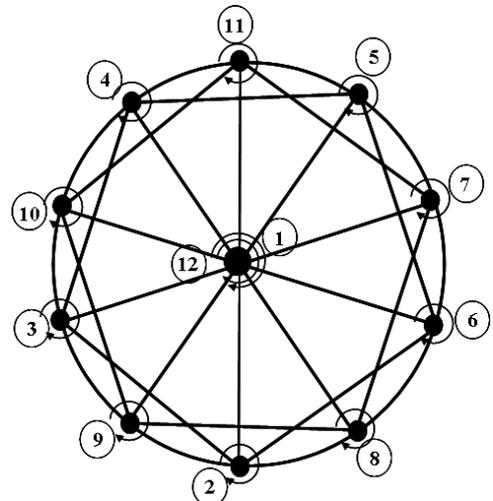

Рис. 2.146. Перестановка $p_5$.

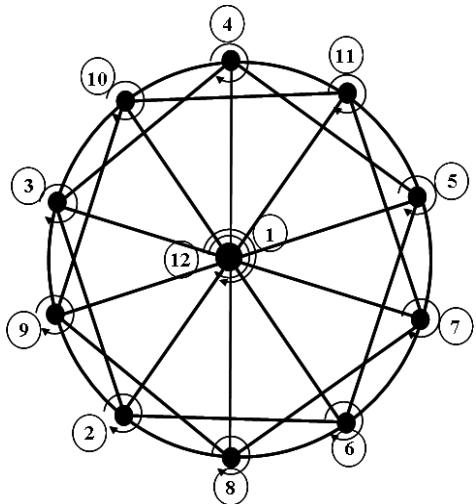

Рис. 2.147. Перестановка p_6.

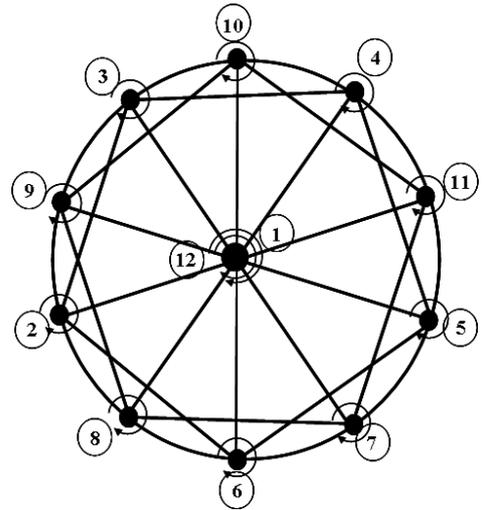

Рис. 2.148. Перестановка p_7.

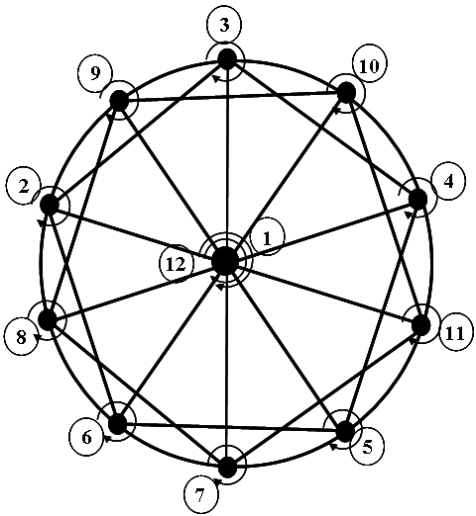

Рис. 2.149. Перестановка p_8.

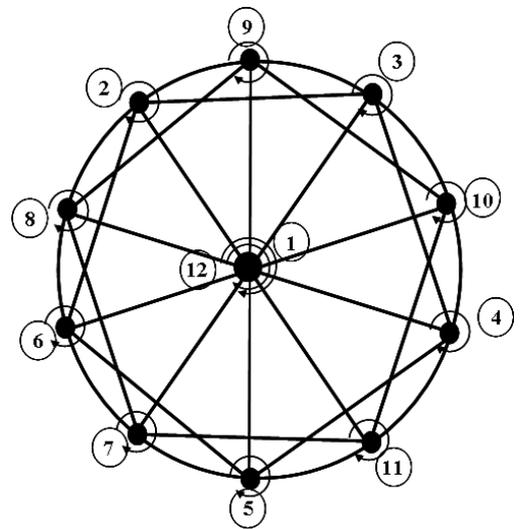

Рис. 2.150. Перестановка p_9.

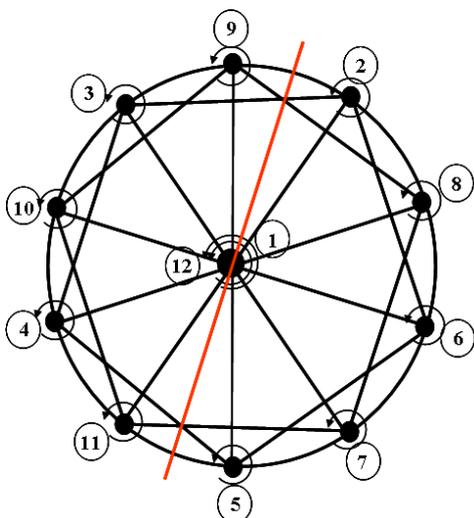

Рис. 2.151. Перестановка p_10.

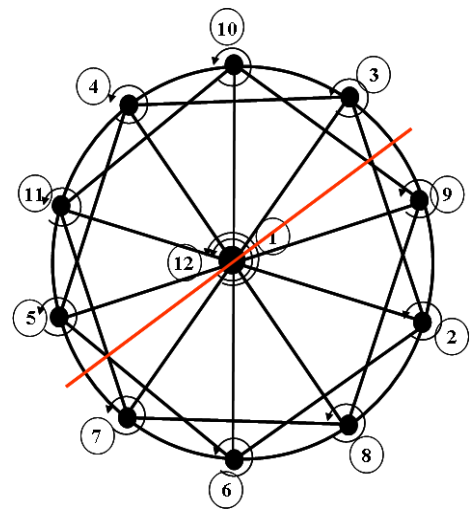

Рис. 2.152. Перестановка p_11.

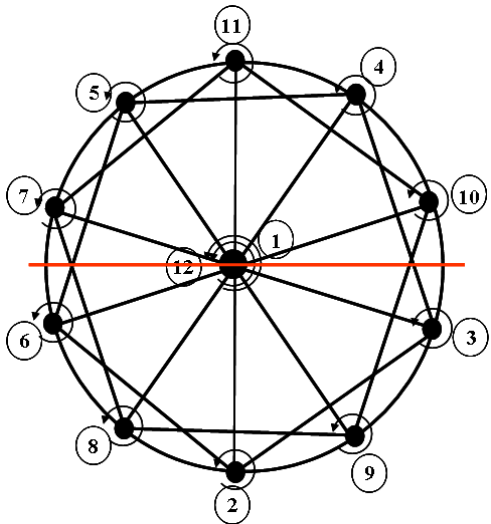

Рис. 2.153. Перестановка p$_{12}$.

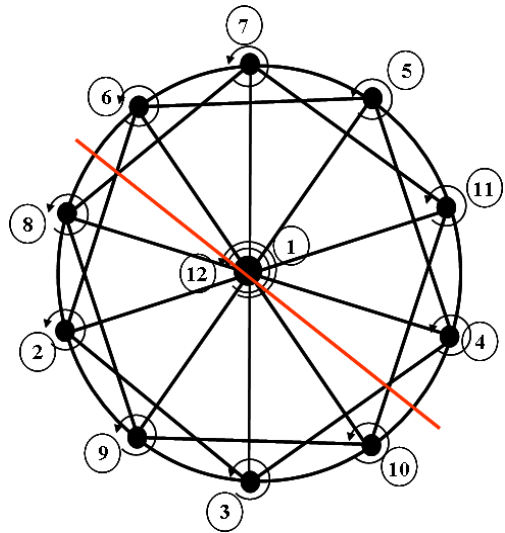

Рис. 2.154. Перестановка p$_{13}$.

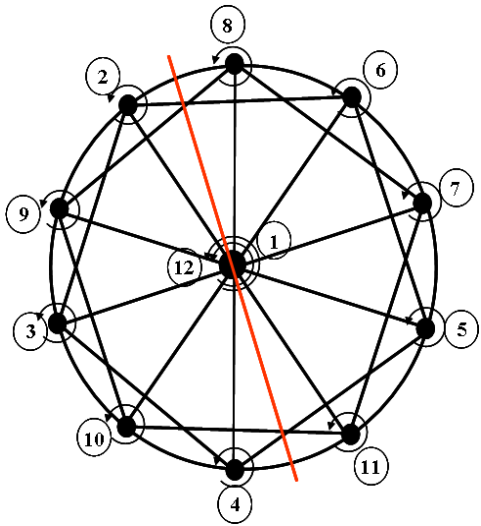

Рис. 2.155. Перестановка p$_{14}$.

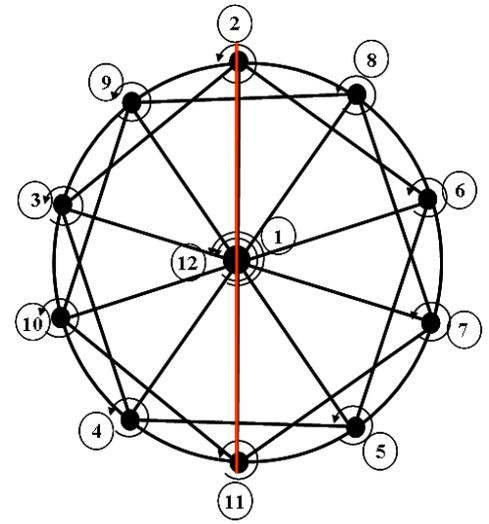

Рис. 2.156. Перестановка p$_{15}$.

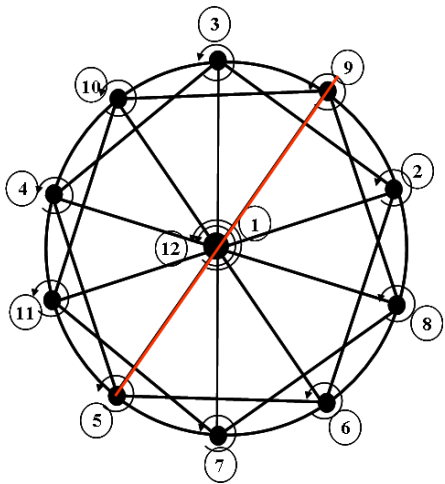

Рис. 2.157. Перестановка p$_{16}$.

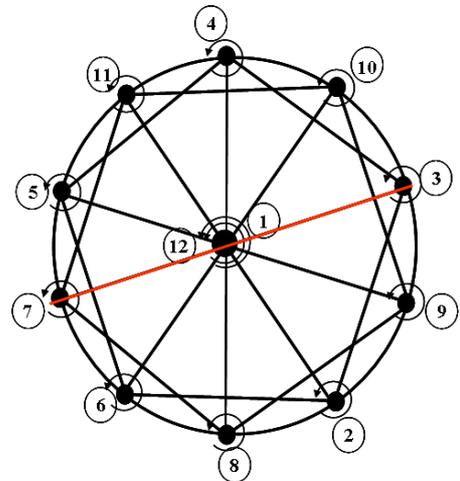

Рис. 2.158. Перестановка p$_{17}$.

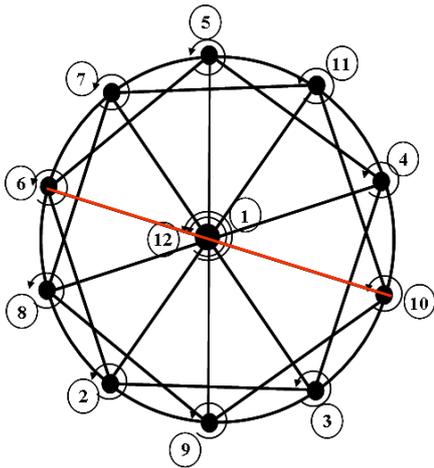 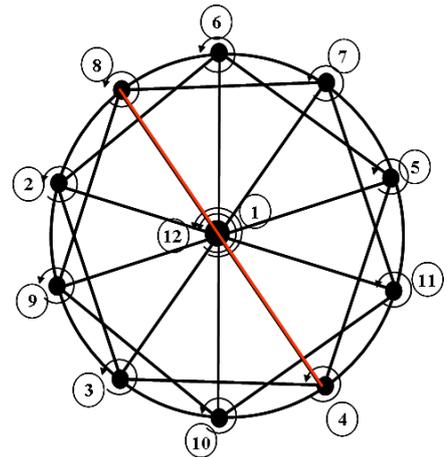

Рис. 2.159. Перестановка $p_{18}$.          Рис. 2.160. Перестановка $p_{19}$.

Так как, количество образующих циклов равно 6, и каждый образующий цикл порождает 20 перестановок, то общее количество перестановок равно 6×20 = 120.

Рассмотрим перестановки для 1-го образующегося цикла графа икосаэдра, приняв за начальное расположение вершин топологический рисунок графа 2.134:

$p_0$ = <1 2 3 4 5 6 7 8 9 10 11 12> = (1)(2)(3)(4)(5)(6)(7)(8)(9)(10)(11)(12);
$p_1$ = <1 8 9 10 11 7 5 6 2 3 4 12> = (1)(8 6 7 5 11 4 10 3 9 2)(12);
$p_2$ = <1 6 2 3 4 5 11 7 8 9 10 12> = (1)(6 5 4 3 2)(11 9 8 7)(12);
$p_3$ = <1 7 8 9 10 11 4 5 6 2 3 12> = (1)(7 4 9 6 11 3 8 5 10 2)(12);
$p_4$ = <1 5 6 2 3 4 10 11 7 8 9 12> = (1)(5 3 6 4 2)(10 8 11 9 7)(12);
$p_5$ = <1 11 7 8 9 10 3 4 5 6 2 12> = (1)(11 2)(7 3)(8 4)(9 5)(10 6)(12);
$p_6$ = <1 4 5 6 2 3 9 10 11 7 8 12> = (1)(4 6 3 5 2)(9 11 8 10 7)(12);
$p_7$ = <1 10 11 7 8 9 2 3 4 5 6 12> = (1)(10 5 8 3 11 6 9 4 7 2)(12);
$p_8$ = <1 3 4 5 6 2 8 9 10 11 7 12> = (1)(3 4 5 6 2)(8 9 10 11 7)(12);
$p_9$ = <1 9 10 11 7 8 6 2 3 4 5 12> = (1)(9 3 10 4 11 5 7 6 8 2)(12);
$p_{10}$ = <1 9 8 7 11 10 4 3 2 6 5 12> = (1)(9 2)(8 3)(7 4)(11 5)(10 6)(12);
$p_{11}$ = <1 10 9 8 7 11 5 4 3 2 6 12> = (1)(10 2)(9 3)(8 4)(7 5)(11 6)(12);
$p_{12}$ = <1 11 10 9 8 7 6 5 4 3 2 12> = (1)(11 2 5 3 6 4)(7 10 8 11 9)(12);
$p_{13}$ = <1 7 11 10 9 8 2 6 5 4 3 12> = (1)(7 2)(11 3)(10 4)(9 5)(8 6)(12);
$p_{14}$ = <1 8 7 11 10 9 3 2 6 5 4 12> = (1)(8 2)(7 3)(11 4)(10 5)(9 6)(12);
$p_{15}$ = <1 2 6 5 4 3 10 9 8 7 11 12> = (1)(2)(6 3)(5 4)(10 7)(9 8)(11)(12);
$p_{16}$ = <1 3 2 6 5 4 11 10 9 8 7 12> = (1)(3 2)(6 4)(5)(11 7)(10 8)(9)(12);
$p_{17}$ = <1 4 3 2 6 5 7 11 10 9 8 12> = (1)(4 2)(3)(6 5)(7)(11 8)(10 9)(12);
$p_{18}$ = <1 5 4 3 2 6 8 7 11 10 9 12> = (1)(5 2)(4 3)(6)(8 7)(11 9)((10)(12);
$p_{19}$ = <1 6 5 4 3 2 9 8 7 11 10 12> = (1)(6 2)(5 3)(4)(9 7)(8)(11 10)(12).

Общий вид таблица Кэли для графа икосаэдра:

$Aut(G_{10}) = $

| $P_{1\,1}$ | $P_{1\,2}$ | $P_{1\,3}$ | $P_{1\,4}$ | $P_{1\,5}$ | $P_{1\,6}$ |
|---|---|---|---|---|---|
| $P_{2\,1}$ | $P_{2\,2}$ | $P_{2\,3}$ | $P_{2\,4}$ | $P_{2\,5}$ | $P_{2\,6}$ |
| $P_{3\,1}$ | $P_{3\,2}$ | $P_{3\,3}$ | $P_{3\,4}$ | $P_{3\,5}$ | $P_{3\,6}$ |
| $P_{4\,1}$ | $P_{4\,2}$ | $P_{4\,3}$ | $P_{4\,4}$ | $P_{4\,5}$ | $P_{4\,6}$ |
| $P_{5\,1}$ | $P_{5\,2}$ | $P_{5\,3}$ | $P_{5\,4}$ | $P_{5\,5}$ | $P_{5\,6}$ |
| $P_{6\,1}$ | $P_{6\,2}$ | $P_{6\,3}$ | $P_{6\,4}$ | $P_{6\,5}$ | $P_{6\,6}$ |

Таблица Кэли для первого блока $P_{1\ 1}$ имеет вид:

|  | $p_0$ | $p_1$ | $p_2$ | $p_3$ | $p_4$ | $p_5$ | $p_6$ | $p_7$ | $p_8$ | $p_9$ | $p_{10}$ | $p_{11}$ | $p_{12}$ | $p_{13}$ | $p_{14}$ | $p_{15}$ | $p_{16}$ | $p_{17}$ | $p_{18}$ | $p_{19}$ |
|---|---|---|---|---|---|---|---|---|---|---|---|---|---|---|---|---|---|---|---|---|
| $p_0$ | $p_0$ | $p_1$ | $p_2$ | $p_3$ | $p_4$ | $p_5$ | $p_6$ | $p_7$ | $p_8$ | $p_9$ | $p_{10}$ | $p_{11}$ | $p_{12}$ | $p_{13}$ | $p_{14}$ | $p_{15}$ | $p_{16}$ | $p_{17}$ | $p_{18}$ | $p_{19}$ |
| $p_1$ | $p_1$ | $p_2$ | $p_3$ | $p_4$ | $p_5$ | $p_6$ | $p_7$ | $p_8$ | $p_9$ | $p_0$ | $p_{16}$ | $p_{17}$ | $p_{18}$ | $p_{19}$ | $p_{15}$ | $p_{10}$ | $p_{11}$ | $p_{12}$ | $p_{13}$ | $p_{14}$ |
| $p_2$ | $p_2$ | $p_3$ | $p_4$ | $p_5$ | $p_6$ | $p_7$ | $p_8$ | $p_9$ | $p_0$ | $p_1$ | $p_{11}$ | $p_{12}$ | $p_{13}$ | $p_{14}$ | $p_{10}$ | $p_{16}$ | $p_{17}$ | $p_{18}$ | $p_{19}$ | $p_{15}$ |
| $p_3$ | $p_3$ | $p_4$ | $p_5$ | $p_6$ | $p_7$ | $p_8$ | $p_9$ | $p_0$ | $p_1$ | $p_2$ | $p_{17}$ | $p_{18}$ | $p_{19}$ | $p_{15}$ | $p_{16}$ | $p_{11}$ | $p_{12}$ | $p_{13}$ | $p_{14}$ | $p_{10}$ |
| $p_4$ | $p_4$ | $p_5$ | $p_6$ | $p_7$ | $p_8$ | $p_9$ | $p_0$ | $p_1$ | $p_2$ | $p_3$ | $p_{12}$ | $p_{13}$ | $p_{14}$ | $p_{10}$ | $p_{11}$ | $p_{17}$ | $p_{18}$ | $p_{19}$ | $p_{15}$ | $p_{16}$ |
| $p_5$ | $p_5$ | $p_6$ | $p_7$ | $p_8$ | $p_9$ | $p_0$ | $p_1$ | $p_2$ | $p_3$ | $p_4$ | $p_{18}$ | $p_{19}$ | $p_{15}$ | $p_{16}$ | $p_{17}$ | $p_{12}$ | $p_{13}$ | $p_{14}$ | $p_{10}$ | $p_{11}$ |
| $p_6$ | $p_6$ | $p_7$ | $p_8$ | $p_9$ | $p_0$ | $p_1$ | $p_2$ | $p_3$ | $p_4$ | $p_5$ | $p_{13}$ | $p_{14}$ | $p_{10}$ | $p_{11}$ | $p_{12}$ | $p_{18}$ | $p_{19}$ | $p_{15}$ | $p_{16}$ | $p_{17}$ |
| $p_7$ | $p_7$ | $p_8$ | $p_9$ | $p_0$ | $p_1$ | $p_2$ | $p_3$ | $p_4$ | $p_5$ | $p_6$ | $p_{19}$ | $p_{15}$ | $p_{16}$ | $p_{17}$ | $p_{18}$ | $p_{13}$ | $p_{14}$ | $p_{10}$ | $p_{11}$ | $p_{12}$ |
| $p_8$ | $p_8$ | $p_9$ | $p_0$ | $p_1$ | $p_2$ | $p_3$ | $p_4$ | $p_5$ | $p_6$ | $p_7$ | $p_{14}$ | $p_{10}$ | $p_{11}$ | $p_{12}$ | $p_{13}$ | $p_{19}$ | $p_{15}$ | $p_{16}$ | $p_{17}$ | $p_{18}$ |
| $p_9$ | $p_9$ | $p_0$ | $p_1$ | $p_2$ | $p_3$ | $p_4$ | $p_5$ | $p_6$ | $p_7$ | $p_8$ | $p_{15}$ | $p_{16}$ | $p_{17}$ | $p_{18}$ | $p_{19}$ | $p_{14}$ | $p_{10}$ | $p_{11}$ | $p_{12}$ | $p_{13}$ |
| $p_{10}$ | $p_{10}$ | $p_{15}$ | $p_{14}$ | $p_{19}$ | $p_{13}$ | $p_{18}$ | $p_{12}$ | $p_{17}$ | $p_{11}$ | $p_{16}$ | $p_0$ | $p_8$ | $p_6$ | $p_4$ | $p_2$ | $p_1$ | $p_9$ | $p_7$ | $p_5$ | $p_3$ |
| $p_{11}$ | $p_{11}$ | $p_{16}$ | $p_{10}$ | $p_{15}$ | $p_{14}$ | $p_{19}$ | $p_{13}$ | $p_{18}$ | $p_{12}$ | $p_{17}$ | $p_2$ | $p_0$ | $p_8$ | $p_6$ | $p_4$ | $p_3$ | $p_1$ | $p_9$ | $p_7$ | $p_5$ |
| $p_{12}$ | $p_{12}$ | $p_{17}$ | $p_{11}$ | $p_{16}$ | $p_{10}$ | $p_{15}$ | $p_{14}$ | $p_{19}$ | $p_{13}$ | $p_{18}$ | $p_4$ | $p_2$ | $p_0$ | $p_8$ | $p_6$ | $p_5$ | $p_3$ | $p_1$ | $p_9$ | $p_7$ |
| $p_{13}$ | $p_{13}$ | $p_{18}$ | $p_{12}$ | $p_{17}$ | $p_{11}$ | $p_{16}$ | $p_{10}$ | $p_{15}$ | $p_{14}$ | $p_{19}$ | $p_6$ | $p_4$ | $p_2$ | $p_0$ | $p_8$ | $p_7$ | $p_5$ | $p_3$ | $p_1$ | $p_9$ |
| $p_{14}$ | $p_{14}$ | $p_{19}$ | $p_{13}$ | $p_{18}$ | $p_{12}$ | $p_{17}$ | $p_{11}$ | $p_{16}$ | $p_{10}$ | $p_{15}$ | $p_8$ | $p_6$ | $p_4$ | $p_2$ | $p_0$ | $p_9$ | $p_7$ | $p_5$ | $p_3$ | $p_1$ |
| $p_{15}$ | $p_{15}$ | $p_{14}$ | $p_{19}$ | $p_{13}$ | $p_{18}$ | $p_{12}$ | $p_{17}$ | $p_{11}$ | $p_{16}$ | $p_{10}$ | $p_9$ | $p_7$ | $p_5$ | $p_3$ | $p_1$ | $p_0$ | $p_8$ | $p_6$ | $p_4$ | $p_2$ |
| $p_{16}$ | $p_{16}$ | $p_{10}$ | $p_{15}$ | $p_{14}$ | $p_{19}$ | $p_{13}$ | $p_{18}$ | $p_{12}$ | $p_{17}$ | $p_{11}$ | $p_1$ | $p_9$ | $p_7$ | $p_5$ | $p_3$ | $p_2$ | $p_0$ | $p_8$ | $p_6$ | $p_4$ |
| $p_{170}$ | $p_{17}$ | $p_{11}$ | $p_{16}$ | $p_{10}$ | $p_{15}$ | $p_{14}$ | $p_{19}$ | $p_{13}$ | $p_{18}$ | $p_{12}$ | $p_3$ | $p_1$ | $p_9$ | $p_7$ | $p_5$ | $p_4$ | $p_2$ | $p_0$ | $p_8$ | $p_6$ |
| $p_{18}$ | $p_{18}$ | $p_{12}$ | $p_{17}$ | $p_{11}$ | $p_{16}$ | $p_{10}$ | $p_{15}$ | $p_{14}$ | $p_{19}$ | $p_{13}$ | $p_5$ | $p_3$ | $p_1$ | $p_9$ | $p_7$ | $p_6$ | $p_4$ | $p_2$ | $p_0$ | $p_8$ |
| $p_{19}$ | $p_{19}$ | $p_{13}$ | $p_{18}$ | $p_{12}$ | $p_{17}$ | $p_{11}$ | $p_{16}$ | $p_{10}$ | $p_{15}$ | $p_{14}$ | $p_7$ | $p_5$ | $p_3$ | $p_1$ | $p_9$ | $p_8$ | $p_6$ | $p_4$ | $p_2$ | $p_0$ |

Граф икосаэдра $G_{10}$ является планарным графом, кольцевая сумма изометрических циклов есть пустое множество $\sum_{i=1}^{20} c_i = \varnothing$.

## 2.8. Рисунок графа и изометрические циклы

Важную роль при построении перестановок, в регулярных графах с равными весами вершин, играют изометрические циклы и топологические рисунки графа. Изометрические циклы порождают образующие циклы. Визуальное воспроизводство топологического рисунка графа позволяет облегчить решение задачи создания множества изоморфных перестановок. В данной главе рассматриваются планарные графы, имеющие нулевое значение функционала Маклейна, для базиса подпространства циклов C графа G[19]. Каждый планарный граф может быть представлен топологическим рисунком.

Будем рассматривать свойства преобразования топологического рисунка графа для построения перестановок.

Известно, что топологический рисунок графа характеризуется вращением его вершин, которое индуцирует простые циклы графа, удовлетворяющие условию теоремы планарности Маклейна.

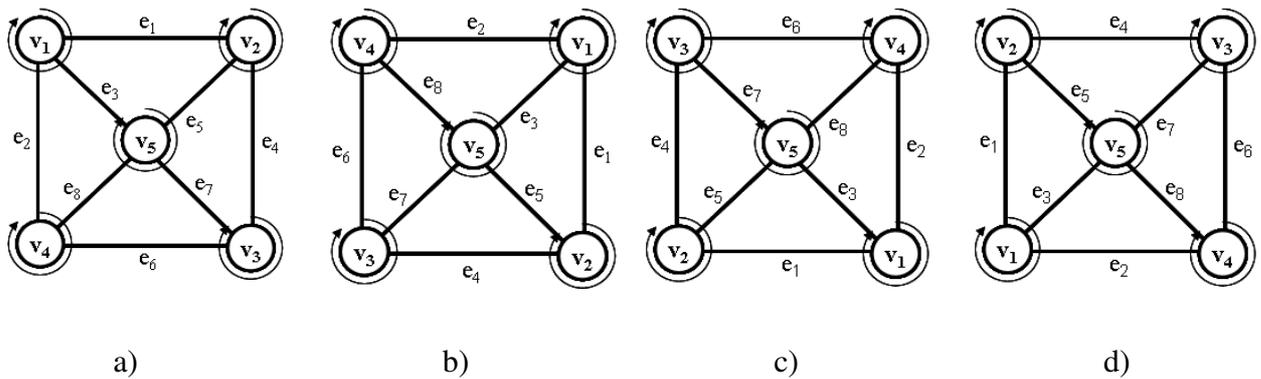

a) b) c) d)

Рис. 2.161. Вращение вершин при поворотах многоугольника.

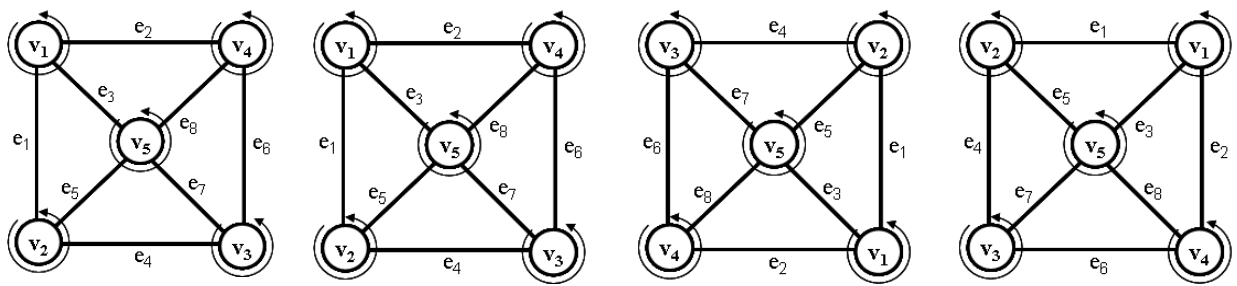

Рис. 2.162. Вращение вершин при осевых отражениях многоугольника.

Составим диаграмму вращения вершин для графа, представленного на рис. 2.161,а:

$v_1$: $v_2$ $v_5$ $v_4$
$v_2$: $v_1$ $v_3$ $v_5$
$v_3$: $v_2$ $v_4$ $v_5$
$v_4$: $v_1$ $v_5$ $v_3$
$v_5$: $v_1$ $v_2$ $v_3$ $v_4$

Для топологических рисунков характеризующих перестановки можно указать следующее правило: вращение вершин не меняет своей структуры при построении перестановок. Вращение изменяет свое направление вращения (по часовой стрелки или против), только при переходе от изображения поворотов относительно центра к изображению осевых отображений.

Существует симметричный топологический рисунок, учитывающий симметрию расположения вершин при поворотах и осевом отображении, допускающий существование нескольких вершин в центре.

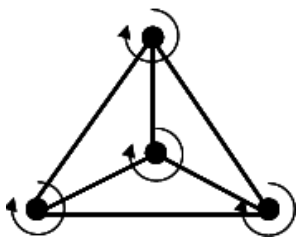

Рис. 2.163. Симметричный топологический рисунок графа тетраэдра.

Так для тетраэдра симметричный топологический рисунок графа имеет вид правильного треугольника с центральной точкой.

Симметричный топологический рисунок графа гектаэдра, представляет собой два связанных четырехугольника, вложенных друг в друга.

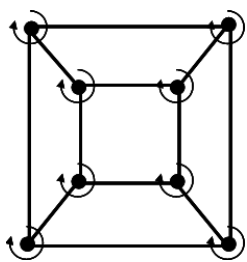

Рис. 2.164. Симметричный топологический рисунок графа гексаэдра.

На рис. 2.166 представлен переход от топологического рисунка графа к симметричному топологическому рисунку графа. Образующий цикл разбивает поверхность на две части внешнюю и внутреннюю. Осевые симметрии можно изобразить только на одной части пространства (как правило во внутренней части), поэтому элементы топологического рисунка графа следует объединить и поместить в одну часть пространства. При этом перемещаемые вершины изменят направление вращения. На рис. 2.166,b вращение вершины **A** представлено в обратном направлении, изменяется также вращение вершин **C,D,E,F**.

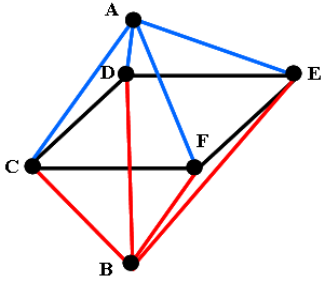
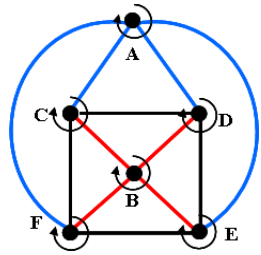
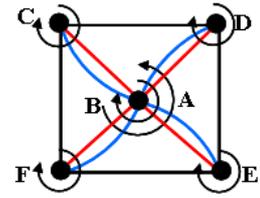

Рис. 2.165. Граф октаэдра.  Рис. 2.166. Переход от топологического рисунка графа к симметричному топологическому рисунку графа октаэдра.

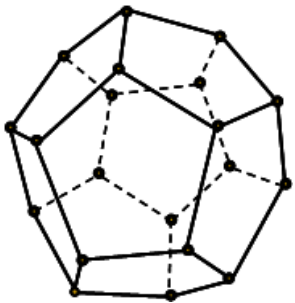
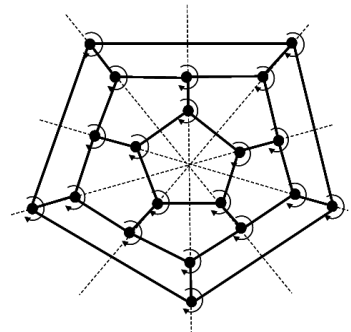

Рис. 2.167. Граф додекаэдра.  Рис. 2.168. Симметричный топологический рисунок графа додекаэдра.

Симметричный рисунок для графа додекаэдра представлен на рис. 2.168.

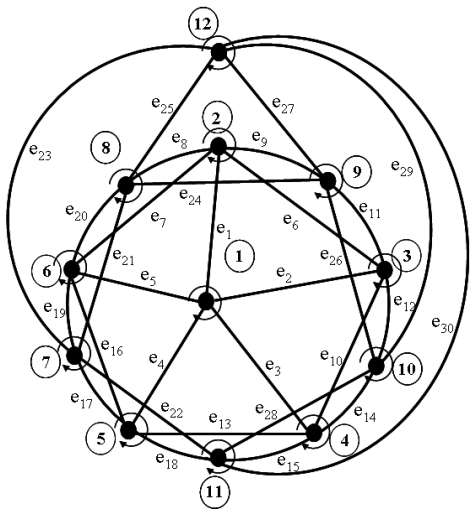
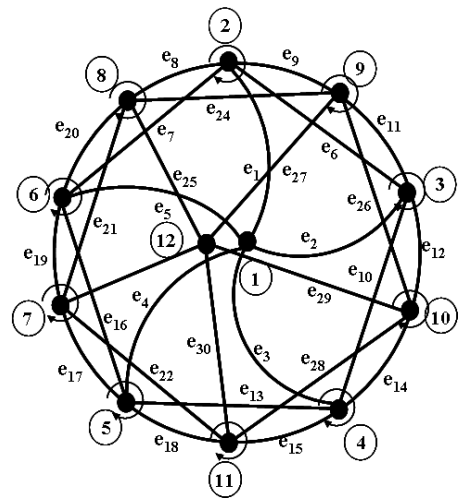

Рис. 2.169. Топологический рисунок графа икосаэдра.  Рис. 2.170. Симметричный топологический рисунок графа икосаэдра.

Образующий цикл разделяет топологический рисунок графа на внешнюю и внутреннюю части. Симметричный топологический рисунок графа икосаэдра представлен на рис. 2.170.

Следует заметить, что диэдральная группа для графа октаэдра может быть записана как диэдральная грппа с двумя центральными точками $D_{4+2}$. Диэдральная группа для додекаэдра может быть записана в виде $D_{5,10,5}$. Диэдральную группу для икосаэдра запишем как $D_{5,5+2}$.

## 2.9. Свойства диэдральных групп

Рассмотрим свойства перестановок в графе тетраэдра. Поставим в соответствие перестановки, полученные топологическим подходом и полученные геометрическим способом:

$p_0 = \langle 1,2,3,4 \rangle = (1)(2)(3)(4) \to \sigma_1$;
$p_1 = \langle 3,1,2,4 \rangle = (1\ 3\ 2)(4) = (3\ 1)(1\ 2)(4) \to \sigma_9$;
$p_2 = \langle 2,3,1,4 \rangle = (1\ 2\ 3)(4) = (2\ 1)(1\ 3)(4) \to \sigma_8$;
$p_3 = \langle 1,3,2,4 \rangle = (1)(2\ 3)(4) \to \sigma_{13}$;
$p_4 = \langle 3,2,1,4 \rangle = (1\ 3)(2)(4) \to \sigma_{21}$;
$p_5 = \langle 2,1,3,4 \rangle = (1\ 2)(3)(4) \to \sigma_{20}$;
$p_6 = \langle 1,2,4,3 \rangle = (1)(2)(3\ 4) \to \sigma_{15}$;
$p_7 = \langle 4,1,2,3 \rangle = (1\ 4\ 3\ 2) = (4\ 1)(1\ 2)(2\ 3) \to \sigma_{17}$;
$p_8 = \langle 2,4,1,3 \rangle = (1\ 2\ 4\ 3) = (2\ 1)(4\ 1)(1\ 3) \to \sigma_{22}$;
$p_9 = \langle 1,4,2,3 \rangle = (1)(2\ 4\ 3) = (1)(4\ 2)(2\ 3) \to \sigma_3$;
$p_{10} = \langle 4,2,1,3 \rangle = (1\ 4\ 3)(2) = (4\ 1)((1\ 3)(2) \to \sigma_5$;
$p_{11} = \langle 2,1,4,3 \rangle = (1\ 2)(3\ 4) \to \sigma_{10}$;
$p_{12} = \langle 1,3,4,2 \rangle = (1)(2\ 3\ 4) = (1)(3\ 2)(4\ 2) \to \sigma_2$;
$p_{13} = \langle 4,1,3,2 \rangle = (1\ 4\ 2)(3) = (4\ 1)(1\ 2)(3) \to \sigma_7$;
$p_{14} = \langle 3,4,1,2 \rangle = (1\ 3)(2\ 4) \to \sigma_{11}$;
$p_{15} = \langle 1,4,3,2 \rangle = (1)(2\ 4)(3) \to \sigma_{14}$;
$p_{16} = \langle 4,3,1,2 \rangle = (1\ 4\ 2\ 3) = (4\ 1)(3\ 2)(1\ 2) \to \sigma_{19}$;
$p_{17} = \langle 3,1,4,2 \rangle = (1\ 3\ 4\ 2) = (3\ 1)(1\ 2)(2\ 4) \to \sigma_{23}$;
$p_{18} = \langle 2,4,3,1 \rangle = (1\ 2\ 4)(3) = (2\ 1)(4\ 1)(3) \to \sigma_6$;
$p_{19} = \langle 3,2,4,1 \rangle = (1\ 3\ 4)(2) = (3\ 1)(4\ 1)(2) \to \sigma_4$;
$p_{20} = \langle 4,3,2,1 \rangle = (1\ 4)(2\ 3) \to \sigma_{12}$;
$p_{21} = \langle 2,3,4,1 \rangle = (1\ 2\ 3\ 4) = (2\ 1)(3\ 1)(4\ 1) \to \sigma_{18}$;
$p_{22} = \langle 3,4,2,1 \rangle = (1\ 3\ 2\ 4) = (3\ 1)(4\ 2)(2\ 1) \to \sigma_{16}$;
$p_{23} = \langle 4,2,3,1 \rangle = (1\ 4)(2)(3) \to \sigma_{24}$.

Здесь символом p обозначаются перестановки, созданные топологическим способом, а символом $\sigma$ перестановки, созданные геометрические способом (см. раздел 2.1). Следуем заметить, что четные и нечетные перестановки, индуцированные образующим циклом, меняют и направление вращения вершин в топологическом рисунке графа. Например, четные перестановки $p_0, p_1, p_2$ индуцированные образующим циклом $c_4$ определяют вращение вершин по часовой стрелки (против часовой стрелки), а нечетные перестановки $p_3, p_4, p_5$ индуцированные тем же образующим циклом $c_4$ определяют вращение вершин против часовой стрелки (по часовой стрелки).

# Комментарии

Во 2-ой главе рассмотрены методы построения группы автоморфизма графа для известных пяти правильных многогранников (платоновых тел). Это планарные регулярные графы с равными весами вершин. Произведен сравнительный анализ топологического и геометрического методов построения групп автоморфизма графа. Показано, что количество перестановок полученные двумя разными методами совпадают.

В основу построения группы автоморфизма графа положен принцип симметрии правильного многоугольника. Показано, что симметрия правильного многоугольника порождается образующим циклом графа, состоящим из множества изометрических циклов графа. Выделение перестановок проводится с помощью топологического рисунка графа и соответствующего симметричного топологического рисунка. Показано, что каждый образующий цикл порождает множество перестановок с мощностью равной его удвоенной длины. По результатам строится множество произведений перестановок в виде таблицы Кэли.

# Глава 3. АВТОМОРФИЗМ РЕГУЛЯРНЫХ ГРАФОВ
## 3.1. Регулярные графы с валентностью равной 3

***Пример 3.1.*** Рассмотрим автоморфизм графа Фрухта. Это один из двух несепарабельных минимальных кубических графов, не имеющих нетривиальных автоморфизмов. Описан Робертом Фрухтом в 1939.

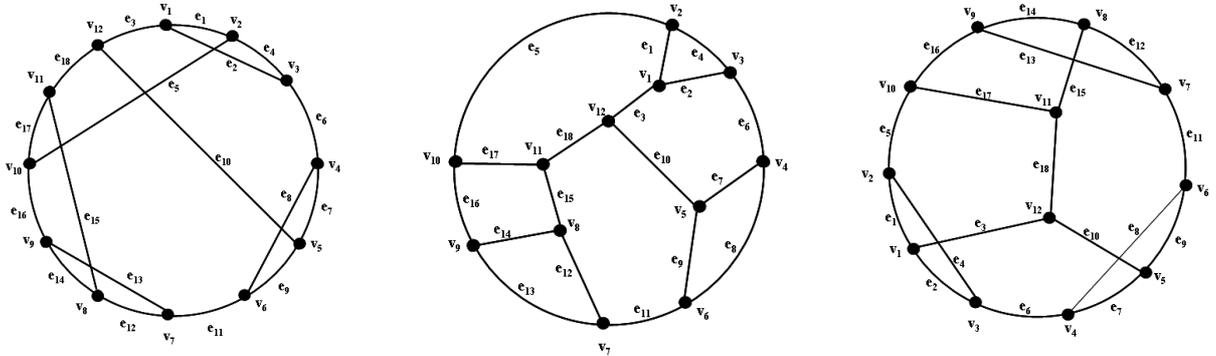

Рис. 3.1. Граф Фрухта.

Количество вершин графа = 12
Количество ребер графа = 18
Количество единичных циклов = 8
Смежность графа:

вершина $v_1$ = {$v_2,v_3,v_{12}$};              вершина $v_2$ = {$v_1,v_3,v_{10}$};
вершина $v_3$ = {$v_1,v_2,v_4$};                 вершина $v_4$ = {$v_3,v_5,v_6$};
вершина $v_5$ = {$v_4,v_6,v_{12}$};              вершина $v_6$ = {$v_4,v_5,v_7$};
вершина $v_7$ = {$v_6,v_8,v_9$};                 вершина $v_8$ = {$v_7,v_9,v_{11}$};
вершина $v_9$ = {$v_7,v_8,v_{10}$};              вершина $v_{10}$ = {$v_2,v_9,v_{11}$};
вершина $v_{11}$ = {$v_8,v_{10},v_{12}$};        вершина $v_{12}$ = {$v_1,v_5,v_{11}$}.

Инцидентность ребер графа:

$e_1 \to (v_1,v_2)$ или $(v_2,v_1)$;             $e_2 \to (v_1,v_3)$ или $(v_3,v_1)$;
$e_3 \to (v_1,v_{12})$ или $(v_{12}.v_1)$;       $e_4 \to (v_2,v_3)$ или $(v_3,v_2)$;
$e_5 \to (v_2,v_{10})$ или $(v_{10},v_2)$;       $e_6 \to (v_3,v_4)$ или $(v_4,v_3)$;
$e_7 \to (v_4,v_5)$ или $(v_5,v_4)$;             $e_8 \to (v_4,v_6)$ или $(v_6,v_4)$;
$e_9 \to (v_5,v_6)$ или $(v_6,v_5)$;             $e_{10} \to (v_5,v_{12})$ или $(v_{12},v_5)$;
$e_{11} \to (v_6,v_7)$ или $(v_7,v_6)$;          $e_{12} \to (v_7,v_8)$ или $(v_8,v_7)$;
$e_{13} \to (v_7,v_9)$ или $(v_9,v_7)$;          $e_{14} \to (v_8,v_9)$ или $(v_9,v_8)$;
$e_{15} \to (v_8,v_{11})$ или $(v_{11},v_8)$;    $e_{16} \to (v_9,v_{10})$ или $(v_{10},v_9)$;
$e_{17} \to (v_{10},v_{11})$ или $(v_{11},v_{10})$;  $e_{18} \to (v_{11},v_{12})$ или $(v_{12},v_{11})$.

Множество изометрических циклов графа:

| Циклы в виде ребер | Циклы в виде вершин |
|---|---|
| $c_1 = \{e_1,e_2,e_4\} \to$ | {$v_1,v_2,v_3$}; |
| $c_2 = \{e_1,e_3,e_5,e_{17},e_{18}\} \to$ | {$v_1,v_2,v_{10},v_{11},v_{12}$}; |
| $c_3 = \{e_2,e_3,e_6,e_7,e_{10}\} \to$ | {$v_1,v_3,v_4,v_5,v_{12}$}; |
| $c_4 = \{e_7,e_8,e_9\} \to$ | {$v_4,v_5,v_6$}; |
| $c_5 = \{e_9,e_{10},e_{11},e_{12},e_{15},e_{18}\} \to$ | {$v_5,v_6,v_7,v_8,v_{11},v_{12}$}; |
| $c_6 = \{e_4,e_5,e_6,e_8,e_{11},e_{13},e_{16}\} \to$ | {$v_2,v_3,v_4,v_6,v_7,v_9,v_{10}$}; |
| $c_7 = \{e_{12},e_{13},e_{14}\} \to$ | {$v_7,v_8,v_9$}; |
| $c_8 = \{e_{14},e_{15},e_{16},e_{17}\} \to$ | {$v_8,v_9,v_{10},v_{11}$}. |

Множество базовых реберных разрезов графа.

$w_0(e_1) = \{e_2,e_3,e_4,e_5\}$;  $w_0(e_2) = \{e_1,e_3,e_4,e_6\}$;
$w_0(e_3) = \{e_1,e_2,e_{10},e_{18}\}$;  $w_0(e_4) = \{e_1,e_2,e_5,e_6\}$;
$w_0(e_5) = \{e_1,e_4,e_{16},e_{17}\}$;  $w_0(e_6) = \{e_2,e_4,e_7,e_8\}$;
$w_0(e_7) = \{e_6,e_8,e_9,e_{10}\}$;  $w_0(e_8) = \{e_6,e_7,e_9,e_{11}\}$;
$w_0(e_9) = \{e_7,e_8,e_{10},e_{11}\}$;  $w_0(e_{10}) = \{e_3,e_7,e_9,e_{18}\}$;
$w_0(e_{11}) = \{e_8,e_9,e_{12},e_{13}\}$;  $w_0(e_{12}) = \{e_{11},e_{13},e_{14},e_{15}\}$;
$w_0(e_{13}) = \{e_{11},e_{12},e_{14},e_{16}\}$;  $w_0(e_{14}) = \{e_{12},e_{13},e_{15},e_{16}\}$;
$w_0(e_{15}) = \{e_{12},e_{14},e_{17},e_{18}\}$;  $w_0(e_{16}) = \{e_5,e_{13},e_{14},e_{17}\}$;
$w_0(e_{17}) = \{e_5,e_{15},e_{16},e_{18}\}$;  $w_0(e_{18}) = \{e_3,e_{10},e_{15},e_{17}\}$.

Кортеж весов ребер : <10,10,14,10,14,12,10,10,10,14,12,10,10,8,12,12,14,16>;
Кортеж весов вершин : <34,34,32,32,34,32,32,30,30,40,42,44>;
Вектор весов ребер : (8,10,10,10,10,10,10,10,10,12,12,12,12,14,14,14,14,16);
Вектор весов вершин : (30,30,32,32,32,32,34,34,34,40,42,44).

### 3.2. Регулярные графы с валентностью равной 4

*Пример 3.2.* Рассмотрим регулярный граф $G_{11}$.

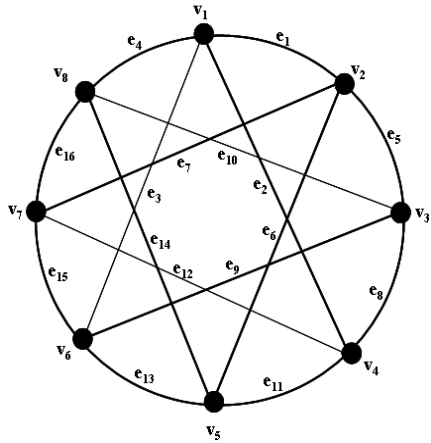

Рис. 3.2. Четырех валентный граф $G_{11}$ на 8 вершин.

Количество вершин в графе = 8.
Количество ребер в графе = 16.

Смежность графа:

вершина $v_1 = \{v_2,v_4,v_6,v_8\}$;  вершина $v_2 = \{v_1,v_3,v_5,v_7\}$;
вершина $v_3 = \{v_2,v_4,v_6,v_8\}$;  вершина $v_4 = \{v_1,v_3,v_5,v_7\}$;
вершина $v_5 = \{v_2,v_4,v_6,v_8\}$;  вершина $v_6 = \{v_1,v_3,v_5,v_7\}$;
вершина $v_7 = \{v_2,v_4,v_6,v_8\}$;  вершина $v_8 = \{v_1,v_3,v_5,v_7\}$.

Совместимость графа:

вершина $v_1 = \{e_1,e_2,e_3,e_4\}$;  вершина $v_2 = \{e_1,e_5,e_6,e_7\}$;
вершина $v_3 = \{e_5,e_8,e_9,e_{10}\}$;  вершина $v_4 = \{e_2,e_8,e_{11},e_{12}\}$;
вершина $v_5 = \{e_6,e_{11},e_{13},e_{14}\}$;  вершина $v_6 = \{e_3,e_9,e_{13},e_{15}\}$;
вершина $v_7 = \{e_7,e_{12},e_{15},e_{16}\}$;  вершина $v_8 = \{e_4,e_{10},e_{14},e_{16}\}$.

Множество базовых реберных разрезов графа.

$w_0(e_1) = \{e_2,e_3,e_4,e_5,e_6,e_7\}$;  $w_0(e_2) = \{e_1,e_3,e_4,e_8,e_{11},e_{12}\}$;
$w_0(e_3) = \{e_1,e_2,e_4,e_9,e_{13},e_{15}\}$;  $w_0(e_4) = \{e_1,e_2,e_3,e_{10},e_{14},e_{16}\}$;
$w_0(e_5) = \{e_1,e_6,e_7,e_8,e_9,e_{10}\}$;  $w_0(e_6) = \{e_1,e_5,e_7,e_{11},e_{13},e_{14}\}$;
$w_0(e_7) = \{e_1,e_5,e_6,e_{12},e_{15},e_{16}\}$;  $w_0(e_8) = \{e_2,e_5,e_9,e_{10},e_{11},e_{12}\}$;
$w_0(e_9) = \{e_3,e_5,e_8,e_{10},e_{13},e_{15}\}$;  $w_0(e_{10}) = \{e_4,e_5,e_8,e_9,e_{14},e_{16}\}$;
$w_0(e_{11}) = \{e_2,e_6,e_8,e_{12},e_{13},e_{14}\}$;  $w_0(e_{12}) = \{e_2,e_7,e_8,e_{11},e_{15},e_{16}\}$;

$w_0(e_{13}) = \{e_3,e_6,e_9,e_{11},e_{14},e_{15}\};$ $\qquad w_0(e_{14}) = \{e_4,e_6,e_{10},e_{11},e_{13},e_{16}\};$
$w_0(e_{15}) = \{e_3,e_7,e_9,e_{12},e_{13},e_{16}\};$ $\qquad w_0(e_{16}) = \{e_4,e_7,e_{10},e_{12},e_{14},e_{15}\}.$

Кортеж весов ребер : <6,6,6,6,6,6,6,6,6,6,6,6,6,6,6,6>;
Кортеж весов вершин : <24,24,24,24,24,24,24,24>;
Вектор весов ребер : (6,6,6,6,6,6,6,6,6,6,6,6,6,6,6,6);
Вектор весов вершин : (24,24,24,24,24,24,24,24).

Изометрические циклы в графе :

$c_1 = \{e_1,e_2,e_5,e_8\};$ $\qquad c_2 = \{e_1,e_2,e_6,e_{11}\};$
$c_3 = \{e_1,e_2,e_7,e_{12}\};$ $\qquad c_4 = \{e_1,e_3,e_5,e_9\};$
$c_5 = \{e_1,e_3,e_6,e_{13}\};$ $\qquad c_6 = \{e_1,e_3,e_7,e_{15}\};$
$c_7 = \{e_1,e_4,e_5,e_{10}\};$ $\qquad c_8 = \{e_1,e_4,e_6,e_{14}\};$
$c_9 = \{e_1,e_4,e_7,e_{16}\};$ $\qquad c_{10} = \{e_2,e_3,e_8,e_9\};$
$c_{11} = \{e_2,e_3,e_{11},e_{13}\};$ $\qquad c_{12} = \{e_2,e_3,e_{12},e_{15}\};$
$c_{13} = \{e_2,e_4,e_8,e_{10}\};$ $\qquad c_{14} = \{e_2,e_4,e_{11},e_{14}\};$
$c_{15} = \{e_2,e_4,e_{12},e_{16}\};$ $\qquad c_{16} = \{e_3,e_4,e_9,e_{10}\};$
$c_{17} = \{e_3,e_4,e_{13},e_{14}\};$ $\qquad c_{18} = \{e_3,e_4,e_{15},e_{16}\};$
$c_{19} = \{e_5,e_6,e_8,e_{11}\};$ $\qquad c_{20} = \{e_5,e_6,e_9,e_{13}\};$
$c_{21} = \{e_5,e_7,e_{10},e_{14}\};$ $\qquad c_{22} = \{e_5,e_7,e_8,e_{12}\};$
$c_{23} = \{e_5,e_7,e_9,e_{15}\};$ $\qquad c_{24} = \{e_5,e_7,e_{10},e_{16}\};$
$c_{25} = \{e_6,e_7,e_{11},e_{12}\};$ $\qquad c_{26} = \{e_6,e_7,e_{13},e_{15}\};$
$c_{27} = \{e_6,e_7,e_{14},e_{16}\};$ $\qquad c_{28} = \{e_8,e_9,e_{11},e_{13}\};$
$c_{29} = \{e_8,e_9,e_{12},e_{15}\};$ $\qquad c_{30} = \{e_8,e_{10},e_{11},e_{14}\};$
$c_{31} = \{e_8,e_{10},e_{12},e_{16}\};$ $\qquad c_{32} = \{e_9,e_{10},e_{13},e_{14}\};$
$c_{33} = \{e_9,e_{10},e_{15},e_{16}\};$ $\qquad c_{34} = \{e_{11},e_{12},e_{13},e_{15}\};$
$c_{35} = \{e_{11},e_{12},e_{14},e_{16}\};$ $\qquad c_{36} = \{e_{13},e_{14},e_{15},e_{16}\}.$

Образующие циклы длиной восемь, строятся как результат кольцевого суммирования изометрических циклов (см. рис. 3.3), используя правило покрытия всех вершин тремя циклами длиной четыре. Количество конфигураций (количество сочетаний) из трех циклов, покрывающих множество всех вершин графа равно 864. Однако множество всех конфигураций состоит из пересекающихся изометрических циклов, количество таких образующих циклов в графе равно 72.

Например, рассмотрим различный состав изометрических циклов для построения образующего цикла $\{e_3,e_4,e_5,e_6,e_8,e_{12},e_{13},e_{16}\}$:

$\{e_3,e_4,e_5,e_6,e_8,e_{12},e_{13},e_{16}\} = c_1 \oplus c_5 \oplus c_{15} = \{e_1,e_2,e_5,e_8\} \oplus \{e_1,e_3,e_6,e_{13}\} \oplus \{e_2,e_4,e_{12},e_{16}\};$
$\{e_3,e_4,e_5,e_6,e_8,e_{12},e_{13},e_{16}\} = c_5 \oplus c_7 \oplus c_{31} = \{e_1,e_3,e_6,e_{13}\} \oplus \{e_1,e_4,e_5,e_{10}\} \oplus \{e_8,e_{10},e_{12},e_{16}\};$
$\{e_3,e_4,e_5,e_6,e_8,e_{12},e_{13},e_{16}\} = c_5 \oplus c_9 \oplus c_{22} = \{e_1,e_3,e_6,e_{13}\} \oplus \{e_1,e_4,e_7,e_{16}\} \oplus \{e_5,e_7,e_8,e_{12}\};$
$\{e_3,e_4,e_5,e_6,e_8,e_{12},e_{13},e_{16}\} = c_{10} \oplus c_{15} \oplus c_{20} = \{e_2,e_3,e_8,e_9\} \oplus \{e_2,e_4,e_{12},e_{16}\} \oplus \{e_5,e_6,e_9,e_{13}\};$
$\{e_3,e_4,e_5,e_6,e_8,e_{12},e_{13},e_{16}\} = c_{11} \oplus c_{15} \oplus c_{19} = \{e_2,e_3,e_{11},e_{13}\} \oplus \{e_2,e_4,e_{12},e_{16}\} \oplus \{e_5,e_6,e_8,e_{11}\};$
$\{e_3,e_4,e_5,e_6,e_8,e_{12},e_{13},e_{16}\} = c_{16} \oplus c_{20} \oplus c_{31} = \{e_3,e_4,e_9,e_{10}\} \oplus \{e_5,e_6,e_9,e_{13}\} \oplus \{e_8,e_{10},e_{12},e_{16}\};$
$\{e_3,e_4,e_5,e_6,e_8,e_{12},e_{13},e_{16}\} = c_{17} \oplus c_{19} \oplus c_{35} = \{e_3,e_4,e_{13},e_{14}\} \oplus \{e_5,e_6,e_8,e_{11}\} \oplus \{e_{11},e_{12},e_{14},e_{16}\};$
$\{e_3,e_4,e_5,e_6,e_8,e_{12},e_{13},e_{16}\} = c_{17} \oplus c_{21} \oplus c_{31} = \{e_3,e_4,e_{13},e_{14}\} \oplus \{e_5,e_7,e_{10},e_{14}\} \oplus \{e_8,e_{10},e_{12},e_{16}\};$
$\{e_3,e_4,e_5,e_6,e_8,e_{12},e_{13},e_{16}\} = c_{17} \oplus c_{22} \oplus c_{27} = \{e_3,e_4,e_{13},e_{14}\} \oplus \{e_5,e_7,e_8,e_{12}\} \oplus \{e_6,e_7,e_{14},e_{16}\};$
$\{e_3,e_4,e_5,e_6,e_8,e_{12},e_{13},e_{16}\} = c_{18} \oplus c_{19} \oplus c_{34} = \{e_3,e_4,e_{15},e_{16}\} \oplus \{e_5,e_6,e_8,e_{11}\} \oplus \{e_{11},e_{12},e_{13},e_{15}\};$
$\{e_3,e_4,e_5,e_6,e_8,e_{12},e_{13},e_{16}\} = c_{18} \oplus c_{20} \oplus c_{29} = \{e_3,e_4,e_{15},e_{16}\} \oplus \{e_5,e_6,e_9,e_{12}\} \oplus \{e_8,e_9,e_{12},e_{15}\};$
$\{e_3,e_4,e_5,e_6,e_8,e_{12},e_{13},e_{16}\} = c_{18} \oplus c_{22} \oplus c_{26} = \{e_3,e_4,e_{15},e_{16}\} \oplus \{e_5,e_7,e_8,e_{12}\} \oplus \{e_6,e_7,e_{13},e_{15}\}.$

Перечислим образующие циклы:

$Q_1=\{e_3,e_4,e_5,e_6,e_8,e_{12},e_{13},e_{16}\}$;
$Q_2=\{e_2,e_4,e_5,e_6,e_8,e_{13},e_{15},e_{16}\}$;
$Q_3=\{e_2,e_3,e_6,e_7,e_8,e_{10},e_{13},e_{16}\}$;
$Q_4=\{e_2,e_3,e_5,e_7,e_8,e_{13},e_{14},e_{16}\}$;
$Q_5=\{e_2,e_3,e_5,e_6,e_{10},e_{12},e_{13},e_{16}\}$;
$Q_6=\{e_2,e_3,e_5,e_6,e_8,e_{14},e_{15},e_{16}\}$;
$Q_7=\{e_3,e_4,e_5,e_7,e_8,e_{11},e_{14},e_{15}\}$;
$Q_8=\{e_2,e_4,e_5,e_7,e_8,e_{13},e_{14},e_{15}\}$;
$Q_9=\{e_2,e_3,e_6,e_7,e_8,e_{10},e_{14},e_{15}\}$;
$Q_{10}=\{e_2,e_3,e_5,e_7,e_{10},e_{11},e_{14},e_{15}\}$;
$Q_{11}=\{e_3,e_4,e_5,e_6,e_8,e_{12},e_{14},e_{15}\}$;
$Q_{12}=\{e_2,e_4,e_6,e_7,e_8,e_9,e_{14},e_{15}\}$;
$Q_{13}=\{e_2,e_4,e_5,e_6,e_9,e_{12},e_{14},e_{15}\}$;
$Q_{14}=\{e_3,e_4,e_5,e_7,e_8,e_{11},e_{13},e_{16}\}$;
$Q_{15}=\{e_2,e_4,e_6,e_7,e_8,e_9,e_{13},e_{16}\}$;
$Q_{16}=\{e_2,e_4,e_5,e_7,e_9,e_{11},e_{13},e_{16}\}$;
$Q_{17}=\{e_1,e_4,e_5,e_8,e_{11},e_{13},e_{15},e_{16}\}$;
$Q_{18}=\{e_1,e_3,e_7,e_8,e_{10},e_{11},e_{13},e_{16}\}$;
$Q_{19}=\{e_1,e_3,e_5,e_{10},e_{11},e_{12},e_{13},e_{16}\}$;
$Q_{20}=\{e_1,e_3,e_5,e_8,e_{12},e_{13},e_{14},e_{16}\}$;
$Q_{21}=\{e_1,e_3,e_5,e_8,e_{11},e_{14},e_{15},e_{16}\}$;
$Q_{22}=\{e_1,e_4,e_5,e_8,e_{12},e_{13},e_{14},e_{15}\}$;
$Q_{23}=\{e_1,e_3,e_6,e_8,e_{10},e_{12},e_{14},e_{15}\}$;
$Q_{24}=\{e_1,e_3,e_5,e_{10},e_{11},e_{12},e_{14},e_{15}\}$;
$Q_{25}=\{e_1,e_4,e_7,e_8,e_9,e_{11},e_{14},e_{15}\}$;
$Q_{26}=\{e_1,e_4,e_5,e_9,e_{11},e_{12},e_{14},e_{15}\}$;
$Q_{27}=\{e_1,e_4,e_6,e_8,e_9,e_{12},e_{13},e_{16}\}$;
$Q_{28}=\{e_1,e_4,e_5,e_9,e_{11},e_{12},e_{13},e_{16}\}$;
$Q_{29}=\{e_1,e_2,e_7,e_8,e_9,e_{13},e_{14},e_{16}\}$;
$Q_{30}=\{e_1,e_2,e_6,e_9,e_{10},e_{12},e_{13},e_{16}\}$;
$Q_{31}=\{e_1,e_2,e_6,e_8,e_{10},e_{13},e_{15},e_{16}\}$;
$Q_{32}=\{e_1,e_2,e_6,e_8,e_9,e_{14},e_{15},e_{16}\}$;
$Q_{33}=\{e_1,e_2,e_7,e_8,e_{10},e_{13},e_{14},e_{15}\}$;
$Q_{34}=\{e_1,e_2,e_6,e_9,e_{10},e_{12},e_{14},e_{15}\}$;
$Q_{35}=\{e_1,e_2,e_7,e_9,e_{10},e_{11},e_{14},e_{15}\}$;
$Q_{36}=\{e_1,e_2,e_7,e_9,e_{10},e_{11},e_{13},e_{16}\}$;
$Q_{37}=\{e_1,e_2,e_5,e_{10},e_{11},e_{13},e_{15},e_{16}\}$;
$Q_{38}=\{e_1,e_2,e_5,e_9,e_{12},e_{13},e_{14},e_{16}\}$;
$Q_{39}=\{e_1,e_2,e_5,e_9,e_{11},e_{14},e_{15},e_{16}\}$;
$Q_{40}=\{e_1,e_2,e_5,e_{10},e_{12},e_{13},e_{14},e_{15}\}$;
$Q_{41}=\{e_3,e_4,e_5,e_6,e_9,e_{11},e_{12},e_{16}\}$;
$Q_{42}=\{e_2,e_4,e_5,e_6,e_9,e_{11},e_{15},e_{16}\}$;
$Q_{43}=\{e_2,e_3,e_6,e_7,e_9,e_{10},e_{11},e_{16}\}$;
$Q_{44}=\{e_2,e_3,e_5,e_7,e_9,e_{11},e_{14},e_{16}\}$;
$Q_{45}=\{e_2,e_3,e_5,e_6,e_{10},e_{11},e_{15},e_{16}\}$;
$Q_{46}=\{e_2,e_3,e_5,e_6,e_9,e_{12},e_{14},e_{16}\}$;
$Q_{47}=\{e_3,e_4,e_6,e_7,e_8,e_{10},e_{11},e_{15}\}$;
$Q_{48}=\{e_2,e_4,e_6,e_7,e_9,e_{10},e_{11},e_{15}\}$;
$Q_{49}=\{e_3,e_4,e_5,e_6,e_{10},e_{11},e_{12},e_{15}\}$;
$Q_{50}=\{e_2,e_4,e_5,e_7,e_{10},e_{11},e_{13},e_{15}\}$;
$Q_{51}=\{e_2,e_4,e_5,e_6,e_{10},e_{12},e_{13},e_{15}\}$;
$Q_{52}=\{e_3,e_4,e_6,e_7,e_8,e_9,e_{11},e_{16}\}$;
$Q_{53}=\{e_1,e_4,e_6,e_8,e_9,e_{11},e_{15},e_{16}\}$;
$Q_{54}=\{e_1,e_3,e_7,e_8,e_9,e_{11},e_{14},e_{16}\}$;
$Q_{55}=\{e_1,e_3,e_6,e_9,e_{10},e_{11},e_{12},e_{16}\}$;
$Q_{56}=\{e_1,e_3,e_6,e_8,e_{10},e_{11},e_{15},e_{16}\}$;
$Q_{57}=\{e_1,e_3,e_6,e_8,e_9,e_{12},e_{14},e_{16}\}$;
$Q_{58}=\{e_1,e_4,e_6,e_9,e_{10},e_{11},e_{12},e_{15}\}$;
$Q_{59}=\{e_1,e_4,e_7,e_8,e_{10},e_{11},e_{13},e_{15}\}$;
$Q_{60}=\{e_1,e_4,e_6,e_8,e_{10},e_{12},e_{13},e_{15}\}$;
$Q_{61}=\{e_3,e_4,e_5,e_7,e_9,e_{11},e_{12},e_{14}\}$;
$Q_{62}=\{e_2,e_4,e_5,e_7,e_9,e_{12},e_{13},e_{14}\}$;
$Q_{63}=\{e_2,e_3,e_6,e_7,e_9,e_{10},e_{12},e_{14}\}$;
$Q_{64}=\{e_2,e_3,e_5,e_7,e_{10},e_{12},e_{13},e_{14}\}$;
$Q_{65}=\{e_3,e_4,e_6,e_7,e_8,e_{10},e_{12},e_{13}\}$;
$Q_{66}=\{e_2,e_4,e_6,e_7,e_9,e_{10},e_{12},e_{13}\}$;
$Q_{67}=\{e_3,e_4,e_5,e_7,e_{10},e_{11},e_{12},e_{13}\}$;
$Q_{68}=\{e_3,e_4,e_6,e_7,e_8,e_9,e_{12},e_{14}\}$;
$Q_{69}=\{e_1,e_4,e_7,e_8,e_9,e_{12},e_{13},e_{14}\}$;
$Q_{70}=\{e_1,e_3,e_7,e_9,e_{10},e_{11},e_{12},e_{14}\}$;
$Q_{71}=\{e_1,e_3,e_7,e_8,e_{10},e_{12},e_{13},e_{14}\}$;
$Q_{72}=\{e_1,e_4,e_7,e_9,e_{10},e_{11},e_{12},e_{13}\}$.

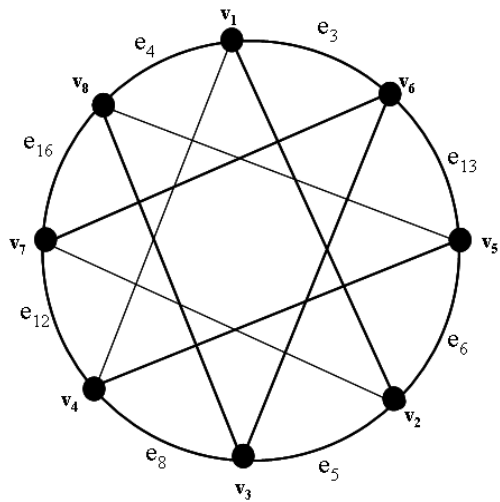

Рис. 3.3. Образующий цикл $Q_1 = \{e_3,e_4,e_5,e_6,e_8,e_{12},e_{13},e_{16}\}$.

Перестановки, порождаемые образующим циклом $Q_1 = \{e_3, e_4, e_5, e_6, e_8, e_{12}, e_{13}, e_{16}\}$:

$p_0$ = < 1 6 5 2 3 4 7 8 > = (1)(2 6 4)(3 5)(7)(8);
$p_1$ = < 8 1 6 5 2 3 4 7 > = (1 8 7 4 5 2)(3 6);
$p_2$ = < 7 8 1 6 5 2 3 4 > = (1 7 3)(2 8 4 6)(5);
$p_3$ = < 4 7 8 1 6 5 2 3 > = (1 4)(2 7)(3 8)(5 6);
$p_4$ = < 3 4 7 8 1 6 5 2 > = (1 3 7 5)(2 4 8)(6);
$p_5$ = < 2 3 4 7 8 1 6 5 > = (1 2 3 4 7 6)(5 8);
$p_6$ = < 5 2 3 4 7 8 1 6 > = (1 5 7)(2)(3)(4)(6 8);
$p_7$ = < 6 5 2 3 4 7 8 1 > = (1 6 7 8)(2 5 4 3);
$p_8$ = < 6 1 8 7 4 3 2 5 > = (1 6 3 8 5 4 7 2);
$p_9$ = < 2 5 6 1 8 7 4 3 > = (1 2 5 8 3 6 7 4);
$p_{10}$ = < 4 3 2 5 6 1 8 7 > = (1 4 5 6)(2 3)(7 8);
$p_{11}$ = < 8 7 4 3 2 5 6 1 > = (1 8)(2 7 6 5)(3 4);
$p_{12}$ = < 1 8 7 4 3 2 5 6 > = (1)(2 8 6)(3 7 5)(4);
$p_{13}$ = < 5 6 1 8 7 4 3 2 > = (1 5 7 3)(2 6 4 8);
$p_{14}$ = < 3 2 5 6 1 8 7 4 > = (1 3 5)(2)(4 6 8)(7);
$p_{15}$ = < 7 4 3 2 5 6 1 8 > = (1 7)(2 4)(3)(5)(6)(8);

Каждый образующий цикл порождает 16 перестановок. В общем случае количество перестановок определится как $72 \times 16 = 1152$.

***Пример 3.3.*** Рассмотрим регулярный граф $G_{13}$.

Смежность графа:

$v_1 = \{v_2, v_4, v_8, v_{10}\};$   $v_2 = \{v_1, v_3, v_5, v_9\};$
$v_3 = \{v_2, v_4, v_6, v_{10}\};$   $v_4 = \{v_1, v_3, v_5, v_7\};$
$v_5 = \{v_2, v_4, v_6, v_8\};$   $v_6 = \{v_3, v_5, v_7, v_9\};$
$v_7 = \{v_4, v_6, v_8, v_{10}\};$   $v_8 = \{v_1, v_5, v_7, v_9\};$
$v_9 = \{v_2, v_6, v_8, v_{10}\};$   $v_{10} = \{v_1, v_3, v_7, v_9\}.$

Центральные разрезы графа:

$s_1 = \{e_1, e_2, e_3, e_4\};$   $s_2 = \{e_1, e_5, e_6, e_7\};$
$s_3 = \{e_5, e_8, e_9, e_{10}\};$   $s_4 = \{e_2, e_8, e_{11}, e_{12}\};$
$s_5 = \{e_6, e_{11}, e_{13}, e_{14}\};$   $s_6 = \{e_9, e_{13}, e_{15}, e_{16}\};$
$s_7 = \{e_{12}, e_{15}, e_{17}, e_{18}\};$   $s_8 = \{e_3, e_{14}, e_{17}, e_{19}\};$
$s_9 = \{e_7, e_{16}, e_{19}, e_{20}\};$   $s_{10} = \{e_4, e_{10}, e_{18}, e_{20}\}.$

Количество вершин в графе = 10.
Количество ребер в графе = 20.

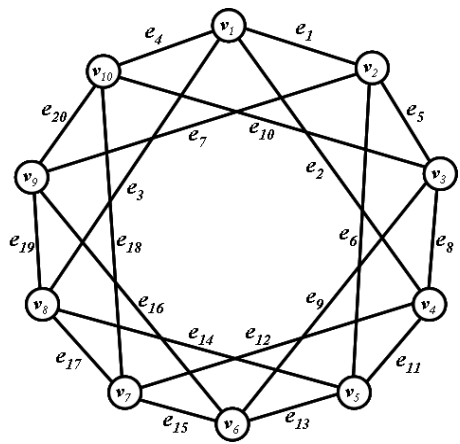

Рис. 3.4. Граф $G_{12}$.

Множество базовых реберных разрезов графа:

$w_0(e_1) = \{e_2,e_3,e_4,e_5,e_6,e_7\}$;   $w_0(e_2) = \{e_1,e_3,e_4,e_8,e_{11},e_{12}\}$;
$w_0(e_3) = \{e_1,e_2,e_4,e_{14},e_{17},e_{19}\}$;   $w_0(e_4) = \{e_1,e_2,e_3,e_{10},e_{18},e_{20}\}$;
$w_0(e_5) = \{e_1,e_6,e_7,e_8,e_9,e_{10}\}$;   $w_0(e_6) = \{e_1,e_5,e_7,e_{11},e_{13},e_{14}\}$;
$w_0(e_7) = \{e_1,e_5,e_6,e_{16},e_{19},e_{20}\}$;   $w_0(e_8) = \{e_2,e_5,e_9,e_{10},e_{11},e_{12}\}$;
$w_0(e_9) = \{e_5,e_8,e_{10},e_{13},e_{15},e_{16}\}$;   $w_0(e_{10}) = \{e_4,e_5,e_8,e_9,e_{18},e_{20}\}$;
$w_0(e_{11}) = \{e_2,e_6,e_8,e_{12},e_{13},e_{14}\}$;   $w_0(e_{12}) = \{e_2,e_8,e_{11},e_{15},e_{17},e_{18}\}$;
$w_0(e_{13}) = \{e_6,e_9,e_{11},e_{14},e_{15},e_{16}\}$;   $w_0(e_{14}) = \{e_3,e_6,e_{11},e_{13},e_{17},e_{19}\}$;
$w_0(e_{15}) = \{e_9,e_{12},e_{13},e_{16},e_{17},e_{18}\}$;   $w_0(e_{16}) = \{e_7,e_9,e_{13},e_{15},e_{19},e_{20}\}$;
$w_0(e_{17}) = \{e_3,e_{12},e_{14},e_{15},e_{18},e_{19}\}$;   $w_0(e_{18}) = \{e_4,e_{10},e_{12},e_{15},e_{17},e_{20}\}$;
$w_0(e_{19}) = \{e_3,e_7,e_{14},e_{16},e_{17},e_{20}\}$;   $w_0(e_{20}) = \{e_4,e_7,e_{10},e_{16},e_{18},e_{19}\}$.

Кортеж весов ребер: $\xi_w(G_{12}) = \langle 12,12,12,12,12,12,12,12,12,12,12,12,12,12,12,12,12,12,12,12 \rangle$;

Кортеж весов вершин: $\zeta_w(G_{12}) = \langle 48,48,48,48,48,48,48,48,48,48 \rangle$;

Вектор весов ребер: $F(\xi_w(G_{12})) = (12,12,12,12,12,12,12,12,12,12,12,12,12,12,12,12,12,12,12,12)$;

Вектор весов вершин: $F(\zeta_w(G_{12})) = \langle 48,48,48,48,48,48,48,48,48,48)$.

Множество изометрических циклов графа:

$c_1 = \{e_1,e_2,e_5,e_8\} \to \{v_1,v_2,v_3,v_4\}$;   $c_2 = \{e_1,e_2,e_6,e_{11}\} \to \{v_1,v_2,v_4,v_5\}$;
$c_3 = \{e_1,e_3,e_6,e_{14}\} \to \{v_1,v_2,v_5,v_8\}$;   $c_4 = \{e_1,e_3,e_7,e_{19}\} \to \{v_1,v_2,v_8,v_9\}$;
$c_5 = \{e_1,e_4,e_5,e_{10}\} \to \{v_1,v_2,v_3,v_{10}\}$;   $c_6 = \{e_1,e_4,e_7,e_{20}\} \to \{v_1,v_2,v_9,v_{10}\}$;
$c_7 = \{e_2,e_3,e_{11},e_{14}\} \to \{v_1,v_4,v_5,v_8\}$;   $c_8 = \{e_2,e_3,e_{12},e_{17}\} \to \{v_1,v_4,v_7,v_8\}$;
$c_9 = \{e_2,e_4,e_8,e_{10}\} \to \{v_1,v_3,v_4,v_{10}\}$;   $c_{10} = \{e_2,e_4,e_{12},e_{18}\} \to \{v_1,v_4,v_7,v_{10}\}$;
$c_{11} = \{e_3,e_4,e_{17},e_{18}\} \to \{v_1,v_7,v_8,v_{10}\}$;   $c_{12} = \{e_3,e_4,e_{19},e_{20}\} \to \{v_1,v_8,v_9,v_{10}\}$;
$c_{13} = \{e_5,e_6,e_8,e_{11}\} \to \{v_2,v_3,v_4,v_5\}$;   $c_{14} = \{e_5,e_6,e_9,e_{13}\} \to \{v_2,v_3,v_5,v_6\}$;
$c_{15} = \{e_5,e_7,e_9,e_{16}\} \to \{v_2,v_3,v_6,v_9\}$;   $c_{16} = \{e_5,e_7,e_{10},e_{20}\} \to \{v_2,v_3,v_9,v_{10}\}$;
$c_{17} = \{e_6,e_7,e_{13},e_{16}\} \to \{v_2,v_5,v_6,v_9\}$;   $c_{18} = \{e_6,e_7,e_{14},e_{19}\} \to \{v_2,v_5,v_8,v_9\}$;
$c_{19} = \{e_8,e_9,e_{11},e_{13}\} \to \{v_3,v_4,v_5,v_6\}$;   $c_{20} = \{e_8,e_9,e_{12},e_{15}\} \to \{v_3,v_4,v_6,v_7\}$;
$c_{21} = \{e_8,e_{10},e_{12},e_{18}\} \to \{v_3,v_4,v_7,v_{10}\}$;   $c_{22} = \{e_9,e_{10},e_{15},e_{18}\} \to \{v_3,v_6,v_7,v_{10}\}$;
$c_{23} = \{e_9,e_{10},e_{16},e_{20}\} \to \{v_3,v_6,v_9,v_{10}\}$;   $c_{24} = \{e_{11},e_{12},e_{13},e_{15}\} \to \{v_4,v_5,v_6,v_7\}$;
$c_{25} = \{e_{11},e_{12},e_{14},e_{17}\} \to \{v_4,v_5,v_7,v_8\}$;   $c_{26} = \{e_{13},e_{14},e_{15},e_{17}\} \to \{v_5,v_6,v_7,v_8\}$;
$c_{27} = \{e_{13},e_{14},e_{16},e_{19}\} \to \{v_5,v_6,v_8,v_9\}$;   $c_{28} = \{e_{15},e_{16},e_{17},e_{19}\} \to \{v_6,v_7,v_8,v_9\}$;
$c_{29} = \{e_{15},e_{16},e_{18},e_{20}\} \to \{v_6,v_7,v_9,v_{10}\}$;   $c_{30} = \{e_{17},e_{18},e_{19},e_{20}\} \to \{v_7,v_8,v_9,v_{10}\}$.

Общее число конфигураций = 1680. Количество образующих циклов = 156.

$Q_1 = \{e_3,e_4,e_6,e_7,e_8,e_9,e_{12},e_{14},e_{16},e_{18}\}$;   $Q_2 = \{e_3,e_4,e_5,e_7,e_8,e_{12},e_{13},e_{14},e_{16},e_{18}\}$;
$Q_3 = \{e_3,e_4,e_5,e_6,e_8,e_{12},e_{13},e_{16},e_{18},e_{19}\}$;   $Q_4 = \{e_3,e_4,e_5,e_6,e_8,e_{12},e_{14},e_{15},e_{16},e_{20}\}$;
$Q_5 = \{e_2,e_4,e_6,e_7,e_8,e_9,e_{14},e_{16},e_{17},e_{18}\}$;   $Q_6 = \{e_2,e_4,e_5,e_7,e_8,e_{13},e_{14},e_{16},e_{17},e_{18}\}$;
$Q_7 = \{e_2,e_4,e_5,e_6,e_8,e_{13},e_{16},e_{17},e_{18},e_{19}\}$;   $Q_8 = \{e_2,e_4,e_5,e_6,e_8,e_{14},e_{15},e_{16},e_{18},e_{19}\}$;
$Q_9 = \{e_2,e_4,e_5,e_6,e_8,e_{14},e_{15},e_{16},e_{17},e_{20}\}$;   $Q_{10} = \{e_2,e_4,e_5,e_6,e_9,e_{12},e_{14},e_{15},e_{19},e_{20}\}$;
$Q_{11} = \{e_2,e_4,e_5,e_6,e_8,e_{13},e_{15},e_{17},e_{19},e_{20}\}$;   $Q_{12} = \{e_2,e_3,e_6,e_7,e_9,e_{10},e_{12},e_{14},e_{16},e_{18}\}$;
$Q_{13} = \{e_2,e_3,e_6,e_7,e_8,e_{10},e_{14},e_{15},e_{16},e_{18}\}$;   $Q_{14} = \{e_2,e_3,e_6,e_7,e_8,e_9,e_{14},e_{15},e_{18},e_{20}\}$;
$Q_{15} = \{e_2,e_3,e_6,e_7,e_9,e_{10},e_{12},e_{14},e_{15},e_{20}\}$;   $Q_{16} = \{e_2,e_3,e_6,e_7,e_8,e_{10},e_{13},e_{15},e_{17},e_{20}\}$;
$Q_{17} = \{e_2,e_3,e_5,e_7,e_{10},e_{12},e_{13},e_{14},e_{16},e_{18}\}$;   $Q_{18} = \{e_2,e_3,e_5,e_7,e_8,e_{13},e_{14},e_{15},e_{18},e_{20}\}$;
$Q_{19} = \{e_2,e_3,e_5,e_6,e_{10},e_{12},e_{14},e_{15},e_{16},e_{20}\}$;   $Q_{20} = \{e_2,e_3,e_5,e_6,e_9,e_{12},e_{14},e_{16},e_{18},e_{20}\}$;
$Q_{21} = \{e_2,e_3,e_5,e_6,e_{10},e_{12},e_{13},e_{16},e_{18},e_{19}\}$;   $Q_{22} = \{e_2,e_3,e_5,e_6,e_8,e_{13},e_{16},e_{17},e_{18},e_{20}\}$;
$Q_{23} = \{e_2,e_3,e_5,e_6,e_8,e_{13},e_{15},e_{18},e_{19},e_{20}\}$;   $Q_{24} = \{e_3,e_4,e_6,e_7,e_8,e_9,e_{12},e_{13},e_{18},e_{19}\}$;
$Q_{25} = \{e_3,e_4,e_5,e_7,e_9,e_{11},e_{12},e_{13},e_{18},e_{19}\}$;   $Q_{26} = \{e_3,e_4,e_5,e_7,e_8,e_{11},e_{13},e_{15},e_{18},e_{19}\}$;
$Q_{27} = \{e_2,e_4,e_6,e_7,e_8,e_9,e_{13},e_{17},e_{18},e_{19}\}$;   $Q_{28} = \{e_2,e_4,e_5,e_7,e_9,e_{11},e_{13},e_{17},e_{18},e_{19}\}$;
$Q_{29} = \{e_2,e_4,e_5,e_7,e_8,e_{13},e_{14},e_{15},e_{18},e_{19}\}$;   $Q_{30} = \{e_2,e_3,e_6,e_7,e_9,e_{10},e_{12},e_{13},e_{18},e_{19}\}$;
$Q_{31} = \{e_2,e_3,e_6,e_7,e_8,e_{10},e_{13},e_{15},e_{18},e_{19}\}$;   $Q_{32} = \{e_2,e_3,e_6,e_7,e_8,e_9,e_{13},e_{17},e_{18},e_{20}\}$;

$Q_{33}=\{e_2,e_3,e_5,e_7,e_{10},e_{11},e_{13},e_{15},e_{18},e_{19}\}$;

$Q_{34}=\{e_2,e_3,e_5,e_7,e_9,e_{11},e_{13},e_{17},e_{18},e_{20}\}$;

$Q_{35}=\{e_3,e_4,e_5,e_7,e_9,e_{11},e_{12},e_{14},e_{15},e_{20}\}$;

$Q_{36}=\{e_3,e_4,e_5,e_7,e_8,e_{12},e_{13},e_{14},e_{15},e_{20}\}$;

$Q_{37}=\{e_3,e_4,e_5,e_7,e_8,e_{11},e_{13},e_{15},e_{17},e_{20}\}$;

$Q_{38}=\{e_3,e_4,e_5,e_7,e_8,e_{11},e_{14},e_{15},e_{16},e_{18}\}$;

$Q_{39}=\{e_3,e_4,e_6,e_7,e_8,e_9,e_{12},e_{13},e_{17},e_{20}\}$;

$Q_{40}=\{e_3,e_4,e_5,e_6,e_8,e_{12},e_{13},e_{16},e_{17},e_{20}\}$;

$Q_{41}=\{e_3,e_4,e_5,e_7,e_9,e_{11},e_{12},e_{13},e_{17},e_{20}\}$;

$Q_{42}=\{e_2,e_4,e_6,e_7,e_8,e_9,e_{14},e_{15},e_{17},e_{20}\}$;

$Q_{43}=\{e_2,e_4,e_5,e_7,e_9,e_{12},e_{13},e_{14},e_{17},e_{20}\}$;

$Q_{44}=\{e_2,e_4,e_5,e_7,e_9,e_{11},e_{14},e_{15},e_{17},e_{20}\}$;

$Q_{45}=\{e_1,e_4,e_7,e_8,e_9,e_{11},e_{14},e_{16},e_{17},e_{18}\}$;

$Q_{46}=\{e_1,e_4,e_5,e_8,e_{11},e_{13},e_{16},e_{17},e_{18},e_{19}\}$;

$Q_{47}=\{e_1,e_4,e_5,e_8,e_{11},e_{14},e_{15},e_{16},e_{18},e_{19}\}$;

$Q_{48}=\{e_1,e_4,e_5,e_8,e_{11},e_{14},e_{15},e_{16},e_{17},e_{20}\}$;

$Q_{49}=\{e_1,e_4,e_5,e_9,e_{11},e_{12},e_{14},e_{15},e_{19},e_{20}\}$;

$Q_{50}=\{e_1,e_4,e_5,e_8,e_{12},e_{13},e_{14},e_{15},e_{19},e_{20}\}$;

$Q_{51}=\{e_1,e_4,e_5,e_8,e_{11},e_{13},e_{15},e_{17},e_{19},e_{20}\}$;

$Q_{52}=\{e_1,e_3,e_7,e_9,e_{10},e_{11},e_{12},e_{14},e_{16},e_{18}\}$;

$Q_{53}=\{e_1,e_3,e_7,e_8,e_{10},e_{11},e_{14},e_{15},e_{16},e_{18}\}$;

$Q_{54}=\{e_1,e_3,e_7,e_8,e_9,e_{11},e_{14},e_{15},e_{18},e_{20}\}$;

$Q_{55}=\{e_1,e_3,e_7,e_9,e_{10},e_{11},e_{12},e_{14},e_{15},e_{20}\}$;

$Q_{56}=\{e_1,e_3,e_7,e_8,e_{10},e_{12},e_{13},e_{14},e_{15},e_{20}\}$;

$Q_{57}=\{e_1,e_3,e_7,e_8,e_{10},e_{11},e_{13},e_{15},e_{17},e_{20}\}$;

$Q_{58}=\{e_1,e_3,e_5,e_{10},e_{11},e_{12},e_{14},e_{15},e_{16},e_{20}\}$;

$Q_{59}=\{e_1,e_3,e_5,e_9,e_{11},e_{12},e_{14},e_{16},e_{18},e_{20}\}$;

$Q_{60}=\{e_1,e_3,e_5,e_{10},e_{11},e_{12},e_{13},e_{16},e_{18},e_{19}\}$;

$Q_{61}=\{e_1,e_3,e_5,e_8,e_{12},e_{13},e_{14},e_{16},e_{18},e_{20}\}$;

$Q_{62}=\{e_1,e_3,e_5,e_8,e_{11},e_{13},e_{16},e_{17},e_{18},e_{20}\}$;

$Q_{63}=\{e_1,e_3,e_5,e_8,e_{11},e_{13},e_{15},e_{18},e_{19},e_{20}\}$;

$Q_{64}=\{e_1,e_4,e_6,e_8,e_9,e_{12},e_{13},e_{17},e_{19},e_{20}\}$;

$Q_{65}=\{e_1,e_4,e_5,e_9,e_{11},e_{12},e_{13},e_{17},e_{19},e_{20}\}$;

$Q_{66}=\{e_1,e_4,e_5,e_8,e_{12},e_{13},e_{14},e_{16},e_{17},e_{20}\}$;

$Q_{67}=\{e_1,e_3,e_6,e_8,e_{10},e_{12},e_{13},e_{16},e_{17},e_{20}\}$;

$Q_{68}=\{e_1,e_3,e_6,e_8,e_9,e_{12},e_{13},e_{18},e_{19},e_{20}\}$;

$Q_{69}=\{e_1,e_3,e_7,e_9,e_{10},e_{11},e_{12},e_{13},e_{17},e_{20}\}$;

$Q_{70}=\{e_1,e_3,e_5,e_{10},e_{11},e_{12},e_{13},e_{16},e_{17},e_{20}\}$;

$Q_{71}=\{e_1,e_3,e_5,e_9,e_{11},e_{12},e_{13},e_{18},e_{19},e_{20}\}$;

$Q_{72}=\{e_1,e_4,e_7,e_8,e_9,e_{12},e_{13},e_{14},e_{18},e_{19}\}$;

$Q_{73}=\{e_1,e_4,e_6,e_8,e_9,e_{12},e_{14},e_{16},e_{18},e_{19}\}$;

$Q_{74}=\{e_1,e_4,e_5,e_9,e_{11},e_{12},e_{14},e_{16},e_{18},e_{19}\}$;

$Q_{75}=\{e_1,e_2,e_7,e_9,e_{10},e_{12},e_{13},e_{14},e_{18},e_{19}\}$;

$Q_{76}=\{e_1,e_2,e_7,e_8,e_{10},e_{13},e_{14},e_{15},e_{18},e_{19}\}$;

$Q_{77}=\{e_1,e_2,e_7,e_8,e_9,e_{13},e_{14},e_{17},e_{18},e_{20}\}$;

$Q_{78}=\{e_1,e_2,e_6,e_9,e_{10},e_{12},e_{14},e_{16},e_{18},e_{19}\}$;

$Q_{79}=\{e_1,e_2,e_6,e_9,e_{10},e_{12},e_{13},e_{17},e_{19},e_{20}\}$;

$Q_{80}=\{e_1,e_2,e_6,e_8,e_{10},e_{14},e_{15},e_{16},e_{18},e_{19}\}$;

$Q_{81}=\{e_1,e_2,e_6,e_8,e_{10},e_{13},e_{16},e_{17},e_{18},e_{19}\}$;

$Q_{82}=\{e_1,e_2,e_6,e_8,e_{10},e_{13},e_{15},e_{17},e_{19},e_{20}\}$;

$Q_{83}=\{e_1,e_2,e_6,e_8,e_{10},e_{14},e_{15},e_{16},e_{17},e_{20}\}$;

$Q_{84}=\{e_1,e_2,e_6,e_8,e_9,e_{14},e_{16},e_{17},e_{18},e_{20}\}$;

$Q_{85}=\{e_1,e_2,e_6,e_8,e_9,e_{14},e_{15},e_{18},e_{19},e_{20}\}$;

$Q_{86}=\{e_1,e_2,e_7,e_9,e_{10},e_{11},e_{14},e_{16},e_{17},e_{18}\}$;

$Q_{87}=\{e_1,e_2,e_7,e_8,e_{10},e_{13},e_{14},e_{16},e_{17},e_{18}\}$;

$Q_{88}=\{e_1,e_2,e_6,e_9,e_{10},e_{12},e_{14},e_{15},e_{19},e_{20}\}$;

$Q_{89}=\{e_1,e_2,e_7,e_9,e_{10},e_{12},e_{13},e_{14},e_{17},e_{20}\}$;

$Q_{90}=\{e_1,e_2,e_7,e_9,e_{10},e_{11},e_{14},e_{15},e_{17},e_{20}\}$;

$Q_{91}=\{e_1,e_2,e_7,e_9,e_{10},e_{11},e_{13},e_{17},e_{18},e_{19}\}$;

$Q_{92}=\{e_1,e_2,e_5,e_{10},e_{11},e_{14},e_{15},e_{16},e_{18},e_{19}\}$;

$Q_{93}=\{e_1,e_2,e_5,e_{10},e_{11},e_{13},e_{16},e_{17},e_{18},e_{19}\}$;

$Q_{94}=\{e_1,e_2,e_5,e_{10},e_{11},e_{13},e_{15},e_{17},e_{19},e_{20}\}$;

$Q_{95}=\{e_1,e_2,e_5,e_{10},e_{12},e_{13},e_{14},e_{16},e_{17},e_{20}\}$;

$Q_{96}=\{e_1,e_2,e_5,e_{10},e_{11},e_{14},e_{15},e_{16},e_{17},e_{20}\}$;

$Q_{97}=\{e_1,e_2,e_5,e_9,e_{12},e_{13},e_{14},e_{18},e_{19},e_{20}\}$;

$Q_{98}=\{e_1,e_2,e_5,e_9,e_{11},e_{14},e_{16},e_{17},e_{18},e_{20}\}$;

$Q_{99}=\{e_1,e_2,e_5,e_9,e_{11},e_{14},e_{15},e_{18},e_{19},e_{20}\}$;

$Q_{100}=\{e_1,e_2,e_5,e_{10},e_{12},e_{13},e_{14},e_{15},e_{19},e_{20}\}$;

$Q_{101}=\{e_3,e_4,e_6,e_7,e_9,e_{10},e_{11},e_{12},e_{15},e_{19}\}$;

$Q_{102}=\{e_3,e_4,e_6,e_7,e_8,e_9,e_{11},e_{15},e_{18},e_{19}\}$;

$Q_{103}=\{e_3,e_4,e_6,e_7,e_8,e_{10},e_{12},e_{13},e_{15},e_{19}\}$;

$Q_{104}=\{e_3,e_4,e_6,e_7,e_8,e_{10},e_{11},e_{15},e_{16},e_{17}\}$;

$Q_{105}=\{e_3,e_4,e_5,e_6,e_9,e_{11},e_{12},e_{16},e_{18},e_{19}\}$;

$Q_{106}=\{e_2,e_4,e_5,e_6,e_9,e_{11},e_{16},e_{17},e_{18},e_{19}\}$;

$Q_{107}=\{e_2,e_4,e_6,e_7,e_9,e_{10},e_{11},e_{15},e_{17},e_{19}\}$;

$Q_{108}=\{e_2,e_3,e_5,e_6,e_{10},e_{11},e_{15},e_{16},e_{18},e_{19}\}$;

$Q_{109}=\{e_2,e_3,e_5,e_6,e_9,e_{11},e_{15},e_{18},e_{19},e_{20}\}$;

$Q_{110}=\{e_2,e_3,e_5,e_6,e_9,e_{11},e_{16},e_{17},e_{18},e_{20}\}$;

$Q_{111}=\{e_2,e_3,e_5,e_6,e_{10},e_{12},e_{13},e_{15},e_{19},e_{20}\}$;

$Q_{112}=\{e_2,e_3,e_5,e_6,e_{10},e_{11},e_{15},e_{16},e_{17},e_{20}\}$;

$Q_{113}=\{e_2,e_3,e_6,e_7,e_9,e_{10},e_{11},e_{16},e_{17},e_{18}\}$;

$Q_{114}=\{e_2,e_3,e_6,e_7,e_9,e_{10},e_{11},e_{15},e_{17},e_{20}\}$;

$Q_{115}=\{e_3,e_4,e_6,e_7,e_9,e_{10},e_{11},e_{12},e_{16},e_{17}\}$;

$Q_{116}=\{e_3,e_4,e_5,e_7,e_{10},e_{11},e_{12},e_{13},e_{16},e_{17}\}$;

$Q_{117}=\{e_3,e_4,e_5,e_6,e_9,e_{11},e_{12},e_{16},e_{17},e_{20}\}$;

$Q_{118}=\{e_3,e_4,e_5,e_6,e_{10},e_{11},e_{12},e_{15},e_{16},e_{19}\}$;

$Q_{119}=\{e_2,e_3,e_5,e_7,e_{10},e_{11},e_{13},e_{16},e_{17},e_{18}\}$;

$Q_{120}=\{e_2,e_4,e_6,e_7,e_9,e_{10},e_{12},e_{14},e_{16},e_{17}\}$;

$Q_{121}=\{e_2,e_4,e_5,e_7,e_{10},e_{12},e_{13},e_{14},e_{16},e_{17}\}$;

$Q_{122}=\{e_2,e_4,e_5,e_7,e_{10},e_{11},e_{14},e_{15},e_{16},e_{17}\}$;

$Q_{123}=\{e_2,e_4,e_5,e_7,e_{10},e_{11},e_{13},e_{15},e_{17},e_{19}\}$;

$Q_{124}=\{e_2,e_4,e_5,e_7,e_9,e_{11},e_{14},e_{15},e_{18},e_{19}\}$;

$Q_{125}=\{e_2,e_4,e_5,e_7,e_{10},e_{12},e_{13},e_{14},e_{15},e_{19}\}$;

$Q_{126}=\{e_2,e_4,e_5,e_6,e_9,e_{11},e_{15},e_{17},e_{19},e_{20}\}$;

$Q_{127}=\{e_2,e_4,e_5,e_6,e_9,e_{12},e_{14},e_{16},e_{17},e_{20}\}$;

$Q_{128}=\{e_2,e_4,e_5,e_6,e_{10},e_{12},e_{14},e_{15},e_{16},e_{19}\}$;

$Q_{129}=\{e_2,e_4,e_5,e_6,e_{10},e_{12},e_{13},e_{16},e_{17},e_{19}\}$;

$Q_{130}=\{e_3,e_4,e_6,e_7,e_8,e_9,e_{11},e_{15},e_{17},e_{20}\}$;

$Q_{131}=\{e_1,e_4,e_6,e_8,e_9,e_{11},e_{15},e_{17},e_{19},e_{20}\}$;

$Q_{132}=\{e_1,e_4,e_6,e_9,e_{10},e_{11},e_{12},e_{16},e_{17},e_{19}\}$;

$Q_{133}=\{e_1,e_3,e_7,e_8,e_{10},e_{11},e_{13},e_{16},e_{17},e_{18}\}$;

$Q_{134}=\{e_1,e_3,e_6,e_8,e_{10},e_{11},e_{15},e_{16},e_{17},e_{20}\}$;

$Q_{135}=\{e_1,e_3,e_6,e_8,e_9,e_{11},e_{16},e_{17},e_{18},e_{20}\}$;

$Q_{136}=\{e_1,e_3,e_6,e_8,e_9,e_{11},e_{15},e_{18},e_{19},e_{20}\}$;

$Q_{137}=\{e_1,e_3,e_6,e_8,e_{10},e_{11},e_{15},e_{16},e_{18},e_{19}\}$; $Q_{138}=\{e_1,e_3,e_6,e_9,e_{10},e_{11},e_{12},e_{16},e_{18},e_{19}\}$;
$Q_{139}=\{e_1,e_3,e_6,e_9,e_{10},e_{11},e_{12},e_{15},e_{19},e_{20}\}$; $Q_{140}=\{e_1,e_3,e_6,e_8,e_{10},e_{12},e_{13},e_{15},e_{19},e_{20}\}$;
$Q_{141}=\{e_1,e_4,e_7,e_8,e_{10},e_{12},e_{13},e_{14},e_{16},e_{17}\}$; $Q_{142}=\{e_1,e_4,e_7,e_8,e_{10},e_{11},e_{14},e_{15},e_{16},e_{17}\}$;
$Q_{143}=\{e_1,e_4,e_7,e_8,e_{10},e_{11},e_{13},e_{15},e_{17},e_{19}\}$; $Q_{144}=\{e_1,e_4,e_7,e_9,e_{10},e_{11},e_{12},e_{14},e_{15},e_{19}\}$;
$Q_{145}=\{e_1,e_4,e_7,e_8,e_9,e_{11},e_{14},e_{15},e_{18},e_{19}\}$; $Q_{146}=\{e_1,e_4,e_7,e_8,e_{10},e_{12},e_{13},e_{14},e_{15},e_{19}\}$;
$Q_{147}=\{e_1,e_4,e_6,e_8,e_9,e_{11},e_{16},e_{17},e_{18},e_{19}\}$; $Q_{148}=\{e_1,e_4,e_6,e_8,e_9,e_{12},e_{14},e_{16},e_{17},e_{20}\}$;
$Q_{149}=\{e_1,e_4,e_6,e_8,e_{10},e_{12},e_{14},e_{15},e_{16},e_{19}\}$; $Q_{150}=\{e_1,e_4,e_6,e_8,e_{10},e_{12},e_{13},e_{16},e_{17},e_{19}\}$;
$Q_{151}=\{e_3,e_4,e_6,e_7,e_8,e_{10},e_{12},e_{14},e_{15},e_{16}\}$; $Q_{152}=\{e_3,e_4,e_5,e_7,e_{10},e_{11},e_{12},e_{14},e_{15},e_{16}\}$;
$Q_{153}=\{e_2,e_3,e_5,e_7,e_9,e_{12},e_{13},e_{14},e_{18},e_{20}\}$; $Q_{154}=\{e_1,e_3,e_7,e_8,e_9,e_{12},e_{13},e_{14},e_{18},e_{20}\}$;
$Q_{155}=\{e_2,e_4,e_6,e_7,e_9,e_{10},e_{12},e_{13},e_{17},e_{19}\}$; $Q_{156}=\{e_1,e_4,e_7,e_9,e_{10},e_{11},e_{12},e_{13},e_{17},e_{19}\}$.

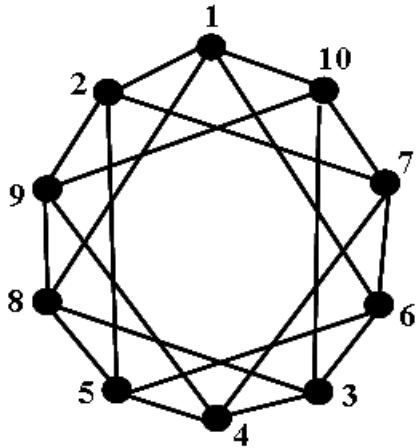 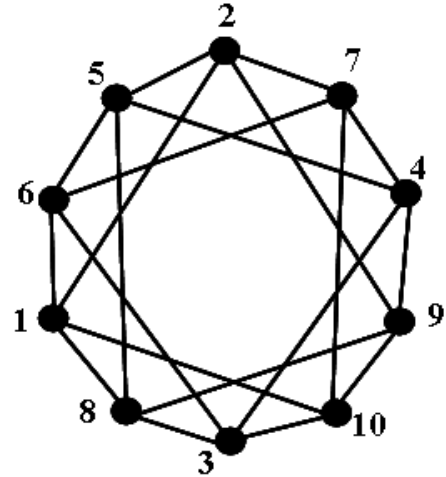

Рис. 3.5. Построение различных образующих циклов.

Образующий цикл, может быть построен как кольцевая сумма четырех изометрических циклов (см. рис. 3.5), используя правило покрытия всех вершин четырьмя циклами длиной четыре. Кольцевая сумма четырех циклов имеет пересечение трех циклов с четвертым циклом по одному элементу.

Например, $Q_{51}=c_1 \oplus c_3 \oplus c_{11} \oplus c_{17}=\{e_1,e_2,e_5,e_8\} \oplus \{e_3,e_4,e_{17},e_{18}\} \oplus \{e_6,e_7,e_{13},e_{16}\} \oplus \{e_1,e_3,e_6,e_{14}\}=$
$= \{e_2,e_4,e_5,e_7,e_8,e_{13},e_{14},e_{16},e_{17},e_{18}\}$. Рассмотрим перестановки индуцируемые циклом $Q_{51}$.

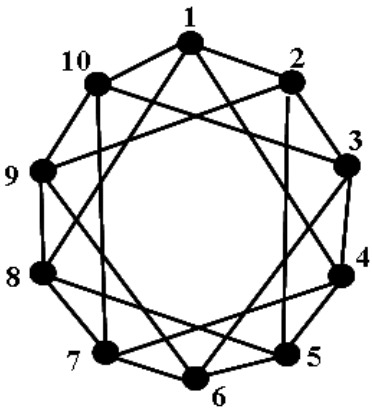 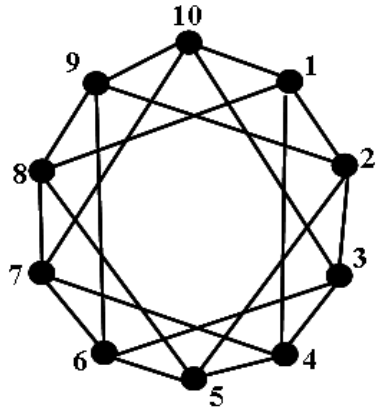 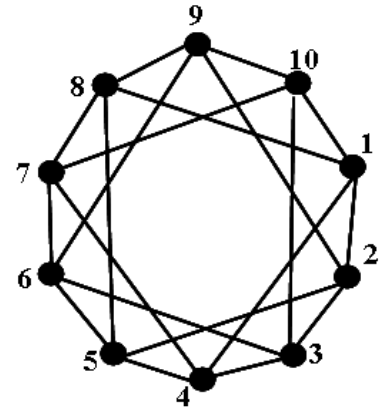

Рис. 3.6. Перестановка $p_0$.   Рис. 3.7. Перестановка $p_1$.   Рис. 3.8. Перестановка $p_2$.

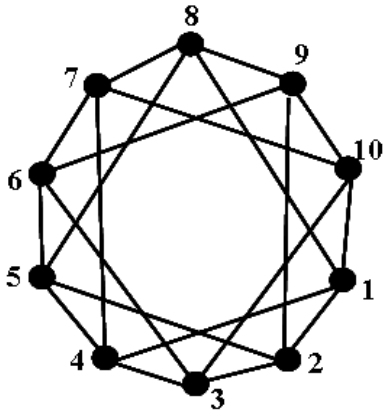
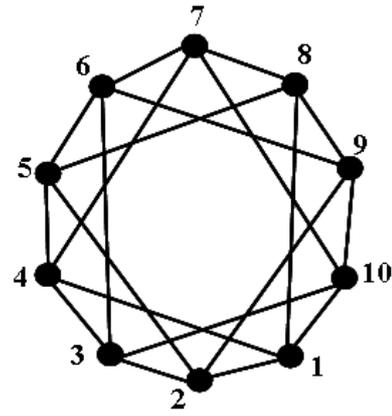
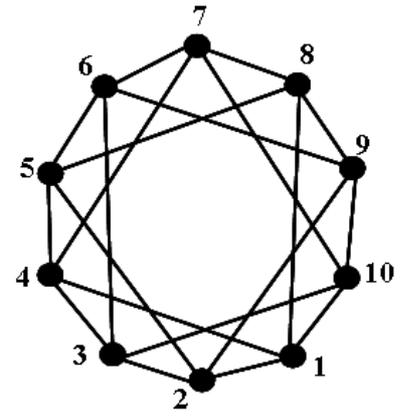

Рис. 3.9. Перестановка $p_3$.   Рис. 3.10. Перестановка $p_4$.   Рис. 3.11. Перестановка $p_5$.

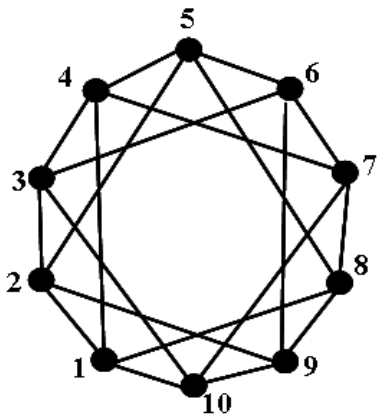
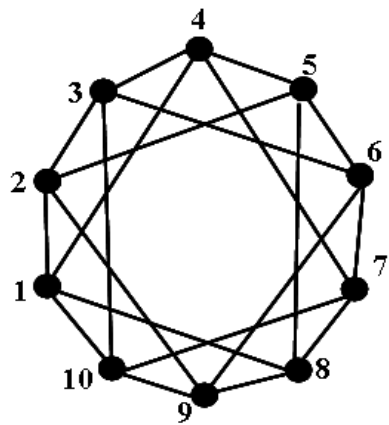
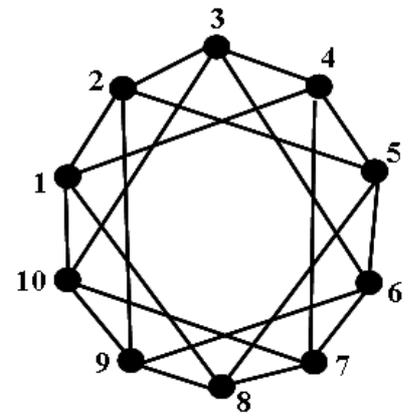

Рис. 3.12. Перестановка $p_6$.   Рис. 3.13. Перестановка $p_7$.   Рис. 3.14. Перестановка $p_8$.

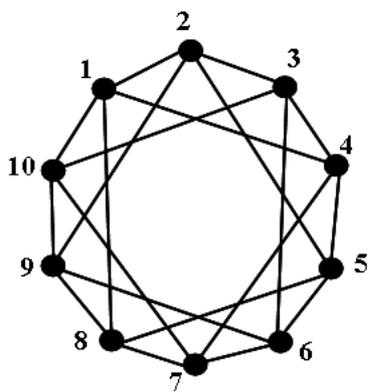
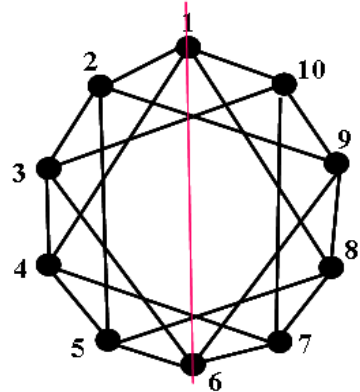
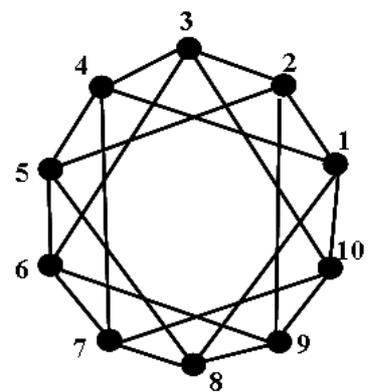

Рис. 3.15. Перестановка $p_9$.   Рис. 3.16. Перестановка $p_{10}$.   Рис. 3.17. Перестановка $p_{11}$.

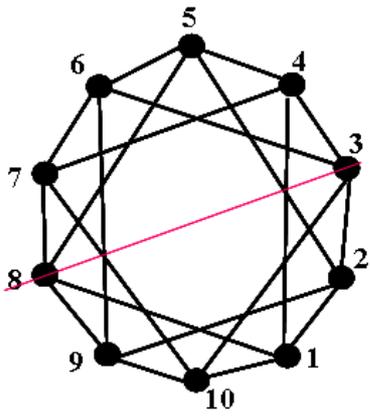
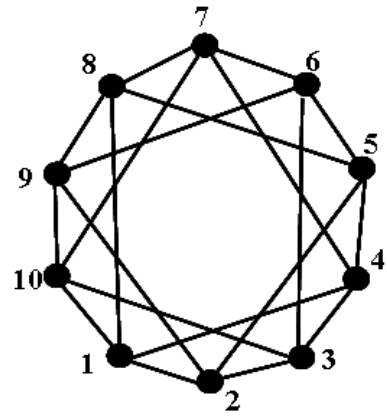
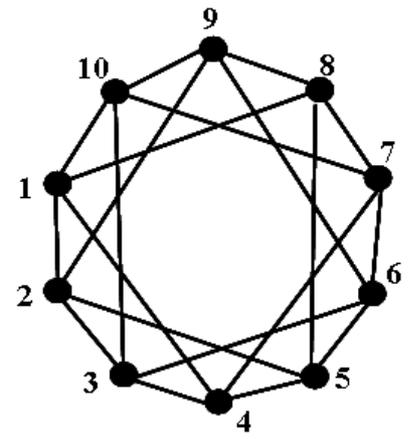

Рис. 3.18. Перестановка p₁₂.   Рис. 3.19. Перестановка p₁₃.   Рис. 3.20. Перестановка p₁₄.

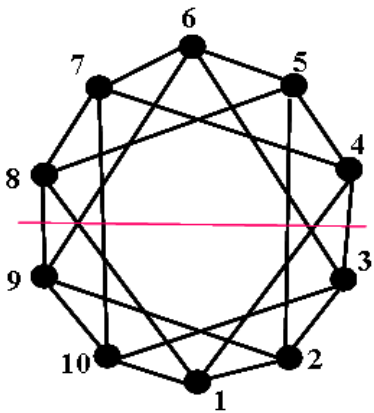
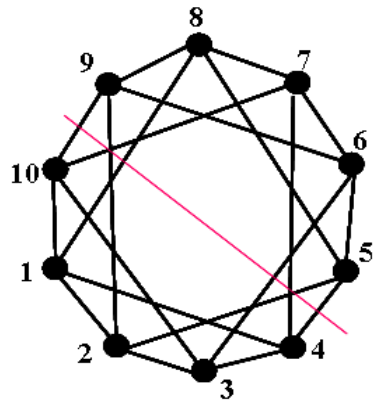
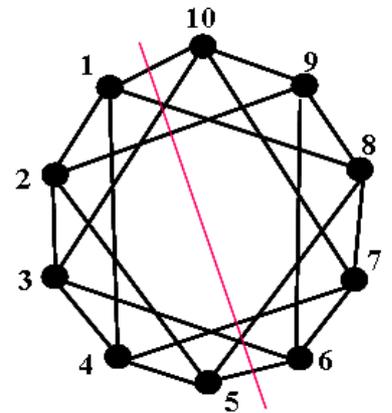

Рис. 3.21. Перестановка p₁₅.   Рис. 3.22. Перестановка p₁₆.   Рис. 3.23. Перестановка p₁₇.

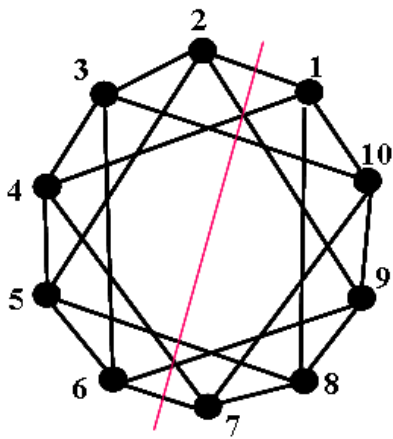
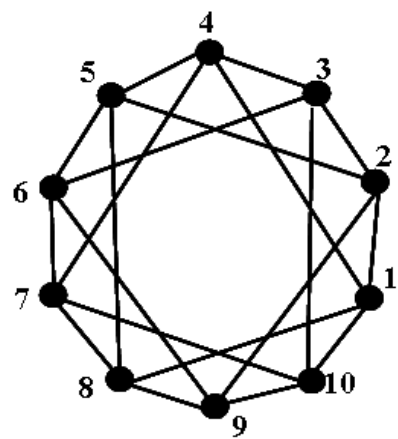

Рис. 3.24. Перестановка p₁₈.   Рис. 3.25. Перестановка p₁₉.

Перестановки для данного образующего цикла представлены на рис. 3.6 – рис. 3.25.

Подмножество изометрических циклов, образует образующий цикл длиной десять. Поэтому общее количество перестановок для данного опорного цикла, можно определить, как удвоенную длину образующего цикла.

$p_0$ = <1 2 3 4 5 6 7 8 9 10> = (1)(2)(3)(4)(5)(6)(7)(8)(9)(10);
$p_1$ = <10 1 2 3 4 5 6 7 8 9> = (1 10 9 8 7 6 5 4 3 2);
$p_2$ = <9 10 1 2 3 4 5 6 7 8> = (1 9 7 5 3)(2 10 8 6 4);
$p_3$ = <8 9 10 1 2 3 4 5 6 7> = (1 8 5 2 9 6 3 10 7 4);
$p_4$ = <7 8 9 10 1 2 3 4 5 6> = (1 7 3 9 5)(2 8 4 10 6);
$p_5$ = <6 7 8 9 10 1 2 3 4 5> = (1 6)(2 7)(3 8)(4 9)(5 10);
$p_6$ = <5 6 7 8 9 10 1 2 3 4> = (1 5 9 3 7)(2 6 10 4 8);
$p_7$ = <4 5 6 7 8 9 10 1 2 3> = (1 4 7 10 3 6 9 2 5 8).
$P_8$ = <3 4 5 6 7 8 9 10 1 2> = (1 3 5 7 9)(2 4 6 8 10);
$P_9$ = <2 3 4 5 6 7 8 9 10 1> = (1 2 3 4 5 6 7 8 9 10);
$P_{10}$ = <1 10 9 8 7 6 5 4 3 2> = (1)(2 10)(3 9)(4 8)(5 7)(6);
$P_{11}$ = <3 2 1 10 9 8 7 6 5 4> = (1 3)(2)(4 10)(5 9)(6 8)(7);
$P_{12}$ = <5 4 3 2 1 10 9 8 7 6> = (1 5)(2 4)(3)(6 10)(7 9)(8);
$P_{13}$ = <7 6 5 4 3 2 1 10 9 8> = (1 7)(2 6)(3 5)(4)(8 10)(9);
$P_{14}$ = <9 8 7 6 5 4 3 2 1 10> = (1 9)(2 8)(3 7)(4 6)(5)(10);
$P_{15}$ = <6 5 4 3 2 1 10 9 8 7> = (1 6)(2 5)(3 4)(7 10)(8 9);
$P_{16}$ = <8 7 6 5 4 3 2 1 10 9> = (1 8)(2 7)(3 6)(4 5)(9 10);
$p_{17}$ = <10 9 8 7 6 5 4 3 2 1> = (1 10)(2 9)(3 8)(4 7)(5 6);
$p_{18}$ = <2 1 10 9 8 7 6 5 4 3> = (1 2)(3 10)(4 9)(5 8)(6 7);
$p_{19}$ = <4 3 2 1 10 9 8 7 6 5> = (1 4)(2 3)(5 10)(6 9)(7 8).

Рассмотрим следующий образующий цикл <$v_1,v_4,v_7,v_{10},v_3,v_6,v_9,v_2,v_5,v_8$>. Определим перестановки для данного опорного цикла.

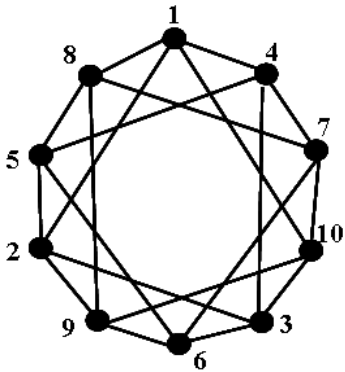 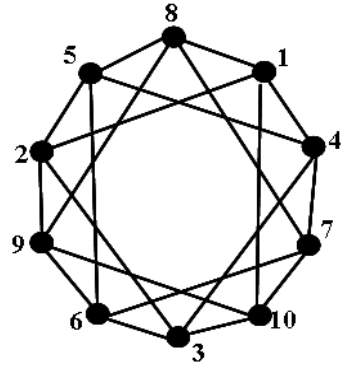 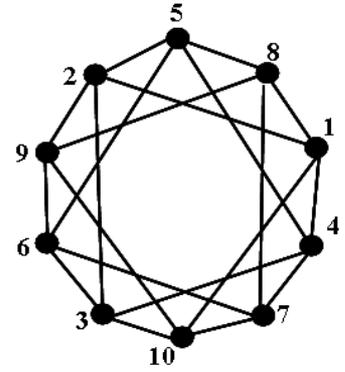

Рис. 3.26. Перестановка $p_{20}$.    Рис. 3.27. Перестановка $p_{21}$.    Рис. 3.28. Перестановка $p_{22}$.

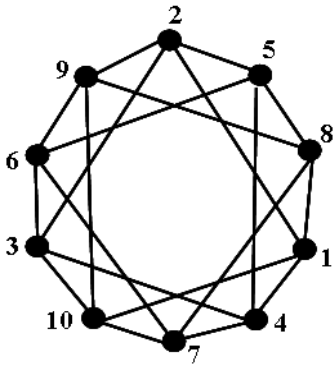

Рис. 3.29. Перестановка p₂₃.

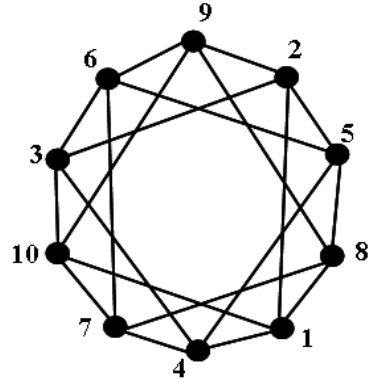

Рис. 3.30. Перестановка p₂₄.

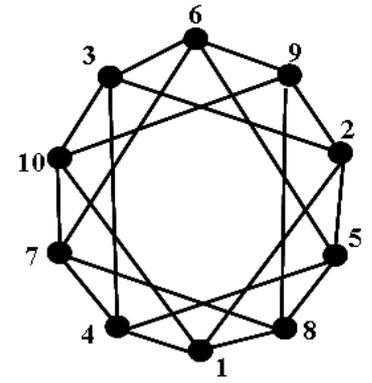

Рис. 3.31. Перестановка p₂₅.

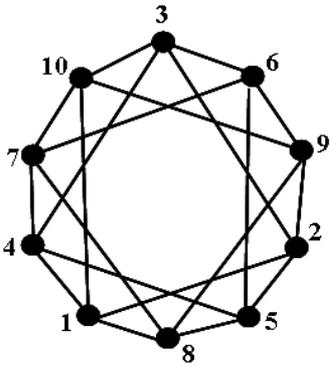

Рис. 3.32. Перестановка p₂₆.

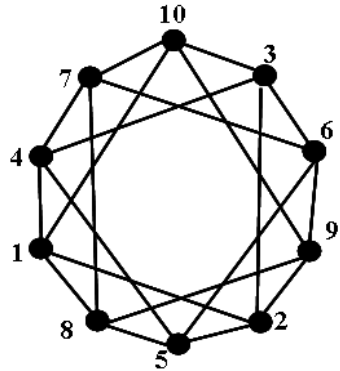

Рис. 3.33. Перестановка p₂₇.

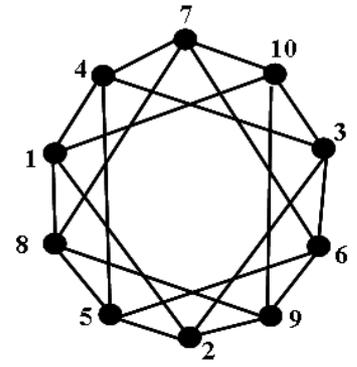

Рис. 3.34. Перестановка p₂₈.

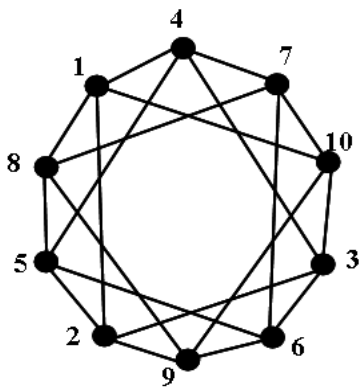

Рис. 3.35. Перестановка p₂₉.

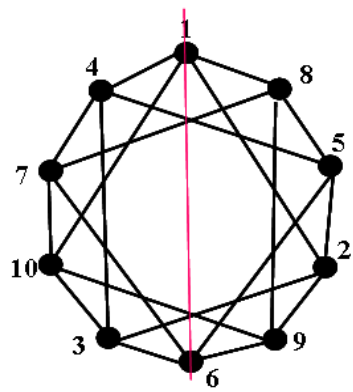

Рис. 3.36. Перестановка p₃₀.

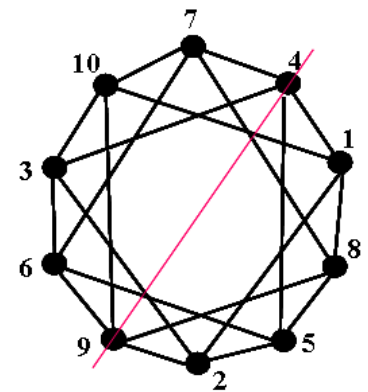

Рис. 3.37. Перестановка p₃₁.

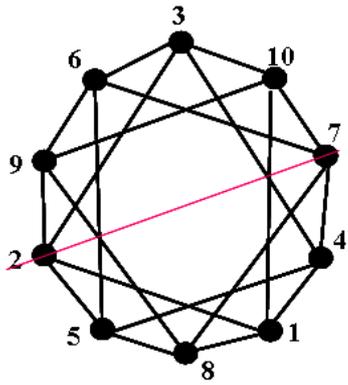 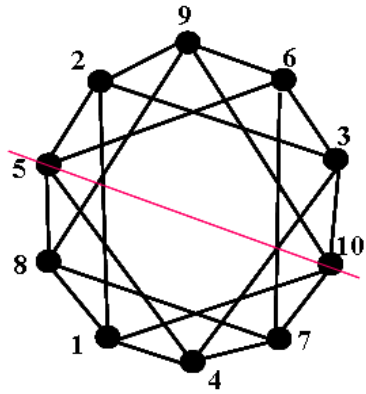 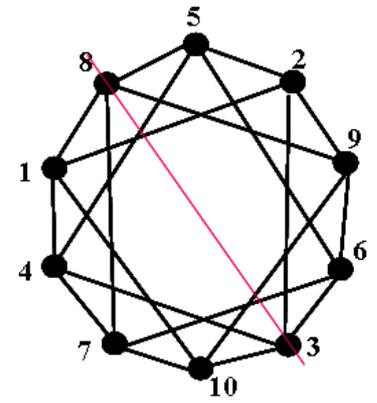

Рис. 3.38. Перестановка p₃₂.    Рис. 3.39. Перестановка p₃₃.    Рис. 3.40. Перестановка p₃₄.

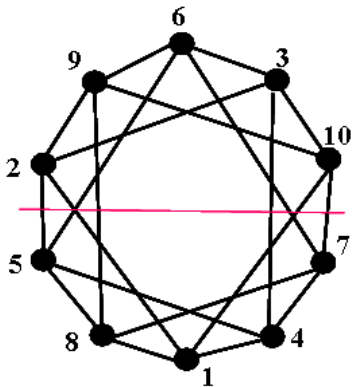 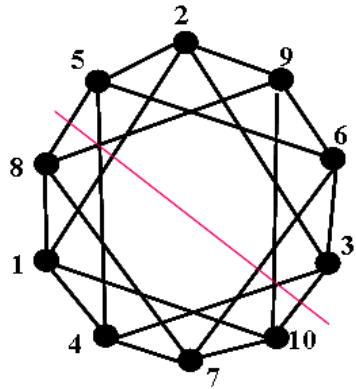 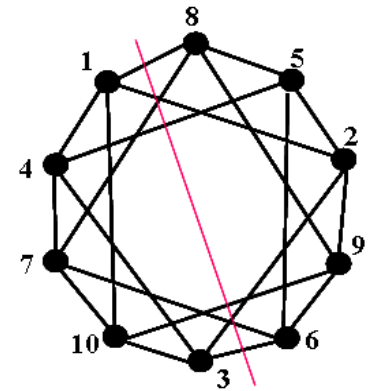

Рис. 3.41. Перестановка p₃₅.    Рис. 3.42. Перестановка p₃₆.    Рис. 3.43. Перестановка p₃₇.

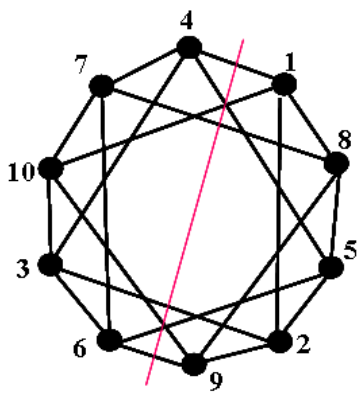 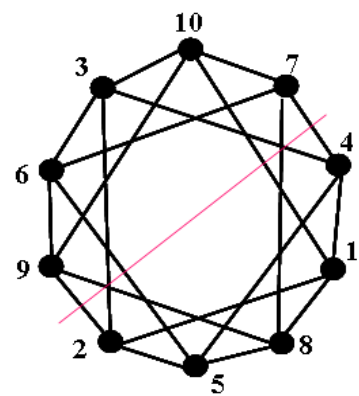

Рис. 3.44. Перестановка p₃₈.    Рис. 3.45. Перестановка p₃₉.

Определим следующие перестановки:

$p_{20}$ = <1 4 7 10 3 6 9 2 5 8> = (1)(2 4 10 8)(3 7 9 5)(6);
$p_{21}$ = <8 1 4 7 10 3 6 9 2 5> = (1 8 9 2)(3 4 7 6)(5 10);
$p_{22}$ = <5 8 1 4 7 10 3 6 9 2> = (1 5 7 3)(2 8 6 10)(4)(9);
$p_{23}$ = <2 5 8 1 4 7 10 3 6 9> = (1 2 5 4)(3 8)(6 7 10 9);
$p_{24}$ = <9 2 5 8 1 4 7 10 3 6> = (1 9 3 5)(2)(4 8 10 6)(7);
$p_{25}$ = <6 9 2 5 8 1 4 7 10 3> = (1 6)(2 9 10 3)(4 5 8 7);

p$_{26}$ = <3 6 9 2 5 8 1 4 7 10> = (1 3 9 7)(2 6 8 4)(5)(10);
p$_{27}$ = <10 3 6 9 2 5 8 1 4 7> = (1 10 7 8)(2 3 6 5)(4 9).
P$_{28}$ = <7 10 3 6 9 2 5 8 1 4> = (1 7 5 9)(2 10 4 6)(3)(8);
P$_{29}$ = <4 7 10 3 6 9 2 5 8 1> = (1 4 3 10)(2 7)(5 6 9 8);
P$_{30}$ = <1 8 5 2 9 6 3 10 7 4> = (1)(2 8 10 4)(3 5 9 7)(6);
P$_{31}$ = <7 4 1 8 5 2 9 6 3 10> = (1 7 9 3)(2 4 8 6)(5)(10);
P$_{32}$ = <3 10 7 4 1 8 5 2 9 6> = (1 5)(2 4)(3)(6 10)(7 9)(8);
P$_{33}$ = <9 6 3 10 7 4 1 8 5 2> = (1 9 5 7)(2 6 4 10)(3)(8);
P$_{34}$ = <5 2 9 6 3 10 7 4 1 8> = (1 5 3 9)(2)(4 6 10 8)(7);
P$_{35}$ = <6 3 10 7 4 1 8 5 2 9> = (1 6)(2 3 10 9)(4 7 8 5);
P$_{36}$ = <2 9 6 3 10 7 4 1 8 5> = (1 2 9 8)(3 6 7 4)(5 10);
P$_{37}$ = <8 5 2 9 6 3 10 7 4 1> = (1 8 7 10)(2 5 6 3)(4 9);
P$_{38}$ = <4 1 8 5 2 9 6 3 10 7> = (1 4 5 2)(3 8)(6 9 10 7);
P$_{39}$ = <10 7 4 1 8 5 2 9 6 3> = (1 10 3 4)(2 7)(5 8 9 6).

Общее количество перестановок графа G$_{12}$ можно определить как умножение количества образующих циклов на количество перестановок в диэдральном цикле D$_{10}$ равное 156×20 = 3120.

*Пример 3.4.* Рассмотрим следующий регулярный граф (см. рис. 3.46).

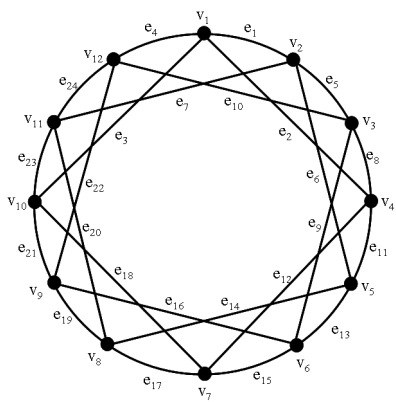

Рис. 4.46. Граф G$_{13}$.

Количество вершин графа = 12.
Количество ребер графа = 24.
Количество изометрических циклов = 27.

Смежность графа:

вершина v$_1$: v$_2$ v$_4$ v$_{10}$ v$_{12}$;           вершина v$_2$: v$_1$ v$_3$ v$_5$ v$_{11}$;
вершина v$_3$: v$_2$ v$_4$ v$_6$ v$_{12}$;              вершина v$_4$: v$_1$ v$_3$ v$_5$ v$_7$;
вершина v$_5$: v$_2$ v$_4$ v$_6$ v$_8$;                 вершина v$_6$: v$_3$ v$_5$ v$_7$ v$_9$;
вершина v$_7$: v$_4$ v$_6$ v$_8$ v$_{10}$;              вершина v$_8$: v$_5$ v$_7$ v$_9$ v$_{11}$;
вершина v$_9$: v$_6$ v$_8$ v$_{10}$ v$_{12}$;           вершина v$_{10}$: v$_1$ v$_7$ v$_9$ v$_{11}$;
вершина v$_{11}$: v$_2$ v$_8$ v$_{10}$ v$_{12}$;        вершина v$_{12}$: v$_1$ v$_3$ v$_9$ v$_{11}$.

Инцидентность графа:

ребро e$_1$: (v$_1$,v$_2$) или (v$_2$,v$_1$);           ребро e$_2$: (v$_1$,v$_4$) или (v$_4$,v$_1$);
ребро e$_3$: (v$_1$,v$_{10}$) или (v$_{10}$,v$_1$);     ребро e$_4$: (v$_1$,v$_{12}$) или (v$_{12}$,v$_1$);
ребро e$_5$: (v$_2$,v$_3$) или (v$_3$,v$_2$);           ребро e$_6$: (v$_2$,v$_5$) или (v$_5$,v$_2$);
ребро e$_7$: (v$_2$,v$_{11}$) или (v$_{11}$,v$_2$);     ребро e$_8$: (v$_3$,v$_4$) или (v$_4$,v$_3$);

ребро $e_9$: $(v_3,v_6)$ или $(v_6,v_3)$;   ребро $e_{10}$: $(v_3,v_{12})$ или $(v_{12},v_3)$;
ребро $e_{11}$: $(v_4,v_5)$ или $(v_5,v_4)$;   ребро $e_{12}$: $(v_4,v_7)$ или $(v_7,v_4)$;
ребро $e_{13}$: $(v_5,v_6)$ или $(v_6,v_5)$;   ребро $e_{14}$: $(v_5,v_8)$ или $(v_8,v_5)$;
ребро $e_{15}$: $(v_6,v_7)$ или $(v_7,v_6)$;   ребро $e_{16}$: $(v_6,v_9)$ или $(v_9,v_6)$;
ребро $e_{17}$: $(v_7,v_8)$ или $(v_8,v_7)$;   ребро $e_{18}$: $(v_7,v_{10})$ или $(v_{10},v_7)$;
ребро $e_{19}$: $(v_8,v_9)$ или $(v_9,v_8)$;   ребро $e_{20}$: $(v_8,v_{11})$ или $(v_{11},v_8)$;
ребро $e_{21}$: $(v_9,v_{10})$ или $(v_{10},v_9)$;   ребро $e_{22}$: $(v_9,v_{12})$ или $(v_{12},v_9)$;
ребро $e_{23}$: $(v_{10},v_{11})$ или $(v_{11},v_{10})$;   ребро $e_{24}$: $(v_{11},v_{12})$ или $(v_{12},v_{11})$.

Множество базовых реберных разрезов графа:

$w_0(e_1) = \{e_2,e_3,e_4,e_5,e_6,e_7\}$;   $w_0(e_2) = \{e_1,e_3,e_4,e_8,e_{11},e_{12}\}$;
$w_0(e_3) = \{e_1,e_2,e_4,e_{18},e_{21},e_{23}\}$;   $w_0(e_4) = \{e_1,e_2,e_3,e_{10},e_{22},e_{24}\}$;
$w_0(e_5) = \{e_1,e_6,e_7,e_8,e_9,e_{10}\}$;   $w_0(e_6) = \{e_1,e_5,e_7,e_{11},e_{13},e_{14}\}$;
$w_0(e_7) = \{e_1,e_5,e_6,e_{20},e_{23},e_{24}\}$;   $w_0(e_8) = \{e_2,e_5,e_9,e_{10},e_{11},e_{12}\}$;
$w_0(e_9) = \{e_5,e_8,e_{10},e_{13},e_{15},e_{16}\}$;   $w_0(e_{10}) = \{e_4,e_5,e_8,e_9,e_{22},e_{24}\}$;
$w_0(e_{11}) = \{e_2,e_6,e_8,e_{12},e_{13},e_{14}\}$;   $w_0(e_{12}) = \{e_2,e_8,e_{11},e_{15},e_{17},e_{18}\}$;
$w_0(e_{13}) = \{e_6,e_9,e_{11},e_{14},e_{15},e_{16}\}$;   $w_0(e_{14}) = \{e_6,e_{11},e_{13},e_{17},e_{19},e_{20}\}$;
$w_0(e_{15}) = \{e_9,e_{12},e_{13},e_{16},e_{17},e_{18}\}$;   $w_0(e_{16}) = \{e_9,e_{13},e_{15},e_{19},e_{21},e_{22}\}$;
$w_0(e_{17}) = \{e_{12},e_{14},e_{15},e_{18},e_{19},e_{20}\}$;   $w_0(e_{18}) = \{e_3,e_{12},e_{15},e_{17},e_{21},e_{23}\}$;
$w_0(e_{19}) = \{e_{14},e_{16},e_{17},e_{20},e_{21},e_{22}\}$;   $w_0(e_{20}) = \{e_7,e_{14},e_{17},e_{19},e_{23},e_{24}\}$;
$w_0(e_{21}) = \{e_3,e_{16},e_{18},e_{19},e_{22},e_{23}\}$;   $w_0(e_{22}) = \{e_4,e_{10},e_{16},e_{19},e_{21},e_{24}\}$;
$w_0(e_{23}) = \{e_3,e_7,e_{18},e_{20},e_{21},e_{24}\}$;   $w_0(e_{24}) = \{e_4,e_7,e_{10},e_{20},e_{22},e_{23}\}$.

Кортеж весов ребер:
<14,16,16,14,14,16,16,14,16,16,14,16,14,16,14,16,14,16,14,16,14,16,14,14>.
Кортеж весов вершин:
<60,60,60,60,60,60,60,60,60,60,60,60>.

Множество изометрических циклов графа:

$c_1 = \{e_1,e_2,e_5,e_8\}$;   $c_2 = \{e_1,e_2,e_6,e_{11}\}$;
$c_3 = \{e_1,e_3,e_7,e_{23}\}$;   $c_4 = \{e_1,e_4,e_5,e_{10}\}$;
$c_5 = \{e_1,e_4,e_7,e_{24}\}$;   $c_6 = \{e_2,e_3,e_{12},e_{18}\}$;
$c_7 = \{e_2,e_4,e_8,e_{10}\}$;   $c_8 = \{e_3,e_4,e_{21},e_{22}\}$;
$c_9 = \{e_3,e_4,e_{23},e_{24}\}$;   $c_{10} = \{e_5,e_6,e_8,e_{11}\}$;
$c_{11} = \{e_5,e_6,e_9,e_{13}\}$;   $c_{12} = \{e_5,e_7,e_{10},e_{24}\}$;
$c_{13} = \{e_6,e_7,e_{14},e_{20}\}$;   $c_{14} = \{e_8,e_9,e_{11},e_{13}\}$;
$c_{15} = \{e_8,e_9,e_{12},e_{15}\}$;   $c_{16} = \{e_9,e_{10},e_{16},e_{22}\}$;
$c_{17} = \{e_{11},e_{12},e_{13},e_{15}\}$;   $c_{18} = \{e_{11},e_{12},e_{14},e_{17}\}$;
$c_{19} = \{e_{13},e_{14},e_{15},e_{17}\}$;   $c_{20} = \{e_{13},e_{14},e_{16},e_{19}\}$;
$c_{21} = \{e_{15},e_{16},e_{17},e_{19}\}$;   $c_{22} = \{e_{15},e_{16},e_{18},e_{21}\}$;
$c_{23} = \{e_{17},e_{18},e_{19},e_{21}\}$;   $c_{24} = \{e_{17},e_{18},e_{20},e_{23}\}$;
$c_{25} = \{e_{19},e_{20},e_{21},e_{23}\}$;   $c_{26} = \{e_{19},e_{20},e_{22},e_{24}\}$;
$c_{27} = \{e_{21},e_{22},e_{23},e_{24}\}$.

Образующий цикл, может быть построен как кольцевая сумма пяти изометрических циклов (см. рис. 3.46), используя правило покрытия всех вершин пятью циклами длиной четыре. Общее число конфигураций = 2229. Количество образующих циклов = 335.

Перечислим образующие циклы:

$Q_1=\{e_2,e_4,e_6,e_7,e_8,e_9,e_{14},e_{15},e_{17},e_{21},e_{22},e_{23}\}$;

$Q_2=\{e_2,e_4,e_6,e_7,e_8,e_9,e_{13},e_{17},e_{18},e_{19},e_{22},e_{23}\}$;
$Q_3=\{e_2,e_4,e_6,e_7,e_8,e_9,e_{13},e_{17},e_{18},e_{20},e_{21},e_{22}\}$;
$Q_4=\{e_2,e_4,e_5,e_6,e_9,e_{12},e_{14},e_{15},e_{20},e_{21},e_{22},e_{23}\}$;
$Q_5=\{e_2,e_4,e_5,e_6,e_8,e_{13},e_{15},e_{17},e_{20},e_{21},e_{22},e_{23}\}$;
$Q_6=\{e_2,e_4,e_5,e_6,e_8,e_{14},e_{15},e_{16},e_{18},e_{20},e_{22},e_{23}\}$;
$Q_7=\{e_2,e_4,e_5,e_7,e_9,e_{12},e_{13},e_{14},e_{17},e_{21},e_{22},e_{23}\}$;
$Q_8=\{e_2,e_4,e_5,e_7,e_9,e_{11},e_{14},e_{15},e_{17},e_{21},e_{22},e_{23}\}$;
$Q_9=\{e_2,e_4,e_5,e_7,e_9,e_{11},e_{13},e_{17},e_{18},e_{19},e_{22},e_{23}\}$;
$Q_{10}=\{e_2,e_4,e_5,e_7,e_9,e_{11},e_{13},e_{17},e_{18},e_{20},e_{21},e_{22}\}$;
$Q_{11}=\{e_2,e_4,e_5,e_7,e_8,e_{13},e_{14},e_{16},e_{17},e_{18},e_{22},e_{23}\}$;
$Q_{12}=\{e_2,e_4,e_5,e_7,e_8,e_{13},e_{14},e_{15},e_{18},e_{19},e_{22},e_{23}\}$;
$Q_{13}=\{e_2,e_4,e_5,e_7,e_8,e_{13},e_{14},e_{15},e_{18},e_{20},e_{21},e_{22}\}$;
$Q_{14}=\{e_2,e_3,e_6,e_7,e_8,e_{10},e_{14},e_{15},e_{16},e_{17},e_{22},e_{23}\}$;
$Q_{15}=\{e_2,e_3,e_6,e_7,e_8,e_{10},e_{13},e_{15},e_{17},e_{19},e_{22},e_{23}\}$;
$Q_{16}=\{e_2,e_3,e_6,e_7,e_8,e_{10},e_{13},e_{16},e_{17},e_{18},e_{20},e_{22}\}$;
$Q_{17}=\{e_2,e_3,e_6,e_7,e_8,e_9,e_{14},e_{15},e_{17},e_{21},e_{22},e_{24}\}$;
$Q_{18}=\{e_2,e_3,e_6,e_7,e_8,e_9,e_{13},e_{17},e_{18},e_{19},e_{22},e_{24}\}$;
$Q_{19}=\{e_2,e_3,e_5,e_6,e_{10},e_{12},e_{14},e_{15},e_{16},e_{20},e_{22},e_{23}\}$;
$Q_{20}=\{e_2,e_3,e_5,e_6,e_9,e_{12},e_{14},e_{15},e_{19},e_{22},e_{23},e_{24}\}$;

……………………………………………………….

$Q_{333}=\{e_3,e_4,e_5,e_6,e_9,e_{11},e_{12},e_{16},e_{17},e_{20},e_{21},e_{24}\}$;
$Q_{334}=\{e_3,e_4,e_5,e_6,e_8,e_{12},e_{13},e_{16},e_{17},e_{20},e_{21},e_{24}\}$;
$Q_{335}=\{e_3,e_4,e_6,e_7,e_8,e_9,e_{12},e_{13},e_{17},e_{19},e_{22},e_{23}\}$;

Рассмотрим образующий цикл $Q_1 = \{e_2,e_4,e_6,e_7,e_8,e_9,e_{14},e_{15},e_{17},e_{21},e_{22},e_{23}\}$ =

= $c_1 \oplus c_3 \oplus c_8 \oplus c_{11} \oplus c_{19}$ = $\{e_1,c_2,c_5,c_8\} \oplus \{e_1,c_3,c_7,c_{23}\} \oplus \{e_3,c_4,c_{21},c_{22}\} \oplus \{e_5,c_6,c_9,c_{13}\} \oplus$

$\oplus \{e_{13},c_{14},c_{15},c_{17}\}$.

Для образующего цикла $Q_1$ составим подгруппу перестановок:

$p_0$ = <1 4 3 6 7 8 5 2 11 10 9 12> = (1)(4 6 8 2)(3)(7 5)(11 9)(10 12);
$p_1$ = <12 1 4 3 6 7 8 5 2 11 10 9> = (12 9 2 1)(4 3)(6 7 8 5)(11 10 9);
$p_2$ = <9 12 1 4 3 6 7 8 5 2 11 10> = (9 5 3 1)(12 10 2)(4)(6)(7)(8)(11);
$p_3$ = <10 9 12 1 4 3 6 7 8 5 2 11> = (10 5 4 1)(12 11 2 9 8 7 6 3);
$p_4$ = <11 10 9 12 1 4 3 6 7 8 5 2> = (11 5 1)(10 8 6 4 12 2)(9 7 3);
$p_5$ = <2 11 10 9 12 1 4 3 6 7 8 5> = (2 11 8 3 10 7 4 9 6 1)(12 5);
$p_6$ = <5 2 11 10 9 12 1 4 3 6 7 8> = (5 9 3 11 7 1)(2)(10 6 12 8 4);
$p_7$ = <8 5 2 11 10 9 12 1 4 3 6 7> = (8 1)(5 10 3 2)(11 6 9 4)(12 7);
$p_8$ = <7 8 5 2 11 10 9 12 1 4 3 6> = (7 9 1)(8 12 6 10 4 2)(5 11 3);
$p_9$ = <6 7 8 5 2 11 10 9 12 1 4 3> = (6 11 4 5 2 7 10 1)(8 9 12 3);
$p_{10}$ = <3 6 7 8 5 2 11 10 9 12 1 4> = (3 7 11 1)(6 2)(8 10 12 4)(5)(9);
$p_{11}$ = <4 3 6 7 8 5 2 11 10 9 12 1> = (4 7 2 3 6 5 8 11 12 1)(10 9);
$p_{12}$ = <1 12 9 10 11 2 5 8 7 6 3 4> = (1)(12 4 10 6 2)(9 7 5 11 3)(8);
$p_{13}$ = <3 4 1 12 9 10 11 2 5 8 7 6> = (3 1)(4 12 6 10 8 2)(9 5)(11 7);
$p_{14}$ = <7 6 3 4 1 12 9 10 11 2 5 8> = (7 9 11 5 1)(6 12 8 10 2)(3)(4);
$p_{15}$ = <5 8 7 6 3 4 1 12 9 10 11 2> = (5 3 7 1)(8 12 2)(6 4)(9)(10)(11);
$p_{16}$ = <11 2 5 8 7 6 3 4 1 12 9 10> = (11 9 1)(2)(5 7 3)(8 4)(6)(12 10);
$p_{17}$ = <9 10 11 2 5 8 7 6 3 4 1 12> = (9 3 11 1)(10 4 2)(5)(8 6)(7)(12);
$p_{18}$ = <4 1 12 9 10 11 2 5 8 7 6 3> = (4 9 8 5 10 7 2 1)(12 3)(11 6);
$p_{19}$ = <6 3 4 1 12 9 10 2 11 5 8 7> = (6 9 11 8 2 3 4 1)(12 7 10 5);
$p_{20}$ = <8 7 6 3 4 1 12 9 11 10 2 5> = (8 9 11 2 7 12 5 4 3 6 1)(10);
$p_{21}$ = <2 5 8 7 6 3 4 1 12 9 10 11> = (2 5 6 3 8 1)(7 4)(12 11 10 9);
$p_{22}$ = <10 11 2 5 8 7 6 3 4 1 12 9> = (10 1)(11 12 9 4 5 8 3 2)(7 6);

p₂₃ = <12 9 10 11 2 5 8 7 6 3 4 1> = (12 1)(9 6 5 2)(10 3)(11 4)(8 7).

Количество подгрупп определяется количеством образующих циклов. Каждая подгруппа состоит из 24 перестановок. Общее количество перестановок можно определить как 335× 24 = 8540.

*Пример 3.5.*

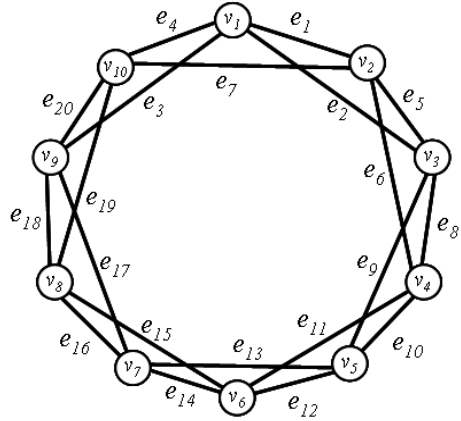

Рис. 3.47. Граф $G_{14}$.

Количество вершин графа = 10. Количество ребер графа = 20
Количество изометрических циклов = 12
Смежность графа $G_{14}$:

$v_1$: {$v_2,v_3,v_9,v_{10}$};    $v_2$: {$v_1,v_3,v_4,v_{10}$};
$v_3$: {$v_1,v_2,v_4,v_5$};    $v_4$: {$v_2,v_3,v_5,v_6$};
$v_5$: {$v_3,v_4,v_6,v_7$};    $v_6$: {$v_4,v_5,v_7,v_8$};
$v_7$: {$v_5,v_6,v_8,v_9$};    $v_8$: {$v_6,v_7,v_9,v_{10}$};
$v_9$: {$v_1,v_7,v_8,v_{10}$};    $v_{10}$: {$v_1,v_2,v_8,v_9$}.

Инцидентность графа $G_{14}$:

$e_1 = (v_1,v_2) \vee (v_2,v_1)$;    $e_2 = (v_1,v_3) \vee (v_3,v_1)$;
$e_3 = (v_1,v_9) \vee (v_9,v_1)$;    $e_4 = (v_1,v_{10}) \vee (v_{10},v_1)$;
$e_5 = (v_2,v_3) \vee (v_3,v_2)$;    $e_6 = (v_2,v_4) \vee (v_4,v_2)$;
$e_7 = (v_2,v_{10}) \vee (v_{10},v_2)$;    $e_8 = (v_3,v_4) \vee (v_4,v_3)$;
$e_9 = (v_3,v_5) \vee (v_5,v_3)$;    $e_{10} = (v_4,v_5) \vee (v_5,v_4)$;
$e_{11} = (v_4,v_6) \vee (v_6,v_4)$;    $e_{12} = (v_5,v_6) \vee (v_6,v_5)$;
$e_{13} = (v_5,v_7) \vee (v_7,v_5)$;    $e_{14} = (v_6,v_7) \vee (v_7,v_6)$;
$e_{15} = (v_6,v_8) \vee (v_8,v_6)$;    $e_{16} = (v_7,v_8) \vee (v_8,v_7)$;
$e_{17} = (v_7,v_9) \vee (v_9,v_7)$;    $e_{18} = (v_8,v_9) \vee (v_9,v_8)$;
$e_{19} = (v_8,v_{10}) \vee (v_{10},v_8)$;    $e_{20} = (v_9,v_{10}) \vee (v_{10},v_9)$.

Множество изометрических циклов графа:

$c_1 = \{e_1,e_2,e_5\} \to \{v_1,v_2,v_3\}$;    $c_2 = \{e_1,e_4,e_7\} \to \{v_1,v_2,v_{10}\}$;
$c_3 = \{e_2,e_3,e_9,e_{13},e_{17}\} \to \{v_1,v_3,v_5,v_7,v_9\}$;    $c_4 = \{e_3,e_4,e_{20}\} \to \{v_1,v_9,v_{10}\}$;
$c_5 = \{e_5,e_6,e_8\} \to \{v_2,v_3,v_4\}$;    $c_6 = \{e_6,e_7,e_{11},e_{15},e_{19}\} \to \{v_2,v_4,v_6,v_8,v_{10}\}$;
$c_7 = \{e_8,e_9,e_{10}\} \to \{v_3,v_4,v_5\}$;    $c_8 = \{e_{10},e_{11},e_{12}\} \to \{v_4,v_5,v_6\}$;
$c_9 = \{e_{12},e_{13},e_{14}\} \to \{v_5,v_6,v_7\}$;    $c_{10} = \{e_{14},e_{15},e_{16}\} \to \{v_6,v_7,v_8\}$;
$c_{11} = \{e_{16},e_{17},e_{18}\} \to \{v_7,v_8,v_9\}$;    $c_{12} = \{e_{18},e_{19},e_{20}\} \to \{v_8,v_9,v_{10}\}$.

В данном графе существует один образующий цикл:

$c_{o1} = c_1 \oplus c_3 \oplus c_4 \oplus c_7 \oplus c_9 \oplus c_{11} = c_2 \oplus c_5 \oplus c_6 \oplus c_8 \oplus c_{10} \oplus c_{12} =$
$= \{e_1, e_4, e_5, e_8, e_{10}, e_{12}, e_{14}, e_{16}, e_{18}, e_{20}\} \to \{v_1, v_2, v_3, v_4, v_5, v_6, v_7, v_8, v_9, v_{10}\}.$

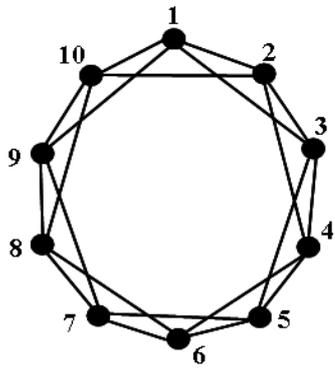 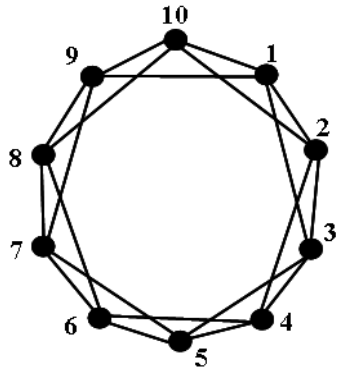 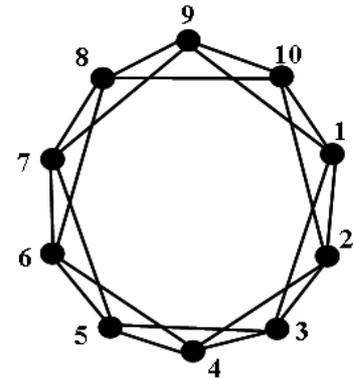

Рис. 3.48. Перестановка $p_0$.  Рис. 3.49. Перестановка $p_1$.  Рис. 3.50. Перестановка $p_2$.

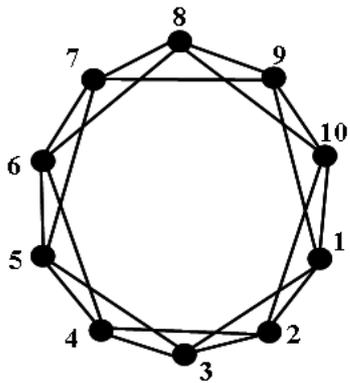 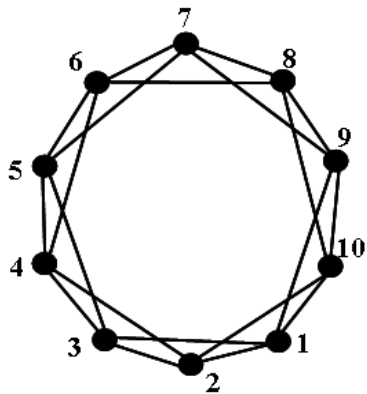 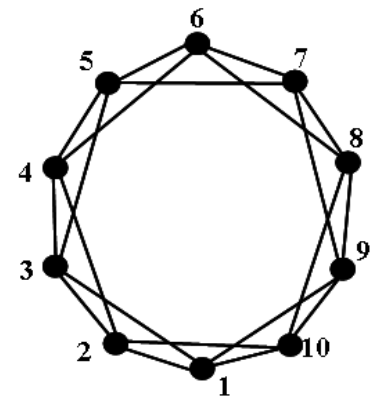

Рис. 3.51. Перестановка $p_3$.  Рис. 3.52. Перестановка $p_4$.  Рис. 3.53. Перестановка $p_5$.

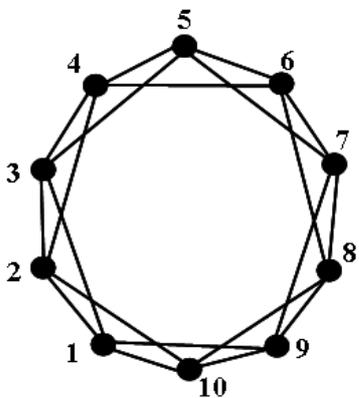 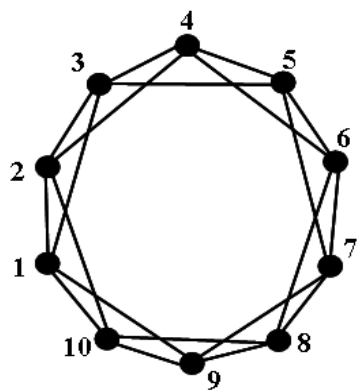 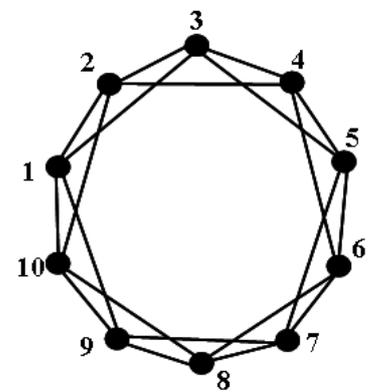

Рис. 3.54. Перестановка $p_6$.  Рис. 3.55. Перестановка $p_7$.  Рис. 3.56. Перестановка $p_8$.

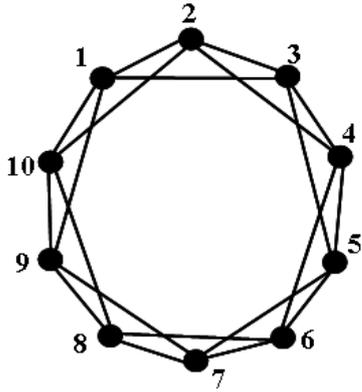 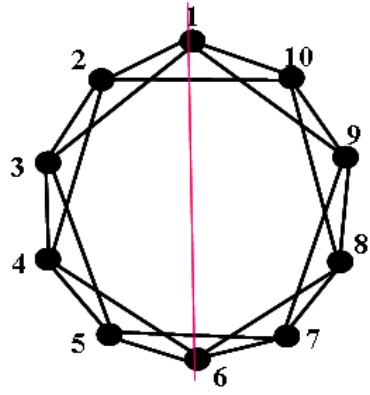 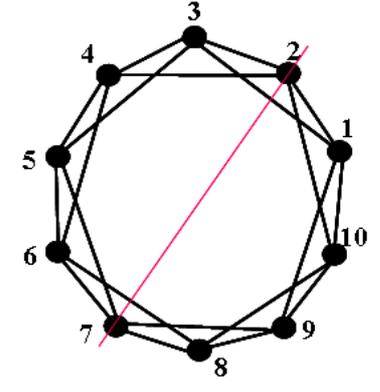

Рис. 3.57. Перестановка p₉.   Рис. 3.58. Перестановка p₁₀.   Рис. 3.59. Перестановка p₁₁.

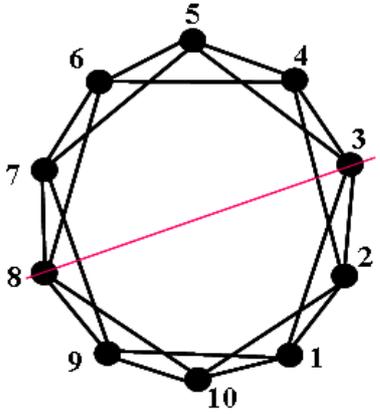 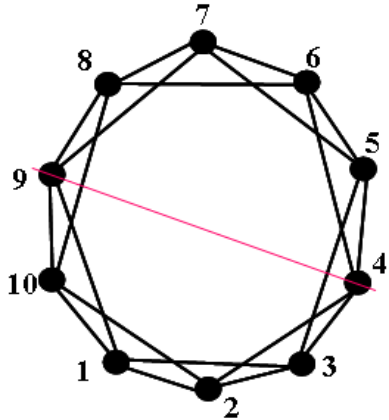 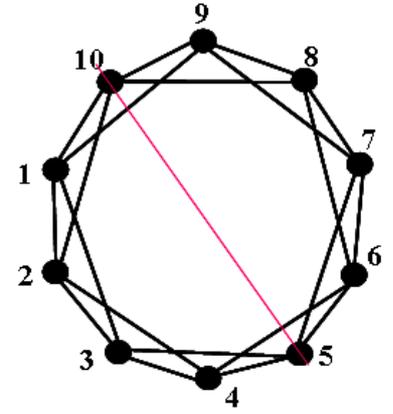

Рис. 3.60. Перестановка p₁₂.   Рис. 3.61. Перестановка p₁₃.   Рис. 3.62. Перестановка p₁₄.

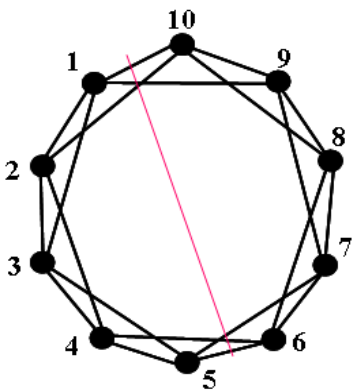 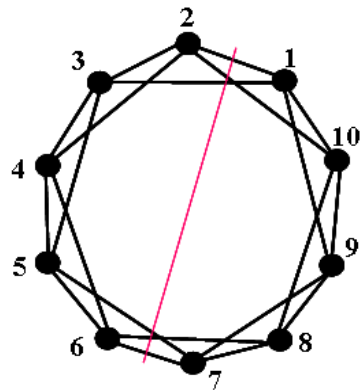 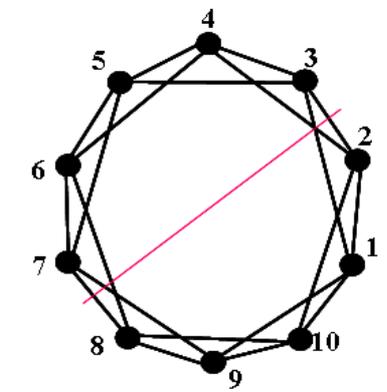

Рис. 3.63. Перестановка p₁₅.   Рис. 3.64. Перестановка p₁₆.   Рис. 3.65. Перестановка p₁₇.

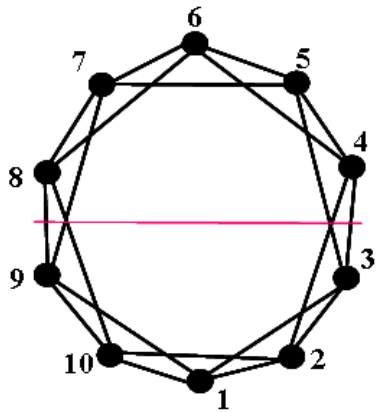
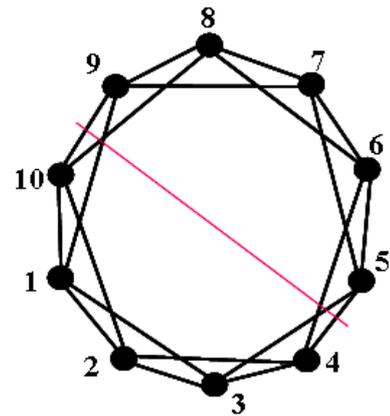

Рис. 3.66. Перестановка $p_{18}$.  Рис. 3.67. Перестановка $p_{19}$.

Определим все перестановки графа $G_{14}$:

$p_0$ = <1 2 3 4 5 6 7 8 9 10> = (1)(2)(3)(4)(5)(6)(7)(8)(9)(10);
$p_1$ = <10 1 2 3 4 5 6 7 8 9> = (1 10 9 8 7 6 5 4 3 2);
$p_2$ = <9 10 1 2 3 4 5 6 7 8> = (1 9 7 5 3)(2 10 8 6 4);
$p_3$ = <8 9 10 1 2 3 4 5 6 7> = (1 8 5 2 9 6 3 10 7 4);
$p_4$ = <7 8 9 10 1 2 3 4 5 6> = (1 7 3 9 5)(2 8 4 10 6);
$p_5$ = <6 7 8 9 10 1 2 3 4 5> = (1 6)(2 7)(3 8)(4 9)(5 10);
$p_6$ = <5 6 7 8 9 10 1 2 3 4> = (1 5 9 3 7)(2 6 10 4 8);
$p_7$ = <4 5 6 7 8 9 10 1 2 3> = (1 4 7 10 3 6 9 2 5 8);
$P_8$ = <3 4 5 6 7 8 9 10 1 2> = (1 3 5 7 9)(2 4 6 8 10);
$P_9$ = <2 3 4 5 6 7 8 9 10 1> = (1 2 3 4 5 6 7 8 9 10);
$P_{10}$ = <1 10 9 8 7 6 5 4 3 2> = (1)(2 10)(3 9)(4 8)(5 7)(6);
$P_{11}$ = <3 2 1 10 9 8 7 6 5 4> = (1 3)(2)(4 10)(5 9)(6 8)(7);
$P_{12}$ = <5 4 3 2 1 10 9 8 7 6> = (1 5)(2 4)(3)(6 10)(7 9)(8);
$P_{13}$ = <7 6 5 4 3 2 1 10 9 8> = (1 7)(2 6)(3 5)(4)(8 10)(9);
$P_{14}$ = <9 8 7 6 5 4 3 2 1 10> = (1 9)(2 8)(3 7)(4 6)(5)(10);
$P_{15}$ = <10 9 8 7 6 5 4 3 2 1> = (1 10)(2 9)(3 8)(4 7)(5 6);
$P_{16}$ = <2 1 10 9 8 7 6 5 4 3> = (1 2)(3 10)(4 9)(5 8)(6 7);
$p_{17}$ = <4 3 2 1 10 9 8 7 6 5> = (1 4)(2 3)(5 10)(6 9)(7 8);
$p_{18}$ = <6 5 4 3 2 1 10 9 8 7> = (1 6)(2 5)(3 4)(7 10)(8 9);
$p_{19}$ = <8 7 6 5 4 3 2 1 10 9> = (1 8)(2 7)(3 6)(4 5)(9 10).

Таблица Кэли в данном случае представляется одним блоком состоящим из 20 строк и 20 столбцов.

Веса вершин векторного инварианта для спектра реберных разрезов графа $G_{14}$ равны между собой $\zeta_w(G_{14})$ =<10×64>. Веса ребер также равны между собой $\xi_w(G_{14})$ = <20×16>. Орбита графа характеризуется множеством вершин образующего цикла lkbyjq 10.

Граф $G_{14}$ является планарным, кольцевая сумма изометрических циклов есть пустое множество $\sum_{i=1}^{12} c_i = \varnothing$.

**Замечание.** Кажется, что образующий цикл можно представить в виде двух связанных изометрических циклов длиной пять (см. рис. 3.68).

Такое представление рисунка графа приводит к нарушению смежности графа. При перестановке строк и столбцов матрица смежностей графа имеет другой вид.

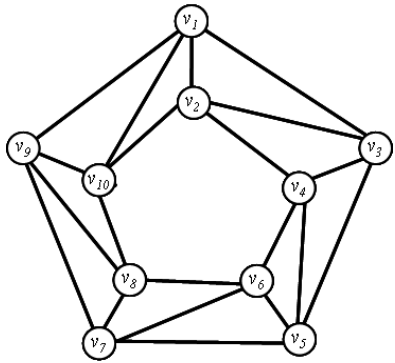 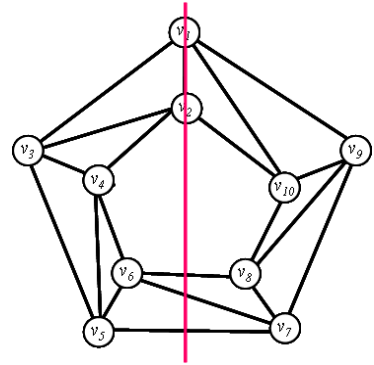

Рис. 3.68. Представление графа двумя циклами.  Рис. 3.69. Симметричное отражение.

### Комментарии

В данной главе, представлены методы определения группы автоморфизмов для регулярных графов с валентностью равной трем и четырем. Для демонстрации метода выбрана только часть регулярных графов. Рассмотрены методы построения образующих циклов для данной группы графов.

## Глава 4. АВТОМОРФИЗМ СИЛЬНО РЕГУЛЯРНЫХ ГРАФОВ

Пусть G — регулярный граф с вершинами и степенью . Говорят, что граф является *сильно регулярным*, если существуют целые числа $\lambda$ и $\mu$ такие, что:

- любые две смежные вершины имеют $\lambda$ общих соседей.
- любые две несмежные вершины имеют $\mu$ общих соседей.

Графы описанного типа иногда обозначаются как $\mathrm{srg}(v,k,\lambda,\mu)$. Четыре параметра $\mathrm{srg}(v,k,\lambda,\mu)$ не являются независимыми и должны удовлетворять следующему условию:

$$\mu(v-k-1) = k(k-\lambda-1) \tag{4.1}$$

Будем рассматривать вопросы построения группы автоморфизмов для сильно регулярных графов.

### 4.1. Граф Петерсена

*Пример 4.1*. Рассмотрим непланарный граф Петерсена (см. рис. 4.1).

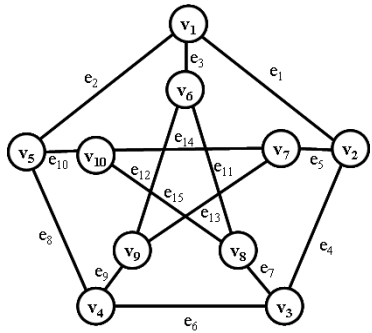

Рис. 4.1. Граф Петерсена $G_{11}$.

Вычислим кортеж весов вершин графа Петерсена.

Кортеж весов вершин: <48,48,48,48,48,48,48,48,48,48>. Таким образом, все вершины имеют равный вес.

Построим множество изометрических циклов графа:

$c_1 = \{e_1,e_2,e_4,e_6,e_8\} \rightarrow \langle v_1,v_2,v_3,v_4,v_5\rangle$;
$c_2 = \{e_1,e_2,e_5,e_{10},e_{14}\} \rightarrow \langle v_2,v_1,v_5,v_{10},v_7\rangle$;
$c_3 = \{e_1,e_3,e_4,e_7,e_{11}\} \rightarrow \langle v_1,v_2,v_3,v_8,v_6\rangle$;
$c_4 = \{e_1,e_3,e_5,e_{12},e_{13}\} \rightarrow = \langle v_2,v_1,v_6,v_9,v_7\rangle$;
$c_5 = \{e_2,e_3,e_{10},e_{11},e_{15}\} \rightarrow \langle v_1,v_5,v_{10},v_8,v_6\rangle$;
$c_6 = \{e_2,e_3,e_8,e_9,e_{12}\} \rightarrow \langle v_5,v_1,v_6,v_9,v_4\rangle$;
$c_7 = \{e_4,e_5,e_6,e_9,e_{13}\} \rightarrow \langle v_2,v_3,v_4,v_9,v_7\rangle$;
$c_8 = \{e_4,e_5,e_7,e_{14},e_{15}\} \rightarrow \langle v_2,v_7,v_{10},v_8,v_3\rangle$;
$c_9 = \{e_6,e_7,e_9,e_{11},e_{12}\} \rightarrow \langle v_3,v_4,v_9,v_6,v_8\rangle$;
$c_{10} = \{e_6,e_7,e_8,e_{10},e_{15}\} \rightarrow \langle v_4,v_3,v_8,v_{10},v_5\rangle$;
$c_{11} = \{e_8,e_9,e_{10},e_{13},e_{14}\} \rightarrow \langle v_5,v_4,v_9,v_7,v_{10}\rangle$;
$c_{12} = \{e_{11},e_{12},e_{13},e_{14},e_{15}\} \rightarrow \langle v_{10},v_7,v_9,v_6,v_8\rangle$.

Будем рассматривать каждый изометрический цикл как образующий. Тогда можно поставить в соответствие опорному циклу симметрический рисунок графа.

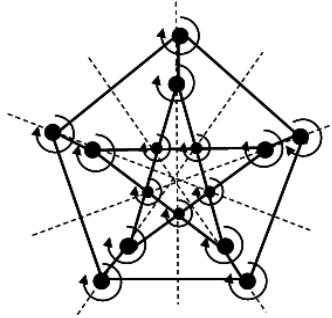

Рис. 4.2 Симметричный рисунок графа Петерсена.

В отличие от симметричных рисунков планарных графов, симметрические рисунки непланарных графов отличаются наличием мнимых вершин характеризующих пересечение рёбер. Но мы их не будем показывать на рисунке

Выберем в качестве образующего цикла – цикл $c_1$ и построим перестановки.

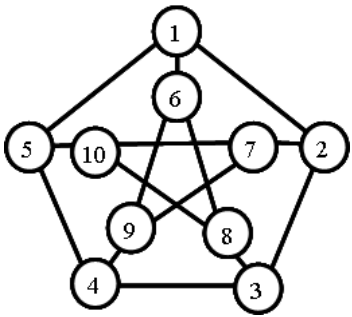 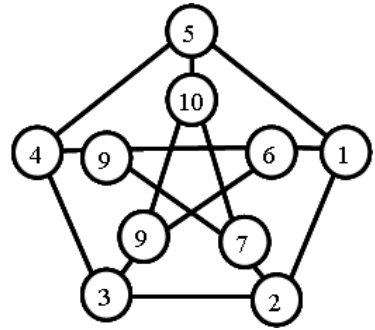

Рис. 4.3. Перестановка $p_0$.          Рис. 4.4. Перестановка $p_1$.

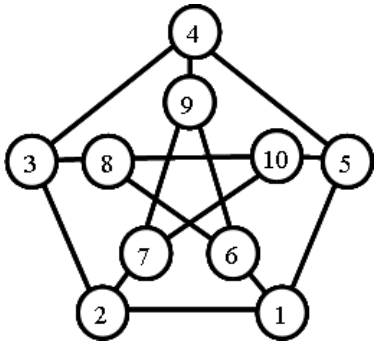 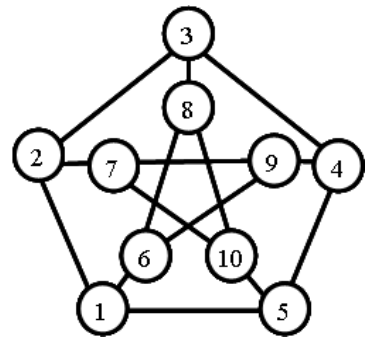

Рис. 4.5. Перестановка $p_2$.          Рис. 4.6. Перестановка $p_3$.

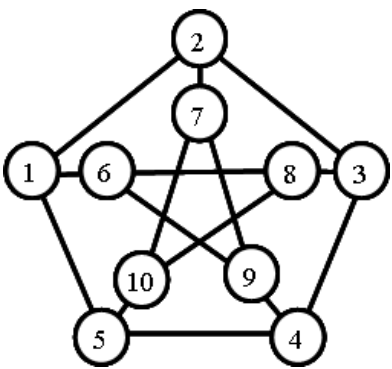 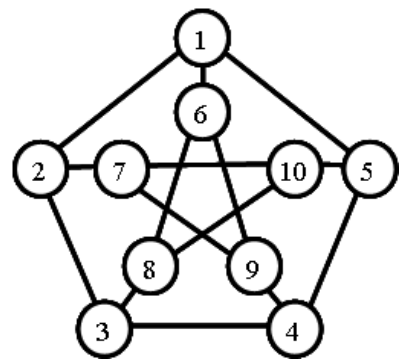

Рис. 4.7. Перестановка $p_4$.          Рис. 4.8. Перестановка $p_5$.

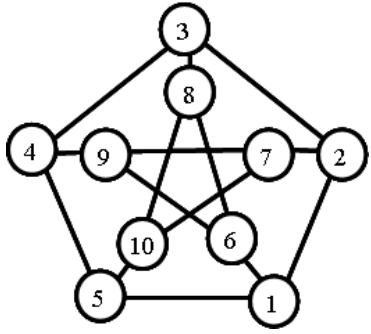 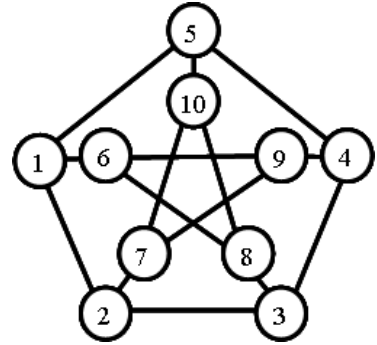

Рис. 4.9. Перестановка $p_6$.          Рис. 4.10. Перестановка $p_7$.

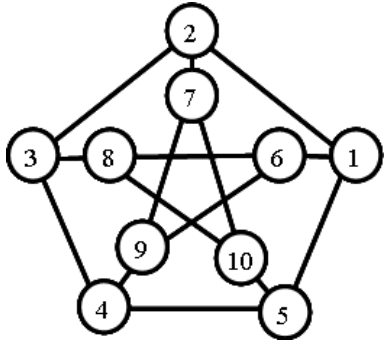 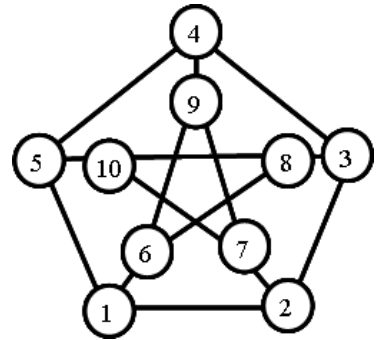

Рис. 4.11. Перестановка $p_8$.          Рис. 4.12. Перестановка $p_9$.

Количество перестановок для образующего цикла равно десяти, количество образующих циклов равно 12. Поэтому множество всех подстановок в графе Петерсена определяется как $12 \times 10 = 120$. Рассмотрим только первые десять перестановок индуцированных циклом $c_1$:

$p_0 = \langle 1,2,3,4,5,6,7,8,9,10 \rangle = (1)(2)(3)(4)(5)(6)(7)(8)(9)(10);$
$p_1 = \langle 5,1,2,3,4,10,6,7,8,9 \rangle = (1\ 5\ 4\ 3\ 2)(6\ 10\ 9\ 8\ 7);$
$p_2 = \langle 4,5,1,2,3,9,10,6,7,8 \rangle = (1\ 4\ 2\ 5\ 3)(6\ 9\ 7\ 10\ 8);$
$p_3 = \langle 3,4,5,1,2,8,9,10,6,7 \rangle = (1\ 3\ 5\ 2\ 4)(6\ 8\ 10\ 7\ 9);$
$p_4 = \langle 2,3,4,5,1,7,8,9,10,6 \rangle = (1\ 2\ 3\ 4\ 5)(6\ 7\ 8\ 9\ 10);$
$p_5 = \langle 1,5,4,3,2,6,10,9,8,7 \rangle = (1)(2\ 5)(3\ 4)(6)(7\ 10)(8\ 9);$
$p_6 = \langle 3,2,1,5,4,8,7,6,10,9 \rangle = (1\ 3)(2)(4\ 5)(6\ 8)(7)(9\ 10);$
$p_7 = \langle 5,4,3,2,1,10,9,8,7,6 \rangle = (1\ 5)(2\ 4)(3)(6\ 10)(7\ 9)(8);$
$p_8 = \langle 2,1,5,4,3,7,6,10,9,8 \rangle = (1\ 2)(3\ 5)(4)(6\ 7)(8\ 10)(9);$
$p_9 = \langle 4,3,2,1,5,9,8,7,6,10 \rangle = (1\ 4)(2\ 3)(5)(6\ 9)(7\ 8)(10).$

Таблица Кэли для графа Петерсена:

$Aut(G_{11}) =$

| $P_{1\,1}$ | $P_{1\,2}$ | $P_{1\,3}$ | $P_{1\,4}$ | … | $P_{1\,8}$ | $P_{1\,9}$ | $P_{1\,10}$ | $P_{1\,11}$ | $P_{1\,12}$ |
|---|---|---|---|---|---|---|---|---|---|
| $P_{2\,1}$ | $P_{2\,2}$ | $P_{2\,3}$ | $P_{2\,4}$ | … | $P_{2\,8}$ | $P_{2\,9}$ | $P_{2\,10}$ | $P_{2\,11}$ | $P_{2\,12}$ |
| $P_{3\,1}$ | $P_{3\,2}$ | $P_{3\,3}$ | $P_{3\,4}$ | … | $P_{3\,8}$ | $P_{3\,9}$ | $P_{3\,10}$ | $P_{3\,11}$ | $P_{3\,12}$ |
| $P_{4\,1}$ | $P_{4\,2}$ | $P_{4\,3}$ | $P_{4\,4}$ | … | $P_{4\,8}$ | $P_{4\,9}$ | $P_{4\,10}$ | $P_{4\,11}$ | $P_{4\,12}$ |
| $P_{5\,1}$ | $P_{5\,2}$ | $P_{5\,3}$ | $P_{5\,4}$ | … | $P_{5\,8}$ | $P_{5\,9}$ | $P_{5\,10}$ | $P_{5\,11}$ | $P_{5\,12}$ |
| … | … | … | … | … | … | … | … | … | … |
| $P_{8\,1}$ | $P_{8\,2}$ | $P_{8\,3}$ | $P_{8\,4}$ | … | $P_{8\,8}$ | $P_{8\,9}$ | $P_{8\,10}$ | $P_{8\,11}$ | $P_{8\,12}$ |
| $P_{9\,1}$ | $P_{9\,2}$ | $P_{9\,3}$ | $P_{9\,4}$ | … | $P_{9\,8}$ | $P_{9\,9}$ | $P_{9\,10}$ | $P_{9\,11}$ | $P_{9\,12}$ |
| $P_{10\,1}$ | $P_{10\,2}$ | $P_{10\,3}$ | $P_{10\,4}$ | … | $P_{10\,8}$ | $P_{10\,9}$ | $P_{10\,10}$ | $P_{10\,11}$ | $P_{10\,12}$ |
| $P_{11\,1}$ | $P_{11\,2}$ | $P_{11\,3}$ | $P_{11\,4}$ | … | $P_{11\,8}$ | $P_{11\,9}$ | $P_{11\,10}$ | $P_{11\,11}$ | $P_{11\,12}$ |
| $P_{12\,1}$ | $P_{12\,2}$ | $P_{12\,3}$ | $P_{12\,4}$ | … | $P_{12\,8}$ | $P_{12\,9}$ | $P_{12\,10}$ | $P_{12\,11}$ | $P_{12\,12}$ |

Часть таблицы умножения перестановок (таблица Кэли) составленных только для образующего цикла $c_1$ (1-ая подгруппа) имеет вид:

$P_{1\,1}$

|  | $p_0$ | $p_1$ | $p_2$ | $p_3$ | $p_4$ | $p_5$ | $p_6$ | $p_7$ | $p_8$ | $p_9$ |
|---|---|---|---|---|---|---|---|---|---|---|
| $p_0$ | $p_0$ | $p_1$ | $p_2$ | $p_3$ | $p_4$ | $p_5$ | $p_6$ | $p_7$ | $p_8$ | $p_9$ |
| $p_1$ | $p_1$ | $p_2$ | $p_3$ | $p_4$ | $p_0$ | $p_8$ | $p_9$ | $p_5$ | $p_6$ | $p_7$ |
| $p_2$ | $p_2$ | $p_3$ | $p_4$ | $p_0$ | $p_1$ | $p_6$ | $p_7$ | $p_8$ | $p_9$ | $p_5$ |
| $p_3$ | $p_3$ | $p_4$ | $p_0$ | $p_1$ | $p_2$ | $p_9$ | $p_5$ | $p_6$ | $p_7$ | $p_8$ |
| $p_4$ | $p_4$ | $p_0$ | $p_1$ | $p_2$ | $p_3$ | $p_7$ | $p_8$ | $p_9$ | $p_5$ | $p_6$ |
| $p_5$ | $p_5$ | $p_7$ | $p_9$ | $p_6$ | $p_8$ | $p_0$ | $p_3$ | $p_1$ | $p_4$ | $p_2$ |
| $p_6$ | $p_6$ | $p_8$ | $p_5$ | $p_7$ | $p_9$ | $p_2$ | $p_0$ | $p_3$ | $p_1$ | $p_4$ |
| $p_7$ | $p_7$ | $p_9$ | $p_6$ | $p_8$ | $p_5$ | $p_4$ | $p_2$ | $p_0$ | $p_3$ | $p_1$ |
| $p_8$ | $p_8$ | $p_5$ | $p_7$ | $p_9$ | $p_6$ | $p_1$ | $p_4$ | $p_2$ | $p_0$ | $p_3$ |
| $p_9$ | $p_9$ | $p_6$ | $p_8$ | $p_5$ | $p_7$ | $p_3$ | $p_1$ | $p_4$ | $p_2$ | $p_0$ |

### 4.2. Сильно регулярный граф Шрикханде

*Пример 4.2*. Будем рассматривать сильно регулярные графы [28,29], в частности граф Шрикханде.

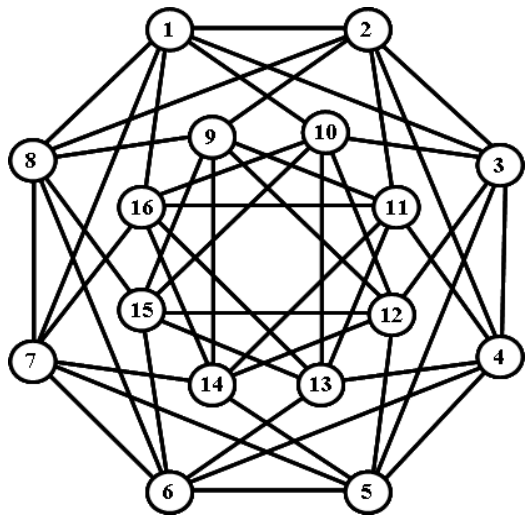

Рис. 4.13. Граф Шрикханде $G_{12}$.

Количество вершин графа = 16.
Количество ребер графа = 48.
Количество изометрических циклов = 44.
Смежность графа:
вершина 1: {$v_2,v_3,v_7,v_8,v_{10},v_{16}$};
вершина 2: {$v_1,v_3,v_4,v_8,v_9,v_{11}$};
вершина 3: {$v_1,v_2,v_4,v_5,v_{10},v_{12}$};
вершина 4: {$v_2,v_3,v_5,v_6,v_{11},v_{13}$};
вершина 5: {$v_3,v_4,v_6,v_7,v_{12},v_{14}$};
вершина 6: {$v_4,v_5,v_7,v_8,v_{13},v_{15}$};
вершина 7: {$v_1,v_5,v_6,v_8,v_{14},v_{16}$};
вершина 8: {$v_1,v_2,v_6,v_7,v_9,v_{15}$};
вершина 9: {$v_2,v_8,v_{11},v_{12},v_{14},v_{15}$};
вершина 10: {$v_1,v_3,v_{12},v_{13},v_{15},v_{16}$};

вершина 11: {$v_2,v_4,v_9,v_{13},v_{14},v_{16}$};
вершина 12: {$v_3,v_5,v_9,v_{10},v_{14},v_{15}$};
вершина 13: {$v_4,v_6,v_{10},v_{11},v_{15},v_{16}$};
вершина 14: {$v_5,v_7,v_9,v_{11},v_{12},v_{16}$};
вершина 15: {$v_6,v_8,v_9,v_{10},v_{12},v_{13}$};
вершина 16: {$v_1,v_7,v_{10},v_{11},v_{13},v_{14}$}.

Инцидентность графа:

$e_1 = (v_1,v_2) \vee (v_2,v_1)$; $\qquad e_2 = (v_1,v_3) \vee (v_3,v_1)$;
$e_3 = (v_1,v_7) \vee (v_7,v_1)$; $\qquad e_4 = (v_1,v_8) \vee (v_8,v_1)$;
$e_5 = (v_1,v_{10}) \vee (v_{10},v_1)$; $\qquad e_6 = (v_1,v_{16}) \vee (v_{16},v_1)$;
$e_7 = (v_2,v_3) \vee (v_3,v_2)$; $\qquad e_8 = (v_2,v_4) \vee (v_4,v_2)$;
$e_9 = (v_2,v_8) \vee (v_8,v_2)$; $\qquad e_{10} = (v_2,v_9) \vee (v_9,v_2)$;
$e_{11} = (v_2,v_{11}) \vee (v_{11},v_2)$; $\qquad e_{12} = (v_3,v_4) \vee (v_4,v_3)$;
$e_{13} = (v_3,v_5) \vee (v_5,v_3)$; $\qquad e_{14} = (v_3,v_{10}) \vee (v_{10},v_3)$;
$e_{15} = (v_3,v_{12}) \vee (v_{12},v_3)$; $\qquad e_{16} = (v_4,v_5) \vee (v_5,v_4)$;
$e_{17} = (v_4,v_6) \vee (v_6,v_4)$; $\qquad e_{18} = (v_4,v_{11}) \vee (v_{11},v_4)$;
$e_{19} = (v_4,v_{13}) \vee (v_{13},v_4)$; $\qquad e_{20} = (v_5,v_6) \vee (v_6,v_5)$;
$e_{21} = (v_5,v_7) \vee (v_7,v_5)$; $\qquad e_{22} = (v_5,v_{12}) \vee (v_{12},v_5)$;
$e_{23} = (v_5,v_{14}) \vee (v_{14},v_5)$; $\qquad e_{24} = (v_6,v_7) \vee (v_7,v_6)$;
$e_{25} = (v_6,v_8) \vee (v_8,v_6)$; $\qquad e_{26} = (v_6,v_{13}) \vee (v_{13},v_6)$;
$e_{27} = (v_6,v_{15}) \vee (v_{15},v_6)$; $\qquad e_{28} = (v_7,v_8) \vee (v_8,v_7)$;
$e_{29} = (v_7,v_{14}) \vee (v_{14},v_7)$; $\qquad e_{30} = (v_7,v_{16}) \vee (v_{16},v_7)$;
$e_{31} = (v_8,v_9) \vee (v_9,v_8)$; $\qquad e_{32} = (v_8,v_{15}) \vee (v_{15},v_8)$;
$e_{33} = (v_9,v_{11}) \vee (v_{11},v_9)$; $\qquad e_{34} = (v_9,v_{12}) \vee (v_{12},v_9)$;
$e_{35} = (v_9,v_{14}) \vee (v_{14},v_9)$; $\qquad e_{36} = (v_9,v_{15}) \vee (v_{15},v_9)$;
$e_{37} = (v_{10},v_{12}) \vee (v_{12},v_{10})$; $\qquad e_{38} = (v_{10},v_{13}) \vee (v_{13},v_{10})$;
$e_{39} = (v_{10},v_{15}) \vee (v_{15},v_{10})$; $\qquad e_{40} = (v_{10},v_{16}) \vee (v_{16},v_{10})$;
$e_{41} = (v_{11},v_{13}) \vee (v_{13},v_{11})$; $\qquad e_{42} = (v_{11},v_{14}) \vee (v_{14},v_{11})$;
$e_{43} = (v_{11},v_{16}) \vee (v_{16},v_{11})$; $\qquad e_{44} = (v_{12},v_{14}) \vee (v_{14},v_{12})$;
$e_{45} = (v_{12},v_{15}) \vee (v_{15},v_{12})$; $\qquad e_{46} = (v_{13},v_{15}) \vee (v_{15},v_{13})$;
$e_{47} = (v_{13},v_{16}) \vee (v_{16},v_{13})$; $\qquad e_{48} = (v_{14},v_{16}) \vee (v_{16},v_{14})$.

Вектор весов ребер : $F(\xi_w((G_2))) = (48 \times 34)$;

Вектор весов вершин : $F(\zeta_w((G_2))) = (16 \times 204)$.

Множество изометрических циклов графа:

$c_1 = \{e_1,e_2,e_7\} \rightarrow \{v_1,v_2,v_3\}$;
$c_2 = \{e_1,e_4,e_9\} \rightarrow \{v_1,v_2,v_8\}$;
$c_3 = \{e_1,e_6,e_{11},e_{43}\} \rightarrow \{v_1,v_2,v_{16},v_{11}\}$;
$c_4 = \{e_2,e_3,e_{13},e_{21}\} \rightarrow \{v_1,v_3,v_7,v_5\}$;
$c_5 = \{e_2,e_5,e_{14}\} \rightarrow \{v_1,v_3,v_{10}\}$;
$c_6 = \{e_3,e_4,e_{28}\} \rightarrow \{v_1,v_7,v_8\}$;
$c_7 = \{e_3,e_6,e_{30}\} \rightarrow \{v_1,v_7,v_{16}\}$;
$c_8 = \{e_4,e_5,e_{32},e_{39}\} \rightarrow \{v_1,v_8,v_{10},v_{15}\}$;
$c_9 = \{e_5,e_6,e_{40}\} \rightarrow \{v_1,v_{10},v_{16}\}$;
$c_{10} = \{e_7,e_8,e_{12}\} \rightarrow \{v_2,v_3,v_4\}$;
$c_{11} = \{e_7,e_{10},e_{15},e_{34}\} \rightarrow \{v_2,v_3,v_9,v_{12}\}$;
$c_{12} = \{e_8,e_9,e_{17},e_{25}\} \rightarrow \{v_2,v_4,v_8,v_6\}$;
$c_{13} = \{e_8,e_{11},e_{18}\} \rightarrow \{v_2,v_4,v_{11}\}$;
$c_{14} = \{e_9,e_{10},e_{31}\} \rightarrow \{v_2,v_8,v_9\}$;

$c_{15} = \{e_{10}, e_{11}, e_{33}\} \rightarrow \{v_2, v_9, v_{11}\}$;
$c_{16} = \{e_{12}, e_{13}, e_{16}\} \rightarrow \{v_3, v_4, v_5\}$;
$c_{17} = \{e_{12}, e_{14}, e_{19}, e_{38}\} \rightarrow \{v_3, v_4, v_{10}, v_{13}\}$;
$c_{18} = \{e_{13}, e_{15}, e_{22}\} \rightarrow \{v_3, v_5, v_{12}\}$;
$c_{19} = \{e_{14}, e_{15}, e_{37}\} \rightarrow \{v_3, v_{10}, v_{12}\}$;
$c_{20} = \{e_{16}, e_{17}, e_{20}\} \rightarrow \{v_4, v_5, v_6\}$;
$c_{21} = \{e_{16}, e_{18}, e_{23}, e_{42}\} \rightarrow \{v_4, v_5, v_{11}, v_{14}\}$;
$c_{22} = \{e_{17}, e_{19}, e_{26}\} \rightarrow \{v_4, v_6, v_{13}\}$;
$c_{23} = \{e_{18}, e_{19}, e_{41}\} \rightarrow \{v_4, v_{11}, v_{13}\}$;
$c_{24} = \{e_{20}, e_{21}, e_{24}\} \rightarrow \{v_5, v_6, v_7\}$;
$c_{25} = \{e_{20}, e_{22}, e_{27}, e_{45}\} \rightarrow \{v_5, v_6, v_{12}, v_{15}\}$;
$c_{26} = \{e_{21}, e_{23}, e_{29}\} \rightarrow \{v_5, v_7, v_{14}\}$;
$c_{27} = \{e_{22}, e_{23}, e_{44}\} \rightarrow \{v_5, v_{12}, v_{14}\}$;
$c_{28} = \{e_{24}, e_{25}, e_{28}\} \rightarrow \{v_6, v_7, v_8\}$;
$c_{29} = \{e_{24}, e_{26}, e_{30}, e_{47}\} \rightarrow \{v_6, v_7, v_{13}, v_{16}\}$;
$c_{30} = \{e_{25}, e_{27}, e_{32}\} \rightarrow \{v_6, v_8, v_{15}\}$;
$c_{31} = \{e_{26}, e_{27}, e_{46}\} \rightarrow \{v_6, v_{13}, v_{15}\}$;
$c_{32} = \{e_{28}, e_{29}, e_{31}, e_{35}\} \rightarrow \{v_7, v_8, v_{14}, v_9\}$;
$c_{33} = \{e_{29}, e_{30}, e_{48}\} \rightarrow \{v_7, v_{14}, v_{16}\}$;
$c_{34} = \{e_{31}, e_{32}, e_{36}\} \rightarrow \{v_8, v_9, v_{15}\}$;
$c_{35} = \{e_{33}, e_{35}, e_{42}\} \rightarrow \{v_9, v_{11}, v_{14}\}$;
$c_{36} = \{e_{33}, e_{36}, e_{41}, e_{46}\} \rightarrow \{v_9, v_{11}, v_{15}, v_{13}\}$;
$c_{37} = \{e_{34}, e_{35}, e_{44}\} \rightarrow \{v_9, v_{12}, v_{14}\}$;
$c_{38} = \{e_{34}, e_{36}, e_{45}\} \rightarrow \{v_9, v_{12}, v_{15}\}$;
$c_{39} = \{e_{37}, e_{39}, e_{45}\} \rightarrow \{v_{10}, v_{12}, v_{15}\}$;
$c_{40} = \{e_{37}, e_{40}, e_{44}, e_{48}\} \rightarrow \{v_{10}, v_{12}, v_{16}, v_{14}\}$;
$c_{41} = \{e_{38}, e_{39}, e_{46}\} \rightarrow \{v_{10}, v_{13}, v_{15}\}$;
$c_{42} = \{e_{38}, e_{40}, e_{47}\} \rightarrow \{v_{10}, v_{13}, v_{16}\}$;
$c_{43} = \{e_{41}, e_{43}, e_{47}\} \rightarrow \{v_{11}, v_{13}, v_{16}\}$;
$c_{44} = \{e_{42}, e_{43}, e_{48}\} \rightarrow \{v_{11}, v_{14}, v_{16}\}$.

Осуществим следующее построение.

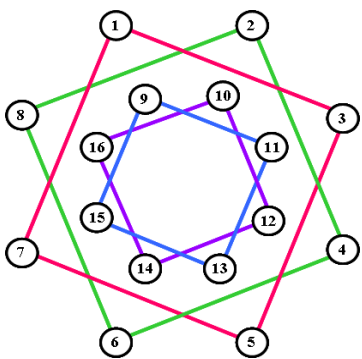
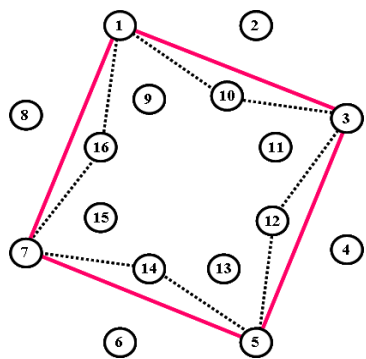
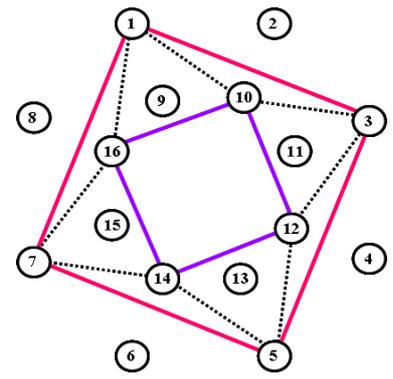

Рис. 4.14. Построение графа.    Рис. 4.15. Первый шаг.    Рис. 4.16. Второй шаг.

Выделим изометрические циклы длиной четыре:

$c_3 = \{e_1, e_6, e_{11}, e_{43}\} \rightarrow \{v_1, v_2, v_{16}, v_{11}\}$;
$c_4 = \{e_2, e_3, e_{13}, e_{21}\} \rightarrow \{v_1, v_3, v_7, v_5\}$;
$c_8 = \{e_4, e_5, e_{32}, e_{39}\} \rightarrow \{v_1, v_8, v_{10}, v_{15}\}$;
$c_{11} = \{e_7, e_{10}, e_{15}, e_{34}\} \rightarrow \{v_2, v_3, v_9, v_{12}\}$;

$c_{12} = \{e_8, e_9, e_{17}, e_{25}\} \rightarrow \{v_2, v_4, v_8, v_6\}$;
$c_{17} = \{e_{12}, e_{14}, e_{19}, e_{38}\} \rightarrow \{v_3, v_4, v_{10}, v_{13}\}$;
$c_{21} = \{e_{16}, e_{18}, e_{23}, e_{42}\} \rightarrow \{v_4, v_5, v_{11}, v_{14}\}$;
$c_{25} = \{e_{20}, e_{22}, e_{27}, e_{45}\} \rightarrow \{v_5, v_6, v_{12}, v_{15}\}$;
$c_{29} = \{e_{24}, e_{26}, e_{30}, e_{47}\} \rightarrow \{v_6, v_7, v_{13}, v_{16}\}$;
$c_{32} = \{e_{28}, e_{29}, e_{31}, e_{35}\} \rightarrow \{v_7, v_8, v_{14}, v_9\}$;
$c_{36} = \{e_{33}, e_{36}, e_{41}, e_{46}\} \rightarrow \{v_9, v_{11}, v_{15}, v_{13}\}$;
$c_{40} = \{e_{37}, e_{40}, e_{44}, e_{48}\} \rightarrow \{v_{10}, v_{12}, v_{16}, v_{14}\}$.

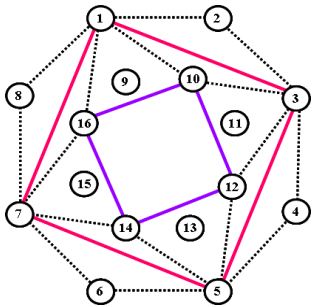 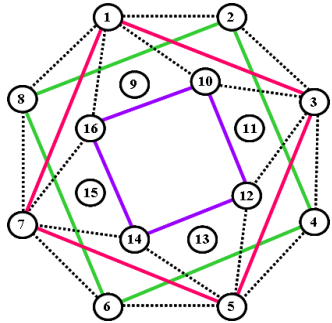 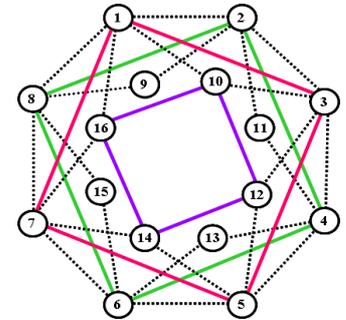

Рис. 4.17. Третий шаг.  Рис. 4.18. Четвертый шаг.  Рис. 4.19. Пятый шаг.

Количество вершин в графе равно 16. Попробуем выбрать подмножество циклов так, чтобы покрыть все множество вершин. Очевидно, что в этом случае нужно выбрать четыре непересекающихся цикла. Проверив все сочетания из 12 элементов по 4, и выберем циклы, объединение которых представляет все множество вершин графа V:

1-е подмножество:

$c_3 = \{e_1, e_6, e_{11}, e_{43}\} \rightarrow \{v_1, v_2, v_{16}, v_{11}\}$;
$c_{17} = \{e_{12}, e_{14}, e_{19}, e_{38}\} \rightarrow \{v_3, v_4, v_{10}, v_{13}\}$;
$c_{25} = \{e_{20}, e_{22}, e_{27}, e_{45}\} \rightarrow \{v_5, v_6, v_{12}, v_{15}\}$;
$c_{32} = \{e_{28}, e_{29}, e_{31}, e_{35}\} \rightarrow \{v_7, v_8, v_{14}, v_9\}$;

2-е подмножество:

$c_4 = \{e_2, e_3, e_{13}, e_{21}\} \rightarrow \{v_1, v_3, v_7, v_5\}$;
$c_{12} = \{e_8, e_9, e_{17}, e_{25}\} \rightarrow \{v_2, v_4, v_8, v_6\}$;
$c_{36} = \{e_{33}, e_{36}, e_{41}, e_{46}\} \rightarrow \{v_9, v_{11}, v_{15}, v_{13}\}$;
$c_{40} = \{e_{37}, e_{40}, e_{44}, e_{48}\} \rightarrow \{v_{10}, v_{12}, v_{16}, v_{14}\}$;

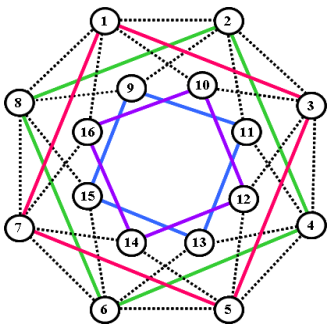 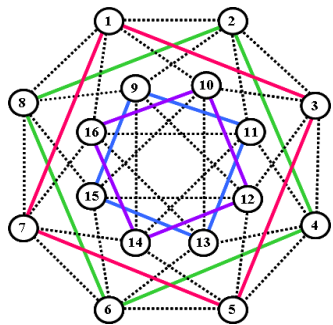 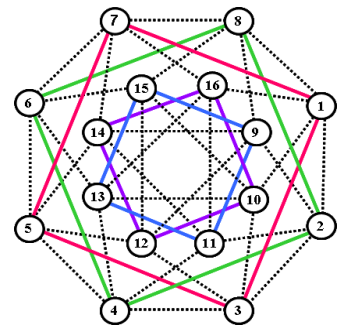

Рис. 4.20. Шестой шаг.  Рис. 4.21. Перестановка $p_0$.  Рис. 4.22. Перестановка $p_1$.

3-е подмножество:

$c_8 = \{e_4, e_5, e_{32}, e_{39}\} \rightarrow \{v_1, v_8, v_{10}, v_{15}\}$;
$c_{11} = \{e_7, e_{10}, e_{15}, e_{34}\} \rightarrow \{v_2, v_3, v_9, v_{12}\}$;
$c_{21} = \{e_{16}, e_{18}, e_{23}, e_{42}\} \rightarrow \{v_4, v_5, v_{11}, v_{14}\}$;
$c_{29} = \{e_{24}, e_{26}, e_{30}, e_{47}\} \rightarrow \{v_6, v_7, v_{13}, v_{16}\}$;

Топологический рисунок графа $G_2$ строится путем расположения на плоскости четырех непересекающихся по вершинам изометрических циклов длиной четыре. На рис. 4.14 представлено расположение изометрических циклов для 2-го подмножества циклов.

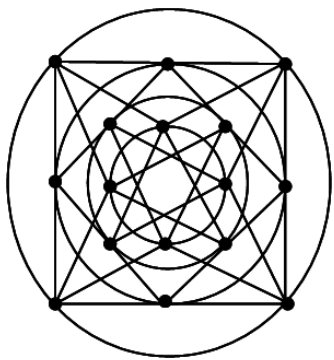

Рис. 4.23. Симметричный топологический рисунок.

На рис. 4.15 представлен результат суммирования цикла $\{v_1, v_3, v_5, v_7\}$ и изометрических циклов длиной три. Рис. 4.16 представляет построение цикла $\{v_{10}, v_{12}, v_{14}, v_{16}\}$. На третьем шаге построения к циклу $\{v_1, v_3, v_5, v_7\}$ присоединяются четыре изометрических цикла длиной три. На рисунках 4.17-4.21 представлены стадии проведения ребер внутри опорного цикла.

Будем рассматривать 2-е подмножество циклов. Распределим циклы по окружностям.

1-ая окружность содержит вершины цикла $\{v_1, v_3, v_5, v_7\}$.
2-ая окружность содержит вершины цикла $\{v_2, v_4, v_6, v_8\}$.
3-ья окружность содержит вершины цикла $\{v_{10}, v_{12}, v_{14}, v_{16}\}$.
4-ая окружность содержит вершины цикла $\{v_9, v_{11}, v_{13}, v_{15}\}$.

Установим, что 4-ая окружность вписана в 3-ью окружность, 3-ья окружность вписана во 2-ую окружность, а 2-ая окружность вписана в 1-ую окружность. Первая окружность характеризует образующий цикл с вершинами $\{v_1, v_3, v_5, v_7\}$. Тогда порядок вписывания окружностей запишем (1 2 3 4).

Для второго множества циклов можно произвести только восемь вариантов симметричного вписывания окружностей с соблюдением смежности (см. рис. 4.24).

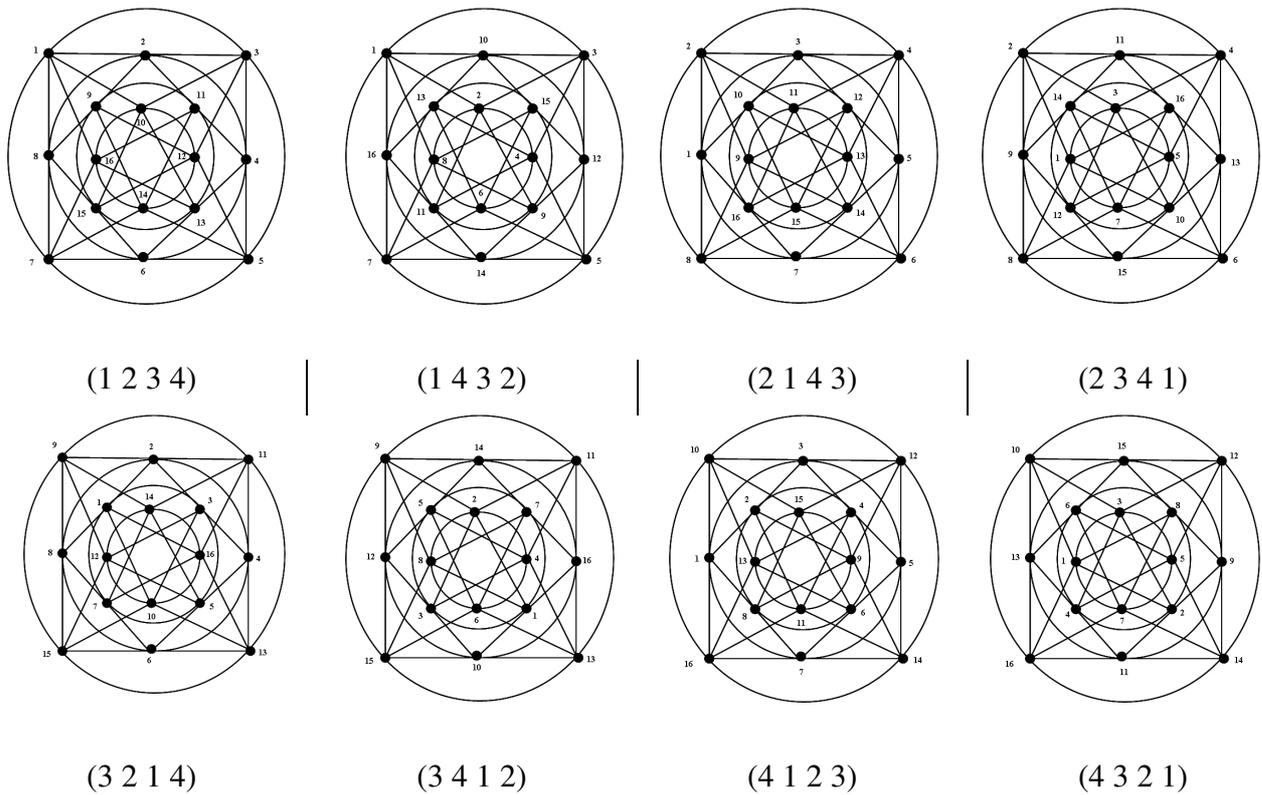

Рис. 4.24. Порядок вписанных окружностей.

Построим следующую последовательность вершин для образующего цикла $\{v_1, v_3, v_5, v_7\}$.

1 3 5 7 1 3 5 7 1 3 5 7 1 3 5 7 → $\{v_1, v_3, v_5, v_7\}$

  2 4 6 8 2 4 6 8 2 4 6 8 2 4 6 → $\{v_2, v_4, v_6, v_8\}$

9 11 13 15 9 11 13 15 9 11 13 15 9 11 13 15 → $\{v_9, v_{11}, v_{13}, v_{15}\}$

 14 16 10 12 14 16 10 12 14 16 10 12 14 16 10 → $\{v_{10}, v_{12}, v_{14}, v_{16}\}$

5 7 1 3 5 7 1 3 5 7 1 3 5 7 1 3 → $\{v_1, v_3, v_5, v_7\}$

Здесь всегда можно определить какие вершины окружают выбранную вершину. Например, 6-ую вершину окружают следующие вершины

   5  7
 4  **6**  8
  13 15

Выберем в качестве образующего цикла $\{v_1, v_3, v_5, v_7\}$ и будем переставлять вершины только в изометрических циклах 2-го подмножества вершин. На рис. 4.22 представлена перестановка $p_1$ графа $G_2$.

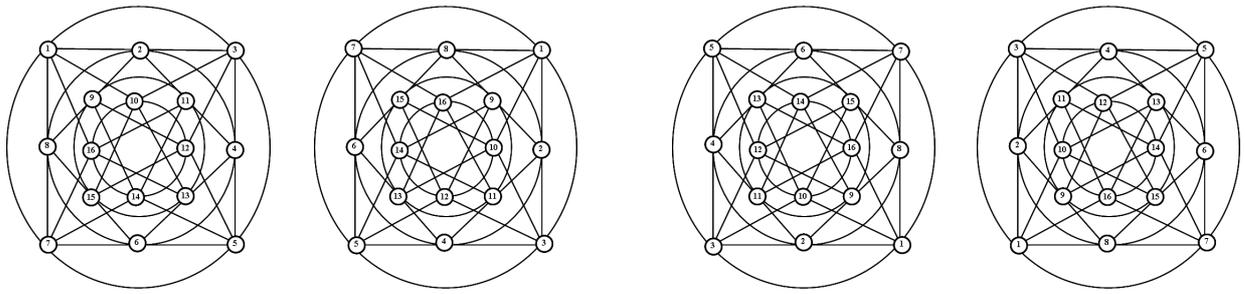

| Перестановка $p_0$. | Перестановка $p_1$. | Перестановка $p_2$. | Перестановка $p_3$. |

Рис. 4.25. Перестановки поворотом.

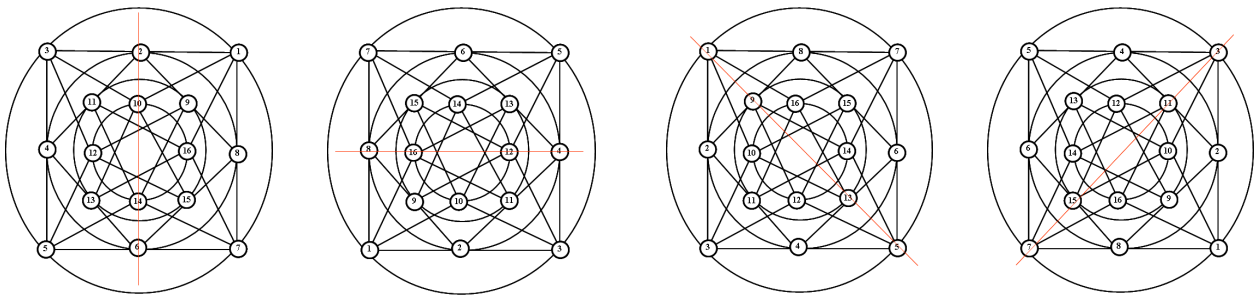

| Перестановка $p_4$. | Перестановка $p_5$. | Перестановка $p_6$. | Перестановка $p_7$. |

Рис. 4.26. Перестановки отображения вокруг осей.

В результате получим следующие перестановки:

$p_0$ = <1 2 3 4 5 6 7 8 9 10 11 12 13 14 15 16> = (1)(2)(3)(4)(5)(6)(7)(8)(9)(10)(11)(12)(13)(14)(15)(16);

$p_1$ = <7 8 1 2 3 4 5 6 15 16 9 10 11 12 13 14> = (1 7 5 3)(2 8 6 4)(9 15 13 11)(10 16 14 12);

$p_2$ = <5 6 7 8 1 2 3 4 13 14 15 16 9 10 11 12> = (1 5)(2 6)(3 7)(4 8)(9 13)(10 14)(11 15)(12 16);

$p_3$ = <3 4 5 6 7 8 1 2 11 12 13 14 15 16 9 10> = (1 3 5 7)(2 4 6 8)(9 11 13 15)(10 12 14 16);

$p_4$ = <3 2 1 8 7 6 5 4 11 10 9 16 15 14 13 12> = (1 3)(2)(4 8)(5 7)(6)(9 11)(10)(12 16)(13 15)(14);

$p_5$ = <7 6 5 4 3 2 1 8 15 14 13 12 11 10 9 16> = (1 7)(2 6)(3 5)(4)(8)(9 15)(10 14)(11 13)(12)(16);

$p_6$ = <1 8 7 6 5 4 3 2 9 16 15 14 13 12 11 10> = (1)(2 8)(3 7)(4 6)(5)(9)(10 16)(11 15)(12 14)(13);

$p_7$ = <5 4 3 2 1 8 7 6 13 12 11 10 9 16 15 14> = (1 5)(2 4)(3)(6 8)(7)(9 13)(10 12)(11)(14 16)(15).

Если в качестве образующих циклов выбрать другие изометрические циклы из 2-го подмножества, то для каждого множества циклов получим 8 вариантов размещения окружностей. И для каждого варианта 8 перестановок. Так как мы имеем 3 подмножества непересекающихся циклов, то граф $G_2$ может порождать 3×8×8=24×8 = 192 перестановок.

Таблица Кэли для графа $G_2$:

$Aut(G_2) =$

| $P_{1\,1}$ | $P_{1\,2}$ | $P_{1\,3}$ | $P_{1\,4}$ | … | $P_{1\,20}$ | $P_{1\,21}$ | $P_{1\,22}$ | $P_{1\,23}$ | $P_{1\,24}$ |
|---|---|---|---|---|---|---|---|---|---|
| $P_{2\,1}$ | $P_{2\,2}$ | $P_{2\,3}$ | $P_{2\,4}$ | … | $P_{2\,20}$ | $P_{2\,21}$ | $P_{2\,22}$ | $P_{2\,23}$ | $P_{2\,24}$ |
| $P_{3\,1}$ | $P_{3\,2}$ | $P_{3\,3}$ | $P_{3\,4}$ | … | $P_{3\,20}$ | $P_{3\,21}$ | $P_{3\,22}$ | $P_{3\,23}$ | $P_{3\,24}$ |
| $P_{4\,1}$ | $P_{4\,2}$ | $P_{4\,3}$ | $P_{4\,4}$ | … | $P_{4\,20}$ | $P_{4\,21}$ | $P_{4\,22}$ | $P_{4\,23}$ | $P_{4\,24}$ |
| $P_{5\,1}$ | $P_{5\,2}$ | $P_{5\,3}$ | $P_{5\,4}$ | … | $P_{5\,20}$ | $P_{5\,21}$ | $P_{5\,22}$ | $P_{5\,23}$ | $P_{5\,24}$ |
| $P_{6\,1}$ | $P_{6\,2}$ | $P_{6\,3}$ | $P_{6\,4}$ | … | $P_{6\,20}$ | $P_{6\,21}$ | $P_{6\,22}$ | $P_{6\,23}$ | $P_{6\,24}$ |
| … | … | … | … | … | … | … | … | … | … |
| $P_{20\,1}$ | $P_{20\,2}$ | $P_{20\,3}$ | $P_{20\,4}$ | … | $P_{20\,20}$ | $P_{20\,21}$ | $P_{20\,22}$ | $P_{20\,23}$ | $P_{20\,24}$ |
| $P_{21\,1}$ | $P_{21\,2}$ | $P_{21\,3}$ | $P_{21\,4}$ | … | $P_{21\,20}$ | $P_{21\,21}$ | $P_{21\,22}$ | $P_{21\,23}$ | $P_{21\,24}$ |
| $P_{22\,1}$ | $P_{22\,2}$ | $P_{22\,3}$ | $P_{22\,4}$ | … | $P_{22\,20}$ | $P_{22\,21}$ | $P_{22\,22}$ | $P_{22\,23}$ | $P_{22\,24}$ |
| $P_{23\,1}$ | $P_{23\,2}$ | $P_{23\,3}$ | $P_{23\,4}$ | … | $P_{23\,20}$ | $P_{23\,21}$ | $P_{23\,22}$ | $P_{23\,23}$ | $P_{23\,24}$ |
| $P_{24\,1}$ | $P_{24\,2}$ | $P_{24\,3}$ | $P_{24\,4}$ | … | $P_{24\,20}$ | $P_{24\,9}$ | $P_{24\,10}$ | $P_{24\,23}$ | $P_{24\,24}$ |

Часть таблицы умножения перестановок (таблица Кэли) составленных только для образующего цикла $c_4$, имеет вид:

$P_{1\,1}$

|     | $p_0$ | $p_1$ | $p_2$ | $p_3$ | $p_4$ | $p_5$ | $p_6$ | $p_7$ |
|---|---|---|---|---|---|---|---|---|
| $p_0$ | $p_0$ | $p_1$ | $p_2$ | $p_3$ | $p_4$ | $p_5$ | $p_6$ | $p_7$ |
| $p_1$ | $p_1$ | $p_2$ | $p_3$ | $p_0$ | $p_7$ | $p_6$ | $p_4$ | $p_5$ |
| $p_2$ | $p_2$ | $p_3$ | $p_0$ | $p_1$ | $p_5$ | $p_4$ | $p_7$ | $p_6$ |
| $p_3$ | $p_3$ | $p_0$ | $p_1$ | $p_2$ | $p_6$ | $p_7$ | $p_5$ | $p_4$ |
| $p_4$ | $p_4$ | $p_6$ | $p_5$ | $p_7$ | $p_0$ | $p_2$ | $p_1$ | $p_3$ |
| $p_5$ | $p_5$ | $p_7$ | $p_4$ | $p_6$ | $p_2$ | $p_0$ | $p_3$ | $p_1$ |
| $p_6$ | $p_6$ | $p_5$ | $p_7$ | $p_4$ | $p_3$ | $p_1$ | $p_0$ | $p_2$ |
| $p_7$ | $p_7$ | $p_4$ | $p_6$ | $p_5$ | $p_1$ | $p_3$ | $p_2$ | $p_0$ |

### 4.3. Сильно регулярный граф (*lattice graph*)

***Пример 4.3.*** Будем рассматривать сильно регулярный граф *(lattice graph)* на 16 вершин представленный на рис. 4.27.

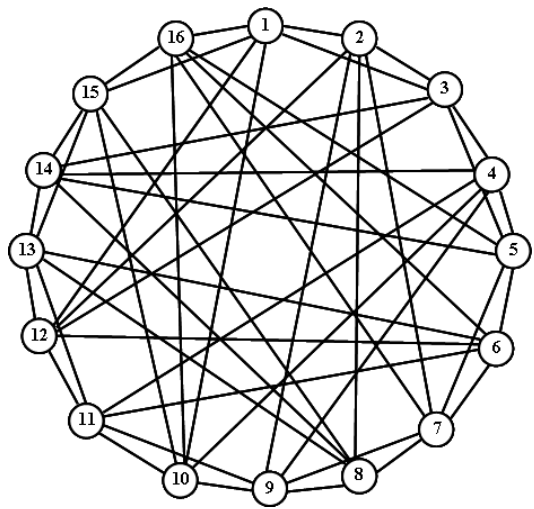

Рис. 4.27. Граф $G_3$ *(lattice graph)*.

Количество вершин в графе = 16. Количество рёбер в графе = 48.

Смежность графа:

вершина  1:  $\{v_2,v_3,v_{10},v_{12},v_{15},v_{16}\}$;
вершина  2:  $\{v_1,v_3,v_7,v_8,v_9,v_{12}\}$;
вершина  3:  $\{v_1,v_2,v_4,v_5,v_{12},v_{14}\}$;
вершина  4:  $\{v_3,v_5,v_9,v_{10},v_{11},v_{14}\}$;
вершина  5:  $\{v_3,v_4,v_6,v_7,v_{14},v_{16}\}$;
вершина  6:  $\{v_5,v_7,v_{11},v_{12},v_{13},v_{16}\}$;
вершина  7:  $\{v_2,v_5,v_6,v_8,v_9,v_{16}\}$;
вершина  8:  $\{v_2,v_7,v_9,v_{13},v_{14},v_{15}\}$;
вершина  9:  $\{v_2,v_4,v_7,v_8,v_{10},v_{11}\}$;
вершина  10:  $\{v_1,v_4,v_9,v_{11},v_{15},v_{16}\}$;
вершина  11:  $\{v_4,v_6,v_9,v_{10},v_{12},v_{13}\}$;
вершина  12:  $\{v_1,v_2,v_3,v_6,v_{11},v_{13}\}$;
вершина  13:  $\{v_6,v_8,v_{11},v_{12},v_{14},v_{15}\}$;
вершина  14:  $\{v_3,v_4,v_5,v_8,v_{13},v_{15}\}$;
вершина  15:  $\{v_1,v_8,v_{10},v_{13},v_{14},v_{16}\}$;
вершина  16:  $\{v_1,v_5,v_6,v_7,v_{10},v_{15}\}$.

Инцидентность графа:

$e_1 = (v_1,v_2) \vee (v_2,v_1)$;
$e_2 = (v_1,v_3) \vee (v_3,v_1)$;
$e_3 = (v_1,v_{10}) \vee (v_{10},v_1)$;
$e_4 = (v_1,v_{12}) \vee (v_{12},v_1)$;
$e_5 = (v_1,v_{15}) \vee (v_{15},v_1)$;
$e_6 = (v_1,v_{16}) \vee (v_{16},v_1)$;
$e_7 = (v_2,v_3) \vee (v_3,v_2)$;
$e_8 = (v_2,v_7) \vee (v_7,v_2)$;
$e_9 = (v_2,v_8) \vee (v_8,v_2)$;
$e_{10} = (v_2,v_9) \vee (v_9,v_2)$;
$e_{11} = (v_2,v_{12}) \vee (v_{12},v_2)$;
$e_{12} = (v_3,v_4) \vee (v_4,v_3)$;
$e_{13} = (v_3,v_5) \vee (v_5,v_3)$;
$e_{14} = (v_3,v_{12}) \vee (v_{12},v_3)$;
$e_{15} = (v_3,v_{14}) \vee (v_{14},v_3)$;
$e_{16} = (v_4,v_5) \vee (v_5,v_4)$;
$e_{17} = (v_4,v_9) \vee (v_9,v_4)$;
$e_{18} = (v_4,v_{10}) \vee (v_{10},v_4)$;
$e_{19} = (v_4,v_{11}) \vee (v_{11},v_4)$;
$e_{20} = (v_4,v_{14}) \vee (v_{14},v_4)$;
$e_{21} = (v_5,v_6) \vee (v_6,v_5)$;
$e_{22} = (v_5,v_7) \vee (v_7,v_5)$;
$e_{23} = (v_5,v_{14}) \vee (v_{14},v_5)$;
$e_{24} = (v_5,v_{16}) \vee (v_{16},v_5)$;
$e_{25} = (v_6,v_7) \vee (v_7,v_6)$;
$e_{26} = (v_6,v_{11}) \vee (v_{11},v_6)$;
$e_{27} = (v_6,v_{12}) \vee (v_{12},v_6)$;
$e_{28} = (v_6,v_{13}) \vee (v_{13},v_6)$;
$e_{29} = (v_6,v_{16}) \vee (v_{16},v_6)$;
$e_{30} = (v_7,v_8) \vee (v_8,v_7)$;
$e_{31} = (v_7,v_9) \vee (v_9,v_7)$;
$e_{32} = (v_7,v_{16}) \vee (v_{16},v_7)$;
$e_{33} = (v_8,v_9) \vee (v_9,v_8)$;
$e_{34} = (v_8,v_{13}) \vee (v_{13},v_8)$;
$e_{35} = (v_8,v_{14}) \vee (v_{14},v_8)$;
$e_{36} = (v_8,v_{15}) \vee (v_{15},v_8)$;
$e_{37} = (v_9,v_{10}) \vee (v_{10},v_9)$;
$e_{38} = (v_9,v_{11}) \vee (v_{11},v_9)$;
$e_{39} = (v_{10},v_{11}) \vee (v_{11},v_{10})$;
$e_{40} = (v_{10},v_{15}) \vee (v_{15},v_{10})$;
$e_{41} = (v_{10},v_{16}) \vee (v_{16},v_{10})$;
$e_{42} = (v_{11},v_{12}) \vee (v_{12},v_{11})$;
$e_{43} = (v_{11},v_{13}) \vee (v_{13},v_{11})$;
$e_{44} = (v_{12},v_{13}) \vee (v_{13},v_{12})$;
$e_{45} = (v_{13},v_{14}) \vee (v_{14},v_{13})$;
$e_{46} = (v_{13},v_{15}) \vee (v_{15},v_{13})$;
$e_{47} = (v_{14},v_{15}) \vee (v_{15}.v_{14})$;
$e_{48} = (v_{15},v_{16}) \vee (v_{16},v_{15})$.

Вектор весов ребер : $48 \times 26$;

Вектор весов вершин : $16 \times 156$.

Множество изометрических циклов графа:

$c_1 = \{e_1,e_2,e_7\} \to \{v_1,v_2,v_3\}$;
$c_2 = \{e_1,e_3,e_{10},e_{37}\} \to \{v_1,v_2,v_9,v_{10}\}$;
$c_3 = \{e_1,e_4,e_{11}\} \to \{v_1,v_2,v_{12}\}$;

$c_4 = \{e_1,e_5,e_9,e_{36}\} \rightarrow \{v_1,v_2,v_8,v_{15}\}$;
$c_5 = \{e_1,e_6,e_8,e_{32}\} \rightarrow \{v_1,v_2,v_7,v_{16}\}$;
$c_6 = \{e_2,e_3,e_{12},e_{18}\} \rightarrow \{v_1,v_3,v_4,v_{10}\}$;
$c_7 = \{e_2,e_4,e_{14}\} \rightarrow \{v_1,v_3,v_{12}\}$;
$c_8 = \{e_2,e_5,e_{15},e_{47}\} \rightarrow \{v_1,v_3,v_{14},v_{15}\}$;
$c_9 = \{e_2,e_6,e_{13},e_{24}\} \rightarrow \{v_1,v_3,v_5,v_{16}\}$;
$c_{10} = \{e_3,e_4,e_{39},e_{42}\} \rightarrow \{v_1,v_{10},v_{11},v_{12}\}$;
$c_{11} = \{e_3,e_5,e_{40}\} \rightarrow \{v_1,v_{10},v_{15}\}$;
$c_{12} = \{e_3,e_6,e_{41}\} \rightarrow \{v_1,v_{10},v_{16}\}$;
$c_{13} = \{e_4,e_5,e_{44},e_{46}\} \rightarrow \{v_1,v_{12},v_{13},v_{15}\}$;
$c_{14} = \{e_4,e_6,e_{27},e_{29}\} \rightarrow \{v_1,v_6,v_{12},v_{16}\}$;
$c_{15} = \{e_5,e_6,e_{48}\} \rightarrow \{v_1,v_{15},v_{16}\}$;
$c_{16} = \{e_7,e_8,e_{13},e_{22}\} \rightarrow \{v_2,v_3,v_5,v_7\}$;
$c_{17} = \{e_7,e_9,e_{15},e_{35}\} \rightarrow \{v_2,v_3,v_8,v_{14}\}$;
$c_{18} = \{e_7,e_{10},e_{12},e_{17}\} \rightarrow \{v_2,v_3,v_4,v_9\}$;
$c_{19} = \{e_7,e_{11},e_{14}\} \rightarrow \{v_2,v_3,v_{12}\}$;
$c_{20} = \{e_8,e_9,e_{30}\} \rightarrow \{v_2,v_7,v_8\}$;
$c_{21} = \{e_8,e_{10},e_{31}\} \rightarrow \{v_2,v_7,v_9\}$;
$c_{22} = \{e_8,e_{11},e_{25},e_{27}\} \rightarrow \{v_2,v_6,v_7,v_{12}\}$;
$c_{23} = \{e_9,e_{10},e_{33}\} \rightarrow \{v_2,v_8,v_9\}$;
$c_{24} = \{e_9,e_{11},e_{34},e_{44}\} \rightarrow \{v_2,v_8,v_{12},v_{13}\}$;
$c_{25} = \{e_{10},e_{11},e_{38},e_{42}\} \rightarrow \{v_2,v_9,v_{11},v_{12}\}$;
$c_{26} = \{e_{12},e_{13},e_{16}\} \rightarrow \{v_3,v_4,v_5\}$;
$c_{27} = \{e_{12},e_{14},e_{19},e_{42}\} \rightarrow \{v_3,v_4,v_{11},v_{12}\}$;
$c_{28} = \{e_{12},e_{15},e_{20}\} \rightarrow \{v_3,v_4,v_{14}\}$;
$c_{29} = \{e_{13},e_{14},e_{21},e_{27}\} \rightarrow \{v_3,v_5,v_6,v_{12}\}$;
$c_{30} = \{e_{13},e_{15},e_{23}\} \rightarrow \{v_3,v_5,v_{14}\}$;
$c_{31} = \{e_{14},e_{15},e_{44},e_{45}\} \rightarrow \{v_3,v_{12},v_{13},v_{14}\}$;
$c_{32} = \{e_{16},e_{17},e_{22},e_{31}\} \rightarrow \{v_4,v_5,v_7,v_9\}$;
$c_{33} = \{e_{16},e_{18},e_{24},e_{41}\} \rightarrow \{v_4,v_5,v_{10},v_{16}\}$;
$c_{34} = \{e_{16},e_{19},e_{21},e_{26}\} \rightarrow \{v_4,v_5,v_6,v_{11}\}$;
$c_{35} = \{e_{16},e_{20},e_{23}\} \rightarrow \{v_4,v_5,v_{14}\}$;
$c_{36} = \{e_{17},e_{18},e_{37}\} \rightarrow \{v_4,v_9,v_{10}\}$;
$c_{37} = \{e_{17},e_{19},e_{38}\} \rightarrow \{v_4,v_9,v_{11}\}$;
$c_{38} = \{e_{17},e_{20},e_{33},e_{35}\} \rightarrow \{v_4,v_8,v_9,v_{14}\}$;
$c_{39} = \{e_{18},e_{19},e_{39}\} \rightarrow \{v_4,v_{10},v_{11}\}$;
$c_{40} = \{e_{18},e_{20},e_{40},e_{47}\} \rightarrow \{v_4,v_{10},v_{14},v_{15}\}$;
$c_{41} = \{e_{19},e_{20},e_{43},e_{45}\} \rightarrow \{v_4,v_{11},v_{13},v_{14}\}$;
$c_{42} = \{e_{21},e_{22},e_{25}\} \rightarrow \{v_5,v_6,v_7\}$;
$c_{43} = \{e_{21},e_{23},e_{28},e_{45}\} \rightarrow \{v_5,v_6,v_{13},v_{14}\}$;
$c_{44} = \{e_{21},e_{24},e_{29}\} \rightarrow \{v_5,v_6,v_{16}\}$;
$c_{45} = \{e_{22},e_{23},e_{30},e_{35}\} \rightarrow \{v_5,v_7,v_8,v_{14}\}$;
$c_{46} = \{e_{22},e_{24},e_{32}\} \rightarrow \{v_5,v_7,v_{16}\}$;
$c_{47} = \{e_{23},e_{24},e_{47},e_{48}\} \rightarrow \{v_5,v_{14},v_{15},v_{16}\}$;
$c_{48} = \{e_{25},e_{26},e_{31},e_{38}\} \rightarrow \{v_6,v_7,v_9,v_{11}\}$;
$c_{49} = \{e_{25},e_{28},e_{30},e_{34}\} \rightarrow \{v_6,v_7,v_8,v_{13}\}$;
$c_{50} = \{e_{25},e_{29},e_{32}\} \rightarrow \{v_6,v_7,v_{16}\}$;
$c_{51} = \{e_{26},e_{27},e_{42}\} \rightarrow \{v_6,v_{11},v_{12}\}$;
$c_{52} = \{e_{26},e_{28},e_{43}\} \rightarrow \{v_6,v_{11},v_{13}\}$;
$c_{53} = \{e_{26},e_{29},e_{39},e_{41}\} \rightarrow \{v_6,v_{10},v_{11},v_{16}\}$;
$c_{54} = \{e_{27},e_{28},e_{44}\} \rightarrow \{v_6,v_{12},v_{13}\}$;

$c_{55} = \{e_{28},e_{29},e_{46},e_{48}\} \rightarrow \{v_6,v_{13},v_{15},v_{16}\}$;
$c_{56} = \{e_{30},e_{31},e_{33}\} \rightarrow \{v_7,v_8,v_9\}$;
$c_{57} = \{e_{30},e_{32},e_{36},e_{48}\} \rightarrow \{v_7,v_8,v_{15},v_{16}\}$;
$c_{58} = \{e_{31},e_{32},e_{37},e_{41}\} \rightarrow \{v_7,v_9,v_{10},v_{16}\}$;
$c_{59} = \{e_{33},e_{34},e_{38},e_{43}\} \rightarrow \{v_8,v_9,v_{11},v_{13}\}$;
$c_{60} = \{e_{33},e_{36},e_{37},e_{40}\} \rightarrow \{v_8,v_9,v_{10},v_{15}\}$;
$c_{61} = \{e_{34},e_{35},e_{45}\} \rightarrow \{v_8,v_{13},v_{14}\}$;
$c_{62} = \{e_{34},e_{36},e_{46}\} \rightarrow \{v_8,v_{13},v_{15}\}$;
$c_{63} = \{e_{35},e_{36},e_{47}\} \rightarrow \{v_8,v_{14},v_{15}\}$;
$c_{64} = \{e_{37},e_{38},e_{39}\} \rightarrow \{v_9,v_{10},v_{11}\}$;
$c_{65} = \{e_{39},e_{40},e_{43},e_{46}\} \rightarrow \{v_{10},v_{11},v_{13},v_{15}\}$;
$c_{66} = \{e_{40},e_{41},e_{48}\} \rightarrow \{v_{10},v_{15},v_{16}\}$;
$c_{67} = \{e_{42},e_{43},e_{44}\} \rightarrow \{v_{11},v_{12},v_{13}\}$;
$c_{68} = \{e_{45},e_{46},e_{47}\} \rightarrow \{v_{13},v_{14},v_{15}\}$.

Выделим изометрические циклы длиной четыре:

$c_2 = \{e_1,e_3,e_{10},e_{37}\} \rightarrow \{v_1,v_2,v_9,v_{10}\}$;
$c_4 = \{e_1,e_5,e_9,e_{36}\} \rightarrow \{v_1,v_2,v_8,v_{15}\}$;
$c_5 = \{e_1,e_6,e_8,e_{32}\} \rightarrow \{v_1,v_2,v_7,v_{16}\}$;
$c_6 = \{e_2,e_3,e_{12},e_{18}\} \rightarrow \{v_1,v_3,v_4,v_{10}\}$;
$c_8 = \{e_2,e_5,e_{15},e_{47}\} \rightarrow \{v_1,v_3,v_{14},v_{15}\}$;
$c_9 = \{e_2,e_6,e_{13},e_{24}\} \rightarrow \{v_1,v_3,v_5,v_{16}\}$;
$c_{10} = \{e_3,e_4,e_{39},e_{42}\} \rightarrow \{v_1,v_{10},v_{11},v_{12}\}$;
$c_{13} = \{e_4,e_5,e_{44},e_{46}\} \rightarrow \{v_1,v_{12},v_{13},v_{15}\}$;
$c_{14} = \{e_4,e_6,e_{27},e_{29}\} \rightarrow \{v_1,v_6,v_{12},v_{16}\}$;
$c_{16} = \{e_7,e_8,e_{13},e_{22}\} \rightarrow \{v_2,v_3,v_5,v_7\}$;
$c_{17} = \{e_7,e_9,e_{15},e_{35}\} \rightarrow \{v_2,v_3,v_8,v_{14}\}$;
$c_{18} = \{e_7,e_{10},e_{12},e_{17}\} \rightarrow \{v_2,v_3,v_4,v_9\}$;
$c_{22} = \{e_8,e_{11},e_{25},e_{27}\} \rightarrow \{v_2,v_6,v_7,v_{12}\}$;
$c_{24} = \{e_9,e_{11},e_{34},e_{44}\} \rightarrow \{v_2,v_8,v_{12},v_{13}\}$;
$c_{25} = \{e_{10},e_{11},e_{38},e_{42}\} \rightarrow \{v_2,v_9,v_{11},v_{12}\}$;
$c_{27} = \{e_{12},e_{14},e_{19},e_{42}\} \rightarrow \{v_3,v_4,v_{11},v_{12}\}$;
$c_{29} = \{e_{13},e_{14},e_{21},e_{27}\} \rightarrow \{v_3,v_5,v_6,v_{12}\}$;
$c_{31} = \{e_{14},e_{15},e_{44},e_{45}\} \rightarrow \{v_3,v_{12},v_{13},v_{14}\}$;
$c_{32} = \{e_{16},e_{17},e_{22},e_{31}\} \rightarrow \{v_4,v_5,v_7,v_9\}$;
$c_{33} = \{e_{16},e_{18},e_{24},e_{41}\} \rightarrow \{v_4,v_5,v_{10},v_{16}\}$;
$c_{34} = \{e_{16},e_{19},e_{21},e_{26}\} \rightarrow \{v_4,v_5,v_6,v_{11}\}$;
$c_{38} = \{e_{17},e_{20},e_{33},e_{35}\} \rightarrow \{v_4,v_8,v_9,v_{14}\}$;
$c_{40} = \{e_{18},e_{20},e_{40},e_{47}\} \rightarrow \{v_4,v_{10},v_{14},v_{15}\}$;
$c_{41} = \{e_{19},e_{20},e_{43},e_{45}\} \rightarrow \{v_4,v_{11},v_{13},v_{14}\}$;
$c_{43} = \{e_{21},e_{23},e_{28},e_{45}\} \rightarrow \{v_5,v_6,v_{13},v_{14}\}$;
$c_{45} = \{e_{22},e_{23},e_{30},e_{35}\} \rightarrow \{v_5,v_7,v_8,v_{14}\}$;
$c_{47} = \{e_{23},e_{24},e_{47},e_{48}\} \rightarrow \{v_5,v_{14},v_{15},v_{16}\}$;
$c_{48} = \{e_{25},e_{26},e_{31},e_{38}\} \rightarrow \{v_6,v_7,v_9,v_{11}\}$;
$c_{49} = \{e_{25},e_{28},e_{30},e_{34}\} \rightarrow \{v_6,v_7,v_8,v_{13}\}$;
$c_{53} = \{e_{26},e_{29},e_{39},e_{41}\} \rightarrow \{v_6,v_{10},v_{11},v_{16}\}$;
$c_{55} = \{e_{28},e_{29},e_{46},e_{48}\} \rightarrow \{v_6,v_{13},v_{15},v_{16}\}$;
$c_{57} = \{e_{30},e_{32},e_{36},e_{48}\} \rightarrow \{v_7,v_8,v_{15},v_{16}\}$;
$c_{58} = \{e_{31},e_{32},e_{37},e_{41}\} \rightarrow \{v_7,v_9,v_{10},v_{16}\}$;
$c_{59} = \{e_{33},e_{34},e_{38},e_{43}\} \rightarrow \{v_8,v_9,v_{11},v_{13}\}$;
$c_{60} = \{e_{33},e_{36},e_{37},e_{40}\} \rightarrow \{v_8,v_9,v_{10},v_{15}\}$;

$c_{65} = \{e_{39},e_{40},e_{43},e_{46}\} \to \{v_{10},v_{11},v_{13},v_{15}\}$.

Количество вершин в графе равно 16. Попробуем выбрать подмножество циклов так, чтобы покрыть все множество вершин. Очевидно, что в этом случае нужно выбрать четыре непересекающихся цикла. Проверив все сочетания из 36 элементов по 4, и выберем циклы, объединение которых представляет все множество вершин графа V:

1-ое подмножество:
$c_2 = \{e_1,e_3,e_{10},e_{37}\} \to \{v_1,v_2,v_9,v_{10}\}$;
$c_{27} = \{e_{12},e_{14},e_{19},e_{42}\} \to \{v_3,v_4,v_{11},v_{12}\}$;
$c_{43} = \{e_{21},e_{23},e_{28},e_{45}\} \to \{v_5,v_6,v_{13},v_{14}\}$;
$c_{57} = \{e_{30},e_{32},e_{36},e_{48}\} \to \{v_7,v_8,v_{15},v_{16}\}$.

2-ое подмножество:
$c_2 = \{e_1,e_3,e_{10},e_{37}\} \to \{v_1,v_2,v_9,v_{10}\}$;
$c_{27} = \{e_{12},e_{14},e_{19},e_{42}\} \to \{v_3,v_4,v_{11},v_{12}\}$;
$c_{45} = \{e_{22},e_{23},e_{30},e_{35}\} \to \{v_5,v_7,v_8,v_{14}\}$;
$c_{55} = \{e_{28},e_{29},e_{46},e_{48}\} \to \{v_6,v_{13},v_{15},v_{16}\}$.

3-ое подмножество:
$c_2 = \{e_1,e_3,e_{10},e_{37}\} \to \{v_1,v_2,v_9,v_{10}\}$;
$c_{27} = \{e_{12},e_{14},e_{19},e_{42}\} \to \{v_3,v_4,v_{11},v_{12}\}$;
$c_{47} = \{e_{23},e_{24},e_{47},e_{48}\} \to \{v_5,v_{14},v_{15},v_{16}\}$;
$c_{49} = \{e_{25},e_{28},e_{30},e_{34}\} \to \{v_6,v_7,v_8,v_{13}\}$.

4-ое подмножество:
$c_2 = \{e_1,e_3,e_{10},e_{37}\} \to \{v_1,v_2,v_9,v_{10}\}$;
$c_{29} = \{e_{13},e_{14},e_{21},e_{27}\} \to \{v_3,v_5,v_6,v_{12}\}$;
$c_{41} = \{e_{19},e_{20},e_{43},e_{45}\} \to \{v_4,v_{11},v_{13},v_{14}\}$;
$c_{57} = \{e_{30},e_{32},e_{36},e_{48}\} \to \{v_7,v_8,v_{15},v_{16}\}$.

5-ое подмножество:
$c_2 = \{e_1,e_3,e_{10},e_{37}\} \to \{v_1,v_2,v_9,v_{10}\}$;
$c_{31} = \{e_{14},e_{15},e_{44},e_{45}\} \to \{v_3,v_{12},v_{13},v_{14}\}$;
$c_{34} = \{e_{16},e_{19},e_{21},e_{26}\} \to \{v_4,v_5,v_6,v_{11}\}$;
$c_{57} = \{e_{30},e_{32},e_{36},e_{48}\} \to \{v_7,v_8,v_{15},v_{16}\}$.

6-е подмножество
$c_4 = \{e_1,e_5,e_9,e_{36}\} \to \{v_1,v_2,v_8,v_{15}\}$;
$c_{27} = \{e_{12},e_{14},e_{19},e_{42}\} \to \{v_3,v_4,v_{11},v_{12}\}$;
$c_{43} = \{e_{21},e_{23},e_{28},e_{45}\} \to \{v_5,v_6,v_{13},v_{14}\}$;
$c_{58} = \{e_{31},e_{32},e_{37},e_{41}\} \to \{v_7,v_9,v_{10},v_{16}\}$.

7-ое подмножество:
$c_4 = \{e_1,e_5,e_9,e_{36}\} \to \{v_1,v_2,v_8,v_{15}\}$;
$c_{29} = \{e_{13},e_{14},e_{21},e_{27}\} \to \{v_3,v_5,v_6,v_{12}\}$;
$c_{41} = \{e_{19},e_{20},e_{43},e_{45}\} \to \{v_4,v_{11},v_{13},v_{14}\}$;
$c_{58} = \{e_{31},e_{32},e_{37},e_{41}\} \to \{v_7,v_9,v_{10},v_{16}\}$.

8-ое подмножество:
$c_4 = \{e_1,e_5,e_9,e_{36}\} \to \{v_1,v_2,v_8,v_{15}\}$;
$c_{31} = \{e_{14},e_{15},e_{44},e_{45}\} \to \{v_3,v_{12},v_{13},v_{14}\}$;
$c_{32} = \{e_{16},e_{17},e_{22},e_{31}\} \to \{v_4,v_5,v_7,v_9\}$;
$c_{53} = \{e_{26},e_{29},e_{39},e_{41}\} \to \{v_6,v_{10},v_{11},v_{16}\}$.

9-е подмножество
$c_4 = \{e_1,e_5,e_9,e_{36}\} \to \{v_1,v_2,v_8,v_{15}\}$;
$c_{31} = \{e_{14},e_{15},e_{44},e_{45}\} \to \{v_3,v_{12},v_{13},v_{14}\}$;
$c_{33} = \{e_{16},e_{18},e_{24},e_{41}\} \to \{v_4,v_5,v_{10},v_{16}\}$;

$c_{48} = \{e_{25}, e_{26}, e_{31}, e_{38}\} \rightarrow \{v_6, v_7, v_9, v_{11}\}$.

10-е подмножество

$c_4 = \{e_1, e_5, e_9, e_{36}\} \rightarrow \{v_1, v_2, v_8, v_{15}\}$;
$c_{31} = \{e_{14}, e_{15}, e_{44}, e_{45}\} \rightarrow \{v_3, v_{12}, v_{13}, v_{14}\}$;
$c_{34} = \{e_{16}, e_{19}, e_{21}, e_{26}\} \rightarrow \{v_4, v_5, v_6, v_{11}\}$;
$c_{58} = \{e_{31}, e_{32}, e_{37}, e_{41}\} \rightarrow \{v_7, v_9, v_{10}, v_{16}\}$.

11-е подмножество

$c_5 = \{e_1, e_6, e_8, e_{32}\} \rightarrow \{v_1, v_2, v_7, v_{16}\}$;
$c_{27} = \{e_{12}, e_{14}, e_{19}, e_{42}\} \rightarrow \{v_3, v_4, v_{11}, v_{12}\}$;
$c_{43} = \{e_{21}, e_{23}, e_{28}, e_{45}\} \rightarrow \{v_5, v_6, v_{13}, v_{14}\}$;
$c_{60} = \{e_{33}, e_{36}, e_{37}, e_{40}\} \rightarrow \{v_8, v_9, v_{10}, v_{15}\}$.

12-е подмножество

$c_5 = \{e_1, e_6, e_8, e_{32}\} \rightarrow \{v_1, v_2, v_7, v_{16}\}$;
$c_{29} = \{e_{13}, e_{14}, e_{21}, e_{27}\} \rightarrow \{v_3, v_5, v_6, v_{12}\}$;
$c_{38} = \{e_{17}, e_{20}, e_{33}, e_{35}\} \rightarrow \{v_4, v_8, v_9, v_{14}\}$;
$c_{65} = \{e_{39}, e_{40}, e_{43}, e_{46}\} \rightarrow \{v_{10}, v_{11}, v_{13}, v_{15}\}$.

13-е подмножество

$c_5 = \{e_1, e_6, e_8, e_{32}\} \rightarrow \{v_1, v_2, v_7, v_{16}\}$;
$c_{29} = \{e_{13}, e_{14}, e_{21}, e_{27}\} \rightarrow \{v_3, v_5, v_6, v_{12}\}$;
$c_{40} = \{e_{18}, e_{20}, e_{40}, e_{47}\} \rightarrow \{v_4, v_{10}, v_{14}, v_{15}\}$;
$c_{59} = \{e_{33}, e_{34}, e_{38}, e_{43}\} \rightarrow \{v_8, v_9, v_{11}, v_{13}\}$.

14-е подмножество

$c_5 = \{e_1, e_6, e_8, e_{32}\} \rightarrow \{v_1, v_2, v_7, v_{16}\}$;
$c_{29} = \{e_{13}, e_{14}, e_{21}, e_{27}\} \rightarrow \{v_3, v_5, v_6, v_{12}\}$;
$c_{41} = \{e_{19}, e_{20}, e_{43}, e_{45}\} \rightarrow \{v_4, v_{11}, v_{13}, v_{14}\}$;
$c_{60} = \{e_{33}, e_{36}, e_{37}, e_{40}\} \rightarrow \{v_8, v_9, v_{10}, v_{15}\}$.

15-е подмножество

$c_5 = \{e_1, e_6, e_8, e_{32}\} \rightarrow \{v_1, v_2, v_7, v_{16}\}$;
$c_{31} = \{e_{14}, e_{15}, e_{44}, e_{45}\} \rightarrow \{v_3, v_{12}, v_{13}, v_{14}\}$;
$c_{34} = \{e_{16}, e_{19}, e_{21}, e_{26}\} \rightarrow \{v_4, v_5, v_6, v_{11}\}$;
$c_{60} = \{e_{33}, e_{36}, e_{37}, e_{40}\} \rightarrow \{v_8, v_9, v_{10}, v_{15}\}$.

16-е подмножество

$c_6 = \{e_2, e_3, e_{12}, e_{18}\} \rightarrow \{v_1, v_3, v_4, v_{10}\}$;
$c_{22} = \{e_8, e_{11}, e_{25}, e_{27}\} \rightarrow \{v_2, v_6, v_7, v_{12}\}$;
$c_{47} = \{e_{23}, e_{24}, e_{47}, e_{48}\} \rightarrow \{v_5, v_{14}, v_{15}, v_{16}\}$;
$c_{59} = \{e_{33}, e_{34}, e_{38}, e_{43}\} \rightarrow \{v_8, v_9, v_{11}, v_{13}\}$.

17-е подмножество

$c_6 = \{e_2, e_3, e_{12}, e_{18}\} \rightarrow \{v_1, v_3, v_4, v_{10}\}$;
$c_{24} = \{e_9, e_{11}, e_{34}, e_{44}\} \rightarrow \{v_2, v_8, v_{12}, v_{13}\}$;
$c_{47} = \{e_{23}, e_{24}, e_{47}, e_{48}\} \rightarrow \{v_5, v_{14}, v_{15}, v_{16}\}$;
$c_{48} = \{e_{25}, e_{26}, e_{31}, e_{38}\} \rightarrow \{v_6, v_7, v_9, v_{11}\}$.

18-е подмножество

$c_6 = \{e_2, e_3, e_{12}, e_{18}\} \rightarrow \{v_1, v_3, v_4, v_{10}\}$;
$c_{25} = \{e_{10}, e_{11}, e_{38}, e_{42}\} \rightarrow \{v_2, v_9, v_{11}, v_{12}\}$;
$c_{43} = \{e_{21}, e_{23}, e_{28}, e_{45}\} \rightarrow \{v_5, v_6, v_{13}, v_{14}\}$;
$c_{57} = \{e_{30}, e_{32}, e_{36}, e_{48}\} \rightarrow \{v_7, v_8, v_{15}, v_{16}\}$.

19-е подмножество

$c_6 = \{e_2, e_3, e_{12}, e_{18}\} \rightarrow \{v_1, v_3, v_4, v_{10}\}$;
$c_{25} = \{e_{10}, e_{11}, e_{38}, e_{42}\} \rightarrow \{v_2, v_9, v_{11}, v_{12}\}$;
$c_{45} = \{e_{22}, e_{23}, e_{30}, e_{35}\} \rightarrow \{v_5, v_7, v_8, v_{14}\}$;
$c_{55} = \{e_{28}, e_{29}, e_{46}, e_{48}\} \rightarrow \{v_6, v_{13}, v_{15}, v_{16}\}$.

20-е подмножество
$c_6 = \{e_2, e_3, e_{12}, e_{18}\} \rightarrow \{v_1, v_3, v_4, v_{10}\}$;
$c_{25} = \{e_{10}, e_{11}, e_{38}, e_{42}\} \rightarrow \{v_2, v_9, v_{11}, v_{12}\}$;
$c_{47} = \{e_{23}, e_{24}, e_{47}, e_{48}\} \rightarrow \{v_5, v_{14}, v_{15}, v_{16}\}$;
$c_{49} = \{e_{25}, e_{28}, e_{30}, e_{34}\} \rightarrow \{v_6, v_7, v_8, v_{13}\}$.

21-е подмножество
$c_8 = \{e_2, e_5, e_{15}, e_{47}\} \rightarrow \{v_1, v_3, v_{14}, v_{15}\}$;
$c_{22} = \{e_8, e_{11}, e_{25}, e_{27}\} \rightarrow \{v_2, v_6, v_7, v_{12}\}$;
$c_{33} = \{e_{16}, e_{18}, e_{24}, e_{41}\} \rightarrow \{v_4, v_5, v_{10}, v_{16}\}$;
$c_{59} = \{e_{33}, e_{34}, e_{38}, e_{43}\} \rightarrow \{v_8, v_9, v_{11}, v_{13}\}$.

22-е подмножество
$c_8 = \{e_2, e_5, e_{15}, e_{47}\} \rightarrow \{v_1, v_3, v_{14}, v_{15}\}$;
$c_{24} = \{e_9, e_{11}, e_{34}, e_{44}\} \rightarrow \{v_2, v_8, v_{12}, v_{13}\}$;
$c_{32} = \{e_{16}, e_{17}, e_{22}, e_{31}\} \rightarrow \{v_4, v_5, v_7, v_9\}$;
$c_{53} = \{e_{26}, e_{29}, e_{39}, e_{41}\} \rightarrow \{v_6, v_{10}, v_{11}, v_{16}\}$.

23-е подмножество
$c_8 = \{e_2, e_5, e_{15}, e_{47}\} \rightarrow \{v_1, v_3, v_{14}, v_{15}\}$;
$c_{24} = \{e_9, e_{11}, e_{34}, e_{44}\} \rightarrow \{v_2, v_8, v_{12}, v_{13}\}$;
$c_{33} = \{e_{16}, e_{18}, e_{24}, e_{41}\} \rightarrow \{v_4, v_5, v_{10}, v_{16}\}$;
$c_{48} = \{e_{25}, e_{26}, e_{31}, e_{38}\} \rightarrow \{v_6, v_7, v_9, v_{11}\}$.

24-е подмножество
$c_8 = \{e_2, e_5, e_{15}, e_{47}\} \rightarrow \{v_1, v_3, v_{14}, v_{15}\}$;
$c_{24} = \{e_9, e_{11}, e_{34}, e_{44}\} \rightarrow \{v_2, v_8, v_{12}, v_{13}\}$;
$c_{34} = \{e_{16}, e_{19}, e_{21}, e_{26}\} \rightarrow \{v_4, v_5, v_6, v_{11}\}$;
$c_{58} = \{e_{31}, e_{32}, e_{37}, e_{41}\} \rightarrow \{v_7, v_9, v_{10}, v_{16}\}$.

25-е подмножество
$c_8 = \{e_2, e_5, e_{15}, e_{47}\} \rightarrow \{v_1, v_3, v_{14}, v_{15}\}$;
$c_{25} = \{e_{10}, e_{11}, e_{38}, e_{42}\} \rightarrow \{v_2, v_9, v_{11}, v_{12}\}$;
$c_{33} = \{e_{16}, e_{18}, e_{24}, e_{41}\} \rightarrow \{v_4, v_5, v_{10}, v_{16}\}$;
$c_{49} = \{e_{25}, e_{28}, e_{30}, e_{34}\} \rightarrow \{v_6, v_7, v_8, v_{13}\}$.

26-е подмножество
$c_9 = \{e_2, e_6, e_{13}, e_{24}\} \rightarrow \{v_1, v_3, v_5, v_{16}\}$;
$c_{22} = \{e_8, e_{11}, e_{25}, e_{27}\} \rightarrow \{v_2, v_6, v_7, v_{12}\}$;
$c_{38} = \{e_{17}, e_{20}, e_{33}, e_{35}\} \rightarrow \{v_4, v_8, v_9, v_{14}\}$;
$c_{65} = \{e_{39}, e_{40}, e_{43}, e_{46}\} \rightarrow \{v_{10}, v_{11}, v_{13}, v_{15}\}$.

27-е подмножество
$c_9 = \{e_2, e_6, e_{13}, e_{24}\} \rightarrow \{v_1, v_3, v_5, v_{16}\}$;
$c_{22} = \{e_8, e_{11}, e_{25}, e_{27}\} \rightarrow \{v_2, v_6, v_7, v_{12}\}$;
$c_{40} = \{e_{18}, e_{20}, e_{40}, e_{47}\} \rightarrow \{v_4, v_{10}, v_{14}, v_{15}\}$;
$c_{59} = \{e_{33}, e_{34}, e_{38}, e_{43}\} \rightarrow \{v_8, v_9, v_{11}, v_{13}\}$.

28-е подмножество
$c_9 = \{e_2, e_6, e_{13}, e_{24}\} \rightarrow \{v_1, v_3, v_5, v_{16}\}$;
$c_{22} = \{e_8, e_{11}, e_{25}, e_{27}\} \rightarrow \{v_2, v_6, v_7, v_{12}\}$;
$c_{41} = \{e_{19}, e_{20}, e_{43}, e_{45}\} \rightarrow \{v_4, v_{11}, v_{13}, v_{14}\}$;
$c_{60} = \{e_{33}, e_{36}, e_{37}, e_{40}\} \rightarrow \{v_8, v_9, v_{10}, v_{15}\}$.

29-е подмножество
$c_9 = \{e_2, e_6, e_{13}, e_{24}\} \rightarrow \{v_1, v_3, v_5, v_{16}\}$;
$c_{24} = \{e_9, e_{11}, e_{34}, e_{44}\} \rightarrow \{v_2, v_8, v_{12}, v_{13}\}$;
$c_{40} = \{e_{18}, e_{20}, e_{40}, e_{47}\} \rightarrow \{v_4, v_{10}, v_{14}, v_{15}\}$;
$c_{48} = \{e_{25}, e_{26}, e_{31}, e_{38}\} \rightarrow \{v_6, v_7, v_9, v_{11}\}$.

30-е подмножество

$c_9 = \{e_2, e_6, e_{13}, e_{24}\} \rightarrow \{v_1, v_3, v_5, v_{16}\};$

$c_{25} = \{e_{10}, e_{11}, e_{38}, e_{42}\} \rightarrow \{v_2, v_9, v_{11}, v_{12}\};$

$c_{40} = \{e_{18}, e_{20}, e_{40}, e_{47}\} \rightarrow \{v_4, v_{10}, v_{14}, v_{15}\};$

$c_{49} = \{e_{25}, e_{28}, e_{30}, e_{34}\} \rightarrow \{v_6, v_7, v_8, v_{13}\}.$

31-е подмножество

$c_{10} = \{e_3, e_4, e_{39}, e_{42}\} \rightarrow \{v_1, v_{10}, v_{11}, v_{12}\};$

$c_{16} = \{e_7, e_8, e_{13}, e_{22}\} \rightarrow \{v_2, v_3, v_5, v_7\};$

$c_{38} = \{e_{17}, e_{20}, e_{33}, e_{35}\} \rightarrow \{v_4, v_8, v_9, v_{14}\};$

$c_{55} = \{e_{28}, e_{29}, e_{46}, e_{48}\} \rightarrow \{v_6, v_{13}, v_{15}, v_{16}\}.$

32-е подмножество

$c_{10} = \{e_3, e_4, e_{39}, e_{42}\} \rightarrow \{v_1, v_{10}, v_{11}, v_{12}\};$

$c_{17} = \{e_7, e_9, e_{15}, e_{35}\} \rightarrow \{v_2, v_3, v_8, v_{14}\};$

$c_{32} = \{e_{16}, e_{17}, e_{22}, e_{31}\} \rightarrow \{v_4, v_5, v_7, v_9\};$

$c_{55} = \{e_{28}, e_{29}, e_{46}, e_{48}\} \rightarrow \{v_6, v_{13}, v_{15}, v_{16}\}.$

33-е подмножество

$c_{10} = \{e_3, e_4, e_{39}, e_{42}\} \rightarrow \{v_1, v_{10}, v_{11}, v_{12}\};$

$c_{18} = \{e_7, e_{10}, e_{12}, e_{17}\} \rightarrow \{v_2, v_3, v_4, v_9\};$

$c_{43} = \{e_{21}, e_{23}, e_{28}, e_{45}\} \rightarrow \{v_5, v_6, v_{13}, v_{14}\};$

$c_{57} = \{e_{30}, e_{32}, e_{36}, e_{48}\} \rightarrow \{v_7, v_8, v_{15}, v_{16}\}.$

34-е подмножество

$c_{10} = \{e_3, e_4, e_{39}, e_{42}\} \rightarrow \{v_1, v_{10}, v_{11}, v_{12}\};$

$c_{18} = \{e_7, e_{10}, e_{12}, e_{17}\} \rightarrow \{v_2, v_3, v_4, v_9\};$

$c_{45} = \{e_{22}, e_{23}, e_{30}, e_{35}\} \rightarrow \{v_5, v_7, v_8, v_{14}\};$

$c_{55} = \{e_{28}, e_{29}, e_{46}, e_{48}\} \rightarrow \{v_6, v_{13}, v_{15}, v_{16}\}.$

35-е подмножество

$c_{10} = \{e_3, e_4, e_{39}, e_{42}\} \rightarrow \{v_1, v_{10}, v_{11}, v_{12}\};$

$c_{18} = \{e_7, e_{10}, e_{12}, e_{17}\} \rightarrow \{v_2, v_3, v_4, v_9\};$

$c_{47} = \{e_{23}, e_{24}, e_{47}, e_{48}\} \rightarrow \{v_5, v_{14}, v_{15}, v_{16}\};$

$c_{49} = \{e_{25}, e_{28}, e_{30}, e_{34}\} \rightarrow \{v_6, v_7, v_8, v_{13}\}.$

36-е подмножество

$c_{13} = \{e_4, e_5, e_{44}, e_{46}\} \rightarrow \{v_1, v_{12}, v_{13}, v_{15}\};$

$c_{16} = \{e_7, e_8, e_{13}, e_{22}\} \rightarrow \{v_2, v_3, v_5, v_7\};$

$c_{38} = \{e_{17}, e_{20}, e_{33}, e_{35}\} \rightarrow \{v_4, v_8, v_9, v_{14}\};$

$c_{53} = \{e_{26}, e_{29}, e_{39}, e_{41}\} \rightarrow \{v_6, v_{10}, v_{11}, v_{16}\}.$

37-е подмножество

$c_{13} = \{e_4, e_5, e_{44}, e_{46}\} \rightarrow \{v_1, v_{12}, v_{13}, v_{15}\};$

$c_{17} = \{e_7, e_9, e_{15}, e_{35}\} \rightarrow \{v_2, v_3, v_8, v_{14}\};$

$c_{32} = \{e_{16}, e_{17}, e_{22}, e_{31}\} \rightarrow \{v_4, v_5, v_7, v_9\};$

$c_{53} = \{e_{26}, e_{29}, e_{39}, e_{41}\} \rightarrow \{v_6, v_{10}, v_{11}, v_{16}\}.$

38-е подмножество

$c_{13} = \{e_4, e_5, e_{44}, e_{46}\} \rightarrow \{v_1, v_{12}, v_{13}, v_{15}\};$

$c_{17} = \{e_7, e_9, e_{15}, e_{35}\} \rightarrow \{v_2, v_3, v_8, v_{14}\};$

$c_{33} = \{e_{16}, e_{18}, e_{24}, e_{41}\} \rightarrow \{v_4, v_5, v_{10}, v_{16}\};$

$c_{48} = \{e_{25}, e_{26}, e_{31}, e_{38}\} \rightarrow \{v_6, v_7, v_9, v_{11}\}.$

39-е подмножество

$c_{13} = \{e_4, e_5, e_{44}, e_{46}\} \rightarrow \{v_1, v_{12}, v_{13}, v_{15}\};$

$c_{17} = \{e_7, e_9, e_{15}, e_{35}\} \rightarrow \{v_2, v_3, v_8, v_{14}\};$

$c_{34} = \{e_{16}, e_{19}, e_{21}, e_{26}\} \rightarrow \{v_4, v_5, v_6, v_{11}\};$

$c_{58} = \{e_{31}, e_{32}, e_{37}, e_{41}\} \rightarrow \{v_7, v_9, v_{10}, v_{16}\}.$

40-е подмножество

$c_{13} = \{e_4, e_5, e_{44}, e_{46}\} \rightarrow \{v_1, v_{12}, v_{13}, v_{15}\};$

$c_{18} = \{e_7, e_{10}, e_{12}, e_{17}\} \to \{v_2, v_3, v_4, v_9\}$;
$c_{48} = \{e_{25}, e_{26}, e_{31}, e_{38}\} \to \{v_6, v_7, v_9, v_{11}\}$;
$c_{53} = \{e_{26}, e_{29}, e_{39}, e_{41}\} \to \{v_6, v_{10}, v_{11}, v_{16}\}$.

41-е подмножество

$c_{14} = \{e_4, e_6, e_{27}, e_{29}\} \to \{v_1, v_6, v_{12}, v_{16}\}$;
$c_{16} = \{e_7, e_8, e_{13}, e_{22}\} \to \{v_2, v_3, v_5, v_7\}$;
$c_{38} = \{e_{17}, e_{20}, e_{33}, e_{35}\} \to \{v_4, v_8, v_9, v_{14}\}$;
$c_{65} = \{e_{39}, e_{40}, e_{43}, e_{46}\} \to \{v_{10}, v_{11}, v_{13}, v_{15}\}$.

42-е подмножество

$c_{14} = \{e_4, e_6, e_{27}, e_{29}\} \to \{v_1, v_6, v_{12}, v_{16}\}$;
$c_{16} = \{e_7, e_8, e_{13}, e_{22}\} \to \{v_2, v_3, v_5, v_7\}$;
$c_{40} = \{e_{18}, e_{20}, e_{40}, e_{47}\} \to \{v_4, v_{10}, v_{14}, v_{15}\}$;
$c_{59} = \{e_{33}, e_{34}, e_{38}, e_{43}\} \to \{v_8, v_9, v_{11}, v_{13}\}$.

43-е подмножество

$c_{14} = \{e_4, e_6, e_{27}, e_{29}\} \to \{v_1, v_6, v_{12}, v_{16}\}$;
$c_{16} = \{e_7, e_8, e_{13}, e_{22}\} \to \{v_2, v_3, v_5, v_7\}$;
$c_{41} = \{e_{19}, e_{20}, e_{43}, e_{45}\} \to \{v_4, v_{11}, v_{13}, v_{14}\}$;
$c_{60} = \{e_{33}, e_{36}, e_{37}, e_{40}\} \to \{v_8, v_9, v_{10}, v_{15}\}$.

44-е подмножество

$c_{14} = \{e_4, e_6, e_{27}, e_{29}\} \to \{v_1, v_6, v_{12}, v_{16}\}$;
$c_{17} = \{e_7, e_9, e_{15}, e_{35}\} \to \{v_2, v_3, v_8, v_{14}\}$;
$c_{32} = \{e_{16}, e_{17}, e_{22}, e_{31}\} \to \{v_4, v_5, v_7, v_9\}$;
$c_{65} = \{e_{39}, e_{40}, e_{43}, e_{46}\} \to \{v_{10}, v_{11}, v_{13}, v_{15}\}$.

45-е подмножество

$c_{14} = \{e_4, e_6, e_{27}, e_{29}\} \to \{v_1, v_6, v_{12}, v_{16}\}$;
$c_{18} = \{e_7, e_{10}, e_{12}, e_{17}\} \to \{v_2, v_3, v_4, v_9\}$;
$c_{45} = \{e_{22}, e_{23}, e_{30}, e_{35}\} \to \{v_5, v_7, v_8, v_{14}\}$;
$c_{65} = \{e_{39}, e_{40}, e_{43}, e_{46}\} \to \{v_{10}, v_{11}, v_{13}, v_{15}\}$.

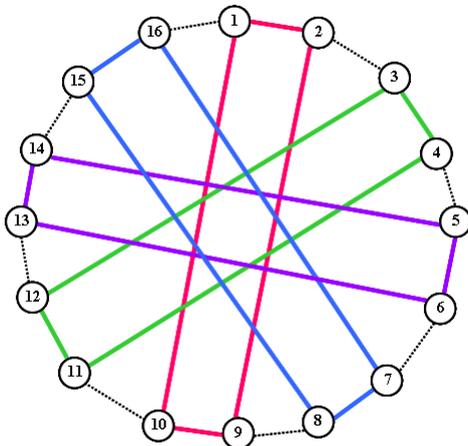

Рис. 4.31. Схематическое расположение четырех непересекающихся циклов.

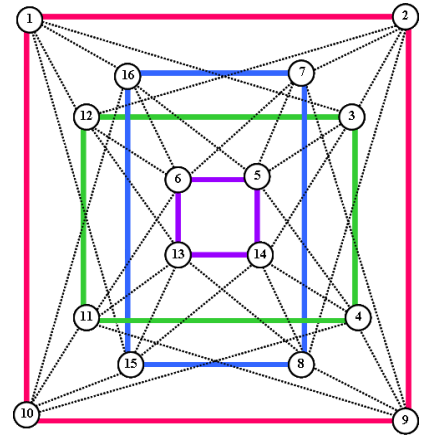

Рис. 4.32. Выбор опорного цикла из подмножества циклов.

Схематическое расположение четырех циклов из подмножества 1, для опорного цикла $\{v_1, v_2, v_9, v_{10}\}$, представлено на рисунке 4.31. Раскрасим циклы 1-го подмножества в разные цвета (см. рис. 4.32):

$c_2 = \{e_1, e_3, e_{10}, e_{37}\} \to \{v_1, v_2, v_9, v_{10}\}$      красный цвет
$c_{27} = \{e_{12}, e_{14}, e_{19}, e_{42}\} \to \{v_3, v_4, v_{11}, v_{12}\}$      зеленый цвет
$c_{43} = \{e_{21}, e_{23}, e_{28}, e_{45}\} \to \{v_5, v_6, v_{13}, v_{14}\}$      фиолетовый цвет

$c_{57} = \{e_{30}, e_{32}, e_{36}, e_{48}\} \rightarrow \{v_7, v_8, v_{15}, v_{16}\}$  синий цвет

Симметричный топологический рисунок графа представлен на рис. 4.33.

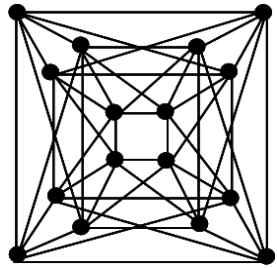

Рис. 4.33. Симметричный топометрический рисунок графа $G_3$.

Если выбрать изометрический цикл $\{v_1, v_2, v_9, v_{10}\}$ в качестве образующего цикла (см. рис. 4.32), то его преобразования индуцируют 8 перестановок. Если выбрать в качестве образующего цикла любой другой изометрический цикл из 1-го подмножества, то его преобразования индуцируют такое же количество перестановок.

Особенность построения симметричного топологического рисунка для образующего цикла $c_2$, состоит в том, что для удовлетворения свойства симметрии приходится объединять синий и зеленый циклы в один цикл (см. рис. 4.35), так как зеленый и синий циклы топологически не различимы. Различие происходит только при геометрическом преобразовании.

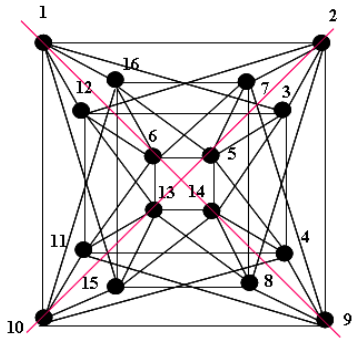
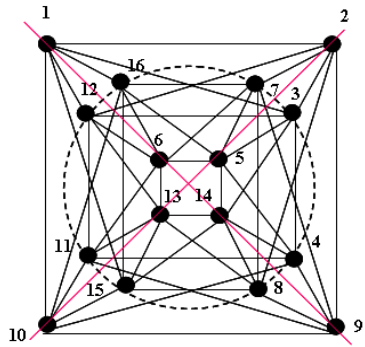
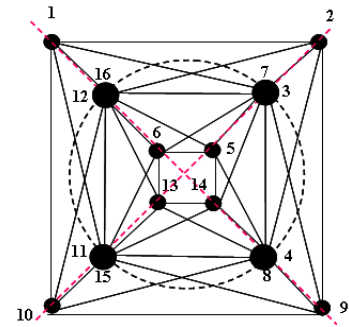

Рис. 4.34. Симметричный топологический рисунок для образующкго цикла $c_2$.

Рис. 4.35. Объединение циклов $\{v_3, v_4, v_{11}, v_{12}\}$ и $\{v_7, v_8, v_{15}, v_{16}\}$ в один цикл.

И тогда количество вариантов выборов образующих циклов в 1-ом подмножестве, также и в других подмножествах непересекающихся циклов, равно 4 (см. рис. 4.36 – 4.39).

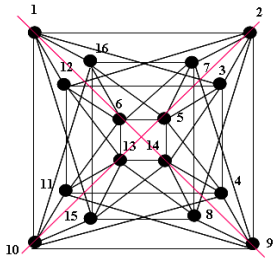 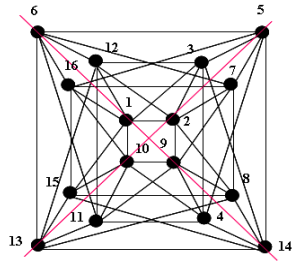 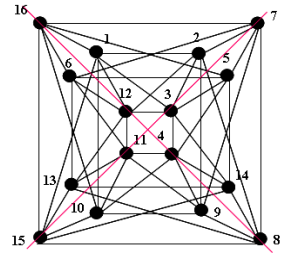 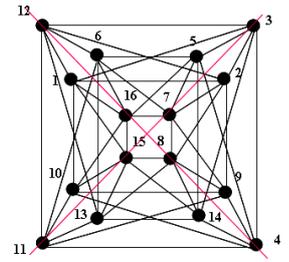

Рис. 4.36. Образующий цикл $c_2$.     Рис. 4.37. Образующий цикл $c_{43}$.     Рис. 4.38. Образующий цикл $c_{57}$.     Рис. 4.39. Образующий цикл $c_{27}$.

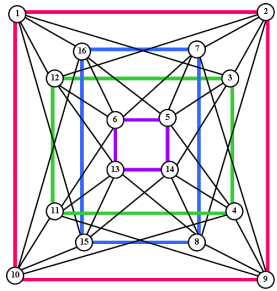 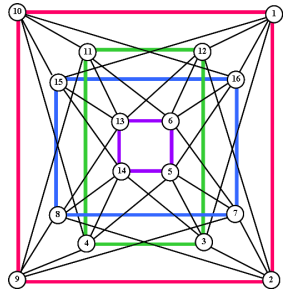 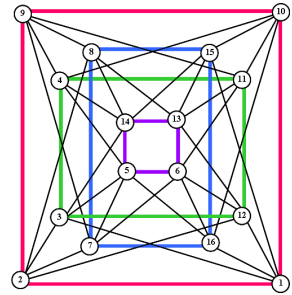 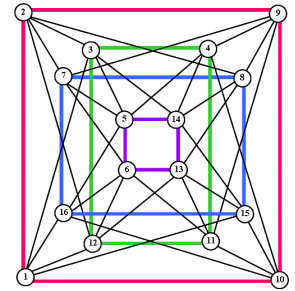

Рис. 4.40. Перестановка $p_0$.     Рис. 4.41. Перестановка $p_1$.     Рис. 4.42. Перестановка $p_2$.     Рис. 4.43. Перестановка $p_3$.

Рассмотрим перестановки, определяемые образующего циклом $\{v_1, v_2, v_9, v_{10}\}$ и порождающие 1-ую подгруппу. На рис. 4.40 – 4.43 представлено геометрическое преобразование, зеленый прямоугольник занимает место синего прямоугольника, а синий прямоугольник занимает место зеленого прямоугольника. Таким образом, для любого подмножества непересекающихся по вершинам циклов, имеем 4×8 = 32 перестановок.

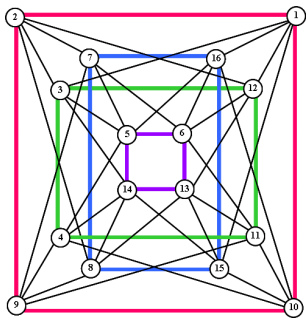 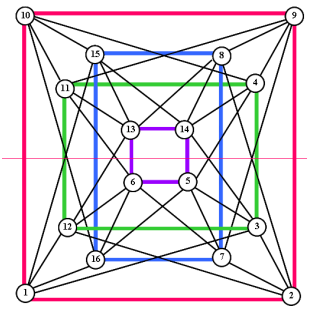 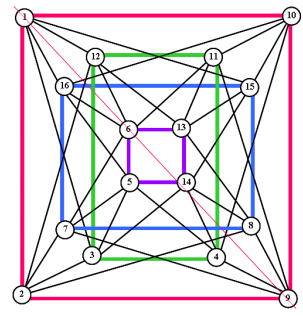 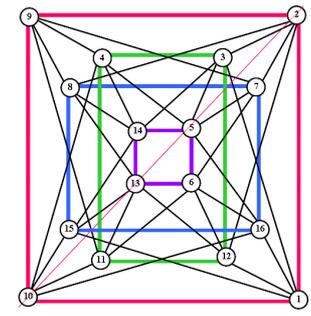

Рис. 4.44. Перестановка $p_4$.     Рис. 4.45. Перестановка $p_5$.     Рис. 4.46. Перестановка $p_6$.     Рис. 4.47. Перестановка $p_7$.

Всего имеем 32 таких подмножеств, поэтому всего перестановок в графе будет 32×32 = 1024.

В результате получим следующие перестановки:

$p_0$ = <1 2 3 4 5 6 7 8 9 10 11 12 13 14 15 16> = (1)(2)(3)(4)(5)(6)(7)(8)(9)(10)(11)(12)(13)(14)(15)(16);

$p_1$ = <7 8 1 2 3 4 5 6 15 16 9 10 11 12 13 14> = (1 7 5 3)(2 8 6 4)(9 15 13 11)(10 16 14 12);

$p_2$ = <5 6 7 8 1 2 3 4 13 14 15 16 9 10 11 12> = (1 5)(2 6)(3 7)(4 8)(9 13)(10 14)(11 15)(12 16);

$p_3$ = <3 4 5 6 7 8 1 2 11 12 13 14 15 16 9 10> = (1 3 5 7)(2 4 6 8)(9 11 13 15)(10 12 14 16);

$p_4$ = <3 2 1 8 7 6 5 4 11 10 9 16 15 14 13 12> = (1 3)(2)(4 8)(5 7)(6)(9 11)(10)(12 16)(13 15)(14);

$p_5$ = <7 6 5 4 3 2 1 8 15 14 13 12 11 10 9 16> = (1 7)(2 6)(3 5)(4)(8)(9 15)(10 14)(11 13)(12)(16);

$p_6$ = <1 8 7 6 5 4 3 2 9 16 15 14 13 12 11 10> = (1)(2 8)(3 7)(4 6)(5)(9)(10 16)(11 15)(12 14)(13);

$p_7$ = <5 4 3 2 1 8 7 6 13 12 11 10 9 16 15 14> = (1 5)(2 4)(3)(6 8)(7)(9 13)(10 12)(11)(14 16)(15).

Таблица Кэли для графа $G_3$:

$Aut(G_2) = $

| $P_{1\,1}$ | $P_{1\,2}$ | $P_{1\,3}$ | $P_{1\,4}$ | … | $P_{1\,125}$ | $P_{1\,126}$ | $P_{1\,127}$ | $P_{1\,128}$ |
|---|---|---|---|---|---|---|---|---|
| $P_{2\,1}$ | $P_{2\,2}$ | $P_{2\,3}$ | $P_{2\,4}$ | … | $P_{2\,125}$ | $P_{2\,126}$ | $P_{2\,127}$ | $P_{2\,128}$ |
| $P_{3\,1}$ | $P_{3\,2}$ | $P_{3\,3}$ | $P_{3\,4}$ | … | $P_{3\,125}$ | $P_{3\,126}$ | $P_{3\,127}$ | $P_{3\,128}$ |
| $P_{4\,1}$ | $P_{4\,2}$ | $P_{4\,3}$ | $P_{4\,4}$ | … | $P_{4\,125}$ | $P_{4\,126}$ | $P_{4\,127}$ | $P_{4\,128}$ |
| $P_{5\,1}$ | $P_{5\,2}$ | $P_{5\,3}$ | $P_{5\,4}$ | … | $P_{5\,125}$ | $P_{5\,126}$ | $P_{5\,127}$ | $P_{5\,128}$ |
| $P_{6\,1}$ | $P_{6\,2}$ | $P_{6\,3}$ | $P_{6\,4}$ | … | $P_{6\,125}$ | $P_{6\,126}$ | $P_{6\,127}$ | $P_{6\,128}$ |
| … | … | … | … | … | … | … | … | … |
| $P_{124\,1}$ | $P_{124\,2}$ | $P_{124\,3}$ | $P_{124\,4}$ | … | $P_{124\,125}$ | $P_{124\,126}$ | $P_{124\,127}$ | $P_{124\,128}$ |
| $P_{125\,1}$ | $P_{125\,2}$ | $P_{125\,3}$ | $P_{125\,4}$ | … | $P_{125\,125}$ | $P_{125\,126}$ | $P_{125\,127}$ | $P_{125\,128}$ |
| $P_{126\,1}$ | $P_{126\,2}$ | $P_{126\,3}$ | $P_{126\,4}$ | … | $P_{126\,125}$ | $P_{126\,126}$ | $P_{126\,127}$ | $P_{126\,128}$ |
| $P_{127\,1}$ | $P_{127\,2}$ | $P_{127\,3}$ | $P_{127\,4}$ | … | $P_{127\,125}$ | $P_{127\,126}$ | $P_{127\,1276}$ | $P_{127\,128}$ |
| $P_{128\,1}$ | $P_{128\,2}$ | $P_{128\,3}$ | $P_{128\,4}$ | … | $P_{128\,125}$ | $P_{128\,126}$ | $P_{128\,127}$ | $P_{128\,128}$ |

Часть таблицы умножения перестановок (таблица Кэли) составленных только для образующего цикла $c_4$ имеет вид:

**Комментарии**

В данной главе, представлены методы определения группы автоморфизмов для сильно регулярных не планарных графов. Рассмотрены методы построения образующих циклов для данной группы графов.

# ВЫВОДЫ

Методы исследования сложных систем по частям [15], позволили выделить несепарабельную часть графа. Это в свою очередь, позволило применить алгебраические методы для решения практических задач [4-6]. Применение цикломатических свойств графа, связывает в единое целое многие задачи теории графов: задачу определения изоморфизма графов, построение топологического рисунка графа, построение группы автоморфизмов, выделение изоморфного подграфа и т.д.

В данной работе рассмотрена связь между векторными инвариантами графа и задачей построения автоморфизма графа. Результаты исследования построения группы автоморфизмов для различных дискретных структур, можно представить в виде:

- представленные методы описания структур характерны только для несепарабельных графовых структур;
- перестановка вершин в графе, может осуществляться только в орбитах имеющих одинаковый вес элементов;
- интегральный и цифровой инварианты графа позволяют определять вес вершины в дискретной графовой структуре;
- *образующий цикл* – это линейная комбинация изометрических циклов графа, вершины которой имеют одинаковый вес и принадлежат одной орбите;
- топологическому рисунку образующего цикла можно поставить в соответствие плоский правильный геометрический многоугольник в качестве обода;
- образующий цикл можно располагать на плоскости в качестве внешнего контура;
- при перемещении вершин по внешнему контуру, порядок последовательности расположения вершин во вращении не меняется, но может меняться направление вращения вершин;
- порядок подгруппы перестановок для образующего цикла определяется его удвоенной длиной;
- размер таблицы Кэли определяется количеством образующего циклов.

В работе представлен алгоритм, основанный на методах алгебры структурных чисел, для определения изоморфизма матриц смежностей при замене строк и столбцов.

Произведен сравнительный анализ геометрического и топологического методов построения перестановок вершин в графе. Рассмотрены основные свойства перестановок.

Данная работа не ставит целью произвести разбиение графов на классы определяющие принципы построения орбит, а ставит целью только показать связь между группой автоморфизмов графа и его интегральным и цифровым инвариантами. Поэтому в качестве иллюстраций рассматриывется незначительное множество графов.

Построение группы автоморфизма графа носит индивидуальный характер, однако выбор совпадающих весов для вершин графа позволяет облегчить создание орбит и опорных циклов. Построение группы автоморфизмов графа, также зависит от выбора топологического рисунка графа основанного на понятии вращения вершин [36].